\documentclass[11pt]{amsart}
\usepackage{enumerate}
\usepackage{amssymb}
\usepackage{bm}
\usepackage[centertags]{amsmath}
\usepackage{amsfonts}
\usepackage{amsthm}
\usepackage{mathrsfs}
\usepackage{mathtools,xfrac}
\usepackage[usenames,dvipsnames,svgnames,table]{xcolor}
\usepackage[utf8]{inputenc}
\usepackage[T1]{fontenc}
\usepackage{enumitem}
\usepackage{hyperref} 
\usepackage[normalem]{ulem}
\usepackage{greek_enumerate}

\makeatletter
\def\l@subsection{\@tocline{2}{0pt}{2.5pc}{2.5pc}{}}%
\makeatother
\DeclareRobustCommand{\SkipTocEntry}[5]{}

\DeclareFontFamily{OT1}{pzc}{}
\DeclareFontShape{OT1}{pzc}{m}{it}{<-> s * [1.10] pzcmi7t}{}
\DeclareMathAlphabet{\mathpzc}{OT1}{pzc}{m}{it}

\DeclareMathOperator\age{age}

\DeclareMathOperator\rank{rank}
\DeclareMathOperator\supp{supp}
\DeclareMathOperator{\dist}{dist}

\numberwithin{equation}{section}

\newtheorem{theorem}{Theorem}[section]
\newtheorem{lemma}[theorem]{Lemma}
\newtheorem{proposition}[theorem]{Proposition}
\newtheorem{corollary}[theorem]{Corollary}

\newtheorem{fact}[theorem]{Fact}
\newtheorem*{claim}{Claim}
\newtheorem{definition}[theorem]{Definition}

\newtheorem{notation}[theorem]{Notation}
\newtheorem{assumption}[theorem]{Assumption}
\newtheorem{example}[theorem]{Example}

\newtheorem{remark}[theorem]{Remark}
\newtheorem{remarks}[theorem]{Remarks}

\newtheorem{problem}{Problem}

\newcommand{\N}{\mathbb{N}}

\newcommand{\T}{\mathcal T}

\makeatletter
\@namedef{subjclassname@2020}{%
  \textup{2020} Mathematics Subject Classification}
\makeatother

\allowdisplaybreaks

\begin{document}
\title[Reflexive Calkin algebras]{Reflexive Calkin algebras}

\author[P. Motakis]{Pavlos Motakis}

\address{Department of Mathematics and Statistics, York University, 4700 Keele Street, Toronto, Ontario, M3J 1P3, Canada}

\email{pmotakis@yorku.ca}

\author[A. Pelczar-Barwacz]{Anna Pelczar-Barwacz}

\address{Jagiellonian University, Faculty of Mathematics and Computer Science, Institute of Mathematics, \L ojasiewicza 6, 30-348 Krak\'ow, Poland}

\email{anna.pelczar@uj.edu.pl}

\thanks{The first author was supported by NSERC Grant RGPIN-2021-03639. The second author was supported by the grant of the National Science Centre (NCN), Poland,
no. UMO-2020/39/B/ST1/01042.}


\begin{abstract}
For a Banach space $X$ denote by $\mathcal{L}(X)$ the algebra of bounded linear operators on $X$, by $\mathcal{K}(X)$ the compact operator ideal on $X$, and by $\mathpzc{Cal}(X) = \mathcal{L}(X)/\mathcal{K}(X)$ the Calkin algebra of $X$. We prove that $\mathpzc{Cal}(X)$ can be an infinite-dimensional reflexive Banach space, even isomorphic to a Hilbert space. More precisely, for every Banach space $U$ with a normalized unconditional basis not having a $c_0$ asymptotic version we construct a Banach space $\mathfrak{X}_U$ and a sequence of mutually annihilating projections $(I_s)_{s=1}^\infty$ on $\mathfrak{X}_U$, i.e., $I_sI_t = 0$, for $s\neq t$, such that $\mathcal{L}(\mathfrak{X}_U) = \mathcal{K}(\mathfrak{X}_U)\oplus[(I_s)_{s=1}^\infty]\oplus\mathbb{C}I$ and $(I_s)_{s=1}^\infty$ is equivalent to $(u_s)_{s=1}^\infty$. In particular, $\mathpzc{Cal}(\mathfrak{X}_U)$ is isomorphic, as a Banach algebra, to the unitization of $U$ with coordinate-wise multiplication. Banach spaces $U$ meeting these criteria include $\ell_p$ and $(\oplus_n\ell_\infty^n)_p$, $1\leq p<\infty$, with their unit vector bases, $L_p$, $1 <p<\infty$, with the Haar system, the asymptotic-$\ell_1$ Tsirelson space and Schlumprecht space with their usual bases, and many others.
\end{abstract}

\maketitle

\tableofcontents

\section{Introduction}

For a Banach space $X$, let $\mathcal{L}(X)$ denote the algebra of bounded linear operators on $X$ and $\mathcal{K}(X)$ denote the compact operator ideal on $X$. We denote by $\mathpzc{Cal}(X) = \mathcal{L}(X)/\mathcal{K}(X)$ the Calkin algebra of $X$. This object was introduced for $X = \ell_2$ by Calkin in 1941 (see \cite{calkin:1941}), who proved that it is simple, thus showing that $\mathcal{L}(\ell_2)$ has a unique non-trivial closed ideal. The Calkin algebra of $\ell_2$ is most famously known for its connections to $K$-theory and Fredholm theory (see \cite{brown:dougla:fillmore:1977}) but also for the independence of ZFC of the existence of outer automorphisms on it (\cite{phillips:weaver:2007} and \cite{farah:2011}). The structural properties of the Calkin algebra of general Banach spaces X were first explored in 1974 in \cite{caradus:pfaffenberger:yood:1974}. Under the mild assumption on a Banach space having the approximation property, $\mathcal{K}(X)$ is the minimum non-trivial closed ideal in $\mathcal{L}(X)$, and hence, for spaces with this property, $\mathpzc{Cal}(X)$ captures the algebraic and analytic structure of $\mathcal{L}(X)$, modulo its minimum non-trivial closed ideal. Therefore, the question of what unital Banach algebras $\mathcal{B}$ are Calkin algebras, i.e., for some Banach space $X$, $\mathcal{B}$ is isomorphic as a Banach algebra to $\mathpzc{Cal}(X)$, is very rudimentary. Nonetheless, it was, until fairly recently, completely inaccessible. This changed with the eventual discovery of powerful Banach spaces construction techniques due to Gowers and Maurey in 1993 (\cite{gowers:maurey:1993} and \cite{gowers:maurey:1997}) and Argyros and Haydon in 2011 (\cite{argyros:haydon:2011}). These relied on earlier landmark Banach spaces such as Tsirelson space (\cite{tsirelson:1974} and \cite{figiel:johnson:1974}), the Bourgain-Delbaen spaces (\cite{bourgain:delbaen:1980}), and Schlumprecht space (\cite{schlumprecht:1991}). We partly revisit the long history behind these techniques in Section \ref{section roadmap}. The famous Argyros-Haydon space, through solving the longstanding scalar-plus-compact problem of Lindenstrauss from \cite{lindenstrauss:1975}, proved that the scalar field is a Calkin algebra. Since then, several more explicit Banach algebras have been shown to be Calkin algebras, e.g., the semigroup algebra of $\mathbb{N}_0$ (\cite{tarbard:2013}), all finite-dimensional semisimple complex algebras (\cite{kania:laustsen:2017}), and all separable commutative $C^*$-algebras (\cite{motakis:puglisi:zisimopoulou:2016} and \cite{motakis:2024}). There are other examples coming from the realm of Banach space theory, such as James space, a specific hereditarily indecomposable Banach algebra, and all non-reflexive Banach spaces with an unconditional basis, albeit with a non-standard multiplication (\cite{motakis:puglisi:tolias:2020}). A question whose answer remained elusive throughout this process is whether $\mathpzc{Cal}(X)$ can possibly be a reflexive and infinite dimensional Banach space. This is a difficult question because, in most operator algebraic contexts, reflexivity implies finite dimensionality. For example, infinite-dimensional $C^*$-algebras are automatically non-reflexive (see, e.g., \cite[page 288]{kadison:ringrose:1983}), and for an infinite dimensional Banach space $X$, $\mathcal{L}(X)$ is always non-reflexive (see \cite{baker:1982}). If, additionally, $X$ has the approximation property, then $\mathcal{K}(X)$ is non-reflexive as well (see \cite{baker:1982}). Therefore, normally, the reflexivity of $\mathpzc{Cal}(X)$ cannot be deduced by the reflexivity of $\mathcal{L}(X)$ or $\mathcal{K}(X)$. Furthermore, for any Banach space $X$ admitting an infinite dimensional Schauder decomposition, $\mathcal{L}(X)$ and $\mathpzc{Cal}(X)$ are non-reflexive. This last Banach-space-theoretic observation drove this research, and it was essential in the second author's paper \cite{pelczar-barwacz:2023} where a reflexive algebra $\mathcal{L}(X)/\mathcal{SS}(X)$ was constructed (see Section \ref{section roadmap}).

We prove that it is possible for $\mathpzc{Cal}(X)$ to be reflexive and infinite dimensional. Before precisely stating our result, we recall some necessary terminology. We use $\N$ to denote the  set of positive integers, i.e. $\N=\{1,2,\dots\}$, whereas we write $\N_0=\{0\}\cup\N$. We use the notation of $\langle A\rangle$ (respectively $[A]$) for the linear span (respectively for the closed linear span) of a subset $A$ of a linear space (respectively a normed space). 
A linear bijection $T$ between two Banach spaces $X$ and $Y$ is called a $C$-isomorphism, for some $C\geq 1$, if $\|T\|\|T^{-1}\|\leq C$. If such $T$ exists, we will say that $X$ and $Y$ are $C$-isomorphic and write $X\simeq^CY$. A Schauder basis of a Banach space $X$ is a sequence $(x_s)_{s=1}^\infty$ in $X$ such that every $x\in X$ admits a representation $x = \sum_{s=1}^\infty a_sx_s$ for a unique scalar sequence $(a_s)_{s=1}^\infty$. A Schauder basis $(x_s)_{s=1}^\infty$ of $X$ has a biorthogonal sequence $(x_s^*)_{s=1}^\infty$ in $X^*$, i.e., $x_{t}^*(\sum_{s=1}^\infty a_sx_s) = a_t$. Also, for $x\in X$, we define its support (with respect to the basis $(x_s)_{s=1}^\infty$) as the set $\mathrm{supp}(x) = \{s\in\N:x^*_s(x)\neq 0\}$ and for $f\in X^*$, we define $\mathrm{supp}(f) = \{s\in\mathbb{N}:f(x_s)\neq 0\}$. We will call a sequence $(y_k)_{k\in I}$ in $X$, indexed over an interval $I\subset\N$, satisfying $\mathrm{\supp}(y_{k-1})<\mathrm{\supp}(y_k)$, $k\in I\setminus \min I$, a block sequence. In $X^*$, block sequences are defined analogously. 
A Schauder basis is called normalized if, for all $s\in\mathbb{N}$, $\|x_s\|=1$. It is called a $C$-unconditional basis, for some $C\geq 1$, if for all scalars $a_1,\ldots,a_n$, $b_1,\ldots,b_n$ such that, for $1\leq s\leq n$, $|a_s|\leq |b_s|$,
\[\Big\|\sum_{s=1}^na_sx_s\Big\|\leq C\Big\|\sum_{s=1}^nb_sx_s\Big\|.\]
A Schauder basis that is $C$-unconditional, for some $C\geq1$, is called an unconditional basis. A sequence in a Banach space is called a Schauder basic sequence if it is a Schauder basis of its closed linear span. Similarly, it is called an unconditional sequence if it is an unconditional basis of its closed linear span. Two commonly indexed
sequences $(y_s)$, $(z_s)$, in potentially different Banach spaces, are called $C$-equivalent, for some $C\geq 1$, if the map $Ty_s = z_s$ extends to a $C$-isomorphism between their closed linear spans. For a Banach space $X$ and closed subspaces $X_1,\ldots,X_n$ of $X$, $X = X_1\oplus\cdots\oplus X_n$ means that for every $x\in X$ there are unique $x_1\in X_1,\ldots,x_n\in X_n$ such that $x = x_1+\cdots+x_n$. Then, by the closed-graph theorem, the maps $x\mapsto x_i$, $1\leq i\leq n$, are bounded linear operators. For $n\in\N$ and $1\leq p\leq \infty$, $\ell_p^n$ denotes $\mathbb{C}^n$ with the $\|\cdot\|_p$-norm.

If a Banach space $U$ has a normalized 1-unconditional basis $(u_s)_{s=1}^\infty$, then the operation of coordinate-wise multiplication
\[\Big(\sum_{s=1}^\infty a_su_s\Big)\Big(\sum_{s=1}^\infty b_su_s\Big) = \sum_{s=1}^\infty a_sb_su_s\]
makes $U$ a non-unital Banach algebra. Since every normalized unconditional basis of a Banach space $U$ can be made 1-unconditional by passing to an equivalent norm, $U$ can always be seen as a Banach algebra with coordinate-wise multiplication.  
Below, we recall the classical notion of an asymptotic version $Z$ of a Banach space $X$ with a basis from \cite[Page 151]{maurey:milman:tomczak-jaegermann:1995} for the special case $Z=c_0$.

\begin{definition}\label{def asymptotic version}
Let $X$ be a Banach space with a Schauder basis $(x_s)_{s=1}^\infty$. We say that $X$ has a $c_0$ asymptotic version if for every $n\in\N$ and $\varepsilon >0$ the following holds:
\begin{align*}
    \text{for every }k_1\in\mathbb{N}\text{ there exists }&y_1\in X\text{ with }\mathrm{supp}(y_1) > k_1,\\
    \text{for every }k_2\in\mathbb{N}\text{ there exists }&y_2\in X\text{ with }\mathrm{supp}(y_2) > k_2,\\
    &\vdots\\
    \text{for every }k_n\in\mathbb{N}\text{ there exists }&y_n\in X\text{ with }\mathrm{supp}(y_n) > k_n,
\end{align*}
such that $(y_i)_{i=1}^n$ is $(1+\varepsilon)$-equivalent to the unit vector basis of $\ell_\infty^n$.
\end{definition}

Our paper is devoted to proving the following.
\begin{theorem}
\label{the theorem}
Let $U$ be a Banach space with a normalized $1$-unconditional basis $(u_s)_{s=1}^\infty$ not having a $c_0$ asymptotic version. Then, there exists a Banach space $\mathfrak{X}_U$ with a Schauder basis and a sequence of infinite-rank norm-one projections $(I_s)_{s=1}^\infty$ on $\mathfrak{X}_U$ such that the following hold.
\begin{enumerate}[label=(\roman*),leftmargin=23pt]
    
    \item The projections $I_s$, $s\in\mathbb{N}$, are mutually annihilating, i.e., for $s\neq t\in\mathbb{N}$, $I_sI_t = 0$.

    \item The sequence $(I_s)_{s=1}^\infty$ in $\mathcal{L}(\mathfrak{X}_U)$ and its canonical image in $\mathpzc{Cal}(\mathfrak{X}_U)$ are 128-equivalent to $(u_s)_{s=1}^\infty$.

    \item $\mathcal{L}(\mathfrak{X}_U) = \mathcal{K}(\mathfrak{X}_U)\oplus [(I_s)_{s=1}^\infty]\oplus \mathbb{C}I$.
    
\end{enumerate}
In particular, $\mathpzc{Cal}(\mathfrak{X}_U)$ is isomorphic, as a Banach algebra, to the unitization of $U$ with coordinate-wise multiplication.
\end{theorem}
Note that because $\mathfrak{X}_U$ has a Schauder basis, it satisfies the approximation property, and thus, $\mathcal{K}(\mathfrak{X}_U)$ is the minimum non-trivial closed ideal in $\mathcal{L}(\mathfrak{X}_U)$.

The simplest examples of Banach spaces $U$ satisfying the assumptions of Theorem \ref{the theorem} are the $\ell_p$ spaces, $1\leq p<\infty$, with their unit vector bases. In particular, the unitization of the infinite-dimensional separable Hilbert space with coordinate-wise multiplication induced by an orthonormal basis is isomorphic as a Banach algebra to the Calkin algebra of some Banach space. The local property of having non-trivial cotype, implies not having a $c_0$ asymptotic version. The $L_p[0,1]$ spaces, for $1<p<\infty$, and more generally, all superreflexive rearrangement invariant spaces on $[0,1]$ have non-trivial cotype and the Haar system as an unconditional basis (see, e.g., \cite[Theorem 2.c.6]{lindenstrauss:tzafriri:1979:partII}). Therefore, by Theorem \ref{the theorem}, their unitizations are Calkin algebras. Due to $L_1$ not having an unconditional basis, Theorem \ref{the theorem} cannot be applied to it, despite it having non-trivial cotype. However, Theorem \ref{the theorem} applies to the Hardy space $H_1$, which has an unconditional basis (see \cite{maurey:1980}). A non-classical Banach space with non-trivial cotype and an unconditional basis is the asymptotic-$\ell_1$ Tsirelson space (\cite{figiel:johnson:1974}). There are also spaces with trivial cotype to which Theorem \ref{the theorem} can be applied, such as $(\oplus_{n=1}^\infty\ell_\infty^n)_p$, for $1\leq p<\infty$, and Schlumprecht space (\cite{schlumprecht:1991}).

Our theorem can also be compared to \cite{argyros:felouzis:2000} and \cite{argyros:raikoftsalis:2012}, where it was shown that every separable reflexive Banach space (e.g., $\ell_2$) is a quotient of a hereditarily indecomposable space, a highly non-classical class of Banach spaces recalled in Section \ref{section roadmap}. Analogously, by our result, many simplistic classical Banach algebras appear as quotient algebras of $\mathcal{L}(X)$ spaces, a class of indirectly defined, and thus more intractable, objects. More generally, although Banach space theory is a mature field with an extensive toolkit of methods for constructing Banach spaces with prescribed properties (see Section \ref{section roadmap}), when it comes to imposing prescribed algebraic and analytic properties on $\mathcal{L}(X)$, there is much left to be desired. Our paper is a contribution to such methods.

\section{Roadmap of ideas and concepts}\label{section roadmap}

The proof of Theorem \ref{the theorem} involves some clear and accessible new ideas, but it also relies on new concepts that are defined with the help of existing heavy machinery. Although this paper is long, it has been written with a focus on clarity of exposition. Initially, we have framed many ideas within the context of some general theory and, at several stages, presented simple examples to which this theory can be applied. Therefore, all the initial principles used in this paper should be accessible to any interested reader with general functional analysis knowledge. We have gathered this material in Part \ref{conceptual part} of the paper. Of course, there is a more technical component based on applying concepts from the theory of hereditarily indecomposable Banach spaces, gathered in Part \ref{HI part}. Prior familiarity with the Gowers-Maurey and Argyros-Haydon constructions (\cite{gowers:maurey:1993} and \cite{argyros:haydon:2011}) makes reading this section more efficient. Although well-established, their techniques have several components in their intuition; an exposition of these is beyond the scope of this paper. However, we have provided detailed proofs of all steps due to substantial required innovations, making this part self-contained. This section provides a roadmap for the ideas and concepts in our paper.

Let us recall some additional standard Banach spaces terminology. A Schauder decomposition  of a Banach space $X$ is a sequence of closed subspaces $\boldsymbol{E} = (E_s)_{s=1}^\infty$ of $X$ such that every $x\in X$ admits a unique representation $x = \sum_{s=1}^\infty x_s$, where, for each $s\in\N$, $x_s\in E_s$. A Schauder decomposition $\boldsymbol{E} = (E_s)_{s=1}^\infty$ of a Banach space $X$ induces a uniformly bounded collection of bounded linear projections $P_J:X\to X$, $J$ is an interval of $\N$, with $P_J(\sum_{s=1}^\infty x_s)=\sum_{s\in J}x_s$. For $x\in X$, we define its support (with respect to the Schauder decomposition $\boldsymbol{E}$) as the set $\mathrm{supp}_{\boldsymbol{E}}(x) = \{s\in\N:P_{\{s\}}x\neq 0\}$ and for $f\in X^*$ we also let $\mathrm{supp}_{\boldsymbol{E}}(f) = \{s\in\mathbb{N}: f\circ P_{\{s\}}\neq 0\}$. We will call a sequence $(y_k)_{k\in I}$, indexed over an interval $I\subset\N$, satisfying $\mathrm{\supp}_{\boldsymbol{E}}(y_{k-1})<\mathrm{\supp}_{\boldsymbol{E}}(y_k)$, $k\in I\setminus \{\min I\}$, an $\boldsymbol{E}$-block sequence. In $X^*$, $\boldsymbol{E}$-block sequences are defined analogously. In some contexts, the natural coordinate system will be a certain Schauder decomposition, and we will skip $\boldsymbol{E}$ in the notation of supports and blocks (specifically in Notation \ref{notation support delta}). If $\boldsymbol{E} = (E_s)_{s=1}^\infty$ is a Schauder decomposition of $X$ such that each $E_s$ is finite-dimensional, then we call it a finite-dimensional decomposition (FDD) of $X$. Similarly, if each $E_s$ is infinite-dimensional, we call it an infinite-dimensional Schauder decomposition. For two Banach spaces $X$, $Y$ and $f\in X^*$, $y\in Y$ we let $f\otimes y:X\to Y$ denote the rank-one bounded linear operator mapping each $x$ to $f(x)y$. Two closed subspaces $Y$ and $Z$ of a Banach space $X$ are said to form a direct sum in $X$ if $Y\cap Z = \{0\}$ and $Y+Z$ is closed, equivalently, the quotient map $Q:X\to X/Z$ restricted on $Y$ is an isomorphism onto its image.

To meaningfully elaborate on the proof of Theorem \ref{the theorem}, we need to briefly revisit some history of construction techniques of Banach spaces, and their relation to the tight control of bounded linear operators, modulo a small ideal. In 1975 (see \cite{lindenstrauss:1975}), Lindenstrauss posed the famous problem of the existence of an infinite dimensional Banach space $X$ with the scalar-plus-compact property, i.e., such that every bounded linear operator $T:X\to X$ is a scalar multiple of the identity plus a compact operator. Note that, for such a space $X$, $\mathpzc{Cal}(X)$ is one-dimensional. For Banach spaces $X$ satisfying the approximation property, $\mathcal{K}(X)$ is the minimum non-trivial closed ideal in $\mathcal{L}(X)$ and thus, in a certain sense, within this class, spaces with the scalar-plus-compact property have the smallest possible algebras of operators. Consequently, a space with the scalar-plus-compact property is said to have very few operators. Although Lindenstrauss's problem was solved affirmatively in 2011 by Argyros and Haydon, this result is a high point of a series of milestones, each addressing an a priori irrelevant problem in the structural theory of Banach spaces. In 1974 (see \cite{tsirelson:1974} or \cite{figiel:johnson:1974}), Tsirelson constructed the first known example of a Banach space not containing an isomorphic copy of $\ell_p$, $1\leq p<\infty$, or $c_0$. Tsirelson space is frequently considered the first non-classical Banach space. Its norm can be described via an implicit formula approximated by a sequence of recursively defined explicit expressions. In 1991 (see \cite{schlumprecht:1991}), Schlumprecht constructed the first example of an arbitrarily distortable Banach space $(X,\|\cdot\|)$, i.e., it has the property that for every $M >0$ there exists an equivalent norm $|\cdot|$ on $X$ such that in every infinite-dimensional subspace $Y$ of $X$ there are vectors $x$, $y$ in $Y$ such that $\|x\| = \|y\|$ and $|x|> M|y|$. The norm of Schlumprecht space is a multi-layered version of the norm of Tsirelson space, and it is the first example in a class of spaces now called mixed-Tsirelson spaces (see Section \ref{section mixed tsirelson}). This class has played a major role in Banach space constructions. Soon after that, ingeniously combining Schlumprecht's technique with an earlier tool of Maurey and Rosenthal from 1977 (see \cite{maurey:rosenthal:1977}), Gowers and Maurey in 1993 (see \cite{gowers:maurey:1993}) solved the unconditional sequence problem by constructing a space with no such sequences. Their space has the additional property that it is hereditarily indecomposable (HI), i.e., no two infinite-dimensional closed subspaces of it form a direct sum. As it was observed by Gowers and Maurey, HI spaces over the complex field have the property that every bounded linear operator on them is a scalar multiple of the identity plus a strictly singular operator. An operator between Banach spaces is called strictly singular if its restriction to any infinite-dimensional subspace of its domain is not an isomorphism onto its image. The collection $\mathcal{SS}(X)$ of strictly singular operators on a Banach space $X$ is an ideal containing $\mathcal{K}(X)$. This ideal came to prominence with the appearance of the Gowers-Maurey space $X_{GM}$, where $\mathcal{L}(X_\mathrm{GM})/\mathcal{SS}(X_\mathrm{GM})$ is one-dimensional. This was the first time an answer to Lindenstrauss's problem appeared within reach. However, it took several more years, until 2011, when the construction of the Argyros-Haydon space $\mathfrak{X}_\mathrm{AH}$ (\cite{argyros:haydon:2011}) led to its solution, i.e., the Calkin algebra $\mathpzc{Cal}(\mathfrak{X}_\mathrm{AH})$ is one-dimensional. The construction of $\mathfrak{X}_\mathrm{AH}$ is based on a combination of the Gowers-Maurey method and an earlier groundbreaking construction technique introduced by Bourgain and Delbaen in 1980 (see \cite{bourgain:delbaen:1980}, page \pageref{argyros-haydon shortcut intro}, and Section \ref{BD section}). It is now well understood that phenomena witnessed modulo the ideal of strictly singular operators on a Gowers-Maurey-type space can frequently also be witnessed modulo the ideal of compact operators on an Argyros-Haydon-type space. In particular, given a unital Banach algebra $\mathcal{B}$, representing it in the form $\mathcal{L}(X)/\mathcal{SS}(X)$, for a Gowers-Maurey-type space $X$, is a first step towards representing $\mathcal{B}$ as the Calkin algebra of some other Argyros-Haydon-type space. A broad viewpoint of this motivates the study of quotient algebras of mixed-Tsirelson and Gowers-Maurey-type spaces and Bourgain-Delbaen and Argyros-Haydon-type spaces. This venture was initiated by Gowers and Maurey in 1997 (see \cite{gowers:maurey:1997}). The quotient algebras of other types of Banach spaces have also been of interest. For example, there is a space $C(K)$ such that $\mathcal{L}(C(K))/\mathcal{SS}(C(K))$ coincides isometrically as a Banach algebra with $C(K)$ (\cite{koszmider:2004} and \cite{plebanek:2004}).

Given a unital Banach algebra $\mathcal{B}$, the problem of finding a Banach space $X$ such that, for $\mathcal{I} = \mathcal{K}(X)$ or $\mathcal{I} = \mathcal{SS}(X)$, $\mathcal{L}(X)/\mathcal{I}$ is isomorphic as a Banach algebra to $\mathcal{B}$ can be broken up into two major components. The first one is the representation of $\mathcal{B}$ in $\mathcal{L}(X)/\mathcal{I}$. This means finding a Banach space $X$ and a subalgebra $\hat{\mathcal{B}}$ of $\mathcal{L}(X)$ such that the closure of its image in the quotient algebra is isomorphic, as a Banach algebra, to $\mathcal{B}$. It is important to note that, when viewed in isolation, this can be achieved almost trivially by taking as $X$, e.g., the unitization of $c_0(\mathcal{B})$. Thus, achieving representation only makes sense in combination with the second major component, the elimination of operators modulo $\mathcal{I}$. This means that $\hat{\mathcal{B}}+\mathcal{I}$ is dense in $\mathcal{L}(X)$, i.e., operators not in the closure of $\hat{\mathcal{B}}+\mathcal{I}$ have been eliminated from $\mathcal{L}(X)$. The latter is achieved by involving Gowers-Maurey or Argyros-Haydon techniques in the construction of $X$. For example, in 1997 (see \cite{gowers:maurey:1997}) Gowers and Maurey constructed a Banach space $X_\mathrm{GM}^{\ell_1}$ such that $\mathcal{L}(X_\mathrm{GM}^{\ell_1})/\mathcal{SS}(X_\mathrm{GM}^{\ell_1})$ is isomorphic as a Banach algebra to the Wiener algebra, i.e., the convolution algebra $\ell_1(\mathbb{Z})$. Their proof leverages on the fact that on the Banach space $\ell_1(\mathbb{Z})$, the closed subalgebra of $\mathcal{L}(\ell_1(\mathbb{Z}))$ generated by the left and right shift operators is the Wiener algebra. On the space $X_\mathrm{GM}^{\ell_1}$, the left and right shift are also bounded and, due to the inherent $\ell_1$-structure of mixed-Tsirelson spaces (specifically Schlumprecht space), they generate the Wiener algebra. The Gowers-Maurey nature of $X_\mathrm{GM}^{\ell_1}$ guarantees the elimination of operators modulo $\mathcal{SS}(X_\mathrm{GM}^{\ell_1})$. Inspired by this proof, in 2013 (see \cite{tarbard:2013}), Tarbard defined an Argyros-Haydon-type space $\mathfrak{X}_\mathrm{Ta}$ with Calkin algebra isomorphic as a Banach algebra to the semigroup algebra of $\N_0$, i.e., the convolution algebra $\ell_1(\N_0)$. On $\mathfrak{X}_\mathrm{Ta}$, a type of right shift operator is bounded, and the unital subalgebra it generates is the semigroup algebra of $\N_0$. The Argyros-Haydon nature of $\mathfrak{X}_\mathrm{Ta}$ provides the elimination of operators modulo $\mathcal{K}(\mathfrak{X}_\mathrm{Ta})$. A recent example of this process was given by the first author in \cite{motakis:2024} where, for every compact metric space $K$, an Argyros-Haydon-type space $\mathfrak{X}_{C(K)}$ is constructed with Calkin algebra isometrically isomorphic as a Banach algebra to $C(K)$. The representation of $C(K)$ in $\mathpzc{Cal}(\mathfrak{X}_{C(K)})$ is based on the following principle. Let $K$ be a compact metric space and $(\kappa_s)_{s=1}^\infty$ be a dense sequence in $K$. Let $U$ be a Banach space with a 1-unconditional basis $(u_s)_{s=1}^\infty$ (e.g., $U = c_0$ or $U = \ell_2$). For $\phi\in C(K)$, let $\hat\phi:U\to U$ be the diagonal operator such that, for $s\in\mathbb{N}$, $\hat\phi (u_s) = \phi(\kappa_s)u_s$. Then, $\widehat C(K) = \{\hat\phi:\phi\in C(K)\}$ is a subalgebra of $\mathcal{L}(U)$ that is isometrically isomorphic as a Banach algebra to $C(K)$. The space $\mathfrak{X}_{C(K)}$ has a (non-unconditional) Schauder basis $(d_\gamma)_{\gamma\in\Gamma}$, where $\Gamma$ is some countable set. For a prechosen dense sequence $(\kappa_\gamma)_{\gamma\in\Gamma}$ in $K$ all entries of which are infinitely repeated, the corresponding representation $\phi\mapsto\hat\phi$ is a well defined, but unbounded, injective linear homomorphism with domain all Lipschitz functions $\phi:K\to\mathbb{C}$ and codomain $\mathcal{L}(\mathfrak{X}_{C(K)})$. By design, when modding by the compact operators, $\widehat{[\,\cdot\,]}:C(K)\to\mathpzc{Cal}(\mathfrak{X}_{C(K)})$ extends to an isometric embedding, in the Banach algebra sense. The Argyros-Haydon nature of $\mathfrak{X}_{C(K)}$ ensures the elimination of operators modulo $\mathcal{K}(\mathfrak{X}_{C(K)})$, i.e., this embedding is onto.

Let us fix a Banach space $U$ with a normalized 1-unconditional basis $(u_s)_{s=1}^\infty$ and analyze the meaning of its presence in $\mathpzc{Cal}(X)$. The lack of a $c_0$ asymptotic version is not yet relevant, so we don't assume it until later. Looking at $U$ as a Banach algebra with coordinate-wise multiplication, $(u_s)_{s=1}^\infty$ is a sequence of mutually annihilating idempotents. Therefore, to detect it in an $\mathpzc{Cal}(X)$ space, we are led to consider sequences of mutually annihilating infinite-rank projections in an $\mathcal{L}(X)$ space. Let us observe that the simplest such collection, a sequence of orthogonal projections on a Hilbert space with pairwise orthogonal infinite dimensional ranges, is equivalent to the unit vector basis of $c_0$, and so is its image in the Calkin algebra. Something similar is true for a sequence of projections $(P_{\{s\}})_{s=1}^\infty$ associated with a Schauder decomposition of a Banach space $X$, which satisfies $\sup_{N\in\N}\|\sum_{s=1}^NP_{\{s\}}\| <\infty$. It then easily follows that $\mathpzc{Cal}(X)$ is non-reflexive, and in particular, if the image of $(P_{\{s\}})_{s=1}^\infty$ in the Calkin algebra is unconditional, it must be equivalent to the unit vector basis of $c_0$. In the second author's recent and highly relevant paper \cite{pelczar-barwacz:2023}, a Gowers-Maurey-type space $X$ was constructed with $\mathcal{L}(X)/\mathcal{SS}(X)$ isomorphic as a Banach algebra to the unitization of the dual of a mixed-Tsirelson space with coordinate-wise multiplication. This mixed-Tsirelson space is a reflexive space with an unconditional basis. The representation is witnessed by a bounded sequence of mutually annihilating projections on $X$ that are not associated to an infinite dimensional Schauder decomposition. The operator structure of this $X$ shares further conceptual similarities via duality to the operator structure of our space $\mathfrak{X}_U$. We briefly elaborate on this relation in Section \ref{section generic sequences of projections} starting in the paragraph before Proposition \ref{dual condition on projections}.

We are now ready to discuss the principle ideas in the representation component of Theorem \ref{the theorem}. For the purpose of exposition and to elucidate our ideas, let us first define a sequence of rank-one projections, which, of course, vanish in the Calkin algebra. At a most basic level, it is based on the following crucial but straightforward identification. Let $(e_s)_{s=1}^\infty$ denote the unit vector basis of $c_0$ and $\mathcal{K}_\mathrm{diag}(c_0,U)$ the subspace of $\mathcal{L}(c_0,U)$ of compact operators $D:c_0\to U$ such that, for $s\in\mathbb{N}$, $De_s$ is a scalar multiple of $u_s$. Then, $\mathcal{K}_\mathrm{diag}(c_0,U) = [(e_s^*\otimes u_s)_{s=1}^\infty]$ and it coincides isometrically with $U$ via the identification $e_s^*\otimes u_s\mapsto u_s$. The caveat is that strictly speaking, $\mathcal{K}_\mathrm{diag}(c_0,U)$ is not a Banach algebra, let alone a subalgebra of an $\mathcal{L}(X)$ space. Instead, let $X = (c_0\oplus U)_\infty$ and, for $s\in\N$, define $x_s = (e_s,u_s)\in X$ and $f_s = (e_s^*,0)\in X^*$. Then, $I_s = f_s\otimes x_s$, $s\in\N$, is a sequence of mutually annihilating idempotents in $\mathcal{L}(X)$ that is isometrically equivalent to $(u_s)_{s=1}^\infty$. In particular, $[(I_s)_{s=1}^\infty]$ is a subalgebra of $\mathcal{L}(X)$ that coincides isometrically as a Banach algebra with $U$. This is proved easily using the identification $\mathcal{K}_\mathrm{diag}(c_0,U) \equiv U$. We put this idea into a general proposition.
\begin{proposition}
Let $U$ be a Banach space with a normalized 1-unconditional basis $(u_s)_{s=1}^\infty$ and $X$ be a Banach space. Let $(x_s)_{s=1}^\infty$ be a sequence in $X$ and $(f_s)_{s=1}^\infty$ be a bounded sequence in $X^*$ such that, for $s,t\in\mathbb{N}$, $f_s(x_t) = \delta_{s,t}$. Put $X_0 = \cap_{s=1}^\infty\mathrm{ker}(f_s)$ and denote $Q:X\to X/X_0$ the quotient map. Assume that
\begin{enumerate}[label=(\alph*)]
    
    \item $(x_s)_{s=1}^\infty$ is equivalent to $(u_s)_{s=1}^\infty$ and

    \item $(Qx_s)_{s=1}^\infty$ is equivalent to the unit vector basis of $c_0$.
\end{enumerate}
 Then, $I_s = f_s\otimes u_s$, $s\in\N$, is a sequence of mutually annihilating rank-one projections in $\mathcal{L}(X)$ that is equivalent to $(u_s)_{s=1}^\infty$. In particular, $[(I_s)_{s=1}^\infty]$ is isomorphic as a Banach algebra to $U$ with coordinate-wise multiplication.
\end{proposition}

Obviously, the above cannot be used directly to study the Calkin algebra of a Banach space because finite rank operators vanish in it. However, it naturally yields a criterion that applies to infinite rank projections. 

\begin{proposition}
\label{toy criterion}
Let $U$ be a Banach space with a normalized 1-unconditional basis $(u_s)_{s=1}^\infty$. Let $X$ be a Banach space and $(I_s)_{s=1}^\infty$ be a bounded sequence of mutually annihilating infinite rank projections in $\mathcal{L}(X)$. For $s\in\N$, put $X_s = I_s(X)$, $X_0 = \cap_{s=1}^\infty\mathrm{ker}(I_s)$, and let $Q:X\to X/X_0$ denote the quotient map. Assume that for every normalized sequence $(x_s)_{s=1}^\infty$, such that $x_s\in X_s$, for all $s\in\N$, the following hold.
\begin{enumerate}[label=(\alph*)]
    
    \item\label{toy criterion a} $(x_s)_{s=1}^\infty$ is equivalent to $(u_s)_{s=1}^\infty$ and

    \item\label{toy criterion b} $(Qx_s)_{s=1}^\infty$ is equivalent to the unit vector basis of $c_0$.
\end{enumerate}
 Then, $(I_s)_{s=1}^\infty$ in $\mathcal{L}(X)$ and its image in $\mathpzc{Cal}(X)$ are equivalent to $(u_s)_{s=1}^\infty$. In particular, $[(I_s)_{s=1}^\infty]$ and $\mathcal{K}(X)$ form a direct sum in $\mathcal{L}(X)$ and $U$ with coordinate-wise multiplication embeds into $\mathpzc{Cal}(X)$ as a Banach algebra.
\end{proposition}

Proposition \ref{toy criterion} is not applicable to the Argyros-Haydon-type space $\mathfrak{X}_U$; assumption \ref{toy criterion a} is too restrictive, and we need to relax it. However, assumption \ref{toy criterion b} will play a central role and it lends some subtle properties to the associated objects. In particular, the roles of $X_0$ and $X_s$, $s\in\N$, are completely different. By \ref{toy criterion b}, $X/X_0$ is isomorphic to $(\oplus_{s=1}^\infty X_s)_{c_0}$ via $x+X_0\mapsto (I_sx)_{s=1}^\infty$. We highlight that unless $U$ is isomorphic to $c_0$, $(X_s)_{s=0}^\infty$ does not form a Schauder decomposition of $X$, i.e.,
there is no bounded projection onto $X_0$ that contains all $X_s$, $s\geq 1$, in its kernel. Assumption \ref{toy criterion b} also yields that if $X$ has the bounded approximation property, the image of $(I_s)_{s=1}^\infty$ in the Calkin algebra is automatically unconditional (Proposition \ref{always unconditional}).

We showcase the above proposition by applying it to a toy example. Let $X$ denote the completion of the space of eventually null sequences
\[
(x_0,x_1,x_2,\ldots)\in c_{00}(\N^2)\times c_{00}(\N)\times c_{00}(\N)\times\cdots
\]
with
\begin{equation}
\label{toy norm}
\big\|(x_0,x_1,x_2,\ldots)\big\| = \max\Big\{\max_{s\in\N}\|x_s\|_\infty,\max_{(i_s)_{s=1}^\infty\in\N^\N}\Big\|\sum_{s=1}^\infty\big(x_s(i_s) - x_0(s,i_s)\big)u_s\Big\|_U\Big\}.
\end{equation}
For $s\in\N$ (but not for $s=0$), denote $I_s:X\to X$ the norm-one linear projection given by
\[I_s(x_0,x_1,\ldots,x_{s-1},x_s,\ldots) = (0,0,\ldots,0,x_s,\ldots).\]
 Then, it is not very hard to show that $(I_s)_{s=1}^\infty$ satisfies the assumptions of Proposition \ref{toy criterion}. In fact, the sequence $(I_s)_{s=1}^\infty$ in $\mathcal{L}(X)$ and its image in $\mathpzc{Cal}(X)$ are isometrically equivalent to $(u_s)_{s=1}^\infty$.

In Section \ref{section generic sequences of projections}, we state Proposition \ref{general condition on projections}, which is a criterion for the representation of $U$ with coordinate-wise multiplication in the Calkin algebra applicable to the Argyros-Haydon-type space $\mathfrak{X}_U$. This criterion is nearly identical to Proposition \ref{toy criterion} in its assumptions and conclusions; the only difference is that assumption \ref{toy criterion a} is replaced with a milder, but slightly less elegant, set of conditions for normalized sequences $(x_s)_{s=1}^\infty$ such that $x_s\in X_s$, for $s\in\N$. All such sequences must satisfy upper $U$-bounds, and sequences that satisfy uniform lower $U$-bounds need to exist asymptotically in a certain sense. Notably, the crucial condition \ref{toy criterion b} remains unchanged. For the purpose of conceptualization, we use this criterion at the end of Section \ref{section generic sequences of projections} to construct, for every sequence of infinite-dimensional separable Banach spaces $(X_s)_{s=1}^\infty$, a Banach space $X$ and a sequence of mutually annihilating projections $(I_s)_{s=1}^\infty$ satisfying the conclusion of Proposition \ref{general condition on projections} and, additionally, for every $s\in\N$, $X_s$ is isomorphic to $I_s(X)$. The norm of this space satisfies a simple formula very similar to \eqref{toy norm}.

The next layer of ideas behind the construction of $\mathfrak{X}_U$ concerns a new type of Bourgain-Delbaen-$\mathscr{L}_\infty$-space  called a bonding. This construction induces a bounded sequence of mutually annihilating projections that automatically satisfy Proposition \ref{toy criterion} \ref{toy criterion b}. In Section \ref{BD section}, we recall in ample detail all the required parts from the theory of Bourgain-Delbaen-$\mathscr{L}_\infty$-spaces. Just for the sake of the introduction, we take a shortcut in defining them. This definition is equivalent to the one presented later but unsuitable as a formal starting point for our approach.\label{argyros-haydon shortcut intro} Let $\Gamma_0\subset\Gamma_1\subset\cdots$ be finite sets and $\Gamma = \cup_{n=0}^\infty\Gamma_n$. For $\gamma\in\Gamma$, we denote $\mathrm{rank}(\gamma)$ the minimum $n\in\N_0$ such that $\gamma\in\Gamma_n$. Consider the standard duality pair $\ell_1(\Gamma)\times\ell_\infty(\Gamma)$ and let $(e_\gamma^*)_{\gamma\in\Gamma}$ denote the unit vector basis of $\ell_1(\Gamma)$. For a collection $(d_\gamma)_{\gamma\in\Gamma}$ in $\ell_\infty(\Gamma)$, the subspace $\mathfrak{X} = [(d_\gamma)_{\gamma\in\Gamma}]$ of $\ell_\infty(\Gamma)$ is called a Bourgain-Delbaen-$C$-$\mathscr{L}_\infty$-space if for all $\gamma,\eta\in\Gamma$ with $\mathrm{rank}(\gamma)\leq\rank(\eta)$, $e_\gamma^*(d_\eta) = \delta_{\gamma,\eta}$ and, for all $n\in\N_0$, the canonical mapping from $\langle\{e_\gamma^*:\gamma\in\Gamma_n\}\rangle$ into the dual of $\langle\{d_\gamma:\gamma\in\Gamma_n\}\rangle$ is a $C$-isomorphism. Then there is a unique collection $(d_\gamma^*)_{\gamma\in\Gamma}$ in $\ell_1(\Gamma)$ such that $(d_\gamma^*,d_\gamma)_{\gamma\in\Gamma}$ is a biorthogonal system. Furthermore, $(d_\gamma)_{\gamma\in\Gamma}$, with any enumeration compatible with the rank, forms a (typically non-unconditional) Schauder basis for $\mathfrak{X}$ and, for all $n\in\N_0$, $\langle\{ d_\gamma^*:\gamma\in\Gamma_n\}\rangle = \langle\{ e_\gamma^*:\gamma\in\Gamma_n\}\rangle$. However, the by far more useful coordinate system is the FDD $(Z_n)_{n=0}^\infty$, where $Z_n = \langle\{d_\gamma:\rank(\gamma) = n\}\rangle$. For every interval $E$ of $\mathbb{N}_0$, we denote by $P_E$ the associated projection onto $\oplus_{n\in E}Z_n$. Starting with the Argyros-Haydon construction, most subsequent Bourgain-Delbaen-$\mathscr{L}_\infty$-spaces satisfy an additional property, namely that for every $\gamma\in\Gamma$ either $e_\gamma^* = d_\gamma^*$ or $e_\gamma^*$ admits an evaluation analysis. The shape of such an analysis may vary depending on the approach, but for the purpose of this paper, this means that
\begin{equation}
\label{intro evaluation analysis}
e_\gamma^* = \sum_{r=1}^a d_{\xi_r}^* + \sum_{r=1}^a\lambda_r e_{\eta_r}^*\circ P_{E_r},
\end{equation}
where $\xi_1,\ldots,\xi_a$, $\eta_1,\ldots,\eta_a$ are in $\Gamma$, $E_1,\ldots,E_a$ are intervals of $\N_0$, and $\lambda_1,\ldots,\lambda_a$ are scalars bounded in absolute value by some uniform $0<\vartheta\leq 1/4$ such that $E_1<\mathrm{rank}(\xi_1)<E_2<\mathrm{rank}(\xi_2)<\cdots$ and $\gamma = \xi_a$. The precise form of this analysis plays a central role in estimating the norm of vectors in $\mathfrak{X}$.

A bonding is a Bourgain-Delbaen-$\mathscr{L}_\infty$-space space with additional restrictions on the form of the evaluation analyses of the functionals $e_\gamma^*$. Section \ref{bonding section} is fully devoted to the definition and study of bondings. For a partition of $\Gamma = \cup_{s=0}^\infty\Gamma^s$ denote, for $s\in\N_0$, $\mathfrak{X}_s = [(d_\gamma)_{\gamma\in\Gamma^s}]$. Also, for $f\in\ell_1(\Gamma)$, denote $\mathrm{supp}_\mathrm{h}(f) = \{s\in\mathbb{N}: f|_{\mathfrak{X}_s}\neq 0\}$. Intentionally, $s=0$ is excluded from the definition of $\mathrm{supp}_\mathrm{h}(f)$ because the spaces $\mathfrak{X}_0$ and $\mathfrak{X}_s$, $s\in\N$, play a completely different role. The space $\mathfrak{X}$ is called a bonding of $(\mathfrak{X}_s)_{s=1}^\infty$ if the following hols. Let $\gamma\in\Gamma$ that admits an evaluation analysis as in \eqref{intro evaluation analysis}.
\begin{enumerate}[label=(\alph*),leftmargin=19pt]

\item\label{bonding intro condition a} If, for some $s\in\N$, $\gamma\in\Gamma^s$ then $\xi_1,\ldots,\xi_a$, $\eta_1,\ldots,\eta_a$ are also in $\Gamma^s$.

\item\label{bonding intro condition b} If $\gamma\in\Gamma^0$, then $\xi_1,\ldots,\xi_a$ are also in $\Gamma^0$ and $\mathrm{supp}_\mathrm{h}(e_{\eta_1}^*\circ P_{E_1}) < \mathrm{supp}_\mathrm{h}(e_{\eta_2}^*\circ P_{E_2})<\cdots$.

\end{enumerate}
In the language of \cite{argyros:motakis:2019}, assumption \ref{bonding intro condition a} means that each $\Gamma^s$ is a self-determined subset of $\Gamma$. There are also Gowers-Maurey analogues of ``\ref{bonding intro condition a} and \ref{bonding intro condition b}'' in \cite{pelczar-barwacz:2023}, but in the Bourgain-Delbaen-$\mathscr{L}_\infty$-setting they are more elementary and have a direct impact: a bonding induces a sequence of mutually annihilating norm-one projections that satisfy Proposition \ref{toy criterion} \ref{toy criterion b}. More specifically, we prove the following in Theorem \ref{bonding general statement}.

\begin{theorem}
Let $\mathfrak{X}$ be a bonding of $(\mathfrak{X}_s)_{s=1}^\infty$ and let $I_s:\langle\{d_\gamma:\gamma\in\Gamma\}\rangle\to \langle\{d_\gamma:\gamma\in\Gamma^s\}\rangle$ denote the coordinate projection, $s\in\N$. Then the $I_s$, $s\in\N$, extend to mutually annihilating norm-one projections on $\mathfrak{X}$ such that $I_s(\mathfrak{X}) = \mathfrak{X}_s$, $s\in\N$, and $\cap_{s=1}^\infty\mathrm{ker}(I_s) = \mathfrak{X}_0$. Furthermore, denoting $Q:\mathfrak{X}\to\mathfrak{X}/\mathfrak{X}_0$ the quotient map, for every normalized sequence $(x_s)_{s=1}^\infty$, such that $x_s\in\mathfrak{X}_s$, $s\in\mathbb{N}$, $(Qx_s)_{s=1}^\infty$ is 2-equivalent to the unit vector basis of $c_0$. In particular, $\mathfrak{X}/\mathfrak{X}_0$ is 2-isomorphic to $(\oplus_{s=1}^\infty\mathfrak{X}_s)_{c_0}$ via $x+\mathfrak{X}_0\mapsto (I_sx)_{s=1}^\infty$ and the image of $(I_s)_{s=1}^\infty$ in $\mathpzc{Cal}(\mathfrak X)$ is unconditional.
\end{theorem}

It is worth commenting that in Bourgain-Delbaen-$\mathscr{L}_\infty$-spaces, the boundedness of non-trivial coordinate projections onto subsets of $\Gamma$ is unusual. For example, constructing any type of non-trivial operator on the original Bourgain-Delbaen-$\mathscr{L}_\infty$-space from \cite{bourgain:delbaen:1980} requires considerable effort (see, e.g., \cite{beanland:mitchell:2010}, \cite{beanland:mitchell:2014}, or \cite[Section 2.5.1]{tarbard:2013}); known non-trivial projections on this space are not coordinate projections.

We wish to represent a specific $U$ with coordinate-wise multiplication in $\mathpzc{Cal}(\mathfrak{X})$. For this, the additional satisfaction of the universal upper $U$ bounds and partial lower $U$ bounds of the earlier discussed criterion (Proposition \ref{general condition on projections}) is necessary. We encode these by imposing $U$-related restrictions on the coefficients $(\lambda_r)_{r=1}^a$ appearing in \eqref{intro evaluation analysis}, for $\gamma\in\Gamma^0$ (but not for $\gamma\in\Gamma^s$, $s\in\N$). These are formally introduced in Definition \ref{upper lower U condition} and called the upper $U$ condition and the partial lower $U$ condition. If both of them are satisfied, we prove in Theorem \ref{bd condition on projections} that the sequence $(I_s)_{s=1}^\infty$ in $\mathcal{L}(\mathfrak{X})$ and its image in $\mathpzc{Cal}(\mathfrak{X})$ are equivalent to $(u_s)_{s=1}^\infty$. Importantly, for $s\in \mathbb{N}$, the only restriction imposed on the evaluation analysis of $e_\gamma^*$, for $\gamma\in\Gamma^s$, is that all its components lie in $\Gamma^s$. This allows, within the boundaries of Bourgain-Delbaen-$\mathscr{L}_\infty$-spaces, almost unlimited freedom in designing each space $\mathfrak{X}_s$. The importance of this is discussed in the next paragraph. At the end of Section \ref{bonding section}, we apply our new tool in a comparatively accessible setting not yet involving Argyros-Haydon techniques; we define an elementary bonding $\mathfrak{B}_U$ to which Theorem \ref{bd condition on projections} is applicable, and therefore, the unitization of $U$ with coordinate-wise multiplication embeds into $\mathpzc{Cal}(\mathfrak{B}_U)$.

The final layer of ideas in the definition of $\mathfrak{X}_U$ concerns the elimination of operators modulo $\mathcal{K}(\mathfrak{X}_U)$ and thus the involvement of Argyros-Haydon techniques. For an $x$ in a bonding $\mathfrak{X}$ we denote $\mathrm{supp}_\mathrm{h}(x) = \{s\in\N: I_sx\neq 0\}$; as always, $s=0$ is not considered. A sequence $y_1,y_2,\ldots$ in $\mathfrak{X}$ is called horizontally block if $\mathrm{supp}_\mathrm{h}(y_1)<\mathrm{supp}_\mathrm{h}(y_2)<\cdots$ and it is block with respect to the standard FDD of $\mathfrak{X}$. The norm structure of linear combinations of horizontally block sequences, i.e., the horizontal structure of $\mathfrak{X}$, is dictated purely by $e_\gamma^*$ for $\gamma\in\Gamma^0$. To maintain upper U bounds for such sequences in $\mathfrak{X}_U$, we are forced to define the set $\Gamma^0$ based on a $q$-convex variant of the Argyros-Haydon construction technique for a value $1\leq q<\infty$ associated with $U$, provided by the fact that $(u_s)_{s=1}^\infty$ has no $c_0$ asymptotic version (see Section \ref{U section}). The partial lower $U$ condition is ensured by including low-complexity $\gamma\in\Gamma^0$ with coefficients $(\lambda_r)_{r=1}^a$ from a collection that describes the norm of $U$. By the Argyros-Haydon nature of $\mathfrak{X}_U$, every operator $T:\mathfrak{X}_U\to\mathfrak{X}_U$ is a scalar multiple of the identity plus a horizontally compact operator. A bounded linear operator $S:\mathfrak{X}_U\to Y$, where $Y$ is some Banach space, is called horizontally compact if for every bounded horizontally block sequence $(y_n)_{n=1}^\infty$ in $\mathfrak{X}_U$, $\lim_n\|Sy_n\| = 0$. The scalar-plus-horizontally-compact property is only the first ingredient in achieving the elimination of operators modulo $\mathcal{K}(\mathfrak{X}_U)$. There are four in total. The second ingredient is that every $\mathfrak{X}_s$, $s\in\N$, satisfies the scalar-plus-compact property. This is achieved by using Argyros-Haydon-techniques in the definition of $\Gamma^s$. The third ingredient is that, for $s\neq t$ in $\N$, the spaces $\mathfrak{X}_s$, $\mathfrak{X}_t$ are totally incomparable, i.e., $\mathcal{L}(\mathfrak{X}_s,\mathfrak{X}_t) = \mathcal{K}(\mathfrak{X}_s,\mathfrak{X}_t)$. The tools for achieving this are found in the Argyros-Haydon paper \cite{argyros:haydon:2011}. The final ingredient is that the horizontal structure of $\mathfrak{X}_U$ is totally incomparable with each $\mathfrak{X}_s$, $s\in\N$. This means that every bounded linear operator $T:\mathfrak{X}_s\to \mathfrak{X}_U$ is a scalar multiple of the inclusion plus a compact operator and that every bounded linear operator $T:\mathfrak{X}_U\to \mathfrak{X}_s$ is horizontally compact. This is made possible by basing each $\mathfrak{X}_s$, $s\in\N$, on a $q$-convex variant of the Argyros-Haydon technique for a value of $q$ matching the one in the definition of $\Gamma^0$. Before proving them, the aforementioned ingredients are used in Section \ref{section calkin algebra XU} to deduce Theorem \ref{the theorem}. All later sections are based on applying HI techniques to prove the four ingredients.

Compared to the standard Argyros-Haydon construction, the $q$-convexity component in the definition of $\mathfrak{X}_U$ substantially impacts the technical aspect of the elimination of operators. Let $p$ denote the conjugate exponent of $q$. We model the definition of $\mathfrak{X}_U$ on a $q$-mixed Tsirelson space $T_q[(\mathcal{A}_{n_j},m_j^{-1})_{j\in\N}]$, where $(m_j)_{j=1}^\infty$ and $(n_j)_{j=1}^\infty$ are sequences of natural numbers that increase sufficiently rapidly.
The norm of this space is defined via a norming set $W_q[(\mathcal{A}_{n_j},m_j^{-1})_{j\in\N}]$. This is the largest subset of $c_{00}(\N)$ such that every $f$ in it is either of the form $\pm e_i^*$, or there are $j\in\N$ and a block sequence $f_1$,\ldots,$f_a$ in $W_q[(\mathcal{A}_{n_j},m_j^{-1})_{j\in\N}]$ such that $a\leq n_j$ and
\[f = \frac{1}{m_j\cdot n_j^{1/p}}\sum_{r=1}^a f_r.\]
For such $f$, we write $\mathrm{weight}(f) = m_j^{-1}$. The space $T_q[(\mathcal{A}_{n_j},m_j^{-1})_{j\in\N}]$ is the completion of $c_{00}(\N)$ with the norm given by $\|x\| = \sup\{|\langle f,x\rangle|:f\in W_q[(\mathcal{A}_{n_j},m_j^{-1})_{j\in\N}]\}$. For more details, see Section \ref{section mixed tsirelson}. This space displays heterogeneity in the behaviour of the norm of linear combinations of normalized block vectors, and this is important in the tight control of bounded linear operators. Certain required estimates in $q$-mixed Tsirelson spaces are proved in Section \ref{section mixed-tsirelson estimates}.

The space $\mathfrak{X}_U$ is modelled on the above as follows. The $\gamma\in\Gamma$ that admit an evaluation analysis can be categorized into ground nodes, simple nodes, and h-nodes. A ground node is a $\gamma\in\Gamma^0$ such that $e_\gamma^*$ is a low-complexity functional providing partial lower $U$ bounds. More precisely,
\[e_\gamma^* = \sum_{r=1}^ad_{\xi_r}^* + \sum_{r=1}^a \lambda_r e_{\eta_r}^*\circ P_{\{\mathrm{rank}(\eta_r)\}}\]
where $\xi_1,\ldots,\xi_r$ are in $\Gamma^0$, $\eta_r\in\Gamma^r$, $r=1,\ldots,a$, and $(\lambda_r)_{r=1}^a$ are from the unit ball of $U^*$ and $a\leq \mathrm{rank}(\eta_1)$. Both simple nodes and h-nodes have evaluation analyses of the form
\[e_\gamma^* = \sum_{r=1}^ad_{\xi_r}^* + \sum_{r=1}^a\frac{1}{m_j\cdot n_j^{1/p}}e_{\eta_r}^*\circ P_{E_r},\]
where $j\in\N$ and $a\leq n_j$, and for such $\gamma$ we write $\mathrm{weight}(\gamma) = m_j^{-1}$. The difference is that simple nodes are in $\Gamma^s$, $s\in\N$, and thus satisfy \ref{bonding intro condition a} of the bonding definition, whereas h-nodes are in $\Gamma^0$ and thus satisfy \ref{bonding intro condition b} of the bonding definition (see page \pageref{bonding intro condition a}). As it is usually the case in this type of construction, whenever $j$ is odd, there are further restrictions on the evaluation analysis imposed by a Maurey-Rosenthal coding function (see \cite{maurey:rosenthal:1977}) on the weights of the $\eta_r$, $1\leq r\leq a$.

A major difference to the standard Argyros-Haydon technique, which uses $q=1$ (and thus, $p=\infty$), lies in the ``residual'' part $\sum_{r=1}^ad_{\xi_r}^*$ of the evaluation analysis of an $e_\gamma^*$ which is incongruent to $q$-convexity. This provides an obstacle in proving the basic inequality, a standard step in HI and related constructions. Let us briefly elaborate on the extent of this. Rapidly increasing sequences (RISs) are a tool that is omnipresent in HI and related constructions. It was introduced no later than 1991 in \cite{schlumprecht:1991} and has since seen limited conceptual modifications. In the Argyros-Haydon space, a bounded linear operator $T$ with domain $\mathfrak{X}_\mathrm{AH}$ is compact if and only if for every RIS $(x_n)_{n=1}^\infty$ in $\mathfrak{X}_\mathrm{AH}$, $\lim_n\|Tx_n\| = 0$. Here, there are two types of RISs: h-RISs concerning the horizontal structure of $\mathfrak{X}_U$ and, for $s\in\N$, $\mathfrak{X}_s$-RISs concerning the internal structure of $\mathfrak{X}_s$. In analogy to the Argyros-Haydon space, a bounded linear operator $T$ with domain $\mathfrak{X}_U$ is horizontally compact if and only if for every h-RIS $(x_n)_{n=1}^\infty$ in $\mathfrak{X}_U$, $\lim_n\|Tx_n\| = 0$. For $s\in\N$, $\mathfrak{X}_s$-RISs in $\mathfrak{X}_s$ similarly characterize compactness of operators with domain $\mathfrak{X}_s$. The role of the basic inequality is to provide high-precision bounds for norm estimates of linear combinations of h-RISs and $\mathfrak{X}_s$-RISs in terms of bounds on linear combinations of basis vectors on some $q$-mixed Tsirelson space. This tool dates at least as far back as 1997 in \cite{argyros:deliyanni:1997}. Due to the standard definition of RISs being an obstacle in proving this inequality, we use a different approach based on the concept of local support. This modification yields a different proof of the basic inequality via an intermediate step called an evaluation refinement. RISs and the basic inequality are treated in Section \ref{RIS section}. A drawback is that these new RISs are unstable. For example, the usual weighted averages of RISs are no longer RISs and thus in Section \ref{dependent sequences} we need to prove a second basic inequality for them. Finally, RISs are ill-suited for studying the structure of arbitrary subspaces of the space $\mathfrak{X}_U$. Fortunately, this does not negatively impact the study of operators with domain the entire space $\mathfrak{X}_U$.

\section{Open problems and directions}

The most obvious question arising from our paper is whether the unitization of all Banach spaces $U$ with an unconditional basis with coordinate-wise multiplication are Calkin algebras. The first such known example was $U = c_0$ from \cite{motakis:puglisi:zisimopoulou:2016}. In \cite[Problem 5]{motakis:puglisi:tolias:2020}, it was asked if any other such spaces existed, which is indeed true, by Theorem \ref{the theorem}. Interestingly, all spaces covered by this theorem are remote from $c_0$ in the asymptotic sense. These two results leave a large gap between them containing such spaces with an unconditional basis as the asymptotic-$c_0$ Tsirelson space from \cite{tsirelson:1974} and the Schreier space from \cite{schreier:1930}. The former is reflexive and asymptotic-$c_0$, and the latter is saturated with subspaces isomorphic to $c_0$.

\begin{problem}
Is the unitization of every Banach space with a normalized unconditional basis with coordinate-wise multiplication a Calkin algebra?
\end{problem}

It is possible that a combination of methods developed in this paper with a transfinite-complexity version of \cite{argyros:motakis:2019} could relax the assumptions on $U$ in Theorem \ref{the theorem} to include spaces with an unconditional basis not containing $c_0$. It is worth recalling that,  by \cite{motakis:puglisi:tolias:2020}, all non-reflexive Banach spaces with an unconditional basis with a non-standard commutative multiplication are Calkin algebras.

Although this is not obvious, our methods show promise in constructing specific types of explicit non-separable Calkin algebras. A Schauder basis $(u_s)_{s=1}^\infty$ of a Banach space $U$ is called 1-subsymmetric if it is 1-unconditional and isometrically equivalent to its subsequences. For a totally ordered set $I$, we let $U(I)$ be the completion of $c_{00}(I)$ with the norm given as follows: for $i_1<\cdots<i_n\in I$ and $a_1,\ldots,a_n\in\mathbb{C}$, if $x = \sum_{s=1}^na_se_{i_s}$ then
\[\|x\|_{U(I)} = \Big\|\sum_{s=1}^na_su_s\Big\|_U.\]
Let $\Delta = \{0,1\}^\N$ with the lexicographical order.

\begin{problem}
Let $U$ be a Banach space with a 1-subsymmetric basis. Is the unitization of $U(\Delta)$ with coordinate-wise multiplication a Calkin algebra? In particular, is the unitization of $\ell_2(\Delta)$ a Calkin algebra?
\end{problem}

A positive answer would yield the first known explicit non-separable Calkin algebra. The relevance of our approach in attacking this problem is that it seems possible to define a type of Bourgain-Delbaen-$\mathscr{L}_\infty$-space called a $\Delta$-bonding of spaces $(\mathfrak{X}_\sigma)_{\sigma\in\Delta}$ such that for $\sigma\neq\sigma'$, $\mathfrak{X}_\sigma\cap\mathfrak{X}_{\sigma'}$ is finite-dimensional.

Let us reiterate a related problem from \cite{motakis:2024}, where it was shown that every separable $C(K)$ space is a Calkin algebra.

\begin{problem}
    Does there exist a non-separable $C(K)$ space that is a Calkin algebra? In particular, are $C(\beta\N)$ and $C(\beta\N\setminus\N)$ Calkin algebras?
\end{problem}

The last part of this question is related to the existence of a Banach space with an unconditional basis on which every bounded linear operator is a diagonal operator plus a compact operator. Among spaces with an unconditional basis, such a space would have the smallest possible algebra of operators modulo the minimum non-trivial closed ideal. Specifically, its Calkin algebra would be $C(\beta\N\setminus\N)$. A reference related to the general problem is \cite{horvath:kania:2021} where, by a pigeonhole argument, for any cardinal $\lambda$ there was shown to exists a $C(K)$ space of density $2^\lambda$ that is not the Calkin algebra of any Banach space of density $\lambda$. Note that such a $C(K)$ could still be the Calkin algebra of a Banach space of larger density.

From an algebraic perspective, the Banach algebra $\mathcal{L}(\mathfrak{X}_U)$ is similar to $\mathcal{L}(X)$, where $X$ is either one of the spaces from \cite{motakis:puglisi:zisimopoulou:2016} or \cite{motakis:2024} with Calkin algebra $C([0,\omega])$ or any of the spaces from \cite{motakis:puglisi:tolias:2020} with Calkin algebras of the form $\mathbb{C}I\oplus \mathcal{K}_\mathrm{diag}(Y)$. For all such $X$, there is a sequence of mutually annihilating infinite rank projections $(I_s)_{s=1}^\infty$ that together with the identity operator $I_0$ and $\mathcal{K}(X)$ generate $\mathcal{L}(X)$. Furthermore, the ideal structure is the same; there is an order-inverting bijection from $\mathcal{P}(\N_0)$ to the non-zero closed two-sided ideals of $\mathcal{L}(X)$ such that $A\mapsto\{T\in\mathcal{L}(X): T\circ I_s\text{ is compact for all }s\in A\}$. In addition, for each of these $X$, the $K_1$-group of $\mathcal{L}(X)$ coincides with the $K_1$-group of $C([0,\omega])$, which is trivial (the former follows, e.g., from \cite[Proposition 4.2]{laustsen:1997} and \cite{novodvorskii:1967}). Although we have not verified this, most likely $K_0(\mathcal{L}(X)) \simeq \mathbb{Z}\oplus K_0(C([0,\omega])) \simeq \mathbb{Z}\oplus C([0,\omega],\mathbb{Z}) \simeq C([0,\omega],\mathbb{Z})$. From an analytic perspective, there is a noteworthy difference concerning the position of the compact operator ideal; for all these $X$ except $\mathfrak{X}_U$, $\mathcal{K}(X)$ is uncomplemented in $\mathcal{L}(X)$.

Although our methods have expanded the list of unital Banach algebras that are Calkin algebras, this list only contains fairly simplistic non-commutative ones. We would thus like to draw attention to the discussion in \cite[Section 9]{motakis:2024} concerning the realization of specific non-commutative Banach algebras as Calkin algebras. As was discussed in Section \ref{section roadmap}, this process comprises a representation and an elimination component. The first one is a non-trivial step, but it can be carried out without any involvement of the Gowers-Maurey and Argyros-Haydon techniques. It would be of great interest to produce explicit operators on Bourgain-Delbaen-$\mathscr{L}_\infty$-spaces $\mathfrak{X}$ (existing ones and ad-hoc ones) such that the unital Banach algebra generated by them in $\mathcal{L}(\mathfrak{X})$ or $\mathpzc{Cal}(\mathfrak{X})$ is a prescribed commutative or non-commutative unital Banach algebra $\mathcal{B}$. References for explicit operator constructions in Bourgain-Delbaen-$\mathscr{L}_\infty$-spaces can be found in \cite{beanland:mitchell:2010}, \cite{beanland:mitchell:2014}, \cite{tarbard:2012}, and \cite{tarbard:2013} (shift operators), \cite{motakis:2024} (diagonal operators), and, of course, in this paper (mutually annihilating projections). It is very likely that such a representation would be amenable to the operator elimination component by involving the Argyros-Haydon technique and would, thus, lead to new examples of explicit Calkin algebras.

\part{Conceptual framework for defining $\mathfrak{X}_U$}
\label{conceptual part}

\section{Sequences of projections in the Calkin algebra}
\label{section generic sequences of projections}

In this section, we provide a criterion that yields the embedding of the unitization of a Banach space $U$ with an unconditional basis $(u_s)_{s=1}^\infty$ into the Calkin algebra of a Banach space $X$. More specifically, this criterion implies that a sequence of mutually annihilating projections $(I_s)_{s=1}^\infty$ on $X$ (i.e.,  $I_sI_t = 0$, for $s\neq t$) is equivalent to $(u_s)_{s=1}^\infty$. This is unusual behaviour; for example, a sequence of orthogonal mutually annihilating projections on a Hilbert space is always isometrically equivalent to the unit vector basis of $c_0$. Some of the conditions in the criterion impose some type of conditionality in the structure of $X$. In particular, the sequence $(I_s)_{s=1}^\infty$ is not associated with an infinite dimensional Schauder decomposition. In page \pageref{dual condition discussion}, we briefly discuss how our criterion relates via duality to a similar one used implicitly in the second author's paper \cite{pelczar-barwacz:2023} where the first known reflexive algebra $\mathcal{L}(X)/\mathcal{SS}(X)$ was constructed. Later in this section, we illustrate our criterion by giving an elementary example of a space $X$ on which it can be applied. The elimination of operators is not applicable to this $X$ as it requires heavy use of tools from the theory of hereditarily indecomposable Banach spaces.

\addtocontents{toc}{\SkipTocEntry}
\subsection*{A general condition on embedding the unitization of $U$ in the Calkin algebra}

The following is a criterion, used in the proof of our main result, for a sequence of uniformly bounded mutually annihilating projections $(I_s)_{s=1}^\infty$ to be equivalent to an unconditional sequence $(u_s)_{s=1}^\infty$.  The criterion is stated in a language that clarifies the properties that need to be satisfied at the level of the underlying Banach space $X$; this spacial interpretation is necessary to determine how to construct such an $X$, which we illustrate in Proposition \ref{simple example of U in Calkin}. Condition \ref{general condition on projections d} easily implies the domination of $(I_s)_{s=1}^\infty$ by $(u_s)_{s=1}^\infty$ and condition \ref{general condition on projections c} implies the domination of $(u_s)_{s=1}^\infty$ by the image of $(I_s)_{s=1}^\infty$ in Calkin algebra.  These conditions are based on the straightforward identification $\mathcal{K}_\mathrm{diag}(c_0,U)\equiv U$ described in Section \ref{section roadmap}.

Let $X$ be a Banach space. Recall that a 1-norming set $W$ for $X$ is a subset of the unit ball of $X^*$ such that, for all $x\in X$, $\|x\| = \sup\{|f(x)|:f\in W\}$. For a subset $S$ of $X^*$ denote $S_\perp = \cap_{f\in S}\mathrm{ker}(f)$, which is a closed subspace of $X$. In particular, if $S$ is finite, then $S_\perp$ is of finite codimension.

\begin{proposition}
\label{general condition on projections}
Let $U$ be a Banach space with a 1-unconditional normalized Schauder basis $(u_s)_{s=1}^\infty$, $X$ be a Banach space, and $(I_s)_{s=1}^\infty$ be a sequence of bounded linear projections on $X$. Denote, for $s\in\N$, $X_s = I_s(X)$, $X_0 = \cap_{s=1}^\infty \mathrm{ker}(I_s)$, and $Q:X\to X/X_0$ the quotient map. Assume that the following hold. 
\begin{enumerate}[leftmargin=19pt,label=(\alph*)]

\item\label{general condition on projections a} The sequence $(I_s)_{s=1}^\infty$ is uniformly bounded, i.e.,  $\sup_{s\geq 1}\|I_s\| = A<\infty$.
    
\item\label{general condition on projections b} The members of $(I_s)_{s=1}^\infty$ are mutually annihilating, i.e., for $s\neq t\in\mathbb{N}$, $I_sI_t = 0$.

\item\label{general condition on projections c} There exists a constant $C$ such that for every  $x_1\in X_1,\ldots,x_N\in X_N$,
    \[\Big\|\sum_{s=1}^NQx_s\Big\| \leq C\max_{1\leq s\leq N}\|x_s\|.\]

\item\label{general condition on projections d} There exists a constant $B$ such that for every $x_1\in X_1,\ldots,x_N\in X_N$,
    \[\Big\|\sum_{s=1}^Nx_s\Big\| \leq B\Big\|\sum_{s=1}^N\|x_s\|u_s\Big\|_U.\]

    \item\label{general condition on projections e} There exist a constant $\beta>0$ and a 1-norming set $W$ for $X$ such that for any finite $S\subset W$ and $a_1,\ldots,a_N\in\mathbb{C}$ there exist $x_1\in X_1\cap S_\perp,\ldots,x_N\in X_N\cap S_\perp$ of norm at most one with 
\[\Big\|\sum_{s=1}^Na_sx_s\Big\| \geq \beta\Big\|\sum_{s=1}^Na_su_s\Big\|_U.\]

\end{enumerate}
Then, for every $a_1,\ldots,a_N\in\mathbb{C}$ and compact operator $K:X\to X$ we have
\begin{equation}
\label{general condition on projections eq1}
\frac{\beta}{2C}\Big\|\sum_{s=1}^Na_su_s\Big\|_U \leq \Big\|\sum_{s=1}^Na_s I_s + K\Big\|\text{ and }\Big\|\sum_{s=1}^Na_s I_s\Big\| \leq AB\Big\|\sum_{s=1}^Na_su_s\Big\|_U.
\end{equation}
Therefore, the sequence $(I_s)_{s=1}^\infty$ in $\mathcal{L}(X)$ and its image in $\mathpzc{Cal}(X)$ are $2ABC/\beta$-equivalent to $(u_s)_{s=1}^\infty$ and, in particular, the spaces  $\mathcal{K}(X)$ and $[(I_s)_{s=1}^\infty]\oplus \mathbb{C}I$ form a direct sum in $\mathcal{L}(X)$.
\end{proposition}

\begin{proof}
We first prove the second part of \eqref{general condition on projections eq1}. Let $T = \sum_{s=1}^Na_sI_s$, fix $z\in X$ with $\|z\| = 1$, and note that, for $s\in\mathbb{N}$, $\|I_sz\| \leq A$. Therefore,
\[\|Tz\| = \Big\|\sum_{s=1}^Na_sI_sz\Big\| \leq B\Big\|\sum_{s=1}^N\|a_sI_sz\|u_s\Big\|_U \leq AB\Big\|\sum_{s=1}^Na_su_s\Big\|_U.\]
For the first part of \eqref{general condition on projections eq1}, fix a compact operator $K:X\to X$. Using \ref{general condition on projections e}, inductively pick, for $i\in\N$, vectors $(x_s^i)_{s=1}^N$ in $X$ of norm at most one and $f_i$ in $W$ such that, putting $S^0 = \emptyset$ and, for $i\in\N$, $S^i = \{f_1,\ldots,f_i\}$, the following are satisfied.
\begin{enumerate}

\item $x_1^i\in X_1\cap S^{i-1}_\perp$, \ldots, $x_N^i\in X_N\cap S^{i-1}_\perp$.

\item $\displaystyle{\Big|f_i\Big(\sum_{s=1}^Na_sx^i_s\Big)\Big| > (1-\varepsilon)\beta\Big\|\sum_{s=1}^Na_su_s\Big\|_U}$.

\end{enumerate}
Using \ref{general condition on projections c}, for $i\in\N$, pick $x_0^i\in X_0$ with
\[\Big\|\sum_{s=1}^nx_s^i - x_0^i\Big\| \leq (1+\varepsilon)C.\]
By the compactness of $K$, we may find an infinite $M\subset \mathbb{N}$ such that, for any $i < j$ in $M$, we have
\[\Big\|K\Big(\sum_{s=1}^N\big(x^i_{s}-x^j_{s}\big)-\big(x_0^i-x_0^j\big)\Big)\Big\| \leq \varepsilon.\]
We then calculate, for $i<j$ in $M$,
\begin{align*}
\|T + K\| &\geq \frac{1}{2(1+\varepsilon)C}\Big\|(T+K)\Big(\sum_{s=1}^N\big(x^i_s-x_s^j\big) - \big(x^i_0-x_0^j\big)\Big)\Big\|\\
&\geq \frac{1}{2(1+\varepsilon)C}\Big(\Big\|\sum_{s=1}^Na_s\big(x^i_s-x_s^j\big)\Big\| -\varepsilon\Big) \geq \frac{1}{2(1+\varepsilon)C}\Big(\Big|f_i\Big(\sum_{s=1}^Na_s\big(x^i_s-x_s^j\big)\Big)\Big| -\varepsilon\Big)\\
&= \frac{1}{2(1+\varepsilon)C}\Big(\Big|f_i\Big(\sum_{s=1}^Na_sx^i_s\Big)\Big| -\varepsilon\Big) \geq \frac{1}{2(1+\varepsilon)C}\Big((1-\varepsilon)\beta\Big\|\sum_{s=1}^Na_su_s\Big\|_U - \varepsilon\Big).
\end{align*}
\end{proof}

\begin{remarks}
\label{variation compact to strictly singular}~
\begin{enumerate}[leftmargin=19pt]

\item In Proposition \ref{general condition on projections}, putting $Z = 
[\cup_{s=0}^\infty X_s]$, then, by  \ref{general condition on projections a}, \ref{general condition on projections b}, and \ref{general condition on projections c}, $Z/X_0$ is $AC$-isomorphic to $(\oplus_{s=1}^\infty X_s)_{c_0}$ via $x+X_0\mapsto (I_sx)_{s=1}^\infty$.

\item In Proposition \ref{general condition on projections} \ref{general condition on projections e}, if we additionally assume that the convex hull of $W$ is norm-dense in the unit ball of $X^*$, then in \eqref{general condition on projections eq1} the factor $1/2$ can be omitted.

\item It is possible to create a variation of Proposition \ref{general condition on projections} in the conclusion of which $\mathcal{K}(X)$ is replaced by $\mathcal{SS}(X)$, the ideal of strictly singular operators and, thus, $\mathpzc{Cal}(X)$ is replaced by $\mathcal{L}(X)/\mathcal{SS}(X)$. To achieve this, the collection of closed subspaces $Y$ of $X$ such that, for all $s\in\N_0$, $Y\cap X_s$ is infinite-dimensional, needs to be incorporated into conditions  \ref{general condition on projections c} and \ref{general condition on projections e}.

\end{enumerate}
\end{remarks}

\begin{remark}
\label{compatibility remark}
The spaces $\mathfrak{X}_U$ constructed in this paper are Bourgain-Delbaen-$\mathscr{L}_\infty$-space and, thus, have an inherent $c_0$ structure. This makes the application of Proposition \ref{general condition on projections} very suitable. In Section \ref{bonding section}, we introduce a general construction process called a bonding of Bourgain-Delbaen-$\mathscr{L}_\infty$-spaces yielding spaces that automatically satisfy conditions Proposition \ref{general condition on projections} \ref{general condition on projections a}, \ref{general condition on projections b}, and \ref{general condition on projections c} via a ``cancellation of coordinates'' (see Proposition \ref{cancellation of coordinates}.) In Section \ref{section Linear combinations of coordinate projections} we formulate conditions on a bonding that guarantee that Proposition \ref{general condition on projections} \ref{general condition on projections d} and \ref{general condition on projections e} are satisfied as well (see Proposition \ref{bonding conditions guaranteeing upper-lower-U}.) We illustrate this  
on page \pageref{elementary BD example with U in Calkin} via a basic example of such a space $\mathfrak{B}_U$ to which this Proposition can be applied.
\end{remark}

\label{dual condition discussion}
A dual version of Proposition \ref{general condition on projections} exists, the strictly singular variation of which was used implicitly by the second author in {\cite{pelczar-barwacz:2023}} to construct the first known reflexive infinite dimensional quotient algebra of the form $\mathcal{L}(X)/\mathcal{SS}(X)$. We do not use or prove this statement; however, we would like to emphasize the way in which the conceptual framework of this paper draws from the earlier work of {\cite{pelczar-barwacz:2023}}. This statement is based on the identification $\mathcal{K}_\mathrm{diag}(U,\ell_1)\equiv [(u_s^*)_{s=1}^\infty]$ which follows by duality and $\mathcal{K}_\mathrm{diag}(c_0,[(u_s^*)_{s=1}^\infty])\equiv [(u_s^*)_{s=1}^\infty]$.  Comparing it to Proposition \ref{general condition on projections},  note that, essentially, the inequality in \ref{general condition on projections d} is reversed to yield \ref{dual general condition on projections beta}, and the upper-$c_0$-and-lower-$U$ conditions from \ref{general condition on projections c} and \ref{general condition on projections e} are replaced with the lower-$\ell_1$-and-upper-$U$ condition in \ref{dual general condition on projections gamma}. Requirement  \ref{general condition on projections a} is no longer necessary as it is implied by condition \ref{dual general condition on projections beta}.

\begin{proposition}[\cite{pelczar-barwacz:2023}]
\label{dual condition on projections}
Let $U$ be a Banach space with a 1-unconditional normalized Schauder basis $(u_s)_{s=1}^\infty$, $X$ be a Banach space, and $(I_s)_{s=1}^\infty$ be a sequence of bounded linear projections on $X$. Denote, for $s\in\N$, $X_s = I_s(X)$ and assume that the following hold. 
\begin{enumerate}[leftmargin=21pt,label=(\greek*),series=dcpList]
    
    \item\label{dual general condition on projections alpha} For every $s\neq t\in\mathbb{N}$, $I_sI_t = 0$.
   
    \item\label{dual general condition on projections beta} There exists a constant $B$ such that for every $N\in\mathbb{N}$ and $x_1\in X_1,\ldots,x_N\in X_N$, $x_0\in\cap_{s=1}^N\mathrm{ker}(I_s)$,
  \[ \Big\|\sum_{s=1}^Nx_s + x_0\Big\| \geq B^{-1}\Big\|\sum_{s=1}^N\|x_s\|u_s\Big\|_U.\]
   
    \item\label{dual general condition on projections gamma} There exist constants $\beta,C>0$ such that for every subspace $Y$ of $X$, with $Y\cap X_s$ infinite dimensional for all $s\in\N$, and $a_1,\ldots,a_N, b_1,\dots,b_N\in\mathbb{C}$ there are vectors $x_1\in X_1\cap Y,\ldots,x_N\in X_N\cap Y$ and $x_0\in \bigcap_{s=1}^N \mathrm{ker}(I_s)\cap Y$ such that 
\begin{equation*}\beta\sum_{s=1}^N |b_s|\leq \Big\|\sum_{s=1}^Nb_sx_s\Big\|\;\text{ and }\; \Big\|\sum_{s=1}^Na_sx_s - x_0\Big\| \leq C\Big\|\sum_{s=1}^Na_su_s\Big\|_U.
\end{equation*}

\end{enumerate}
Then the sequence $(I_s)_{s=1}^\infty$ in $\mathcal{L}(X)$ and its image in $\mathcal{L}(X)/\mathcal{SS}(X)$ are $BC/\beta$-equivalent to $(u^*_s)_{s=1}^\infty$ and, in particular, the spaces $[(I_s)_{s=1}^\infty]\oplus \mathbb{C}I$ and $\mathcal{SS}(X)$ form a direct sum in $\mathcal{L}(X)$.
\end{proposition}

\begin{remark}

The space from \cite{pelczar-barwacz:2023} is based on a Schlumprecht-type construction \cite{schlumprecht:1991}; specifically, it has a mixed-Tsirelson unconditional frame \cite{argyros:deliyanni:1997} incompatible with Proposition \ref{general condition on projections}. This still leaves Proposition \ref{dual condition on projections},  based on the identification $\mathcal{K}_\mathrm{diag}(U,\ell_1)\equiv U^*$, as an available tool. Although it is successfully implemented in \cite{pelczar-barwacz:2023}, contrary to the $\mathscr{L}_\infty$-case, its application is not as natural (see Remark \ref{compatibility remark}), and achieving condition \ref{dual general condition on projections gamma} requires the involvement of Gowers-Maurey-type construction techniques, i.e., techniques from the theory of hereditarily indecomposable spaces.
\end{remark}

\addtocontents{toc}{\SkipTocEntry}
\subsection*{An example of a space with the unitization of $U$ in its Calkin algebra}

For the purpose of conceptualization, given a sequence of infinite-dimensional separable Banach spaces $(X_s)_{s=1}^\infty$ we provide a simple example of a Banach space $X$ that contains isomorphic copies of all $X_s$ and to which Proposition \ref{general condition on projections} can be applied. Therefore, such $X$ contains the unitization of $U$ in its Calkin algebra. The construction in this subsection is not used in the sequel but it provides some intuition. Because this $X$ is so simply defined (i.e., without the use of methods from HI spaces), the representation of the unitization of $U$ in $\mathpzc{Cal}(X)$ is not onto.

\begin{proposition}
\label{simple example of U in Calkin}
Let $U$ be a Banach space with a $1$-unconditional normalized basis $(u_s)_{s=1}^\infty$ and $(X_s)_{s=1}^\infty$ be a sequence of separable infinite-dimensional Banach spaces. Then, there exists a separable Banach space $X$ and a sequence of norm-one linear projections $(I_s)_{s=1}^\infty$ on $X$ such that the following hold.
\begin{enumerate}[label=(\roman*)]

\item $I_sI_t = 0$, for $s\neq t$.

\item\label{simple example of U in Calkin item 2} The sequence $(I_s)_{s=1}^\infty$ in $\mathcal{L}(X)$ and its image in $\mathpzc{Cal}(X)$ are  $3$-equivalent to $(u_s)_{s=1}^\infty$.

\item For every $s\in\mathbb{N}$, $I_s(X)$ is isometrically isomorphic to $X_s$.

\item The spaces $\mathcal{K}(X)$ and $[(I_s)_{s=1}^\infty]\oplus\mathbb{C}I$ form a direct sum in $\mathcal{L}(X)$.

\end{enumerate}
\end{proposition}

\begin{proof}
For a simpler proof, we will use a result due to Pe\l czynski from \cite{pelczynski:1976}, namely that for any separable Banach space $X$ and $\varepsilon>0$, there exists a biorthogonal system $(e_i,e^*_i)_{i=1}^\infty$ in $X\times X^*$ such that $\sup_{i\in\N}\|e_i\|\|e_i^*\|\leq 1+\varepsilon$, $[\{e_i:i\in\N\}] = X$, and the set $W = B_{X^*}\cap\langle\{e^*_i:i\in\N\}\rangle$ is a 1-norming set for $X$. That $W$ can be chosen 1-norming is not stated explicitly but follows from \cite[Proof of Theorem 1, page 296]{pelczynski:1976} starting with a countable 1-norming set $\{x_i^*:i\in\N\}$ and property (5) stated inside that proof, which implies $\langle\{e^*_i:i\in\N\}\rangle\supset\langle\{x^*_i:i\in\N\}\rangle$.

For every $s\in\mathbb{N}$ and $\varepsilon>0$, fix a biorthogonal system $(x_{s,i},f_{s,i})_{i=1}^\infty$ in $X_s\times X_s^*$ such that, for $i\in\N$, $\|x_{s,i}\| = 1$ and $\|f_{s,i}\|\leq 1+\varepsilon$, $X_s = [\{x_{s,i}:i\in\N\}]$, and $W_s = B_{X_s^*}\cap \langle\{f_{s,i}:i\in\N\}\rangle$ is a 1-norming set for $X_s$. Put, for $s,i\in\N$, $g_{s,i} = \|f_{s,i}\|^{-1}f_{s,i}$ and denote by $(e_{s,i},e_{s,i}^*)_{i,s\in\mathbb{N}^2}$ the standard biorthogonal system of $c_{00}(\mathbb{N}^2)$.

Let $Z$ denote the vector subspace of $c_{00}(\mathbb{N}^2)\times \Pi_{s=1}^\infty X_s$ of all sequences $(x_0,x_1,\ldots,x_s,\ldots)$ that are eventually zero and let $X$ be the completion of $Z$ under the norm described as follows: if $x = (x_0,x_1,\ldots,x_s,\ldots)$ then
\[\|x\| = \max\Big\{\max_{s\geq 1}\|x_s\|,\sup_{i\in\mathbb{N}}\Big\|\sum_{s=1}^\infty\big(g_{s,i}(x_s) + e_{s,i}^*(x_0)\big)u_s\Big\|_U\Big\}.\]

Denote by $X_0$ the closure of the canonical image of $c_{00}(\mathbb{N}^2)$ in $X$. It is not difficult to verify that, for $s\geq 1$, the canonical embedding of $X_s$ in $X$ is an isometry and the map $(x_0,x_1,\ldots,x_{s-1},x_s,x_{s+1},\ldots)\mapsto (0,0,\ldots,0,x_s,0,\ldots)$ defines a projection $I_s:X\to X_s$ with $\|I_s\| = 1$. Therefore, assumptions \ref{general condition on projections a}  and \ref{general condition on projections b} of Proposition \ref{general condition on projections} are satisfied for $A = 1$ and it is straightforward to check that \ref{general condition on projections d} is satisfied for $B=1$ as well. For context, let us point out that $X_0$ is isometric to $c_0(U)$ and $X_0$ might not be complemented in $X$.

We verify \ref{general condition on projections c} for $C=1$ in the canonical image of the dense subspaces $\langle\{x_{s,i}:i\in\N\}\rangle$ of $X_s$, $s\in\N$. For $x_s = \sum_{i=1}^ma_{s,i}x_{s,i}$, $1\leq s\leq N$, putting $x_0 = \sum_{s=1}^N\sum_{i=1}^ma_{s,i}\|f_{s,i}\|^{-1}e_{s,i}$ we ``cancel the coordinates'' of $\sum_{s=1}^nx_s$ to directly obtain
\[\Big\|\sum_{s=1}^nx_s - x_0\Big\| = \max_{1\leq s\leq N}\|x_s\|.\]

We verify property \ref{general condition on projections e} for $\beta = 1/(1+\varepsilon)$. Naturally identify, for $s\in\N$, $W_s$ with a subset of $X^*$ and, for $i\in\N$, $e_{s,i}^*$ with a functional in $X^*$. We define $W = (\cup_{s=0}^\infty W_s)\cup (\cup_{i=1}^\infty W^{(i)})$, where, for $i\in\N$,
\[W^{(i)} = \Big\{\sum_{s=1}^N\lambda_s(g_{s,i}+e_{s,i}^*):N\in\N\text{ and }\Big\|\sum_{s=1}^N\lambda_su_s^*\Big\|_{U^*}\leq 1\Big\}.\]
By the formula for the norm of $X$, this is a 1-norming set. Note that, for every $f\in W$ and $s\in\N$, the sequence $(f(x_{s,i}))_{i\in\N}$ is eventually null. To conclude the proof, fix a finite $S\subset W$ and $a_1,\ldots,a_N\in\mathbb{C}$. Then, for every $i\in\N$, we have
\[\Big\|\sum_{s=1}^Na_sx_{s,i}\Big\| \geq \frac{1}{(1+\varepsilon)}\Big\|\sum_{s=1}^Na_su_s\Big\|_U\]
and, for $1\leq s\leq N$, $(x_{s,i})_{i=1}^\infty$  is eventually in $S_\perp$. The conclusion follows if $0<\varepsilon<1/2$.
\end{proof}

\begin{remarks}\phantom{A}~
\begin{enumerate}[label=(\roman*),leftmargin=19pt]

\item In Proposition \ref{simple example of U in Calkin}, if we denote $X_0 = \cap_{s=1}^\infty\mathrm{ker}(I_s)$, then $X/X_0$ is isometrically isomorphic to $(\oplus_{s=1}^\infty X_s)_{c_0}$.

\item If, for each $s\in\N$, $X_s$ has a weakly null bimonotone Schauder basis $(x_{s,i})_{i=1}^\infty$ (or, more generally, a weakly null 1-norming Auerbach basis), we may improve Proposition \ref{simple example of U in Calkin} \ref{simple example of U in Calkin item 2} so that the sequence $(I_s)_{s=1}^\infty$ in $\mathcal{L}(X)$ and it image in $\mathpzc{Cal}(X)$ are isometrically equivalent to $(u_s)_{s=1}^\infty$.

\end{enumerate}
\end{remarks}

We close this section with a general observation about Banach spaces admitting a sequence of projections satisfying assumptions \ref{general condition on projections a}, \ref{general condition on projections b}, and \ref{general condition on projections c} of Proposition \ref{general condition on projections}. Although we do not use it explicitly, it reveals a limitation of the methods developed in the paper, namely bondings may not be suitable to achieve as Calkin algebras spaces with conditional submultiplicative Schauder bases, such as those from \cite{motakis:puglisi:tolias:2020}.

A Banach space $X$ is said to have the bounded approximation property if, for some $\lambda>0$, for every finite subset $S$ of $X$ and $\varepsilon>0$ there exists a finite rank linear operator $F:X\to X$, with $\|F\|\leq \lambda$, such that, for $x\in S$, $\|Fx - x\| < \varepsilon$.

\begin{proposition}
\label{always unconditional}
Let $X$ be a Banach space and $(I_s)_{s=1}^\infty$ be a sequence of bounded linear projections  on $X$ satisfying assumptions \ref{general condition on projections a}, \ref{general condition on projections b}, and \ref{general condition on projections c} of Proposition \ref{general condition on projections}. Then, the sequence $(I_s)_{s=1}^\infty$ in $\mathcal{L}(X)$ is unconditional. If, furthermore, $X$ has the bounded approximation property then, for $M = \{s\in\N:I_s\text{ is of infinite rank}\}$, the image of $(I_s)_{s\in M}$ in $\mathpzc{Cal}(X)$ is unconditional.
\end{proposition}

\begin{proof}
In the general case, fix scalars $a_1,\ldots,a_N$, $b_1,\ldots,b_N$ with $|a_s|\leq |b_s|$, for $1\leq s\leq N$. Take $x\in X$ with $\|x\|\leq 1$ and, for $C'>C$, pick $x_0\in X_0$ such that, if $y = \sum_{s=1}^N(a_s/b_s)I_sx - x_0$, then $\|y\| \leq C'A$ (we followed the usual convention $0/0 = 0$). It follows that $C'A\|\sum_{s=1}^Nb_sI_s\| \geq \|\sum_{s=1}^Nb_sI_s(y)\| = \|\sum_{s=1}^Na_sI_s(x)\|$ and, in conclusion, $(I_s)_{s=1}^\infty$ is $CA$-unconditional.

If we additionally assume that $X$ has the $\lambda$-approximaton property, for $T\in\mathcal{L}(X)$, let
\[|T| = \sup\Big\{\liminf_i\liminf_j\|Tx_i - Tx_j\|:(x_i)_{i=1}^\infty\text{ is a sequence in }B_X\Big\}\]
which defines a seminorm on $\mathcal{L}(X)$ such that $(1+\lambda)^{-2}\mathrm{dist}(T,\mathcal{K}(X))\leq |T|\leq 2\mathrm{dist}(T,\mathcal{K}(X))$. The inequality $|T|\leq 2\mathrm{dist}(T,\mathcal{K}(X))$ is easily shown by observing that, for $K\in\mathcal{K}(X)$, $|K| = 0$. To show the other inequality, inductively define a sequence $(z_i)_{i=1}^\infty$ in $B_X$ and a sequence of finite rank operators $F_i:X\to X$, with $\|F_i\|\leq \lambda$, for $i\in\N$, such that the following hold:
\begin{enumerate}[label=(\roman*)]

    \item $\liminf_i\|(I-F_i)T(I-F_i)z_i\| = \liminf_i\|(I-F_i)T(I-F_i)\|$.

    \item For $i<j$, $\|(I - F_j)T(I-F_i)z_i\| \leq 1/j$.
    
\end{enumerate}
We obtain
\[
\begin{split}
|T| &\geq(1+\lambda)^{-1}\liminf_i\liminf_j\|T(I-F_i)x_i-T(I-F_j)x_j\|\\
& \geq (1+\lambda)^{-2}\liminf_i\liminf_j\|(I-F_j)T(I-F_i)z_i-(I-F_j)T(I-F_j)z_j\|\\
& = (1+\lambda)^{-2}\liminf_j\|(I-F_j)T(I-F_j)z_j\| = (1+\lambda)^{-2}\liminf_j\|(I-F_j)T(I-F_j)\|\\
&\geq (1+\lambda)^{-2}\mathrm{dist}(T,\mathcal{K}(X)).
\end{split}
\]
One then repeats, for $|\cdot|$, the argument used to prove the unconditionality of $(I_s)_{s=1}^\infty$ in $\mathcal{L}(X)$ to obtain that its image in $\mathpzc{Cal}(X)$ is $2(1+\lambda)^2AC$-unconditional.
\end{proof}

\section{The Bourgain-Delbaen-$\mathscr{L}_\infty$ construction method}
\label{BD section}
A separable $\mathscr{L}_\infty$-space is a Banach space $X$ for which there exists a sequence $Y_0\subset Y_1\subset\cdots$ of finite dimensional subspaces of $X$ with union dense in $X$ such that $\sup_n\mathrm{dist}(Y_n,\ell_\infty^{\mathrm{dim}(Y_n)}) <\infty$. This is a class of Banach spaces introduced by Lindenstauss and Pe\l czynski in {\cite{lindenstrauss:pelczynski:1968}} while studying absolutely summing operators on $L_p$ spaces. To construct the first example of a separable RNP space not embedding in a dual space, Bourgain and Delbaen invented in {\cite{bourgain:delbaen:1980}} the synonymous construction method of defining non-classical $\mathscr{L}_\infty$-spaces. In {\cite{argyros:haydon:2011}} Argyros and Haydon used this method, in combination with the Gowers-Maurey construction from {\cite{gowers:maurey:1993}}, to construct the first known example of a space with the scalar-plus-compact property. In this section, we discuss the theory behind directly defining Bourgain-Delbaen $\mathscr{L}_\infty$-spaces and how some of their quotients define further such spaces. We follow the language from {\cite{argyros:gasparis:motakis:2016}} and {\cite{argyros:motakis:2014}} that discusses Bourgain-Delbaen spaces as abstract objects based on the notation of {\cite{argyros:haydon:2011}}. In this paper, the abstract language serves two main purposes. The first one is to obtain the space $\mathfrak{X}_U$ as a quotient of a cruder Bourgain-Delbaen-$\mathscr{L}_\infty$-space $\bar{\mathfrak{X}}$. This approach has notational advantages; it has been explicitly used in {\cite{argyros:motakis:2019}}, {\cite{manoussakis:pelczar-barwacz:swietek:2017}}, and {\cite{motakis:2024}}. A similar scheme was used implicitly earlier by Tarbard in {\cite{tarbard:2013}}. The second purpose is more substantial; we develop the necessary language and define a class of Bourgain-Delbaen-$\mathscr{L}_\infty$-spaces, called bondings, and formulate appropriate conditions under which such spaces satisfy the assumptions of Proposition \ref{general condition on projections}.

\subsection{The general Bourgain-Delbaen scheme}
\label{Bourgain-Delbaen general section}
We discuss, in a general language, the Bourgain-Delbaen method of constructing $\mathscr{L}_\infty$-spaces.

\addtocontents{toc}{\SkipTocEntry}
\subsection*{Abstract Bourgain-Delbaen spaces}
There exist several approaches to the abstract presentation of Bourgain-Delbaen-$\mathscr{L}_\infty$-spaces, e.g., in {\cite{argyros:haydon:2011}} as spaces resulting from duality. Ours is to a define Bourgain-Delbaen-$\mathscr{L}_\infty$-space $\mathfrak{X}_{(\Gamma_n,i_n)}$ as a direct limit of finite-dimensional spaces $\ell_\infty(\Gamma_n)$, $n\in\mathbb{N}_0$, and compatible linear maps $i_{m,n}:\ell_\infty(\Gamma_m)\to\ell_\infty(\Gamma_n)$, $0\leq m\leq n$, satisfying specific additional assumptions. In practice, this is carried out via a recursive algorithm that explicitly yields the $n$'th element as a closed expression based on the preceding elements. 

\begin{notation}
As usual, for a set $\Gamma$, denote
\[\ell_\infty(\Gamma) = \big\{x:\Gamma\to\mathbb{C}\text{ bounded}\big\}\]
endowed with the norm $\|x\|_\infty = \sup_{\gamma\in\Gamma}|x(\gamma)|$. For $\gamma\in\Gamma$ we denote by $e_\gamma^*:\ell_\infty(\Gamma)\to\mathbb{C}$ the usual evaluation functional given by $e_\gamma^*(x) = x(\gamma)$. We make the convention $\ell_\infty(\emptyset) = \{0\}$.

Let $\Gamma_0\subset\Gamma_1$ be finite sets.
\begin{enumerate}[label=(\arabic*),leftmargin=19pt]
    
    \item Denote by $r_{\Gamma_0}:\ell_\infty(\Gamma_1)\to\ell_\infty(\Gamma_0)$ the usual restriction map given by $r_{\Gamma_0}(x) = x|_{\Gamma_0}$, i.e., $(r_{\Gamma_0}(x))(\gamma) = x(\gamma)$, for $x\in\ell_\infty(\Gamma_1)$ and $\gamma\in\Gamma_0$.

    \item Any linear right inverse $i:\ell_\infty(\Gamma_0)\to\ell_\infty(\Gamma_1)$ of $r_{\Gamma_0}$ is called a linear extension map. That is, $(i(x))(\gamma) = x(\gamma)$, for every $x\in\ell_\infty(\Gamma_0)$ and $\gamma\in\Gamma_0$.

\end{enumerate}
When $\Gamma_0 = \emptyset$ and, thus, $\ell_\infty(\Gamma_0) = \{0\}$, we declare that the zero map is both a restriction map and a linear extension map.

Let $\Gamma_0\subset\Gamma_1\subset\cdots\subset\Gamma_n\subset\cdots$ be finite sets and let, for each $n\in\mathbb{N}$, $i_{n-1,n}:\ell_\infty(\Gamma_{n-1})\to\ell_\infty(\Gamma_n)$ be a linear extension map. We will refer to such a pair $(\Gamma_n)_{n=0}^\infty$, $(i_{n-1,n})_{n=1}^\infty$ as a linear extension scheme.

\begin{enumerate}[label=(\arabic*),leftmargin=19pt]\setcounter{enumi}{2}

\item For $m\leq n\in\mathbb{N}_0$ denote by $i_{m,n}:\ell_\infty(\Gamma_m)\to\ell_\infty(\Gamma_n)$ the linear extension map given by $i_{m,n} = i_{n-1,n}\circ i_{n-2,n-1}\circ\cdots \circ i_{m,m+1}$. Make the convention $i_{n,n} = id:\ell_\infty(\Gamma_n)\to\ell_\infty(\Gamma_n)$.

\item Denote $\Gamma = \cup_{n=0}^\infty\Gamma_n$ and for each $n\in\mathbb{N}_0$ define $i_n:\ell_\infty(\Gamma_n)\to\mathbb{C}^\Gamma$ given by $(i_n(x))(\gamma) = \lim_{N\to\infty}(i_{n,N}(x))(\gamma)$ for every $x\in\ell_\infty(\Gamma_n)$ and $\gamma\in\Gamma$.

\end{enumerate}
\end{notation}

\begin{assumption}
\label{general BD extension assumption}
Let $(\Gamma_n)_{n=0}^\infty$, $(i_{n-1,n})_{n=1}^\infty$ be a a linear extension scheme and assume that there exists $C<\infty$ with $\sup_{m\leq n}\|i_{m,n}\| \leq C$.
\end{assumption}

\begin{definition}
\label{BD definition}
A linear extension scheme $(\Gamma_n)_{n=0}^\infty$, $(i_{n-1,n})_{n=1}^\infty$ that satisfies Assumption \ref{general BD extension assumption}, for some finite $C$, is called a Bourgain-Delbaen scheme. Then, for every $n\in\mathbb{N}_0$, the operator $i_n:\ell_\infty(\Gamma_n)\to\ell_\infty(\Gamma)$ is well defined and $\|i_n\| \leq C$. Setting $Y_n = i_n(\ell_\infty(\Gamma_n))$, $Y_0\subset Y_1\subset\cdots\subset Y_n\subset\cdots$ and each $Y_n$ is $C$-isomorphic to $\ell_\infty(\Gamma_n)$. Therefore, the space $\mathfrak{X}_{(\Gamma_n,i_n)} = \overline{\cup_{n=0}^\infty Y_n}$ is a $C$-$\mathscr{L}_\infty$-subspace of $\ell_\infty(\Gamma)$ called a Bourgain-Delbaen-$C$-$\mathscr{L}_\infty$-space.
\end{definition}

When we call a space a Bourgain-Delbaen-$\mathscr{L}_\infty$-space, we mean that it is a Bourgain-Delbaen-$C$-$\mathscr{L}_\infty$-space, for some $C\geq 1$.

Despite its simplicity, as will see in Sections \ref{structures in BD} and \ref{section self-determined}, this abstract definition comes with some structures that allow the study of such a space. It was proved in {\cite{argyros:gasparis:motakis:2016}} that every separable $\mathscr{L}_\infty$-space is isomorphic to some Bourgain-Delbaen-$\mathscr{L}_\infty$-space.

\begin{example}
For $n\in\N_0$ take $\Gamma_n = \{0,\ldots,n\}$ and for $n\in\N$ and $x\in\ell_\infty(\Gamma_{n-1})$ let $i_{n-1,n}(x) = x^\frown x(n-1)$, i.e., extend $x$ by repeating its last entry. Then, it is easy to see that $(\Gamma_n)_{n=0}^\infty$, $(i_{n-1,n})_{n=1}^\infty$ is a Bourgain-Delbaen scheme and that the resulting space $\mathfrak{X}_{(\Gamma_n,i_n)}$ is isometrically isomorphic to $c$, the space of convergent sequences with the $\|\cdot\|_\infty$-norm.
\end{example}

\begin{remark}
\label{equivalent BD schemes}
From an inductive construction viewpoint, the following equivalent formulation of a linear extension scheme is particularly useful.
\begin{enumerate}[leftmargin=19pt,label=(\arabic*)]

\item\label{equivalent BD schemes functionals} Let $\Delta_0,\Delta_1,\ldots,\Delta_n,\ldots$ be a sequence of disjoint finite sets and for each $n\in\mathbb{N}_0$ denote $\Gamma_n = \cup_{m=0}^n\Delta_m$. Assume that for each $n\in\mathbb{N}$ and $\gamma\in\Delta_n$ we are given a linear functional $c_\gamma^*:\ell_\infty(\Gamma_{n-1})\to\mathbb{C}$. Then, for each $n\in\mathbb{N}$, the linear map $i_{n-1,n}:\ell_\infty(\Gamma_{n-1})\to\ell_\infty(\Gamma_n)$ given by
\[\big(i_{n-1,n}(x)\big)(\gamma) = \left\{
	\begin{array}{ll}
		x(\gamma)  & \mbox{if } \gamma\in\Gamma_{n-1}, \\
		c_\gamma^*(x) & \mbox{if } \gamma\in\Delta_n
	\end{array}
\right.\]
defines a linear extension map. Therefore, $(\Delta_n)_{n=0}^\infty$, $((c_\gamma^*)_{\gamma\in\Delta_n})_{n=1}^\infty$ induces a linear extension scheme $(\Gamma_n)_{n=0}^\infty$, $(i_{n-1,n})_{n=1}^\infty$.

\item Conversely, let $(\Gamma_n)_{n=0}^\infty$, $(i_{n-1,n})_{n=1}^\infty$ be a given linear extension scheme. Denote $\Delta_0 = \Gamma_0$ and for $n\in\mathbb{N}$ denote $\Delta_n = \Gamma_n\setminus\Gamma_{n-1}$. Let $n\in\mathbb{N}$ and for $\gamma\in\Delta_n$ denote by $c_\gamma^*:\ell_\infty(\Gamma_{n-1})\to\mathbb{C}$ the linear functional given by $c_\gamma^*(x) = e_\gamma^*(i_{n-1,n}(x))$. Then, the pair $(\Delta_n)_{n=0}^\infty$, $((c_\gamma^*)_{\gamma\in\Delta_n})_{n=1}^\infty$ induces the given linear extension scheme $(\Gamma_n)_{n=0}^\infty$, $(i_{n-1,n})_{n=1}^\infty$.

\end{enumerate}
In light of this equivalence, we also refer to a pair $(\Delta_n)_{n=0}^\infty$, $((c_\gamma^*)_{\gamma\in\Delta_n})_{n=1}^\infty$ as in \ref{equivalent BD schemes functionals} as a linear extension scheme and, if additionally Assumption \ref{general BD extension assumption} is satisfied, as a Bourgain-Delbaen scheme. This formulation is the one used in practice to define Bourgain-Delbaen-$\mathscr{L}_\infty$-spaces.

For obvious reasons, the functionals $c_\gamma^*$ are called extension functionals.
\end{remark}

\begin{example}
Define $\Delta_0 = \{0,1\}$ and, for $n\in\N$, let $\Delta_n = \{k/2^n: 1\leq k<2^n,\;k\text{ odd}\}$. For $n\in\N$ and $\gamma = k/2^n\in\Delta_n$ define $c_\gamma^*:\ell_\infty(\Gamma_{n-1})\to\mathbb{C}$ as follows. Let $\gamma^- = (k-1)/2^{n}\in\Delta_{n-1}$ and $\gamma_n^+ = (k+1)/2^{n}\in\Delta_{n-1}$, such that $\gamma = (1/2)(\gamma^-+\gamma^+)$, and put $c_\gamma^*(x) = (1/2)(x(\gamma^-)+x(\gamma^+))$. It is easy to see that $(\Delta_n)_{n=0}^\infty$, $((c_\gamma^*)_{\gamma\in\Delta_n})_{n=1}^\infty$ defines a Bourgain-Delbaen scheme. By naturally identifying each $\ell_\infty(\Gamma_n)$ with a space of piecewise linear continuous functions on $[0,1]$ it follows that the space $\mathfrak{X}_{(\Gamma_n,i_n)}$ is isometrically isomorphic to $C[0,1]$.
\end{example}

\addtocontents{toc}{\SkipTocEntry}
\subsection*{Achieving Assumption \ref{general BD extension assumption}}
Given a linear extension scheme  $(\Gamma_n)_{n=0}^\infty$, $(i_{n-1,n})_{n=1}^\infty$ (or, equivalently, $(\Delta_n)_{n=0}^\infty$, $((c_\gamma^*)_{\gamma\in\Delta_n})_{n=1}^\infty$) there is a special condition, introduced in {\cite{bourgain:delbaen:1980}}, that, when imposed on the extension functionals $c_\gamma^*$, automatically yields the boundedness Assumption \ref{general BD extension assumption}. Although this is standard, our approach is rather formal, which allows us to develop a general theory of bondings in Section \ref{bonding section} and then construct the space $\mathfrak{X}_U$ within this framework. We first introduce some required terminology.

\begin{notation}
Let $(\Gamma_n)_{n=0}^\infty$, $(i_{n-1,n})_{n=1}^\infty$ be a linear extension scheme and $\gamma\in\Gamma$. Let $\mathrm{rank}(\gamma)$ denote the unique $n\in\mathbb{N}_0$ such that $\gamma\in\Delta_n$.
\end{notation}

\begin{notation}
Let $(\Gamma_n)_{n=0}^\infty$, $(i_{n-1,n})_{n=1}^\infty$ be a linear extension scheme and $n\in\mathbb{N}_0$.
\begin{enumerate}[label=(\arabic*),leftmargin=19pt]

\item For every $0\leq k\leq n$ denote by $P^{(n)}_k:\ell_\infty(\Gamma_n)\to\ell_\infty(\Gamma_n)$ the linear projection $i_{k,n}\circ r_{\Gamma_k}$, i.e., for $x\in\ell_\infty(\Gamma)$, $P_k^{(n)}(x) = i_{k,n}(x|_{\Gamma_k})$.

\item For $E = \{k,k+1,\ldots,m\}\subset\{0,\ldots,n\}$ denote by $P_E^{(n)}:\ell_\infty(\Gamma_n)\to\ell_\infty(\Gamma_n)$ the linear projection $P_E^{(n)} = P^{(n)}_m - P^{(n)}_{k-1}$, where $P^{(n)}_{-1} = 0$.

\end{enumerate}
\end{notation}

Recall that a partial function $f$ with domain a set $X$ is a function $f$ with domain $\mathrm{dom}(f)\subset X$ and codomain some set $Y$. For $x\in X\setminus \mathrm{dom}(f)$ we say that $f(x)$ is undefined.

\begin{assumption}
\label{BD construction assumption}
Let $(\Gamma_n)_{n=0}^\infty$, $(i_{n-1,n})_{n=1}^\infty$ be a linear extension scheme. Assume that there exists $1<\vartheta < 1/2$ and partial functions $\xi(\cdot),E(\cdot),\lambda(\cdot),\eta(\cdot)$ with domain $\Gamma$ such that, for every $n\in\mathbb{N}_0$, every $\gamma\in\Gamma$, with $\rank(\gamma) = n$, satisfies one of the following:
\begin{enumerate}[label=(0),leftmargin=19pt]
    
    \item\label{BD construction assumption 0} $\xi(\gamma)$, $E(\gamma)$, $\lambda(\gamma)$, and $\eta(\gamma)$ are undefined and, if $n\geq 1$, $c_\gamma^* = 0$.
    
\end{enumerate}
\begin{enumerate}[label=(\alph*),leftmargin=19pt]

    \item\label{BD construction assumption 1}
    $\xi(\gamma)$ is undefined, $E = E(\gamma)$ is an interval of $\{0,\ldots,n-1\}$, $\lambda = \lambda(\gamma)$ is in $\mathbb{C}$ with $|\lambda|\leq\vartheta$, $\eta = \eta(\gamma)$ is in $\Gamma_{n-1}$, and
    \[c_\gamma^* = \lambda e_\eta^*\circ P_E^{(n-1)}.\]
    \item\label{BD construction assumption 2} $\xi = \xi(\gamma)$ is in $\Gamma_{n-1}\cap\mathrm{dom}(E(\cdot))$, $E = E(\gamma)$ is an interval of $\{\rank(\xi)+1,\ldots,n-1\}$, $\lambda=\lambda(\gamma)$ is in $\mathbb{C}$ with $|\lambda|\leq\vartheta$, $\eta = \eta(\gamma)$ is in $\Gamma_{n-1}$ with $\rank(\xi) < \rank(\eta)$, and
    \[c_\gamma^* = e^*_\xi + \lambda e_\eta^*\circ P_E^{(n-1)}.\]
    
\end{enumerate}
\end{assumption}

\begin{remark}
Necessarily, $\mathrm{dom}(\xi(\cdot)) \subset \mathrm{dom}(E(\cdot)) = \mathrm{dom}(\lambda(\cdot)) = \mathrm{dom}(\eta(\cdot))$. By definition, the partial functions $\xi(\cdot)$, $E(\cdot)$, $\lambda(\cdot)$, $\eta(\cdot)$ encode the form of the extension functionals $(c^*_\gamma)_{\gamma\in\Gamma}$.
\end{remark}

\begin{proposition}[{\cite[Lemma 4.1]{bourgain:delbaen:1980} or \cite[Theorem 3.4]{argyros:haydon:2011}}]
\label{bd construction proposition}
Let $(\Gamma_n)_{n=0}^\infty$, $(i_{n-1,n})_{n=1}^\infty$ be a linear extension scheme that satisfies Assumption \ref{BD construction assumption}, for some $0<\vartheta<1/2$. Then $\sup_{m\leq n}\|i_{m,n}\| \leq 1/(1-2\vartheta)$ and, therefore, for $C=1/(1-2\vartheta)$, $\mathfrak{X}_{(\Gamma_n,i_n)}$ is a Bourgain-Delbaen-$C$-$\mathscr{L}_\infty$-space.
\end{proposition}

The following notation can be defined easily now but will only become useful later when we define the evaluation analysis of a $\gamma$ in Section \ref{section evaluation analysis}.

\begin{remark}
\label{gamma tuples}
Let $\mathfrak{X}_{(\Gamma_n,i_n)}$ be a Bourgain-Delbaen-$\mathscr{L}_\infty$-space that satisfies Assumption \ref{BD construction assumption}. For $\gamma\in\Gamma$ that satisfies Assumption  \ref{BD construction assumption} \ref{BD construction assumption 1} or  \ref{BD construction assumption 2} there exists a unique $a\in\N$, denoted $\mathrm{age}(\gamma)$, and $a$-tuples $\vec\xi(\gamma) = (\xi_1,\ldots,\xi_a)$, $\vec E(\gamma) = (E_1,\ldots,E_a)$, $\vec\lambda(\gamma) = (\lambda_1,\ldots,\lambda_a)$, and $\vec\eta(\gamma) = (\eta_1,\ldots,\eta_a)$ such that:
\begin{enumerate}[label=(\roman*)]

\item $\xi_a =\gamma$, for $1\leq r< a$, $\xi_r = \xi(\xi_{r+1})$, and $\xi(\xi_1)$ is undefined,

\item  for $1\leq r\leq a$, $E_r = E(\xi_r)$, $\lambda_r = \lambda(\xi_r)$, and $\eta_r = \eta(\xi_r)$.

\end{enumerate}
Then, $\vec\xi(\cdot)$, $\vec E(\cdot)$, $\vec \lambda(\cdot)$, $\vec\eta(\cdot)$ are partially defined functions with domain $\Gamma$, collectively called the analysis of $\gamma$.

We complete the definition of the function $\mathrm{age}:\Gamma\to\N_0$ letting $\age(\gamma)=0$ for any $\gamma\in\Gamma$ satisfying Assumption  \ref{BD construction assumption} \ref{BD construction assumption 0}.

Although we do not exploit this explicitly in this paper, we may imbue $\Gamma$ with a tree structure in which $\gamma'\preceq\gamma$ whenever $\vec\xi(\gamma')$ is an initial segment of $\vec\xi(\gamma)$. For example, if $\vec\xi(\gamma) = (\xi_1,\ldots,\xi_a)$ and, for some $1\leq r\leq a$, we take $\gamma' = \xi_r$ then $\vec\xi(\gamma') = (\xi_1,\ldots,\xi_r)$, which is an initial segment of $\vec\xi(\gamma)$. This may be an infinitely rooted tree with roots all $\gamma\in\Gamma$ satisfying Assumption  \ref{BD construction assumption} \ref{BD construction assumption 0} or  \ref{BD construction assumption 1}. This is why members $\gamma$ of $\Gamma$ are sometimes referred to as nodes.
\end{remark}

At this point, we can conceptualize the Bourgain-Delbaen method via a specific simple algorithmic scheme that is very similar to the one given by Bourgain and Delbaen in {\cite{bourgain:delbaen:1980}}. We emphasize that verifying Assumption \ref{BD construction assumption} for some $n\in\mathbb{N}$ and $0<\vartheta<1/2$ depends only on the objects $(\Gamma_m)_{m=0}^n$ and $(i_{m-1,m})_{m=1}^n$ (or equivalently $(\Delta_m)_{m=0}^n$, $((c_\gamma^*)_{\gamma\in\Delta_m})_{m=1}^n$).

\begin{example}
\label{schur bd space}
Fix $0<\vartheta<1/2$. Define $\Delta_0 = \{(0)\}$ and, assuming that $\Delta_0,\Delta_1,\ldots,\Delta_{n-1}$ and $(c_\gamma^*)_{\gamma\in\Delta_m}$, $1\leq m\leq n-1$ have been defined, define $\Delta_n = \Delta_n^{0}\cup\Delta_n^\mathrm{a}\cup\Delta_n^\mathrm{b}$ and $(c_\gamma^*)_{\gamma\in\Delta_n}$ as follows. Let $\Delta_n^{0} = \{(n)\}$ and, for $\gamma = (n)$, leave $\xi(\gamma)$, $E(\gamma)$, $\lambda(\gamma)$, and $\eta(\gamma)$ undefined and put $c_\gamma^* = 0$.

Let
\[\Delta_n^\mathrm{a} = \Big\{(n,E,\vartheta,\eta): \eta\in\Gamma_{n-1},E\subset\{0,\ldots,n-1\}\text{ interval}\Big\}\]
and, for $\gamma =  (n,E,\vartheta,\eta)\in\Delta_{n}^\mathrm{a}$ leave $\xi(\gamma)$ undefined, put $E(\gamma) = E$, $\eta(\gamma) = \eta$, $\lambda(\gamma) = \vartheta$, and $c_\gamma^* = \vartheta e_\eta^*\circ P^{(n-1)}_E$.

Let
\begin{align*}
\Delta_n^\mathrm{b} = \Big\{(n,\xi,E,\vartheta,\eta): &\;\xi\in\Gamma_{n-2}\cap\mathrm{dom}(E(\cdot)),\;E\subset\{\rank(\xi)+1,\ldots,n-1\}\text{ interval},\\
&\;\eta\in\Gamma_{n-1}\text{ with }\rank(\xi)<\rank(\eta)\Big\}
\end{align*}
and, for $\gamma = (n,\xi,E,\vartheta,\eta)\in\Delta_n^\mathrm{b}$, put $\xi(\gamma) = \xi$, $E(\gamma) = E$, $\lambda(\gamma) = \lambda$, $\eta(\gamma) = \vartheta$, and $c_\gamma^* = e_\xi^* + \vartheta e_\eta^*\circ P^{(n-1)}_E$.
\end{example}

In this example, each element $\gamma$ of $\Gamma$ is a tuple, the first coordinate of which is always $\rank(\gamma)$ and the remaining coordinates encode the form of the extension functional $c_\gamma^*$ and the partial functions $\xi = \xi(\gamma)$, $E = E(\gamma) $, $\lambda = \lambda(\gamma)$, and $\eta = \eta(\gamma)$ which are either defined or undefined as required. Here, $\lambda(\gamma) = \vartheta$ whenever it is defined. The boundedness assumption is satisfied by design for $C = 1/(1-2\vartheta)$; therefore, this algorithm defines a Bourgain-Delbaen-$C$-$\mathscr{L}_\infty$-space. 

\begin{remark}
We do not assume that all Bourgain-Delbaen-$\mathscr{L}_\infty$-spaces satisfy Assumption \ref{BD construction assumption}; there are several existing variations of it and it would be cumbersome to state and handle a version that covers them all. Furthermore, the basic structural properties follow from Definition \ref{BD definition}.
\end{remark}

\addtocontents{toc}{\SkipTocEntry}
\subsection*{The structure of general Bourgain-Delbaen-$\mathscr{L}_\infty$-spaces}
\label{structures in BD}
We discuss some structures that are available in every Bourgain-Delbaen-$\mathscr{L}_\infty$-space, namely a finite dimensional Schauder decomposition and a Schauder basis.

\begin{remark}
\label{extensions defined on big space}
Let $\mathfrak{X}_{(\Gamma_n,i_n)}$ be a Bourgain-Delbaen-$\mathscr{L}_\infty$-space. Formally, for each $\gamma\in\Gamma$, the extension functional $c_\gamma^*$ is only defined when $\rank(\gamma) = n\geq 1$ and its domain is the space $\ell_\infty(\Gamma_{n-1})$. We can naturally extend this to define for each $\gamma\in\Gamma$ an extension functional $c_\gamma^*:\ell_\infty(\Gamma)\to\mathbb{C}$ as follows.
\begin{enumerate}[label=(\arabic*)]
    
    \item If $\rank(\gamma) = 0$ then $c_\gamma^* = 0$; this is not a true extension functional.

\item If $\rank(\gamma) = n\geq 1$ then $c_\gamma^*(x) = c_\gamma^*(x|_{\Gamma_{n-1}})$.

\end{enumerate}
\end{remark}

\begin{notation}\label{notation d gamma}
Let $\mathfrak{X}_{(\Gamma_n,i_n)}$ be a Bourgain-Delbaen-$\mathscr{L}_\infty$-space and $\gamma\in\Gamma$.
\begin{enumerate}[label=(\arabic*)]

\item Denote $d_\gamma^* = e_\gamma^* - c_\gamma^*:\ell_\infty(\Gamma)\to\mathbb{C}$ and,

\item if $\rank(\gamma) = n$, denote $d_\gamma = i_n(e_\gamma)$.

\end{enumerate}
\end{notation}

Note that, by definition, for all $\gamma\in\Gamma$, $\|d_\gamma\|_{\ell_\infty(\Gamma)}\leq C$.

\begin{proposition}[{\cite[Proposition 2.17]{argyros:gasparis:motakis:2016}}]
\label{remark general BD spaces functionals}
Let $\mathfrak{X}_{(\Gamma_n,i_n)}$ be a Bourgain-Delbaen-$\mathscr{L}_\infty$-space.
\begin{enumerate}[label=(\roman*),leftmargin=21pt]

\item The sequence $(d_\gamma,d_\gamma^*)_{\gamma\in\Gamma}$ forms a biorthogonal system in $\ell_\infty(\Gamma)\times\ell_\infty(\Gamma)^*$. 

\item In $\ell_\infty(\Gamma)^*$, for every $n\in\mathbb{N}_0$, $\langle\{d^*_\gamma:\gamma\in\Gamma_n\}\rangle = \langle\{e^*_\gamma:\gamma\in\Gamma_n\}\rangle$ and, in particular, $\langle\{d^*_\gamma:\gamma\in\Gamma\}\rangle = \langle\{e^*_\gamma:\gamma\in\Gamma\}\rangle$.

\end{enumerate}
\end{proposition}

\begin{proposition}[{\cite[Remark 2.6 and Proposition 2.13]{argyros:gasparis:motakis:2016}}]
\label{dual equivalence}
Let $\mathfrak{X}_{(\Gamma_n,i_n)}$ be a Bourgain-Delbaen-$C$-$\mathscr{L}_\infty$-space.
\begin{enumerate}[label=(\roman*),leftmargin=21pt]

\item For each $n\in\mathbb{N}_0$, $(d_\gamma)_{\gamma\in\Delta_n}$ is $C$-equivalent to the unit vector basis of $\ell_\infty(\Delta_n)$.

\item For each $n\in\N_0$ and scalar coefficients $(\lambda_\gamma)_{\gamma\in\Gamma_n}$,
\[\frac{1}{C}\Big\|\sum_{\gamma\in\Gamma_n}\lambda_\gamma e^*_\gamma\Big\|_{\ell_\infty(\Gamma)^*}\leq \Big\|\sum_{\gamma\in\Gamma_n}\lambda_\gamma e^*_\gamma\Big\|_{\mathfrak{X}_{(\Gamma_n,i_n)}^*} \leq \Big\|\sum_{\gamma\in\Gamma_n}\lambda_\gamma e^*_\gamma\Big\|_{\ell_\infty(\Gamma)^*}.\]
That is, linear combinations of $(e^*_\gamma)_{\gamma\in\Gamma}$, and, thus, also of $(d^*_\gamma)_{\gamma\in\Gamma}$, restricted on the subspace $\mathfrak{X}_{(\Gamma_n,i_n)}$, have equivalent norm estimates to those obtained on their entire domain $\ell_\infty(\Gamma)$.

\end{enumerate}
\end{proposition}

\begin{remark}[{\cite[Remark 2.10]{argyros:gasparis:motakis:2016}}]
\label{basis enumeration remark}
We do not use this explicitly, but we mention for context that any enumeration $(\gamma_i)_{i\in\mathbb{N}}$ of $\Gamma$ with the property $i\leq j$ implies $\rank(\gamma_i)\leq\rank(\gamma_j)$ makes $(d_{\gamma})_{\gamma\in\Gamma}$ a Schauder basis of $\mathfrak{X}_{(\Gamma_n,i_n)}$ with associated biorthogonal sequence $(d_\gamma^*)_{\gamma\in\Gamma}$ (viewed in $\mathfrak{X}_{(\Gamma_n,i_n)}^*$). The most commonly used coordinate system in $\mathfrak{X}_{(\Gamma_n,i_n)}$ is the FDD introduced below.
\end{remark}

\begin{notation}
\label{basis notation}
Let $\mathfrak{X}_{(\Gamma_n,i_n)}$ be a Bourgain-Delbaen-$\mathscr{L}_\infty$-space.
\begin{enumerate}[label=(\arabic*),leftmargin=19pt]

\item\label{basis notation 1}
For $n\in\mathbb{N}_0$ denote by $P_n:\ell_\infty(\Gamma)\to\ell_\infty(\Gamma)$ the linear projection $i_n\circ r_{\Gamma_n}$, i.e., for $x\in\ell_\infty(\Gamma)$, $P_n(x) = i_n(x|_{\Gamma_n})$.

\item For an interval $E = [m,n]$ of $\mathbb{N}_0$ denote by  $P_E:\ell_\infty(\Gamma)\to \ell_\infty(\Gamma)$ the linear projection $P_n - P_{m-1}$ (where $P_{-1} = 0$).

\item For an interval $E =[m,\infty)$ of $\N_0$, denote $P_E = I - P_{m-1}$.

\end{enumerate}
\end{notation}

\begin{remark}
\label{drop finite notation}
Let $n\in\mathbb{N}_0$ and $E$ be an interval of $\{0,\ldots,n\}$. The similar choice of notation for $P^{(n)}_E$ and $P_E$ is deliberate; if $\gamma\in\Gamma_n$ and $x\in\ell_\infty(\Gamma)$, then $e_\gamma^*\circ P_E^{(n)}(x|_{\Gamma_n}) = e_\gamma^*\circ P_E(x)$. This simplifies the formulae of extension functionals. For example, if $n\in\N$, $\gamma\in\Delta_n$, and $c_\gamma^*:\ell_\infty(\Gamma_{n-1})\to\mathbb{C}$ satisfies a formula as in Assumption \ref{BD construction assumption}
\[c_\gamma^* = \lambda e_\eta^*\circ P^{(n-1)}_E\text{ or }c_\gamma^* = e_\xi^* + \lambda e_\eta^*\circ P^{(n-1)}_E\]
then the natural extension $c_\gamma^*:\ell_\infty(\Gamma)\to \mathbb{C}$ given by Remark \ref{extensions defined on big space} satisfies the same formula
\[c_\gamma^* = \lambda e_\eta^*\circ P_E\text{ or }c_\gamma^* = e_\xi^* + \lambda e_\eta^*\circ P_E.\]
\end{remark}

\begin{proposition}[{\cite[Remark 2.6, Proposition 2.8, and Remark 2.9]{argyros:gasparis:motakis:2016}}]
\label{remark general BD spaces}
Let $\mathfrak{X}_{(\Gamma_n,i_n)}$ be a Bourgain-Delbaen-$\mathscr{L}_\infty$-space.
\begin{enumerate}[label=(\roman*),leftmargin=23pt]

\item\label{remark general BD spaces basis of FDD component}  For $n\in\mathbb{N}$, $(d_\gamma)_{\gamma\in\Delta_n}$ is $C$-equivalent to the unit vector basis of $\ell_\infty(\Delta_n)$ and, thus, $Z_n = \langle\{d_\gamma:\gamma\in\Delta_n\}\rangle$ is $C$-isomorphic to $\ell_\infty(\Delta_n)$.

\item The projections $(P_n)_{n=0}^\infty$, restricted on $\mathfrak{X}_{(\Gamma_n,i_n)}$, define an FDD of $\mathfrak{X}_{(\Gamma_n,i_n)}$ as follows.
\begin{enumerate}[label=(\greek*)]

\item For every $n\in\N_0$, $P_n(\mathfrak{X}_{(\Gamma_n,i_n)}) = i_n(\ell_\infty(\Gamma_n)) = Y_n$ (see Definition \ref{BD definition}) and,

\item for every $n, m\in\mathbb{N}_0$, $P_nP_m = P_{n\wedge m}$.

\end{enumerate}
Therefore, for each $n\in\N_0$, $Z_n = (P_n-P_{n-1})(\mathfrak{X}_{(\Gamma_n,i_n)})$ and $(Z_n)_{n=0}^\infty$ forms a FDD of $\mathfrak{X}_{(\Gamma_n,i_n)}$.

\item\label{remark general BD spaces norm of d-gamma-star}  For any finite interval $E$ of $\mathbb{N}_0$, $P_E = \sum_{n\in E}\sum_{\gamma\in\Delta_n}d_\gamma^*\otimes d_\gamma$. In particular, for every $\gamma\in\Gamma$, $d_\gamma^* = e^*_\gamma\circ P_{\{\rank(\gamma)\}} = e^*_\gamma\circ P_{[\rank(\gamma),\infty)}$.

\end{enumerate}
\end{proposition}

\begin{notation}\label{notation support delta}
For a Bourgain-Delbaen-$\mathscr{L}_\infty$-space $\mathfrak{X}_{(\Gamma_n,i_n)}$ and $x\in\mathfrak{X}_{(\Gamma_n,i_n)}$ we let
\[\supp(x)=\{n\in\N_0: d^*_\gamma(x)\neq 0 \text{ for some }\gamma\in \Delta_n\},\] 
that is, $\supp (x)$ denotes the support of $x$ with respect to the FDD $(Z_n)_{n=0}^\infty$. 
We also let $\mathrm{range}(x)$ be the smallest interval  in $\N_0$ containing $\supp(x)$. A sequence $(x_n)_n$ in $\mathfrak{X}_{(\Gamma_n,i_n)}$ is called block, if it is block with respect to the FDD $(Z_n)_{n=0}^\infty$, i.e. if $\supp (x_n)<\supp(x_{n+1})$ for all $n$. Analogously, for any $b^*\in \mathfrak X_{(\Gamma_n,i_n)}^*$, we let 
\[\supp(b^*)=\{n\in\N_0: b^*(d_\gamma)\neq 0 \text{ for some }\gamma\in \Delta_n\},\]
i.e. $\supp(b^*)$ is the support of the functional $b^*$ with respect to the FDD $(Z_n)_{n=0}^\infty$ (see Section \ref{section roadmap}). 
\end{notation}

\addtocontents{toc}{\SkipTocEntry}
\subsection*{The evaluation analysis of nodes}\label{section evaluation analysis}

One of the most important concepts in studying Argyros-Haydon-type Bourgain-Delbaen spaces is the evaluation analysis of functionals $e_\gamma^*$ that first appeared in \cite{argyros:haydon:2011}. This is a formula that draws a connection between such spaces and mixed-Tsirelson (and similar) spaces, on which certain Bourgain-Delbaen spaces are modelled.

\begin{proposition}\label{evaluation analysis}
Let $(\Gamma_n)_{n=0}^\infty$, $(i_{n-1,n})_{n=1}^\infty$ be a linear extension scheme that satisfies Assumption \ref{BD construction assumption}, for some $0<\vartheta<1/2$. Then, for any $\gamma\in\Gamma$, either $e_\gamma^* = d_\gamma^*$ or
\begin{equation}
\label{evaluation analysis equation}
e_\gamma^* = \sum_{r=1}^a d_{\xi_r}^* + \sum_{r=1}^a \lambda_r e_{\eta_r}^*\circ  P_{E_r},
\end{equation}
where $\mathrm{age}(\gamma) = a$, $\vec\xi(\gamma) = (\xi_1,\ldots,\xi_a)$, $\vec E(\gamma) = (E_1,\ldots,E_a)$, $\vec\lambda(\gamma) = (\lambda_1,\ldots,\lambda_a)$, and $\vec\eta(\gamma) = (\eta_1,\ldots,\eta_a)$ were defined in Remark \ref{gamma tuples}. The expression \eqref{evaluation analysis equation} is called the evaluation analysis of $\gamma$. Furthermore,
\begin{enumerate}[label=(\roman*),leftmargin=23pt]

\item $\xi_1,\ldots,\xi_a\in\Gamma$ with $\rank(\xi_1)<\cdots<\rank(\xi_a)$ and $\xi_a = \gamma$,

\item $E_1$ is an interval of $[0,\rank(\xi_1))$ and $E_r$ is an interval of $(\rank(\xi_{r-1}),\rank(\xi_r))$, $r=2,\ldots,a$,

\item $\lambda_r$ is a scalar with $|\lambda_r|\leq\vartheta$, for $r=1,\ldots,a$, and

\item $\eta_r\in\Gamma$ with $\rank(\eta_r) < \rank(\xi_r)$, for $r=1,\ldots,a$.

\end{enumerate} 
\end{proposition}

\begin{remark}
\label{evaluation analysis analysis}
The evaluation analysis is the main tool used to compute the norm of linear combinations of block vectors in Bourgain-Delbaen-$\mathscr{L}_\infty$-spaces satisfying Assumption \ref{BD construction assumption}. The main components of this analysis are the following.
\begin{enumerate}[label=(\arabic*),leftmargin=21pt]

\item The sequence of coefficients $\vec \lambda(\gamma) = (\lambda_1,\ldots,\lambda_a)$. Its shape impacts the resulting space greatly. For example, in the Argyros-Haydon space with the scalar-plus-compact property (\cite{argyros:haydon:2011}), they are modelled after a mixed-Tsirelson space, a Schlumprecht-type space (see \cite{schlumprecht:1991}).

\item The relative positions of $\vec\xi(\gamma) = (\xi_1,\ldots,\xi_a)$ and $\vec\eta(\gamma) = (\eta_1,\ldots,\eta_a)$ in $\Gamma$. For example, in the Argyros-Haydon space, a Maurey-Rosenthal-type coding function (see \cite{maurey:rosenthal:1977}) sometimes imposes restrictions on the members of $\vec \eta(\gamma)$. In this paper, we introduce a type of Bourgain-Delbaen scheme called a bonding, where a different kind of control of on the relative position of $\vec\xi(\gamma)$ and $\vec\eta(\gamma)$ yields sequences of projections automatically satisfying Proposition \ref{general condition on projections} \ref{general condition on projections a}, \ref{general condition on projections b}, and \ref{general condition on projections c}  (see Sections \ref{section self-determined} and \ref{bonding section}.)

\end{enumerate}
To some extent, the restrictions imposed on the intervals $E(\gamma) = (E_1,\ldots,E_a)$ also play a role.
\end{remark}

\begin{proof}[Proof of Proposition \ref{evaluation analysis}]
If $\mathrm{rank}(\gamma)=0$, then $\gamma$ must satisfy Assumption \ref{BD construction assumption} \ref{BD construction assumption 0}. Hence, $c_\gamma^* = 0$, i.e., $e^*_\gamma=d^*_\gamma$ and the conclusion holds. For $\gamma$ satisfying Assumption \ref{BD construction assumption} \ref{BD construction assumption 1} or \ref{BD construction assumption 2}, we perform an induction on $\mathrm{age}(\gamma)$, i.e., the $a$ for which $\vec\xi(\gamma) = (\xi_1,\ldots,\xi_a)$, $\vec E(\gamma) = (E_1,\ldots,E_a)$, $\vec\lambda(\gamma) = (\lambda_1,\ldots,\lambda_a)$, and $\vec\eta(\gamma) = (\eta_1,\ldots,\eta_a)$ are defined.

If $a = 1$, then $\xi_1 = \gamma$, $\xi(\gamma)$ is undefined, $E_1 = E(\gamma)$, $\lambda_1 = \lambda(\gamma)$, and $\eta_1 = \eta(\gamma)$. Therefore, $\gamma$ satisfies Assumption \ref{BD construction assumption} \ref{BD construction assumption 1} and, thus, $c_\gamma^* = \lambda_1 e^*_{\gamma_1}\circ P_{E_1}$, where $E_1$ is an interval of $[0,\mathrm{rank}(\xi_1))$, $\lambda_1$ is a scalar with $|\lambda_1|\leq\vartheta$, and $\eta_1\in\Gamma$ with $\mathrm{rank}(\eta_1)<\mathrm{rank}(\xi_1)$. We obtain
\[e_\gamma^* = d_\gamma^* + c_\gamma^* = d_{\xi_1}^* + \lambda_1e_{\eta_1}^*\circ P_{E_1}.\]

For the inductive step, take $\gamma$ with $\mathrm{age}(\gamma) = a > 1$ and assume that the conclusion holds for any $\gamma'$ with $\mathrm{age}(\gamma') < a$. Write $\vec\xi(\gamma) = (\xi_1,\ldots,\xi_{a-1},\xi_a)$, $\vec E(\gamma) = (E_1,\ldots,E_{a-1},E_a)$, $\vec\lambda(\gamma) = (\lambda_1,\ldots,\lambda_{a-1},\lambda_a)$, and $\vec\eta(\gamma) = (\eta_1,\ldots,\eta_{a-1},\eta_a)$. Because $ \xi_{\gamma} = \xi(\xi_a) = \xi_{a-1}$, $E(\gamma) = E_a$, $\lambda(\gamma) = \lambda_a$, and $\eta(\lambda) = \eta_a$, $\gamma$ satisfies Assumption \ref{BD construction assumption} \ref{BD construction assumption 2} and, thus, $\mathrm{rank}(\xi_{a-1})<\mathrm{rank}(\gamma)$, $E_a$ is an interval of $(\mathrm{rank}(\xi_{a-1}),\mathrm{rank}(\xi_a))$, $\lambda_a$ is a scalar with $|\lambda_a|\leq\vartheta$, $\eta_a\in\Gamma$ with $\mathrm{rank}(\eta_a) < \mathrm{rank}(\xi_a)$, and
\[e_\gamma^* = d_\gamma^* + c_\gamma^* = d_{\xi_a}^* + e_{\xi_{a-1}}^* + \lambda_ae_{\eta_a}^*\circ P_{E_a}.\]
Because $\vec\xi (\xi_{a-1}) = (\xi_1,\ldots,\xi_{a-1})$, the inductive hypothesis yields the conclusion. 
\end{proof}

\begin{example}
In the extension scheme  of Example \ref{schur bd space} each $\gamma\in\Gamma$ has evaluation analysis
\[ e_\gamma^* = \sum_{r=1}^a d_{\xi_r}^* + \vartheta\sum_{r=1}^a  e_{\eta_r}^*\circ  P_{E_r}\]
with the associated tuples $\vec\xi(\gamma)$, $\vec E(\gamma)$, $\vec\lambda(\gamma)$, and $\vec\eta(\gamma)$ satisfying  $(i)$-$(iv)$ of Proposition \ref{evaluation analysis}. 
\end{example}

\begin{remark}
\label{remark basis vector norm}
In a Bourgain-Delbaen-C-$\mathscr{L}_\infty$-space $\mathfrak{X}_{(\Gamma_n,i_n)}$, for $\xi\in\Gamma$, we have $\|d_\xi\|_{\ell_\infty(\Gamma)}\leq C$. In spaces satisfying Assumption \ref{BD construction assumption}, a simple application of the evaluation analysis improves this to $\|d_\xi\|_{\ell_\infty(\Gamma)} = 1$. Indeed, for $\gamma\in\Gamma$, either $e_\gamma^*(d_\xi) = 0$ or $e_\gamma^*(d_\xi) = d_{\xi_r}^*(d_\xi)$ or $e_\gamma^*(d_\xi) = \lambda_re_{\eta_r}^*(d_\xi)$, for some $1\leq r\leq a$. An induction on the age of $\gamma$ yields the promised bound. Similarly, and more generally, for every $n\in\N$, $(d_\xi)_{\xi\in\Delta_n}$ is isometrically equivalent to the unit vector basis of $\ell_\infty(\Delta_n)$.
\end{remark}

\subsection{Self-determined subsets of $\Gamma$}
\label{section self-determined}
We recall a method of defining Bourgain-Delbaen $\mathscr{L}_\infty$-spaces with more elaborate properties, as quotients of $\mathscr{L}_\infty$-spaces of more regular type. Athough it was first formalized in  {\cite{argyros:motakis:2019}} it had been used implicitly before, e.g., by Tarbard in {\cite{tarbard:2013}}. This method has also been used in {\cite{manoussakis:pelczar-barwacz:swietek:2017} and \cite{motakis:2024}}. In all aforementioned applications, this concept was used for notational convenience. In the present paper, self-determined sets play a more special role as we also exploit their properties to define in a Bourgain-Delbaen space complemented subspaces that satisfy the assumptions of Proposition \ref{general condition on projections} (see Section \ref{bonding section}).

The following definition originates from {\cite[Definition 1.4]{argyros:motakis:2019}}. The equivalence to the statements given below follows from {\cite[Proposition 1.5]{argyros:motakis:2019}}.

\begin{definition}
\label{sd definition}
Let $\mathfrak{X}_{(\Gamma_n,i_n)}$ be a Bourgain-Delbaen space and $\Gamma'\subset\Gamma$. For $n\in\mathbb{N}_0$ denote $\Gamma'_n = \Gamma\cap\Gamma_n$. We call $\Gamma'$ a self-determined subset of $\Gamma$ provided  one of the following equivalent conditions is satisfied.
\begin{enumerate}[label=(\alph*)]

\item\label{sd definition2} $\langle\{d_\gamma^*:\gamma\in\Gamma'\}\rangle = \langle\{e_\gamma^*:\gamma\in\Gamma'\}\rangle$.

\item\label{sd definition complement} For all $\eta\in\Gamma'$ and $\gamma\in\Gamma\setminus\Gamma'$, $e_\eta^*(d_\gamma) = 0$.

\item\label{sd definition algorithm}  For every $n\in\mathbb{N}$ and $\gamma\in\Gamma'_n$ the extension functional $c_\gamma^*:\ell_\infty(\Gamma_{n-1})\to\mathbb{C}$ satisfies
\[c_\gamma^*\in\langle\big\{e_\xi^*\circ P^{(n-1)}_{[k,m]}:\xi\in\Gamma_{n-1}'\text{ and }0\leq m\leq k<n\big\}\rangle,\]
following the convention $\langle\emptyset\rangle = \{0\}$.
\end{enumerate}
\end{definition}
In other words, a self determined subset of $\Gamma$ is a $\Gamma'\subset\Gamma$ such that for every $\gamma\in\Gamma'$ the extension functional $c_\gamma^*$ satisfies a formula that only uses members of $\Gamma'$ of smaller rank. 

\begin{notation}
Let $\mathfrak{X}_{(\Gamma_n,i_n)}$ be a Bourgain-Delbaen-$\mathscr{L}_\infty$-space and $\Gamma'\subset\Gamma$ be self-determined.
\begin{enumerate}[label=(\arabic*),leftmargin=19pt]

\item For $n\in\mathbb{N}_0$ denote $\Gamma_n' = \Gamma_n\cap\Gamma'$ and $\Delta_n' = \Delta_n\cap \Gamma'$.

\item For $n\geq 1$ define $c_\gamma^{\prime*}:\ell_\infty(\Gamma_{n-1}')\to\mathbb{C}$ given by $c_\gamma^{\prime*}(x) = c_\gamma^*(x)$, where $x\in\ell_\infty(\Gamma_{n-1}')$ is naturally identified with a vector in $\ell_\infty(\Gamma_{n-1})$ by putting zero in the missing coordinates.

\end{enumerate}
We denote by $(\Gamma_n')_{n=0}^\infty$, $(i_{n-1,n}')_{n=1}^\infty$ the corresponding linear extension scheme.
\end{notation}

\begin{proposition} \label{general s-d bd space} 
{\cite[Proposition 1.9 and Proposition 1.12]{argyros:motakis:2019}} Let $\mathfrak{X}_{(\Gamma_n,i_n)}$ be a Bourgain-Delbaen-$C$-$\mathscr{L}_\infty$-space and $\Gamma'\subset\Gamma$ be self-determined.
\begin{enumerate}[label=(\roman*),leftmargin=19pt]

\item\label{general s-d bd space1}
The extension scheme $(\Gamma_n')_{n=0}^\infty$, $(i_{n-1,n}')_{n=1}^\infty$ satisfies the boundedness Assumption \ref{general BD extension assumption} with constant $C$, and therefore it defines a Bourgain-Delbaen-$C$-$\mathscr{L}_\infty$-subspace $\mathfrak{X}_{(\Gamma_n',i_n')}$ of $\ell_\infty(\Gamma')$.

\item\label{general quotient}
If we denote by $r_{\Gamma'}:\ell_\infty(\Gamma)\to\ell_\infty(\Gamma')$ the natural restriction map, then $r_{\Gamma'}(\mathfrak{X}_{(\Gamma_n,i_n)}) = \mathfrak{X}_{(\Gamma_n',i_n')}$. In particular, if $R = r_{\Gamma'}|_{\mathfrak{X}_{(\Gamma_n,i_n)}}$, then $R:\mathfrak{X}_{(\Gamma_n,i_n)}\to\mathfrak{X}_{(\Gamma_n',i_n')}$ is an onto map of norm at most one with  $\mathrm{ker}(R) = [\{d_\gamma:\gamma\in\Gamma\setminus\Gamma'\}]$, which is also a $C$-$\mathscr{L}_\infty$-space.

\end{enumerate}
 \end{proposition}

The most important feature of a self-determined subset $\Gamma'$ of $\Gamma$ is that $\mathfrak{X}_{(\Gamma_n',i_n')}$ inherits the formulae of extension functionals from $\mathfrak{X}_{(\Gamma_n,i_n)}$. To elaborate on this let us first denote by $i_{m,n}'$, $c_\gamma^{\prime*}$, $P_E^{\prime(n)}$, $P_E'$, $d_\gamma^{\prime*}$ the objects from Sections \ref{Bourgain-Delbaen general section} and \ref{structures in BD} associated to $\mathfrak{X}_{(\Gamma_n',i_n')}$. By naturally identifying both $\mathfrak{X}_{(\Gamma_n,i_n)}$ and $\mathfrak{X}_{(\Gamma_n',i_n')}$ with subspaces of $\ell_\infty(\Gamma)$, we may use the common notation $e_\gamma^*$ for the coordinate evaluation functionals for both $\mathfrak{X}_{(\Gamma_n,i_n)}$ and $\mathfrak{X}_{(\Gamma_n',i_n')}$.

\begin{fact}[Preservation of formulae, {\cite[Proposition 1.13]{argyros:motakis:2019}}] 
\label{preservation of formulae}
Let $\Gamma'\subset\Gamma$ be self-determined and let $\gamma\in\Gamma_n'$, and thus, by Definition \ref{sd definition} \ref{sd definition algorithm}, the linear extension functional $c_\gamma^*:\ell_\infty(\Gamma_n)\to\mathbb{C}$ satisfies a formula    
\[c_\gamma^* = \sum_{\xi\in \Gamma_{n-1}'}\sum_{0\leq k\leq m<n}\lambda_{\xi,E}e_\xi^*\circ P^{(n-1)}_{[k,m]}.\]
Then, the linear extension functional $c_\gamma^{\prime*}:\ell_\infty(\Gamma_n')\to\mathbb{C}$ satisfies the same formula 
\[c_\gamma^{\prime*} = \sum_{\xi\in \Gamma_{n-1}'}\sum_{0\leq k\leq m<n}\lambda_{\xi,E}e_\xi^*\circ P^{\prime(n-1)}_{[k,m]}.\]
\end{fact}

The convenience offered by self-determined sets is not easily demonstrated without introducing some complexity. For the purpose of conceptualization, let us briefly mention a toy example. Here, the partial functions $\xi(\cdot)$, $E(\cdot)$, $\lambda(\cdot)$, and $\eta(\cdot)$ are defined in the obvious way, as in Example \ref{schur bd space}, using the coordinates of each $\gamma$.

\begin{example}
\label{toy self-determined}
Fix $0<\vartheta<1/2$ and a sequence $(\vartheta_j)_{j=1}^\infty$ in $(0,\vartheta]$ such that $\vartheta_1 =\vartheta$. Define $\Delta_0 = \{(0,0)\}$ and, assuming that $\Delta_0,\Delta_1,\ldots,\Delta_{n-1}$ and $(c_\gamma^*)_{\gamma\in\Delta_m}$, $1\leq m\leq n-1$ have been defined, define $\Delta_n = \Delta_n^{0}\cup\Delta_n^\mathrm{a}\cup\Delta_n^\mathrm{b}$ and $(c_\gamma^*)_{\gamma\in\Delta_n}$ as follows. Let $\Delta_n^{0} = \{(n)\}$ and for $\gamma\in\Delta_n^{0}$, $c_\gamma^* = 0$. Let
\[\Delta_n^\mathrm{a} = \Big\{(n,E,\vartheta_j,\eta): E\subset\{0,\ldots,n-1\}\text{ interval}, 1\leq j\leq n, \eta\in\Gamma_{n-1}\Big\}\]
and for $\gamma =  (n,E,\vartheta_j,\eta)\in\Delta_{n}^\mathrm{a}$ put $c_\gamma^* = \vartheta_j e_\eta^*\circ P^{(n-1)}_E$. Let
\begin{align*}
    \Delta_n^\mathrm{b} = \Big\{(n,\xi,E,\vartheta_j,\eta):&\;\xi\in\Gamma_{n-2}\cap\mathrm{dom}(E(\cdot)),E\subset\{\rank(\xi)+1,\ldots,n-1\}\text{ interval},\\
    & \; \vartheta_j = \lambda(\xi),\eta\in\Gamma_{n-1}\text{ with }\rank(\xi)<\rank(\eta)\Big\}
\end{align*}
and for $\gamma = (n,\xi,E,\vartheta_j,\eta)\in\Delta_n^\mathrm{b}$ put $c_\gamma^* = e_\xi^* + \vartheta_j e_\eta^*\circ P^{(n-1)}_E$.

Then, $\mathfrak{X}_{(\Gamma_n,i_n)}$ is a Bourgain-Delbaen-$C$-$\mathscr{L}_\infty$-space, for $C = 1/(1-2\vartheta)$. For each $\gamma\in\Gamma$, either $e_\gamma^* = d_\gamma^*$ or  it has evaluation analysis
\[ e_\gamma^* = \sum_{r=1}^a d_{\xi_r}^* + \vartheta_j\sum_{r=1}^a e_{\eta_r}^*\circ  P_{E_r},\]
where, $\vartheta_j = \lambda(\gamma)$.

We choose a self determined subset of $\Gamma$ as follows. Fix $n_0\in\mathbb{N}_0$, $\Lambda\subset\mathbb{N}$, and a function $\sigma:\Gamma\to\Lambda$. Let $\Delta_0'=\cdots=\Delta_{n_0-1}'=\emptyset$ and $\Delta_{n_0}' = \Delta_{n_0}^0 = \{(n_0)\}$. For $n > n_0$ define $\Delta_n' = \Delta_n^0\cup\Delta_n^{\mathrm{a}\prime}\cup\Delta_n^{\mathrm{b}\prime}$ where
\begin{align*}
\Delta_n^{\mathrm{a}\prime} &= \Big\{(n,E,\vartheta_j,\eta)\in\Delta_n^\mathrm{a}: \eta\in\Gamma_{n-1}',j\in\Lambda\Big\}\text{ and}\\
\Delta_n^{\mathrm{b}\prime} &= \Big\{(n,\xi,E,\vartheta_j,\eta)\in\Delta_n^\mathrm{b}: \xi,\eta\in\Gamma_{n-1}',\lambda(\eta) = \vartheta_{\sigma(\xi)}\Big\}.
\end{align*}
\end{example}
 
The set $\Gamma'$ in the above example is self-determined by design. Note that the first non-empty $\Delta_{n_0}^\prime$ only contains $\gamma$ with $c_\gamma^* = 0$, otherwise Definition \ref{sd definition} \ref{sd definition algorithm} would be violated. For every $\gamma\in\Gamma'$, either $e_\gamma^* = d_\gamma^*$ or it has evaluation analysis
\[e_\gamma^* = \sum_{r=1}^a d_{\xi_r}^* + \vartheta_j\sum_{r=1}^a e_{\eta_r}^*\circ  P_{E_r},\]
where $j\in\Lambda$, $\vartheta_j = \lambda(\gamma)$ and, for $2\leq r\leq a$, $\lambda(\eta_r) = \vartheta_{\sigma(\xi_{r-1})}$. The $\xi_r$ and $\eta_r$ are in $\Gamma'$ and, furthermore, have rank at least $n_0$. If, in particular, $\sigma$ is an injection then, for $2\leq r\leq a$, $\lambda(\eta_r)$ determines $\xi_{r-1}$. By the preservation of formulae, the properties of the evaluation analysis are inherited by the space $\mathfrak{X}_{(\Gamma_n',i_n')}$. Intuitively, the space $\mathfrak{X}_{(\Gamma_n',i_n')}$ is the result of imposing additional restrictions to the algorithmic process that defines the extension coordinates and functionals of $\mathfrak{X}_{(\Gamma_n,i_n)}$. In fact, $\mathfrak{X}_{(\Gamma_n',i_n')}$ is isomorphic to a quotient of $\mathfrak{X}_{(\Gamma_n,i_n)}$.

\section{Bonding of a sequence of Bourgain-Delbaen-$\mathscr{L}_\infty$-spaces}
\label{bonding section}
This is one of the most important sections of this paper in terms of explaining the ingredients and how they are used to prove the properties of the space $\mathfrak{X}_U$. We introduce a  construction scheme called a bonding of a sequence of Bourgain-Delbaen-$\mathscr{L}_\infty$-spaces. This scheme produces a Bourgain-Delbaen-$\mathscr{L}_\infty$-space $\mathfrak{X}_{(\Gamma_n,i_n)}$ satisfying Assumption \ref{BD construction assumption} and a sequence of disjoint self-determined subsets $(\Gamma^s)_{s=1}^\infty$ of $\Gamma$ such that, for $s\in\N$, the coordinate projection $I_s$ on $\mathfrak{X}_s = [(d_\gamma)_{\gamma\in\Gamma^s}]$ is bounded and $\mathfrak{X}_s\simeq \mathfrak{X}_{(\Gamma^s_n,i^s_n)}$. The scheme allows complete freedom on the internal structure of the subsets $\Gamma^s$, $s\in\N$, but imposes restrictions on the structure of $\Gamma^0 = \Gamma\setminus\cup_{s=1}^\infty\Gamma^s$. As a result, Proposition \ref{general condition on projections} \ref{general condition on projections a}, \ref{general condition on projections b}, and \ref{general condition on projections c} are satisfied automatically. In particular, if $\mathfrak{X}_0 = [(d_\gamma)_{\gamma\in\Gamma^0}]$, then $\mathfrak{X}_{(\Gamma_n,i_n)}/\mathfrak{X}_0\simeq (\oplus_{s=1}^\infty\mathfrak{X}_s)_{c_0}$  (see Theorem \ref{bonding general statement}). 
In Section \ref{section Linear combinations of coordinate projections}, we formulate additional restrictions on a bonding scheme that, when imposed on $\vec\lambda(\gamma)$, for $\gamma\in\Gamma^0$, the remaining conditions of Proposition \ref{general condition on projections} are satisfied. 
Crucially, the bonding scheme allows the incorporation of Argyros-Haydon-type structure and the space $\mathfrak{X}_U$ defined later is also a bonding.

\subsection{The general construction of a bonding}

We first observe an obvious condition that when imposed on a self-determined subset $\Gamma'$ of $\Gamma$, the coordinate projection $I_{\Gamma'}$ onto $\mathfrak{X}' = [(d_\gamma)_{\gamma\in\Gamma'}]$ is bounded.

\begin{notation}
\label{coordinate pro}
Let $\mathfrak{X}_{(\Gamma_n,i_n)}$ be a Bourgain-Delbaen-$\mathscr{L}_\infty$-space and $\Gamma'\subset\Gamma$. We denote by $I_{\Gamma'}:\langle\{d_\gamma:\gamma\in\Gamma\}\rangle\to\langle\{d_\gamma:\gamma\in\Gamma\}\rangle$ the linear projection given by
\[
I_{\Gamma'}(d_\gamma) = \left\{
	\begin{array}{ll}
		d_\gamma  & \mbox{if } \gamma\in\Gamma' \\
		0 & \mbox{if } \gamma\in\Gamma\setminus\Gamma'.
	\end{array}
\right.
\]
We call $I_{\Gamma'}$ the coordinate projection onto $\Gamma'$.
\end{notation}

\begin{remark}
\label{how to defined composition with unbounded}
For $\Gamma'\subset\Gamma$, $I_{\Gamma'}$ need not always extend to a bounded linear operator on $\mathfrak{X}_{(\Gamma_n,i_n)}$ however, for $\gamma\in\Gamma$ with $\mathrm{rank}(\gamma) = n$, we may always define the bounded linear functional $e_\gamma^*\circ I_{\Gamma'}:\ell_\infty(\Gamma)\to\mathbb{C}$ with
\[e_\gamma^*\circ I_{\Gamma'} \coloneqq \sum_{\eta\in\Gamma'}e_\gamma^*(d_\eta)d_\eta^* = \sum_{\eta\in\Gamma_n'}e_\gamma^*(d_\eta)d_\eta^* .\]
The notation $e_\gamma^*\circ I_{\Gamma'}$ is appropriate; for every $x\in\langle\{d_\eta:\eta\in\Gamma\}\rangle$ we have $e_\gamma^*(I_{\Gamma'}x) = (e_\gamma^*\circ I_{\Gamma'})(x)$. Then, clearly, $I_{\Gamma'}$  admits a bounded extension to $\mathfrak{X}_{(\Gamma_n,i_n)}$ if an only if $\sup_{\gamma\in\Gamma}\|e_\gamma^*\circ I_{\Gamma'}\|_{\mathfrak{X}_{(\Gamma_n,i_n)}^*} = A <\infty$ and then $\|I_{\Gamma'}\|_{\mathcal{L}(\mathfrak{X}_{(\Gamma_n,i_n)})} =  A$.
\end{remark}

\begin{notation}
\label{pre-bonding notation}
Let $\mathfrak{X}_{(\Gamma_n,i_n)}$ be a Bourgain-Delbaen-$\mathscr{L}_\infty$-space and $(\Gamma^s)_{s=1}^\infty$ be a sequence of pairwise disjoint subsets of $\Gamma$.
\begin{enumerate}

\item For $s\in\N$ denote $I_s = I_{\Gamma^s}$ and $\mathfrak{X}_s = [(d_\gamma)_{\gamma\in\Gamma^s}]$.

\item Denote $\Gamma^0 = \Gamma\setminus\cup_{s=1}^\infty\Gamma^s$ and $\mathfrak{X}_0 = [(d_\gamma)_{\gamma\in\Gamma^0}]$.

\item For a linear functional $f:\langle\{d_\gamma:\gamma\in\Gamma\}\rangle\to\mathbb{C}$ denote
\[\mathrm{supp}_\mathrm{h}(f) = \{s\in\N:f(d_\gamma) \neq 0,\text{ for some }\gamma\in\Gamma^s\}.\]

\end{enumerate}
\end{notation}

\begin{remark}
We call the set $\mathrm{supp}_\mathrm{h}(f)$ the horizontal support of $f$ and it measures the interaction of $f$ with the subspaces $\mathfrak{X}_s$, $s\in\N$, but not with $\mathfrak{X}_0$. In particular, for $\gamma\in\Gamma^0$, $\mathrm{supp}_\mathrm{h}(d_\gamma^*) = \emptyset$.
\end{remark}

\begin{definition}
\label{bonding definition}
Let $\mathfrak{X}_{(\Gamma_n,i_n)}$ be a Bourgain-Delbaen-$\mathscr{L}_\infty$-space and $(\Gamma^s)_{s=1}^\infty$ be a sequence of  disjoint subsets of $\Gamma$. We call $\mathfrak{X}_{(\Gamma_n,i_n)}$ a bonding of $(\mathfrak{X}_s)_{s=1}^\infty$ if the following conditions are met.
\begin{enumerate}[label=(\alph*),leftmargin=19pt]

\item\label{bonding definition1} Assumption \ref{BD construction assumption} is satisfied, for some $0<\vartheta\leq1/4$. In particular, by Proposition \ref{evaluation analysis}, for $\gamma\in\Gamma$, either $e_\gamma^* = d_\gamma^*$ or $e_\gamma^*$ admits an evaluation analysis
\[e_\gamma^* = \sum_{r=1}^ad_{\xi_r}^* + \sum_{r=1}^a\lambda_re_{\eta_r}^*\circ P_{E_r}.\]

\item\label{bonding definition2} For $s\in\N$ and $\gamma\in\Gamma^s$ either $e_\gamma^* = d_\gamma^*$  or in the evaluation analysis of $e_\gamma^*$, for $1\leq r\leq a$, $\xi_r,\eta_r\in\Gamma^s$. In particular, $\Gamma^s$ is a self-determined subset of $\Gamma$.

\item\label{bonding definition3} For $\gamma\in\Gamma^0$ either $e_\gamma^* = d_\gamma^*$ or in the evaluation analysis of $e_\gamma^*$, for $1\leq r\leq a$, $\xi_r\in\Gamma^0$ and, for $2\leq r\leq a$, $\max(\mathrm{supp}_\mathrm{h}(e_{\eta_{r-1}}^*\circ P_{E_{r-1}})) < \min(\mathrm{supp}_\mathrm{h}(e_{\eta_r}^*\circ P_{E_r}))$.
\end{enumerate}
\end{definition}

The most general result about bondings is the following.

\begin{theorem}
\label{bonding general statement}
Let $\mathfrak{X}_{(\Gamma_n,i_n)}$ be a Bourgain-Delbaen-$C$-$\mathscr{L}_\infty$-space and $(\Gamma^s)_{s=1}^\infty$ be a sequence of  disjoint subsets of $\Gamma$. If $\mathfrak{X}_{(\Gamma_n,i_n)}$ is a bonding of $(\mathfrak{X}_s)_{s=1}^\infty$ then the following hold.
\begin{enumerate}[label=(\roman*),leftmargin=21pt]

\item\label{bonding general statement 1} For $s\in\N$, $I_s$ extends to a norm-one linear projection $I_s:\mathfrak{X}_{(\Gamma_n,i_n)}\to \mathfrak{X}_{(\Gamma_n,i_n)}$ and its image $\mathfrak{X}_s$ is isometrically isomorphic to $\mathfrak{X}_{(\Gamma_n^s,i_n^s)}$.

\item\label{bonding general statement 2} For every $x_1\in\mathfrak{X}_1,\ldots,x_N\in\mathfrak{X}_N$ there exists $x_0\in\mathfrak{X}_0$ such that
\[\Big\|\sum_{s=1}^Nx_s - x_0\Big\|_{\ell_\infty(\Gamma)} \leq C\max_{1\leq s\leq N}\|x_s\|_{\ell_\infty(\Gamma)}.\]
In particular, $\mathfrak{X}_{(\Gamma_n,i_n)}/\mathfrak{X}_0$ is $C$-isomorphic to $(\oplus_{s=1}^\infty\mathfrak{X}_s)_{c_0}$ via $x+\mathfrak{X}_0\mapsto (I_sx)_{s=1}^\infty$.

\end{enumerate}

\end{theorem}

The second assertion implies that under appropriate assumptions (e.g., if $c_0$ does not embed in $\mathfrak{X}_{(\Gamma_n,i_n)}$) $\mathfrak{X}_0$ is uncomplemented in $\mathfrak{X}_{(\Gamma_n,i_n)}$.

We will prove the above theorem in pieces (see Proposition \ref{projection onto sd set} and Remark \ref{better estimate remark} for \ref{bonding general statement 1} and Proposition \ref{cancellation of coordinates} for \ref{bonding general statement 2} below), but first we discuss the bonding definition.

\begin{remarks}
Some analysis of the definition of a bonding is due on the basis of Remark \ref{evaluation analysis analysis}.
\begin{enumerate}[label=(\arabic*),leftmargin=21pt]

\item For each $s$, the only restriction on $\Gamma^s$ is that it is a self-determined set. There is virtually complete freedom in the choice of coefficients $\vec\lambda(\gamma) = (\lambda_1,\ldots,\lambda_a)$, for $\gamma\in\Gamma^s$. Therefore, for each $s\in\N$, the subspace $\mathfrak{X}_s$ can  be chosen to be isomorphic to any Bourgain-Delbaen-$\mathscr{L}_\infty$-space satisfying Assumption \ref{BD construction assumption} (for some $0<\vartheta\leq 1/4$ that is uniform in $s\in\N$.) In this paper, the $\mathfrak{X}_s$, $s\geq 1$, will be strongly incomparable Argyros-Haydon-like spaces. In particular, for $s\in\N$, $\mathfrak{X}_s$ has the scalar-plus-compact property and, for $s\neq t\in\N$, $\mathcal{L}(\mathfrak{X}_s,\mathfrak{X}_t) = \mathcal{K}(\mathfrak{X}_s,\mathfrak{X}_t)$.

\item For $\gamma\in\Gamma^0$ there is a meaningful restriction on the positions of $\vec\xi(\gamma) = (\xi_1,\ldots,\xi_a)$, $\vec\eta(\gamma) = (\eta_1,\ldots,\eta_a)$ that forces the functional $e_\gamma^*$ to interact with the spaces $\mathfrak{X}_s$, $s\in\N$, in an ordered manner that automatically yields Proposition \ref{general condition on projections} \ref{general condition on projections a}, \ref{general condition on projections b}, and, as we will see later, \ref{general condition on projections c}. However, the coefficients $\vec\lambda(\gamma) = (\lambda_1,\ldots,\lambda_a)$ can be chosen freely. This allows
\begin{enumerate}[label=(\greek*)]

\item the achievement of Proposition \ref{general condition on projections} \ref{general condition on projections d} and \ref{general condition on projections e} via imposing $U$ related restrictions on $\vec\lambda(\gamma) = (\lambda_1,\ldots,\lambda_a)$ and

\item the control of the ``horizontal'' behaviour of bounded linear operators on $\mathfrak{X}_{(\Gamma_n,i_n)}$ by incorporating Argyros-Haydon-like restrictions on $\vec\lambda(\gamma) = (\lambda_1,\ldots,\lambda_a)$.  In particular, the space $\mathfrak{X}_U$ has the scalar-plus-horizontally-compact property defined in Section \ref{section calkin algebra XU}.

\end{enumerate}
\end{enumerate}
\end{remarks}

\addtocontents{toc}
{\SkipTocEntry}
\subsection*{Boundedness of coordinate projections in Bourgain-Delbaen-$\mathscr{L}_\infty$-spaces}
\label{bd coordinate projections}
Here, we prove Theorem \ref{bonding general statement} \ref{bonding general statement 1}. Actually, we prove a slightly more general statement about coordinate projections $I_{\Gamma'}$ in Bourgain-Delbaen spaces $\mathfrak{X}_{(\Gamma_n,i_n)}$ that may be of independent interest.

This is the main assumption that yields the boundedness of $I_{\Gamma'}$.
 
\begin{assumption}
\label{restrictions witnessed on interval}
Let  $\mathfrak{X}_{(\Gamma_n,i_n)}$ be a Bourgain-Delbaen-$\mathscr{L}_\infty$-space and $\Gamma'\subset\Gamma$. Assume that for every $\gamma\in\Gamma$ there exists an interval $E$ of $\N_0$ (that may also be empty) such that $e_\gamma^*\circ I_{\Gamma'} = e_\gamma^*\circ P_E$.
\end{assumption}
 
Observe that the equality $e_\gamma^*\circ I_{\Gamma'} = e_\gamma^*\circ P_{E}$ is valid on the domain $\mathfrak{X}_{(\Gamma_n,i_n)}$ if and only if it is valid on the domain $\ell_\infty(\Gamma)$ because both functionals are finite linear combinations of $(e_\eta^*)_{\eta\in\Gamma}$.
 
 \begin{remark}
Assumption \ref{restrictions witnessed on interval} can be reformulated in terms of the basis $(d_\eta)_{\eta\in\Gamma}$, with a suitable fixed enumeration of $\Gamma$ (as described in Remark \ref{basis enumeration remark}) as follows: the support of $e_\gamma^*$ intersected with  $\Gamma'$, is equal to the support of $e_\gamma^*$ intersected with some interval $F\subset\Gamma$. The notions of support and of being an interval are regarded with the fixed enumeration of $\Gamma$. The crucial fact that $F$ is an interval, instead of an arbitrary subset of $\Gamma$,  automatically yields the boundedness of $I_{\Gamma'}$. An analogous approach was previously used by the second author in {\cite[Remark 2.2]{pelczar-barwacz:2023}} for the boundedness of coordinate projections in a Gowers-Maurey setting.
\end{remark}

 \begin{proposition}
 \label{projection onto sd set}
Let $\mathfrak{X}_{(\Gamma_n,i_n)}$ be a Bourgain-Delbaen-$C$-$\mathscr{L}_\infty$-space and $\Gamma'\subset\Gamma$.
\begin{enumerate}[label=(\roman*),leftmargin=19pt]

\item\label{projection onto sd set1} If $\Gamma'$ satisfies Assumption \ref{restrictions witnessed on interval} then $I_{\Gamma'}$ extends to a bounded linear projection of norm at most $2C$ onto $\mathfrak{X}' = [ (d_\gamma)_{\gamma\in\Gamma'}]$.

\item\label{projection onto sd set2} If, additionally to satisfying Assumption \ref{restrictions witnessed on interval}, $\Gamma'$ is self-determined then $\mathfrak{X}'$ is $2C$-isomorphic to $\mathfrak{X}_{(\Gamma_n',i_n')}$.

\end{enumerate}
\end{proposition}

\begin{proof}
Because for every interval $E$ of $\mathbb{N}_0$ we have $\|P_E\|_{\mathcal{L}(\ell_\infty(\Gamma))}\leq 2C$ it is immediate that for all $\gamma\in\Gamma$ we have $\|e^*_\gamma\circ I_{\Gamma'}\|_{\ell_\infty(\Gamma)^*} \leq 2C$ and consequently $I_{\Gamma'}$ extends to a bounded linear projection on $\mathfrak{X}_{(\Gamma_n,i_n)}$ of norm at most $2C$.

For the second assertion, let $r_{\Gamma'}:\ell_\infty(\Gamma)\to\ell_\infty(\Gamma')$ denote the usual restriction operator. By Proposition \ref{general s-d bd space} \ref{general quotient}, $R = r|_{\Gamma'}:\mathfrak{X}_{(\Gamma_n,i_n)}\to\mathfrak{X}_{(\Gamma_n',i_n')}$ has norm one, it is onto, and $R(\mathfrak{X}')$ is dense in $\mathfrak{X}_{(\Gamma_n',i_n')}$. We will prove that $R|_{\mathfrak{X}'}$ is $1/(2C)$-bounded below and, to that end, we fix $x\in\mathfrak{X}'$. We will show that for an arbitrary fixed $\gamma\in\Gamma$ we have $\|Rx\|_{\ell_\infty(\Gamma')} \geq 1/(2C)|e_\gamma^*(x)|$. By Remark \ref{how to defined composition with unbounded}, $e_\gamma^*\circ I_{\Gamma'}\in\langle\{d_\eta^*:\eta\in\Gamma'\}\rangle$ and, by self-determinedness, $\langle\{d_\eta^*:\eta\in\Gamma'\}\rangle = \langle\{e_\eta^*:\eta\in\Gamma'\}\rangle$. We may, therefore, express $e_\gamma^*\circ I_{\Gamma'}$ as a finite  linear combination $\sum_{\eta\in\Gamma'}\lambda_\eta e_\eta^*$. On one hand
\begin{equation}
\label{projection onto sd set eq1}
|e_\gamma^*(x)| = \big|\big(e_\gamma^*\circ I_{\Gamma'}\big)(x)\big| = \Big|\sum_{\eta\in\Gamma'}\lambda_\eta x(\eta)\Big| \leq \Big(\sum_{\eta\in\Gamma'}|\lambda_\eta|\Big) \|Rx\|_{\ell_\infty(\Gamma')}.
\end{equation}
On the other hand, there exists an interval ${E}$ of $\mathbb{N}_0$ such that $e_\gamma^*\circ I_{\Gamma'} = e_\gamma^*\circ P_{E}$ and $e_\gamma^*\circ P_{E}:\ell_\infty(\Gamma)\to\mathbb{C}$ has norm at most $2C$. Therefore,
\[\sum_{\eta\in\Gamma'}|\lambda_\eta| = \big\|e_\gamma^*\circ I_{\Gamma'}\big\|_{\ell_\infty^*(\Gamma)} \leq 2C.\]
\end{proof}

\begin{remark}
\label{better estimate remark}
In Proposition \ref{projection onto sd set} \ref{projection onto sd set1}, if we replace Assumption \ref{restrictions witnessed on interval} with the assumption $\sup_{\gamma\in\Gamma}\|e_\gamma^*\circ I_{\Gamma'}\|_{\ell^*_\infty(\Gamma)} \leq A$ then we can deduce $\|I_{\Gamma'}\|_{\mathcal{L}(\mathfrak{X}_{(\Gamma_n,i_n)})} \leq A$. If, furthermore, $\Gamma'$ is self-determined then $\mathfrak{X}_{(\Gamma_n',i_n')}\simeq^{A}\mathfrak{X}'$. If, in particular, $A=1$ then $\mathfrak{X}'$ is 1-complemented and isometric to $\mathfrak{X}_{(\Gamma_n',i_n')}$
\end{remark}

\begin{proof}[Proof of Theorem \ref{bonding general statement} \ref{bonding general statement 1}]
Although we could use Proposition \ref{projection onto sd set}, instead we show that, for $s\in\N$, $\sup_{\gamma\in\Gamma}\|e_\gamma^*\circ I_s\|_{\ell_\infty(\Gamma)^*} = 1$. By Remark \ref{better estimate remark}, the conclusion will follow. First, let $\gamma\in\Gamma\setminus\Gamma^0$. If $\gamma\in\Gamma^s$ then, because $\Gamma^s$ is self-determined, $e_\gamma^*\circ I_s = e_\gamma^*$ and, thus, $\|e_\gamma^*\circ I_s\|_{\ell_\infty(\Gamma)^*} = 1$. Otherwise, if $\gamma\in\Gamma^t$, for $t\neq s$, then, by Definition \ref{sd definition} \ref{sd definition complement}, $e_\gamma^*\circ I_s = 0$.

We treat the case $\gamma\in\Gamma^0$ by induction on $\mathrm{rank}(\gamma)$. If $\mathrm{rank}(\gamma) = 0$ then $e_\gamma^* = d_\gamma^*$ and, thus, for every $s\in\N$, $e_\gamma^*\circ I_s = 0$. For the inductive step, let $\gamma$ have rank $n$ and write its evaluation analysis
\[e_\gamma^* = \sum_{r=1}^ad_{\xi_r}^* + \sum_{r=1}^a\lambda_re_{\eta_r}^*\circ P_{E_r}\]
satisfying the conclusion of Definition \ref{bonding definition} \ref{bonding definition3}. For $s\in\N$, there is at most one $1\leq r\leq a$ with $s\in\mathrm{supp}_\mathrm{h}(e_{\eta_r}^*\circ P_{E_r})$ and, thus, $e_\gamma^*\circ I_s = \lambda_re_{\eta_r}^*\circ P_{E_r}\circ I_s = \lambda_re_{\eta_r}^*\circ I_s\circ P_{E_r}$. We, therefore, have, by the inductive assumption,
\[\|e_\gamma^*\circ I_s\|_{\ell_\infty(\Gamma)^*} \leq \vartheta \|e_{\eta_r}^*\circ I_s\|_{\ell_\infty(\Gamma)^*}\|P_{E_r}\|_{\mathcal{L}(\ell_\infty(\Gamma))} \leq \vartheta\frac{2}{1-2\vartheta} \leq 1,\]
because $\vartheta\leq 1/4$.
\end{proof}

\addtocontents{toc}
{\SkipTocEntry}
\subsection*{Cancellation of coordinates}
We prove Theorem \ref{bonding general statement} \ref{bonding general statement 2}. We actually prove a slightly more general statement about disjoint sequences of self-determined sets in general Bourgain-Delbaen spaces.

\begin{proposition}[Cancellation of coordinates]
\label{cancellation of coordinates}
Let $\mathfrak{X}_{(\Gamma_n,i_n)}$ be a Bourgain-Delbaen-$C$-$\mathscr{L}_\infty$-space, $(\Gamma^{s})_{s=1}^\infty$ be a sequence of pairwise disjoint self-determined subsets of $\Gamma$.
\begin{enumerate}[leftmargin=19pt,label=(\roman*)]

\item Let $N\in\mathbb{N}$ and $x_1\in\mathfrak{X}_1,\ldots,x_N\in\mathfrak{X}_N$ with $\supp(x_s)$ finite, for $1\leq s\leq N$. Then, there exists $x_0\in\mathfrak{X}_0$ with finite support and  $\min\supp(x_0) \geq \min(\cup_{s=1}^N\supp(x_s))$ such that
\[\Big\|\sum_{s=1}^Nx_s - x_0\Big\|_{\ell_\infty(\Gamma)} \leq C\max_{1\leq s\leq N}\|x_s\|_{\ell_\infty(\Gamma)}.\]

\item The quotient $\mathfrak{X}_{(\Gamma_n,i_n)}/\mathfrak{X}_0$ is $C$-isomorphic to $(\oplus_{s=1}^\infty\mathfrak{X}_{(\Gamma_n^s,i^s_n)})_{c_0}$. If, additionally, $\mathfrak{X}_{(\Gamma_n,i_n)}$ is a bonding of $(\mathfrak{X}_s)_{s=1}^\infty$, then $\mathfrak{X}_{(\Gamma_n,i_n)}/\mathfrak{X}_0$ is $C$-isomorphic to $(\oplus_{s=1}^\infty\mathfrak{X}_s)_{c_0}$ via $x+\mathfrak{X}_0\mapsto (I_sx)_{s=1}^\infty$.
\end{enumerate}
\end{proposition}

Theorem \ref{bonding general statement} \ref{bonding general statement 2} is included in the above statement. To prove Proposition \ref{cancellation of coordinates}, we first state and prove a lemma.

\begin{lemma}
\label{self-determinacy complement representations}
Let $\mathfrak{X}_{(\Gamma_n,i_n)}$ be a Bourgain-Delbaen-$\mathscr{L}_\infty$-space, $\Gamma'$ be a self-determined subset of $\Gamma$, and put $\mathfrak{X}_0 = [\{d_\gamma:\gamma\in\Gamma\setminus\Gamma'\}]$. Then, for every $m\leq n\in\mathbb{N}_0$ and $\gamma_0\in\Gamma\setminus\Gamma'$ with $\rank(\gamma_0) = m$, viewing $e_{\gamma_0}$ as a member of $\ell_\infty(\Gamma_n)$, we have that $i_n(e_{\gamma_0})\in\mathfrak{X}_0$ and $\supp(i_n(e_{\gamma_0})) \subset \{m,\ldots,n\}$.
\end{lemma}

\begin{proof}
For $k\in\mathbb{N}$ denote $\Delta^{0}_k = \Delta_k\setminus\Gamma'$. We will inductively prove the more precise statement
\[i_n(e_{\gamma_0}) = d_{\gamma_0} - \sum_{k=m+1}^n\sum_{\gamma\in\Delta_k^{0}}\lambda_\gamma d_\gamma.\]

This is proved by induction of $n-m = 0,1,\ldots$. In the base case, i.e., when $m = n$, by definition, $i_n(e_{\gamma_0}) = d_{\gamma_0}$ and, thus, the conclusion holds. Let $\ell\geq 0$ and assume that for every $m\leq n$ with $n-m = \ell$ the conclusion holds. Let $m\leq n$ such that $n - m = \ell+1$ and take $\gamma_0\in\Delta_m^{0}$. First, viewing $e_{\gamma_0}$ as a member of $\Gamma_{n-1}$, apply the inductive hypothesis to write
\begin{equation}
\label{self-determinacy complement representations eq1}
i_{n-1}(e_{\gamma_0}) = d_{\gamma_0} - \sum_{k=m+1}^{n-1}\sum_{\gamma\in\Delta_k^{0}}\lambda_\gamma d_\gamma.
\end{equation}
Put $y = i_{n-1}(e_{\gamma_0})$. We will prove that, viewing $e_{\gamma_0}$ as a member of $\ell_\infty(\Gamma_n)$,
\[i_n(e_{\gamma_0}) = y - \sum_{\gamma\in\Delta_n^{0}}e_\gamma^*(y)d_\gamma,\]
which will complete the proof.

Because $r_{\Gamma_n}:\langle\{d_\gamma:\gamma\in\Gamma_n\}\rangle\to\ell_\infty(\Gamma_n)$ is a bijection it suffices to show that for every $\xi\in\Gamma_n$,
\[e_\xi^*\Big(y - \sum_{\gamma\in\Delta_n^{0}}e_\gamma^*(y)d_\gamma\Big) = \delta_{\xi,\gamma_0}.\]

We write $\Gamma_n = \Delta_n^{0}\cup(\Delta_n\setminus\Delta_n^{0})\cup \Gamma_{n-1}$ and examine separately the case in which $\xi$ is in each of these sets. If $\xi\in\Delta_n^{0}$ then $e_\xi^*(\sum_{\gamma\in\Delta_n^{0}}e_\gamma^*(y)d_\gamma) = e_\xi^*(y)$ which yields the desired result. If $\xi\in\Delta_n\setminus\Delta_n^{0} \subset\Gamma'$ we obviously have $e_\xi^*(\sum_{\gamma\in\Delta_n^0}e_\gamma^*(y)d_\gamma) = 0$. Also, by applying Definition \ref{sd definition} \ref{sd definition complement} to \eqref{self-determinacy complement representations eq1}, we immediately obtain $e_\xi^*(y) = 0$. Finally, if $\xi\in\Gamma_{n-1}$ then obviously $e_\xi^*(\sum_{\gamma\in\Delta_n^{0}}e_\gamma^*(y)d_\gamma) = 0$ and $e_\xi^*(y) = e_\xi^*(e_{\gamma_0}) = \delta_{\xi,\gamma_0}$.
\end{proof}

\begin{proof}[Proof of Proposition \ref{cancellation of coordinates}]
The set $\Gamma' = \cup_{s=1}^\infty\Gamma^{s}$,  as a union of self-determined sets, is self-determined and, therefore satisfies the conclusion of Lemma \ref{self-determinacy complement representations}. Let now $x_1\in\mathfrak{X}_1,\ldots,x_N\in\mathfrak{X}_N$ with $\supp(x_s)$ finite and put $L = \min(\cup_{s=1}^N\supp(x_s))$. For $M\in\mathbb{N}$ sufficiently large, each $x_s$ is of the form
\[x_s = \sum_{k=L}^M\sum_{\gamma\in\Delta_k^{s}}d_\gamma^*(x_s)d_\gamma,\]
where $\Delta_k^{s} = \Gamma^{s}\cap\Delta_k$, for $0\leq s\leq N$ and $L\leq k\leq M$. Put $z = \sum_{s=1}^Nx_s$ and
\[x_0 = \sum_{k=L}^M\sum_{\gamma\in\Delta_k^{0}}e_\gamma^*(z)i_M(e_\gamma),\]
which, by Lemma \ref{self-determinacy complement representations}, is in $\mathfrak{X}_0$ and $\min\supp(x_0) \geq L$.

Because $\|z-x_0\|_{\ell_\infty(\Gamma)} = \|i_Mr_{\Gamma_M}(z-x_0)\|_{\ell_\infty(\Gamma)} \leq C\|r_{\Gamma_M}(z-x_0)\|_{\ell_\infty(\Gamma_M)}$, it suffices to show $|e_\gamma^*(z-x_0)|\leq \max_{1\leq s\leq N}\|r_{\Gamma^s}(x_s)\|_{\ell_\infty(\Gamma^s)}$, for $\gamma\in\Gamma_M$. We take $\gamma\in\Gamma_M$ and examine two cases. If $\gamma\in\Gamma^{0}$, then, obviously, $e_\gamma^*(z-x_0)=0$. If $\gamma\in\Gamma^{s}$, for some $s\in\mathbb{N}$, then because $\Gamma^{s}$ is self-determined, we may apply to it Definition \ref{sd definition} \ref{sd definition complement}, i.e., $|e_\gamma^*(z-x_0)| = |e_\gamma^*(x_s)| \leq \|x_s\|_{\ell_\infty(\Gamma)}$ (where we made the convention $x_s =0$ for $s>N$.)

To prove the second assertion first note that, actually, we have showed that for $x\in\mathfrak{X}_{(\Gamma_n,i_n)}$ and $x_0\in\mathfrak{X}_0$, $\|x+x_0\|_{\ell_\infty(\Gamma)} \leq C\max_{1\leq s\leq N}\|r_{\Gamma^s}(x)\|_{\ell_\infty(\Gamma^s)}$. By Definition \ref{sd definition} \ref{sd definition complement} the map $j:\mathfrak{X}_{(\Gamma_n,i_n)}/\mathfrak{X}_0\to(\oplus_{s=1}^\infty\mathfrak{X}_{(\Gamma_n^s,i^s_n)})_{c_0}$, given by $j(x+\mathfrak{X}_0) =  (r_{\Gamma^s}(x))_{s=1}^\infty\in$, is well defined and has norm at most $C$. Definition \ref{sd definition} \ref{sd definition complement} also yields that $j$ is non-contractive. Proposition \ref{general s-d bd space} \ref{general quotient} also implies that $j$ has dense range. In particular, $\mathfrak{X}_{(\Gamma_n,i_n)}/\mathfrak{X}_0$ is $C$-isomorphic to $(\oplus_{s=1}^\infty\mathfrak{X}_{(\Gamma_n^s,i^s_n)})_{c_0}$.

To conclude the proof, recall that in a bonding each $\Gamma^s$ satisfies Assumption \ref{restrictions witnessed on interval} then, for $s\in\N$, by the proof of Proposition \ref{projection onto sd set} \ref{projection onto sd set2}, $I_s = (r_{\Gamma^s}|_{\mathfrak{X}_s})^{-1}r_{\Gamma^s}$ and 
\[\|(r_{\Gamma^s}|_{\mathfrak{X}_s})^{-1}\| \leq \sup_{\gamma\in\Gamma}\|e_\gamma^*\circ I_s\|_{\ell_\infty(\Gamma)^*} = 1.\]
\end{proof}

\begin{remark}
\label{bonding unconditional projections}
Let $\mathfrak{X}_{(\Gamma_n,i_n)}$ be a Bourgain-Delbaen-$\mathscr{L}_\infty$-space and $(\Gamma^s)_{s=1}^\infty$ be a sequence of pairwise disjoint self-determined subsets, each of which satisfies Assumption \ref{restrictions witnessed on interval}, e.g., if $\mathfrak{X}_{(\Gamma_n,i_n)}$ is a bonding of $(\mathfrak{X}_s)_{s=1}^\infty$. By Proposition \ref{cancellation of coordinates} and Proposition \ref{always unconditional}, the sequence $(I_s)_{s=1}^\infty$ in $\mathcal{L}(\mathfrak{X}_{(\Gamma_n,i_n)})$ and its image in $\mathpzc{Cal}(\mathfrak{X}_{(\Gamma_n,i_n)})$ are unconditional.
\end{remark}

\subsection{The sequence of coordinate projections on a bonding. }
\label{section Linear combinations of coordinate projections}
By Remark \ref{bonding unconditional projections}, in a bonding $\mathfrak{X}_{(\Gamma_n,i_n)}$, the sequence $(I_s)_{s=1}^\infty$ in $\mathcal{L}(\mathfrak{X}_{(\Gamma_n,i_n)})$ and its image in $\mathpzc{Cal}(\mathfrak{X}_{(\Gamma_n,i_n)})$ are always unconditional. In this section, we discuss necessary conditions that when imposed on the coefficients $\vec\lambda(\gamma) = (\lambda_1,\ldots,\lambda_a)$, for $\gamma\in\Gamma^0$, then the sequence $(I_s)_{s=1}^\infty$ in $\mathcal{L}(\mathfrak{X}_{(\Gamma_n,i_n)})$ and it image in $\mathpzc{Cal}(\mathfrak{X}_{(\Gamma_n,i_n)})$ are equivalent to a prechosen unconditional sequence $(u_s)_{s=1}^\infty$. 
In the final part of this section, we utilize the developed theory by giving, for an arbitrary space $U$ with a 1-unconditional normalized Schauder basis $(u_s)_{s=1}^\infty$, a simple example of a bonding $\mathfrak{B}_U$ to which this theory is applied.

\begin{definition}
\label{upper lower U condition}
Let $U$ be a Banach space with a $1$-unconditional normalized basis $(u_s)_{s=1}^\infty$. Let $\mathfrak{X}_{(\Gamma_n,i_n)}$ be a Bourgain-Delbaen-$C$-$\mathscr{L}_\infty$-space, $(\Gamma^{s})_{s=1}^\infty$ be a sequence of disjoint self-determined subsets of $\Gamma$. Assume that $\mathfrak{X}_{(\Gamma_n,i_n)}$ is a bonding of $(\mathfrak{X}_s)_{s=1}^\infty$ and recall that for every $\gamma\in\Gamma^0$, either $e_\gamma^* = d_\gamma^*$ or $e_\gamma^*$ admits an evaluation analysis
\[e_\gamma^* = \sum_{r=1}^ad_{\xi_r}^* + \sum_{r=1}^a\lambda_re_{\eta_r}^*\circ P_{E_r}\]
such that, for $1\leq r\leq a$, $\xi_r\in\Gamma^0$ and the sets $(\mathrm{supp}_\mathrm{h}(e_{\eta_r}^*\circ P_{E_r}))_{r=1}^a$ are successive subsets of $\N$.

\begin{enumerate}[leftmargin=19pt,label=(\alph*)]

\item\label{upper lower U condition upper} We say that $\mathfrak{X}_{(\Gamma_n,i_n)}$ satisfies the upper $U$ condition if, for every $\gamma\in\Gamma^0$, either $e_\gamma^* = d_\gamma^*$ or the sequences of coefficients $\vec\lambda(\gamma) = (\lambda_1,\ldots,\lambda_a)$  and nodes $\vec\eta(\gamma)=(\eta_1,\dots,\eta_a)$ associated to $\gamma$  have the property that for all $f_r$ in the unit ball of $U^*$ with $\mathrm{supp}(f_r)\subset \mathrm{supp}_\mathrm{h}(e_{\eta_r}^*\circ P_{E_r})$, for $1\leq r\leq a$, we have
\[\Big\|\sum_{r=1}^a\lambda_rf_r\Big\|_{U^*}\leq 1.\]

\item\label{upper lower U condition lower} We say that $\mathfrak{X}_{(\Gamma_n,i_n)}$ satisfies the partial lower $U$ condition, for some constant $\beta>0$, if for every $x = \sum_{s=1}^Na_su_s\in U$ and $m\geq 0$, there exists $\gamma\in\Gamma^0$  with $\vec\lambda(\gamma)= (\lambda_1,\ldots,\lambda_a)$, $\vec\eta(\gamma)= (\eta_1,\ldots,\eta_a)$, and $\vec E(\gamma) = (E_1,\ldots,E_a)$  satisfying the following:
\begin{enumerate}[label=(\greek*)]

\item $\displaystyle{\Big|\sum^a_{s=1}\lambda_sa_s\Big| = \Big|\Big(\sum_{s=1}^a\lambda_su_s^*\Big)(x)\Big| \geq \beta\|x\|_U}$,

\item for $1\leq r\leq a$, $\eta_r\in\Gamma^r$ and $\mathrm{rank}(\eta_r)\in E_r$, and

\item $\mathrm{rank}(\eta_1) \geq m$.
    
\end{enumerate}

\end{enumerate}
\end{definition}

\begin{proposition}
\label{bonding conditions guaranteeing upper-lower-U}
Let $U$ be a Banach space with a $1$-unconditional normalized basis $(u_s)_{s=1}^\infty$. Let $\mathfrak{X}_{(\Gamma_n,i_n)}$ be a Bourgain-Delbaen-$C$-$\mathscr{L}_\infty$-space, $(\Gamma^{s})_{s=1}^\infty$ be a sequence of disjoint self-determined subsets of $\Gamma$. Assume that $\mathfrak{X}_{(\Gamma_n,i_n)}$ is a bonding of $(\mathfrak{X}_s)_{s=1}^\infty$.
\begin{enumerate}[leftmargin=19pt,label=(\roman*)]

\item\label{bonding conditions guaranteeing upper-lower-U upper} If $\mathfrak{X}_{(\Gamma_n,i_n)}$ satisfies the upper $U$ condition then, for every $x_1,\in\mathfrak{X}_1,\ldots,x_N\in\mathfrak{X}_N$,
\[\Big\|\sum_{s=1}^Nx_s\Big\|_{\ell_\infty(\Gamma)} \leq 2C\Big\|\sum_{s=1}^N\|x_s\|_{\ell_\infty(\Gamma)}u_s\Big\|_U.\]

\item\label{bonding conditions guaranteeing upper-lower-U lower} If $\mathfrak{X}_{(\Gamma_n,i_n)}$ satisfies the partial lower $U$ condition, for some constant $\beta>0$, then, for every $x = \sum_{s=1}^Na_su_s\in U$ and $m\geq 0$, there exist basis vectors $d_{\eta_1}\in\mathfrak{X}_1,\ldots,d_{\eta_N}\in\mathfrak{X}_N$ such that, for $1\leq r\leq a$, $\min\supp(d_{\eta_r})\geq m$ and
\[\Big\|\sum_{s=1}^Na_sd_{\eta_s}\Big\|_{\ell_\infty(\Gamma)} \geq \beta \Big\|\sum_{s=1}^Na_su_s\Big\|_U.\]

\end{enumerate}
\end{proposition}

\begin{proof}
The second assertion follows trivially by taking $\gamma\in\Gamma^0$ given by Definition \ref{upper lower U condition} \ref{upper lower U condition lower} and the basis vectors $d_{\eta_1},\ldots,d_{\eta_r}$, where $\eta_1,\ldots,\eta_r$ appear in the evaluation analysis of $\gamma$.

For the first assertion, we prove, by induction on $\mathrm{rank}(\gamma)$, the following statement. For every $x_1\in\mathfrak{X}_1,\ldots,x_N\in\mathfrak{X}_N$ there exist intervals $F_1,\ldots,F_N$ of $\N_0$ such that
\[\Big|e_\gamma^*\Big(\sum_{s=1}^Nx_s\Big)\Big| \leq \Big\|\sum_{s=1}^N\|P_{F_s}(x_s)\|_{\ell_\infty(\Gamma)}u_s\Big\|_U.\]
The conclusion then follows by the 1-unconditionality of $(u_s)_{s=1}^\infty$. We proceed with the inductive proof. If $e_\gamma^* = d_\gamma^*$, e.g., when $\mathrm{rank}(\gamma)=0$, the desired bound follows easily. Let $\gamma\in\Gamma$ such that the conclusion holds for all $\gamma'\in\Gamma$ with smaller rank. If $\gamma \in\Gamma^{s_0}$, for some $s_0\in\N$, Definition \ref{sd definition} \ref{sd definition complement} implies $e_\gamma^*(\sum_{s=1}^Nx_s) = e_\gamma^*(x_{s_0})$ (making the convention $x_s = 0$, for $s>N$) which directly gives the desired bound. If $\gamma\in\Gamma^0$ admits an evaluation analysis
\[e_\gamma^* = \sum_{r=1}^ad_{\xi_r}^* + \sum_{r=1}^a\lambda_re_{\eta_r}^*\circ P_{E_r},\]
by the inductive hypothesis, for $1\leq r\leq a$, putting $A_r = \mathrm{supp}_\mathrm{h}(e_{\eta_r}\circ P_{E_r})$, for some intervals $(F_s^r)_{s\in A_r}$ of $\N$ and norm-one $f_r\in U^*$ with $\mathrm{supp}(f_r) = A_r$,
\begin{align*}
\Big|e_{\gamma}^*\Big(\sum_{s=1}^Nx_s\Big)\Big| & = \Big|\sum_{r=1}^a\lambda_re_{\eta_r}^*\Big(\sum_{s\in A_r}P_{E_r}(x_s)\Big)\Big| \leq \sum_{r=1}^a|\lambda_r|\Big\|\sum_{s\in A_r}\big\|P_{F^r_s\cap E_r}(x_s)\big\|_{\ell_\infty(\Gamma)}u_s\Big\|_U\\
\span = \Big|\sum_{r=1}^a\lambda_r f_r\Big(\sum_{s\in A_r}\big\|P_{F^r_s\cap E_r}(x_s)\big\|_{\ell_\infty(\Gamma)}u_s\Big)\Big| = \Big|\Big(\sum_{r=1}^a\lambda_rf_r\Big)\Big(
\sum_{s\in A_r}\big\|P_{F^r_s\cap E_r}(x_s)\big\|_{\ell_\infty(\Gamma)}u_s\Big)\Big|\\
&\leq \Big\|\sum_{s=1}^N\big\|P_{F_s}(x_s)\big\|_{\ell_\infty(\Gamma)}u_s\Big\|_U,
\end{align*}
for appropriate intervals $F_1,\ldots,F_N$ of $\N$ (by the upper $U$ condition). 
\end{proof}

\begin{theorem}
\label{bd condition on projections}
Let $U$ be a Banach space with a $1$-unconditional normalized basis $(u_s)_{s=1}^\infty$. Let $\mathfrak{X}_{(\Gamma_n,i_n)}$ be a Bourgain-Delbaen-$C$-$\mathscr{L}_\infty$-space, $(\Gamma^{s})_{s=1}^\infty$ be a sequence of disjoint self-determined subsets of $\Gamma$. Assume that $\mathfrak{X}_{(\Gamma_n,i_n)}$ is a bonding of $(\mathfrak{X}_s)_{s=1}^\infty$ that satisfies the upper $U$ condition and the partial lower $U$ condition, for some $\beta>0$. Then, for every $a_1,\ldots,a_N\in\mathbb{C}$ and compact operator $K:\mathfrak{X}_{(\Gamma_n,i_n)}\to \mathfrak{X}_{(\Gamma_n,i_n)}$ we have
\begin{equation}
\label{operator equivalence}
\frac{\beta}{2C}\Big\|\sum_{s=1}^Na_su_s\Big\|_U \leq \Big\|\sum_{s=1}^Na_s I_{s} + K\Big\|\text{ and }\Big\|\sum_{s=1}^Na_s I_{s}\Big\| \leq 2C\Big\|\sum_{s=1}^Na_su_s\Big\|_U.
\end{equation}
Therefore, the sequence $(I_s)_{s=1}^\infty$ in $\mathcal{L}(\mathfrak{X}_{(\Gamma_n,i_n)})$ and its image in $\mathpzc{Cal}(\mathfrak{X}_{(\Gamma_n,i_n)})$ are $4C^2/\beta$-equivalent to $(u_s)_{s=1}^\infty$ and, in particular, $[(I_s)_{s=1}^\infty]\oplus \mathbb{C}I$ and $\mathcal{K}(\mathfrak{X}_{(\Gamma_n,i_n)})$ form a direct sum in $\mathcal{L}(\mathfrak{X}_{(\Gamma_n,i_n)})$.
\end{theorem}

\begin{proof}
We verify that the assumptions of Proposition \ref{general condition on projections} are satisfied for the sequence of projections $(I_s)_{s=1}^\infty$. Condition \ref{general condition on projections a} follows from Theorem \ref{bonding general statement} \ref{bonding general statement 1} whereas condition \ref{general condition on projections b} is an immediate consequence of the definition of the projection $I_s = I_{\Gamma^s}$, $s\in\N$, (Notation \ref{coordinate pro}) and the disjointness of the sets $\Gamma^s$, $s\in\N$. Condition \ref{general condition on projections c} follows from Theorem \ref{bonding general statement} \ref{bonding general statement 2} whereas \ref{general condition on projections d} follows from Proposition \ref{bonding conditions guaranteeing upper-lower-U} \ref{bonding conditions guaranteeing upper-lower-U upper}. Finally, condition \ref{general condition on projections e} follows from Proposition \ref{bonding conditions guaranteeing upper-lower-U} \ref{bonding conditions guaranteeing upper-lower-U lower}. Take the 1-norming set $W = \{e_\gamma^*:\gamma\in\Gamma\}$ and note that for a finite $S\subset\Gamma$, any $x\in\mathfrak{X}_{(\Gamma_n,i_n)}$ with $\min\mathrm{supp}(x) > \max_{\gamma\in S}\mathrm{rank}(\gamma)$, $x\in S_\perp$. Also, recall that, by Remark \ref{remark basis vector norm}, the basis vectors $d_\gamma$, $\gamma\in\Gamma$, have norm one.
\end{proof}

\addtocontents{toc}{\SkipTocEntry}
\subsection*{A simple bonding with the unitization of $U$ in its Calkin algebra}
\label{elementary BD example with U in Calkin} 
For any Banach $U$ space with a 1-unconditional and normalized basis $(u_s)_{s=1}^\infty$, we define a bonding $\mathfrak{B}_U$ satisfying assumptions of Theorem \ref{bd condition on projections}. Thus,  the unitization of $U$ embeds in $\mathpzc{Cal}(\mathfrak{B}_U)$. Within the context of Bourgain-Delbaen-$\mathscr{L}_\infty$ spaces, the construction of $\mathfrak{B}_U$ is concise and elegant. It serves as a semi-classical model for the space $\mathfrak{X}_U$ (see Section \ref{xu}), where, for appropriate $U$, the embedding of the unitization of $U$ into $\mathpzc{Cal}(\mathfrak{X}_U)$ is onto. This will be achieved by imposing Argyros-Haydon-like restrictions to the coefficients $\vec\lambda(\gamma) =(\lambda_1,\ldots,\lambda_a)$ and nodes $\vec\eta(\gamma) = (\eta_1,\ldots,\eta_a)$ associated to each $\gamma\in\Gamma$.

We fix a Banach space $U$ with a $1$-unconditional normalized basis $(u_s)_{s=1}^\infty$. We will define a bonding satisfying Assumption \ref{BD construction assumption} with $\vartheta=\frac{1}{4}$ and the assumptions of Theorem \ref{bd condition on projections} for $\beta = 1/8$. We fix an increasing sequence of finite sets $D_0\subset D_1\subset\cdots\subset D_r\subset\cdots$ with each $D_r$ a $1/2^{r+1}$-net of the complex unit disk. We denote

\[\mathcal{B} = \Bigg\{(\lambda_1,\dots,\lambda_a):a\in\mathbb{N},\;\lambda_r\in \frac{1}{4}D_r\;\text{ for }1\leq r\leq a,\text{ and }\Big\|\sum_{r=1}^a\lambda_ru_r^*\Big\|_{U^*}\leq 1\Bigg\}.\]

\begin{remark}
\label{pre two norming}
The subset
\[\left\{\sum_{r=1}^a\lambda_ru^*_r:\;a\in\mathbb{N}, (\lambda_r)_{r=1}^a\in\mathcal{B}\right\}\]
of the unit ball of $U^*$ is (1/8)-norming for $U$. That is, for every $x=\sum_{s=1}^\infty \mu_s u_s$ in $U$ there exists $(\lambda_1,\ldots,\lambda_a)\in\mathcal{B}$ such that $|\sum_{r=1}^a\lambda_r\mu_r| \geq (1/8)\|x\|$. This is because the $\lambda_r$ are chosen in $(1/4)D_r$ and each $D_r$ is a $1/2^{r+1}$-net of the unit disk.
\end{remark}

By induction on $n\in\N_0$, we define sets $\Delta_n$ together with extension functionals $(c^*_\gamma)_{\gamma\in\Delta_n}$. Each $\Delta_n$ will be a disjoint union of finite sets $\Delta_n^s$, $s=0,1,\ldots,n$. For each $\gamma$ in $\Delta_n$, we denote by $\mathrm{index}(\gamma)$ the unique $0\leq s\leq n$ for which $\gamma\in\Delta_n^s$. Eventually, for $s\in\N_0$, $\Gamma^s = \{\gamma\in\Gamma:\mathrm{index}(\gamma) = s\}$. In this example, for $s\neq 0$, all $\gamma\in\Gamma^s$ will satisfy Assumption \ref{BD construction assumption} \ref{BD construction assumption 0} and, thus, $c_\gamma^* = 0$. However, for $\gamma\in\Gamma^0$, we will impose restrictions on $\vec\lambda(\gamma)$.

Let $\Delta_0= \Delta_0^0 = \{(0)\}$. Assuming that we have defined $\Delta_0,\ldots,\Delta_{n-1}$ let $\Delta_n = \Delta_n^0\cup\Delta_n^1\cup\cdots\cup\Delta_n^n$, where, for $1\leq s\leq n$, $\Delta_n^s = \{(n,s)\}$ and for $\gamma = (n,s)$ we put $c_\gamma^* = 0$, leaving $\xi(\gamma)$, $E(\gamma)$, $\lambda(\gamma)$, and $\eta(\gamma)$ undefined. The set $\Delta_n^0$ is the disjoint union of two sets $\Delta_n^{0,a}$ and $\Delta_n^{0,b}$ defined as follows. Let 
\[\Delta_n^{0,a}=\Big\{ (n,\lambda,\eta): (\lambda)\in\mathcal{B},\eta\in\Gamma_{n-1}\text{ with }\mathrm{index}(\eta) = 1\Big\}\]
and for $\gamma = (n,\lambda,\eta) \in\Delta_n^{0,a}$ leave $\xi(\gamma)$ undefined and put $E(\gamma) = \{\mathrm{rank}(\eta)\}$, $\lambda(\gamma) = \lambda$, $\eta(\gamma) = \gamma$, and
\[c^*_\gamma=\lambda e^*_\eta\circ P^{(n-1)}_{\{\mathrm{rank}(\eta)\}}.\]
Let 
\begin{align*}
\Delta_n^{0,b}=\Big\{ (n,\xi,\lambda,\eta):&\; \xi\in \Gamma_{n-2}\cap\mathrm{dom}(E(\cdot))\text{ with }\mathrm{index}(\xi) = 0,\; \vec\lambda(\xi)^\frown (\lambda)\in\mathcal{B},\\
&\eta\in\Gamma_{n-1}\text{ with }\mathrm{index}(\eta) = \mathrm{age}(\xi)+1\text{ and }\mathrm{rank}(\xi)<\mathrm{rank}(\eta)\Big\}
\end{align*}
 and for $\gamma=(n,\xi,\lambda,\eta)\in\Delta_n^{0,b}$ define $\xi(\gamma) = \xi$, $E(\gamma) = \{\mathrm{rank}(\eta)\}$, $\lambda(\gamma) = \lambda$, $\eta(\gamma) = \eta$, and
\[c^*_\gamma=e^*_\xi+\lambda e^*_\eta\circ P^{(n-1)}_{\{\mathrm{rank}(\eta)\}}.\]

Trivially, for $s\geq 1$, the set $\Gamma^s$ is a self-determined subset of $\Gamma$ because the extension functional $c_\gamma^*$ associated with every $\gamma\in\Gamma^s$ is always zero. In fact, this easily yields that the space $\mathfrak{X}_{(\Gamma^s_n, i_n^s)}$ is naturally isometrically isomorphic to $c_0(\Gamma^s)$.

 Note that any node $\gamma\in\Gamma^0$, with $\mathrm{rank}(\gamma)\geq 1$, has evaluation analysis
\begin{equation}
\label{e-a BU-space}
e_\gamma^* =\sum_{r=1}^a d_{\xi_r}^* + \sum_{r=1}^a \lambda_r e_{\eta_r}^*\circ P_{\{\mathrm{rank}(\eta_r)\}} = \sum_{r=1}^a d_{\xi_r}^* + \sum_{r=1}^a \lambda_r d_{\eta_r}^*,
\end{equation}
where $\xi_r \in\Gamma^0$ and $\eta_r \in\Gamma^r$, for $1\leq r\leq a$, and $\vec \lambda(\gamma) = (\lambda_1,\dots,\lambda_a)\in\mathcal{B}$.  For any $s\geq 1$ and $\gamma\in\Gamma^s$, because $c_\gamma^* = 0$, $e^*_\gamma=d^*_\gamma$. This easily yields that $\mathfrak{B}_U = \mathfrak{X}_{(\Gamma_n,i_n)}$ is a bonding that satisfies the upper $U$ condition. Indeed, let $f_1,\ldots,f_r$ be functionals in the unit ball of $U^*$ such that, for $1\leq r\leq a$, $\mathrm{supp}(f_r)\subset\mathrm{supp}_\mathrm{h}(e_{\eta_r}^*\circ P_{\{\mathrm{rank}(\eta_r)\}}) = \{r\}$. Then, for $1\leq r\leq a$, $f_r = z_ru_r^*$, where $z_r$ is a scalar of magnitude at most one. Therefore, by the 1-unconditionality of $(u_s^*)_{s=1}^\infty$,
\[\Big\|\sum_{r=1}^a\lambda_rf_r\Big\|_{U^*} \leq \Big\|\sum_{r=1}^a\lambda_ru^*_r\Big\|_{U^*}\leq 1,\]
because $\vec \lambda(\gamma) = (\lambda_1,\dots,\lambda_a)\in\mathcal{B}$.

The following statement encompasses the main property of the space $\mathfrak{B}_U$.

\begin{proposition}    
\label{BU projections are U}
The space $\mathfrak{B}_U$ satisfies the following.
\begin{enumerate}[label=(\roman*),leftmargin=19pt]

\item For $s\in\N$, the image $\mathfrak{X}_s$ of each projection $I_s$ is isometrically isomorphic to $c_0$.

\item The sequence $(I_s)_{s=1}^\infty$ in $\mathcal{L}(\mathfrak{B}_U)$ and its image in $\mathpzc{Cal}(\mathfrak{B}_U)$ are 128-equivalent to $(u_s)_{s=1}^\infty$. In particular, the spaces $[(I_s)_{s=1}^\infty]\oplus \mathbb{C}I$ and $\mathcal{K}(\mathfrak{B}_U)$ form a direct sum in $\mathcal{L}(\mathfrak{B}_U)$.

\end{enumerate}
\end{proposition}

\begin{proof}
We only need to verify that $\mathfrak{B}_U$ satisfies the partial lower $U$ condition, for $\beta = 1/8$. Let $x=\sum_{s=1}^N\mu_su_s$ and $m\in\N$. Take $(\lambda_1,\dots,\lambda_N)\in\mathcal{B}$ with $\|u\| \leq 8\Big|\Big(\sum_{s=1}^N\lambda_su_s^*\Big)(u)\Big|=8\Big|\sum_{s=1}^N\mu_s\lambda_s\Big|$. 
We construct, inductively on $s$, nodes $(\xi_s)_{s=1}^N\subset\Gamma^0$  and simple nodes $(\eta_s)_{s=1}^N$ such that
 $\mathrm{rank}(\eta_1) \geq m$ and, for $1\leq r\leq N$, $e^*_{\xi_{r}}$ has evaluation analysis
\[e_{\xi_r}^* = \sum_{s=1}^{r} d_{\xi_s}^* + \sum_{s=1}^r \lambda_s d_{\eta_s}^*.\]
Then, $\gamma = \xi_N$ will witness the partial lower $U$ property.

We pick $\eta_1\in \Gamma^1$ arbitrarily with $\mathrm{rank}(\eta_1)>m$ and let $\xi_1=(n, \lambda_1,\eta_1)$, where $n=\mathrm{rank}(\eta_1)+1$. Assume we have chosen $\xi_1,\dots,\xi_s$ and $\eta_1,\dots,\eta_s$ for some $1\leq s<n$. Then  pick $\eta_{s+1}\in\Gamma^{s+1}$ with  $\mathrm{rank}(\eta_{s+1})>\mathrm{rank}(\xi_s)$. Then, also $\mathrm{rank}(\eta_{s+1})>\mathrm{rank}(\eta_s)$. Let $\xi_{s+1}=(n,\xi_s, \lambda_{s+1}, \eta_{s+1})$, where $n=\mathrm{rank}
(\eta_{s+1})+1$, which completes the inductive construction.  
\end{proof}

\begin{remark}
Although in this example, for simplicity, the component spaces $\mathfrak{B}_s$, $s\geq1$, are all isomorphic to $c_0$, this not a necessary restriction; in fact, we could have chosen any sequence of Bourgain-Delbaen-$\mathscr{L}_\infty$-spaces satisfying Assumption \ref{BD construction assumption}, for $0<\vartheta\leq 1/4$, (or, if we slightly modify the definition of a bonding, by {\cite{argyros:gasparis:motakis:2016}}, any sequence of $C$-$\mathscr{L}_\infty$-spaces).
\end{remark}


\section{The space $U$}

\label{U section}

In order to define $\mathfrak{X}_U$, for the Banach space $U$ with an unconditional basis $(u_n)_{n=1}^\infty$ that has no $c_0$ asymptotic version (see Definition \ref{def asymptotic version}), we need to determine a parameter $1\leq q<\infty$, depending on $U$. The parameter is tied to the aforementioned asymptotic property of $U$ and how it is characterized in terms of a lower-$q$ estimate with respect to an FDD defined by a decomposition of the unconditional basis $(u_s)_{s=1}^\infty$ into finite subsets.
The lower $q$-estimate of $U$ together with a $q$-convex implementation of the Argyros-Haydon technique allows the upper $U$ condition (Definition \ref{upper lower U condition} \ref{upper lower U condition upper}) to be satisfied in the resulting bonding $\mathfrak{X}_U$. The upper $U$ condition's only purpose is to aid in the representation of $U$ in the Calkin algebra of a bonding (see Proposition \ref{bonding conditions guaranteeing upper-lower-U} and Theorem \ref{bd condition on projections}).

Let us recall some standard terminology.
Given a Banach space $X$ with a Schauder basis $(x_s)_{s=1}^\infty$, an FDD $\boldsymbol{E} = (E_r)_{r=1}^\infty$ of $X$ is called a blocking of $(x_s)_{s=1}^\infty$ if there exists a partition of $\N$ into a sequence of consecutive intervals $(J_r)_{r=1}^\infty$ such that, for each $r\in\N$, $E_r = \langle\{x_s:s\in J_r\}\rangle$. 

The following is a result seemingly known to experts, but we could not find a direct reference to it.

\begin{theorem}
\label{no asymptotic space c0 equivalent}
Let $U$ be a Banach space with an unconditional basis $(u_s)_{s=1}^\infty$ for which $c_0$ is not an asymptotic version. Then, there exist a blocking of $(u_s)_{s=1}^\infty$ into an FDD $\boldsymbol{\Upsilon} = (U_r)_{r=1}^\infty$, $1\leq q<\infty$, and $\theta_U>0$ such that for any $\boldsymbol{\Upsilon}$-block sequence $(x_k)_{k=1}^n$ in $U$
\[\Big\|\sum_{k=1}^nx_k\Big\|_U \geq \theta_U\Big(\sum_{k=1}^n\|x_k\|_U^q\Big)^{1/q}.\]
That is, $U$ satisfies lower $q$-estimates for sequences of $\boldsymbol{\Upsilon}$-block vectors with constant $\theta_U$.
\end{theorem}

The proof of the above theorem is similar to results proved in {\cite[Section 3]{knaust:odell:schlumprecht:1999}}. For the sake of completeness, we provide a proof in the Appendix (see Section \ref{asymptotics section}). The following result comprises a standard dualization leading to an upper $p$-estimate in $U^*$. For completeness, we provide a proof.

\begin{proposition}\label{conjugate estimate U}
Let $U$ be a Banach space with an unconditional basis $(u_s)_{s=1}^\infty$ for which $c_0$ is not an asymptotic version and let the blocking $\boldsymbol{\Upsilon}$ of $(u_s)_{s=1}^\infty$, $1\leq q<\infty$, and $\theta_U>0$ be given by Theorem \ref{no asymptotic space c0 equivalent}.
For the conjugate exponent $p$ of $q$, the space $U^*$ satisfies upper $p$-estimates for $\boldsymbol{\Upsilon}$-block sequences with constant $\theta_U^{-1}$, i.e. for any $\boldsymbol{\Upsilon}$-block sequence $(f_k)_{k=1}^n$ in $U^*$,
\[\Big\|\sum_{k=1}^nf_k\Big\|_{U^*} \leq \theta_U^{-1}\Big(\sum_{k=1}^n\|f_k\|_{U^*}^p\Big)^{1/p}.\]
\end{proposition}

\begin{proof}
Let $f_1,\ldots,f_n$ be $\boldsymbol{\Upsilon}$-block. Put $r_0 = 0$ and, for $1\leq k\leq n$, $r_k = \max\supp_{\boldsymbol{\Upsilon}}(f_k)$. Let $F_1 = \cup_{r=1}^{r_1}J_r$, for $1<k< n$, $F_k = \cup_{r=r_{k-1}+1}^{r_k}J_r$, and $F_n = \cup_{r=r_{n-1}+1}^\infty J_r$. For any $x\in U$, the sequence $x|_{F_1} = \sum_{s\in F_1}u^*_s(x)u_s,\ldots,x|_{F_n} = \sum_{s\in F_n}u^*_s(x)u_s$ is $\boldsymbol{\Upsilon}$-block and $x|_{F_1}+\cdots+x|_{F_n} = x$, therefore,
\[\Big|\Big(\sum_{k=1}^nf_k\Big)(x)\Big|\leq \sum_{k=1}^n|f_k(x)|\leq \Big(\sum_{k=1}^n\|f_k\|_{U^*}^p\Big)^{1/p}\Big(\sum_{k=1}^n\|x|_{F_k}\|_U^q\Big)^{1/q}\leq \theta_U^{-1}\Big(\sum_{k=1}^n\|f_k\|_{U^*}^p\Big)^{1/p}\|x\|.\]
\end{proof}

For the rest of the paper, we fix a Banach space $U$ with a normalized 1-unconditional basis $(u_s)_{s=1}^\infty$ not having a $c_0$ asymptotic version. We fix the associated parameters $1\leq q<\infty$ and its conjugate exponent $1< p\leq\infty$, $\theta_U>0$, and the blocking $\boldsymbol{\Upsilon} = (U_r)_{r=1}^\infty$ provided by Theorem \ref{no asymptotic space c0 equivalent}. Actually, our proofs of the properties of the space $\mathfrak{X}_U$ rely on the upper $p$-estimates for $\boldsymbol{\Upsilon}$-block sequences with constant $\theta_U^{-1}$ provided by Proposition \ref{conjugate estimate U}. We record this in the following.

\begin{assumption}
\label{dual p assumption}
$U$ is a Banach space with a 1-unconditional basis $(u_s)_{s=1}^\infty$ such that there exist $1< p\leq \infty$, $\theta_U>0$, and a blocking $\boldsymbol{\Upsilon} = (U_r)_{r=1}^\infty$ of $(u_s)_{s=1}^\infty$ such that $U^*$ satisfies upper $p$-estimates for $\boldsymbol{\Upsilon}$-block sequences with constant $\theta_U^{-1}$.
\end{assumption}

We also fix the following notation to use in the definition of $\mathfrak{X}_U$.

\begin{notation}\label{upsilon blocking}
Let the blocking $\boldsymbol{\Upsilon} = (U_r)_{r=1}^\infty$ of $(u_s)_{s=1}^\infty$ be given by the partition $(J_r)_{r=1}^\infty$ of $\N$ into consecutive intervals, i.e., $U_r = \langle\{u_s:s\in J_r\}\rangle$, for $r\in\N$. For $A,B\subset\N$ we write $A\ll B$ provided $\max\{r\in\N: J_r\cap A\neq\emptyset\}<\min\{r\in\N: J_r\cap B\neq \emptyset\}$. For $n\in\N$ and $A\subset\N$, we write $n\ll A$ to mean $\{n\}\ll A$.  We say that a sequence $(A_k)_{k=1}^n$ of subsets of $\N$ is $\boldsymbol{\Upsilon}$-block, if $A_1\ll \dots\ll A_n$. 
\end{notation}

\section{$q$-mixed-Tsirelson spaces}\label{section mixed tsirelson}

Mixed-Tsirelson spaces were introduced in \cite{argyros:deliyanni:1997}, and they have their origins in Schlumprecht space from \cite{schlumprecht:1991}. They serve as a fundamental model for the Gowers-Maurey construction from \cite{gowers:maurey:1993} and all subsequent HI, and related, constructions, such as the Argyros-Haydon one from \cite{argyros:haydon:2011}. Moreover, Mixed-Tsirelson spaces are used directly, via the basic inequality (see, e.g., \cite{argyros:tolias:2004}) to perform computations, which yield bounds for linear combinations of special sequences of vectors that are core to proving HI or scalar-plus-compact properties in the aforementioned constructions. We follow this paradigm by modelling the spaces $\mathfrak{X}_U$ on versions of mixed-Tsirelson spaces that satisfy a type of $q$-convexity. This extra property is crucial in achieving the $U$-related conditions of Definition \ref{upper lower U condition} in a way that is compatible with the mixed-Tsirelson model.

Recall that $1<p\leq \infty$ and its conjugate exponent $1\leq q<\infty$ have been fixed. Our approach for defining mixed-Tsirelson space follows the notation of \cite{argyros:todorcevic:2005}. We focus on a special subclass of them, which we call $q$-mixed-Tsirelson spaces. These spaces are very closely related to a similar class, the $q$-spaces, introduced by Manoussakis in \cite{manoussakis:2001} and further studied in \cite{manoussakis:pelczar:2011}.

\begin{definition}\label{q-convex mT}
For sequences $(n_j)_{j=1}^\infty$ in  $\N$ and $(\theta_j)_{j=1}^\infty$ in $(0,1]$, denote by $W_q[(\mathcal{A}_{n_j},\theta_j)_{j\in\N}]$ the smallest subset $W$ of $c_{00}(\mathbb{N})$ satisfying the following properties:
\begin{enumerate}[leftmargin=21pt]
    \item $W$ contains the unit vector basis $(e^*_n)_{n=1}^\infty$.
    \item $W$ is closed under restrictions to subsets of integers and changes of signs, i.e., for any $\sum_{n}a_ne_n\in W$, $E\subset\N$, and $(\varepsilon_n)_{n=1}^\infty$ in $\{-1,0,1\}$, also $\sum_{n\in E}\varepsilon_na_ne^*_n$ is in $W$. 
    \item $W$ is closed under taking specific weighted averages, namely, for any $j\in\N$ and $f_1<\dots<f_d$ in $W$ with $d\leq n_j$, also,
    \[\theta_jn_j^{-1/p}(f_1+\dots+f_d)\]
    is in $W$.
\end{enumerate}
Then, any $f\in W_q[(\mathcal{A}_{n_j},\theta_j)_{j\in\N}]$ is either of the form $f = \pm e_n$, for some $n\in\N$, or $f = \theta_jn_j^{-1/p}(f_1+\dots+f_d)$, for some $j\in\N$ and $f_1<\cdots<f_d\in W_q[(\mathcal{A}_{n_j},\theta_j)_{j\in\N}]$, $d\leq n_j$. In the latter case, we will write $\mathrm{weight}(f) = \theta_j$.

Let $\|\cdot\|$ be the norm on $c_{00}(\mathbb{N})$ with the norming set $W_q[(\mathcal{A}_{n_j},\theta_j)_{j\in\N}]$, i.e., for $x\in c_{00}(\mathbb{N})$,
\[\|x\|=\sup\Big\{|f(x)|:f\in W_q[(\mathcal{A}_{n_j},\theta_j)_{j\in\N}]\Big\}.\]
The $q$-mixed-Tsirelson space $T_q[(\mathcal{A}_{n_j},\theta_j)_{j\in\N}]$ is defined as the completion of $(c_{00}(\mathbb{N}), \|\cdot\|)$.    
\end{definition}

\begin{remark}
It is well known (see, e.g., the proof of \cite[Proposition 1.1]{argyros:deliyanni:1997}) and it follows readily that $(e_n)_{n=1}^\infty$ is an unconditional Schauder basis for $T_q[(\mathcal{A}_{n_j},\theta_j)_{j\in\N}]$. 
\end{remark}

The above construction, for $q=1$, yields a classical mixed-Tsirelson space, as they were defined in {\cite{argyros:todorcevic:2005}}, and if, additionally, $\theta_j = \log(j+1)$ and $n_j = j$ then it yields Schlumprecht space. In this space, the unit vector basis of $\ell^n_1$ is well-isomorphic to block vectors in every block subspace, which is a consequence of Krivine's theorem (see \cite{krivine:1976}) and $\theta_jn_j^r\to\infty$, for all $r>0$. This hereditary $\ell_1$ finite block representation is one of the main ingredients in proving that Schlumprecht space is arbitrarily distortable and that the Gowers-Maurey space is HI. In the case $q>1$, achieving the hereditary $\ell_q$ finite block representation of $\ell_q^n$'s is more subtle. This was investigated in detail in \cite{manoussakis:2001} and \cite{manoussakis:pelczar:2011}. In the latter, such a representation was achieved via a $q$-convexification of the classical mixed-Tsirelson spaces from \cite{argyros:deliyanni:1997}. Our construction does not depend on hereditary $\ell_q$ finite block representation. It demands uniform upper $\ell_q$-bounds for all vectors, but lower bounds are only required for ``flat'' linear combinations of vectors. This is why our definition of $q$-mixed-Tsirelson spaces only contains $\theta_jn_j^{1/p}$-weighted sums of functionals.

We fix, for the rest of the paper, a pair of integers sequences satisfying the following.

\begin{assumption}\label{assumption integers} Let $(m_j)_{j=1}^\infty$, $(n_j)_{j=1}^\infty$ be sequences of integers such that the following hold.
\begin{enumerate}[label=(\alph*)]
    \item\label{1 assumption integers} $m_1\geq 16\cdot\theta_U^{-1}$,
    \item\label{2 assumption integers} $m_{j+1}\geq 4m_j^2$, $j\in\N$,
    \item\label{3 assumption integers} $n_1\geq m_1^2$, 
    \item\label{4 assumption integers} $n_{j+1}\geq (16n_j)^{q\cdot\log_2(m_{j+1})}$, $j\in\N$. 
\end{enumerate}
\end{assumption}
There are many variations of Assumption \ref{assumption integers} depending on the context, for example, the factor $\theta_U^{-1}$ in \ref{1 assumption integers} is specific to our construction. For appropriate $\theta$, we will consider spaces of the form $T_q[(\mathcal{A}_{n_j},\theta m_j^{-1})_{j\in\N}]$. The Bourgain-Delbaen construction of the space $\mathfrak{X}_U$ will be modelled on such a space, and it will also be used later (see Section \ref{RIS section}) as a tool for estimating the norm of linear combinations of special sequences of vectors.

\section{The definition of the space $\mathfrak{X}_U$}
In this section, we give the definition of the Bourgain-Delbaen-$\mathscr{L}_\infty$-space $\mathfrak{X}_U$. This is achieved in two steps. We first define a general Bourgain-Delbaen-$\mathscr{L}_\infty$-space $\bar{\mathfrak{X}} = \mathfrak{X}_{(\bar\Gamma_n,\bar i_n)}$. The definition of this space is rather ``unrestrained'' and mostly resembles the original Bourgain-Delbaen space from {\cite{bourgain:delbaen:1980}} with the Schur property. It involves some additional notational ingredients that facilitate the next step of the process, namely selecting an appropriate self-determined subset $\Gamma = \cup_{s=0}^\infty\Gamma^s$ of $\bar\Gamma$ to obtain $\mathfrak{X}_U = \mathfrak{X}_{(\Gamma_n,i_n)}$ as a quotient of $\bar{\mathfrak{X}}$. The space $\mathfrak{X}_U$ is a bonding satisfying the upper and partial lower $U$ conditions and, importantly, also incorporates Argyros-Haydon ingredients, both in $\Gamma^0$ governing its horizontal structure and in the sets $\Gamma^s$, $s\in\N$, governing structure of the component spaces $\mathfrak{X}_s$, $s\in\N$.

\subsection{The space $\bar{\mathfrak{X}}$}
We define the ``generic'' Bourgain-Delbaen-$\mathscr{L}_\infty$-space $\bar{\mathfrak{X}} = \mathfrak{X}_{(\bar\Gamma_n,\bar i_n)}$. As this is not the main space, we use ``bar''-notation for all the objects associated to this space, e.g., $\bar \Gamma_n$, $\bar \Delta_n$, $\bar c_\gamma^*$, $\bar P_E^{(n-1)}$, $\bar d_\gamma$, and $\bar d_\gamma^*$. For coordinate functionals $e_\gamma^*$ we use the same notation for all spaces. Eventually, $\mathfrak{X}_U$ will be defined by choosing a self-determined subset $\Gamma\subset\bar\Gamma$.

Recall that we have fixed the following parameters to use in the definition of $\bar{\mathfrak{X}}$.
\begin{enumerate}[label=(\alph*),leftmargin=19pt]
\label{parameters bar X}
\item A Banach space $U$ with a 1-unconditional basis, $1< p\leq \infty$, and a blocking $\boldsymbol{\Upsilon} = (U_r)_{r=1}^\infty$ satisfying Assumption \ref{dual p assumption}. We also denote $q$ the conjugate exponent of $p$.

\item A pair of sequences $(m_j)_{j=1}^\infty, (n_j)_{j=1}^\infty$ of positive integers satisfying Assumption \ref{assumption integers}. We also let $m_0=1$, but do not define an integer $n_0$.

\item We also fix an increasing sequence of of finite sets $D_0\subset D_1\subset\cdots\subset D_n\subset\cdots$, with each $D_n$ an $1/2^{n+1}$-net of the complex unit disk containing $\{1\}\cup\{\pm 4n_j^{-1/p}: 1\leq j\leq n\}$.     

\end{enumerate}

We design the space $\bar{\mathfrak{X}}$ to satisfy Assumption \ref{BD construction assumption}, following similar notation as in Example \ref{schur bd space}, Example \ref{toy self-determined}, and the example $\mathfrak{B}_U$ from page \pageref{elementary BD example with U in Calkin}. In addition, we will recursively define an auxiliary function $\mathrm{weight}:\bar\Gamma\to\{m_j^{-1}:j\in\N_0\}\cup\{0\}$, related to the coefficient partial function $\lambda(\cdot)$. Put $\bar\Delta_0 = \{(0)\}$ and for $\gamma = (0)$ define $\mathrm{weight}(\gamma) = 0$ and leave $\xi(\gamma)$, $E(\gamma)$, $\lambda(\gamma)$ and $\eta(\gamma)$ undefined. All subsequent sets $\bar\Delta_n$ will be defined as finite collections of tuples $\gamma$ and each such tuple will encode the corresponding extension functional $\bar c_\gamma^*$. Assume $\bar\Delta_0,\bar\Delta_1,\ldots,\bar\Delta_{n-1}$ and the function $\mathrm{weight}:\bar\Gamma_{n-1}\to\{m_j^{-1}:j\in\mathbb{N}_0\}\cup\{0\}$ have been defined. We will define $\bar\Delta_n = \{(n)\}\cup \bar\Delta_n^\mathrm{a}\cup \bar\Delta_n^\mathrm{b}$. For $\gamma = (n)$, we put $c_\gamma^* = 0$, $\mathrm{weight}(\gamma) = 0$, and leave $\xi(\gamma)$, $E(\gamma)$, $\lambda(\gamma)$ and $\eta(\gamma)$ undefined. We let
\begin{equation}
\label{generic tuples 0}
\bar\Delta_n^\mathrm{a} = \Big\{(n,m_{j}^{-1},E,\lambda,\eta):E\subset\{0,\ldots,n-1\}\text{ interval}, 0\leq j\leq n, \lambda\in\frac{1}{4}D_n, \eta\in\Gamma_{n-1}\Big\}.
\end{equation}
For such a tuple $\gamma = (n,m_{j}^{-1},E,\lambda,\eta)$ define $\mathrm{weight}(\gamma) = m_{j}^{-1}$, $E(\gamma) = E$, $\lambda(\gamma) = m_j^{-1}\lambda$, $\eta(\gamma) = \eta$, and leave $\xi(\gamma)$ undefined. Finally, put
\[\bar c_\gamma^* = \frac{1}{m_j}\lambda e_\eta^*\circ \bar P_E^{(n-1)}.\]
Next, we let
\begin{equation}
\label{generic tuples 1}
\begin{split}
    \bar\Delta_n^\mathrm{b} = \Big\{(n,m_j^{-1},\xi,E,\lambda,\eta):&\;\xi\in\Gamma_{n-2}\cap\mathrm{dom}(E(\cdot)),E\subset\{\rank(\xi)+1,\ldots,n-1\}\text{ interval},\\
    &m^{-1}_j = \mathrm{weight}(\xi),\;\lambda\in\frac{1}{4}D_n,\eta\in\Gamma_{n-1}\text{ with }\rank(\xi)<\rank(\eta)\Big\}.
\end{split}
\end{equation}
For such a tuple $\gamma = (n,m_{j}^{-1},\xi,E,\lambda,\eta)$ define $\mathrm{weight}(\gamma) = m_{j}^{-1}$, $\xi(\gamma)=\xi$, $E(\gamma) = E$, $\lambda(\gamma) = m_j^{-1}\lambda$, and $\eta(\gamma) = \eta$. Finally, put
\[\bar c_\gamma^* = e^*_{\xi} + \frac{1}{m_j}\lambda e_\eta^*\circ \bar P_E^{(n-1)}.\]

By design, the extension scheme $(\bar\Gamma_n)_{n=0}^\infty$, $(\bar i_{n-1,n})_{n=1}^\infty$ satisfies Assumption \ref{BD construction assumption} for $\vartheta = 1/4$ and, therefore, defines a Bourgain-Delbaen-$2$-$\mathscr{L}_\infty$-space $\bar{\mathfrak{X}} = \mathfrak{X}_{(\bar\Gamma_n,\bar i_n)}$.

\begin{remark}
Note that in \eqref{generic tuples 1} there is no restriction on the age of $\xi$ (defined in Remark \ref{gamma tuples}), as in the Argyros-Haydon construction \cite{argyros:haydon:2011}, and for any $\gamma\in\Gamma$, either $e_\gamma^* = d_\gamma^*$ or
\begin{equation}
e_\gamma^* = \sum_{r=1}^a d_{\xi_r}^* + \sum_{r=1}^a \frac{1}{m_j}\lambda_r e_{\eta_r}^*\circ  P_{E_r},
\end{equation}
where $\mathrm{weight}(\gamma) = m_j^{-1}$, $\vec\lambda(\gamma) = (m_j^{-1}\lambda_1,\ldots,m_j^{-1}\lambda_a)$, $\mathrm{age}(\gamma) = a$, $\vec\xi(\gamma) = (\xi_1,\ldots,\xi_a)$, $\vec E(\gamma) = (E_1,\ldots,E_a)$, and $\vec\eta(\gamma) = (\eta_1,\ldots,\eta_a)$.
\end{remark}

\subsection{The space $\mathfrak{X}_U$} \label{xu}

We shall define $\mathfrak{X}_U$ by choosing an appropriate self-determined subset $\Gamma$ of $\bar\Gamma$. There are disjoint self-determined subsets $(\Gamma^s)_{s=1}^\infty$ such that, following Notation \ref{pre-bonding notation}, the space $\mathfrak{X}_U = \mathfrak{X}_{(\Gamma_n,i_n)}$ is a bonding of $(\mathfrak{X}_s)_{s=1}^\infty$ and it addresses the following main components.

\begin{enumerate}[label=(\arabic*),leftmargin=19pt]

\item The set $\Gamma^0$ is modeled after a $q$-convex variant of the Argyros-Haydon $\mathscr{L}_\infty$-space. There are three specific considerations in its design.
\begin{enumerate}[label=(\roman*)]

\item The bonding $\mathfrak{X}_U$ satisfies the upper $U$ condition. This is the initial reason for using a $q$-convex variant of an Argyros-Haydon space and including the relation ``$\ll$'' in the definition of nodes. 

\item The bonding $\mathfrak{X}_U$ satisfies the partial lower $U$ condition. This is ensured by the inclusion in $\Gamma^0$ of low-complexity ``ground functionals'' that are subject to a Schreier admissibility restriction and carry coefficients from the unit ball of $U^*$.

\item The horizontal block structure of $\mathfrak{X}_U$ is very tight; this enables the control of the behavior of bounded linear operators on horizontally block sequences (cf. Definition \ref{horizontal vector support definition}), and it is achieved via the Argyros-Haydon-type construction.

\end{enumerate}

\item The spaces $\mathfrak{X}_s$, $s\in\N$, are $q$-convex variants of Argyros-Haydon $\mathscr{L}_\infty$-spaces, such that the following are satisfied.
\begin{enumerate}[label=(\roman*),leftmargin=19pt]

\item Each of them has the scalar-plus-compact property and they are pairwise very incomparable, meaning that, for $s\neq t\in\N$, any bounded linear operator $T:\mathfrak{X}_s\to\mathfrak{X}_t$ is compact (see {\cite[Section 10.2]{argyros:haydon:2011}}).

\item The structure of each $\mathfrak{X}_s$, $s\in\N$, is very incomparable to the horizontal one. The implementation of this requires that the component spaces share similar $q$-convexity properties as the horizontal structure, for the same value of $q$.

\end{enumerate}

\end{enumerate}

In addition to the parameters gathered before defining $\bar{\mathfrak{X}}$ (see page \pageref{parameters bar X}), we introduce additional ones used in defining the self-determined subset $\Gamma$ of $\bar\Gamma$.
\begin{enumerate}[leftmargin=19pt,label=(\alph*)]
\setcounter{enumi}{3}

\item We fix a partition $\N_0=\bigcup_{s=0}^\infty L_s$ into  pairwise disjoint sets $(L_s)_{s=0}^\infty$, such that $0\in L_0$ and, for each $s\in\N_0$,  $L_s\cap 4\N$, $L_s\cap (4\N-2)$ and $L_s\cap (2\N-1)$ are infinite.

\item We define a function $\mathrm{index}:\bar\Gamma\to\N_0$ by letting, for $\gamma\in\bar\Gamma$,
\[\mathrm{index}(\gamma) =
\left\{
	\begin{array}{ll}
		s & \mbox{if } \mathrm{weight}(\gamma) = m_j^{-1}\neq 0\text{ and }j\in L_s, \\
		\rank(\gamma) & \mbox{if } \mathrm{weight}(\gamma) = 0.
	\end{array}
\right.
\]

An analogous index function appeared in the example in Section \ref{elementary BD example with U in Calkin}.

\item For $s\in\mathbb{N}_0$ we let $\bar\Gamma^s=\{\gamma\in\bar\Gamma: \mathrm{index} (\gamma)=s\}$.

The goal is to choose $\Gamma^0\subset \bar\Gamma^0$, $\Gamma^1\subset\bar\Gamma^1$,\ldots  with $\Gamma=\bigcup_{s=0}^\infty\Gamma^s$. Note that each $\bar\Gamma^s$ contains a unique $\gamma$ with $c_\gamma^* = 0$, namely $\gamma = (s,0)$. This will be included in $\Gamma^s$, and it can be interpreted as the ``seed'' from which other elements of $\Gamma^s$ will proliferate (see the final paragraph in Example \ref{toy self-determined}). 

\item We fix an injective coding function $\sigma: \bar\Gamma\to \N$ such that, for any $\gamma\in\bar\Gamma$, we have $\sigma(\gamma)>\ \mathrm{ rank}(\gamma)$ 
and, if $\mathrm{index}(\gamma) = s$, then $\sigma(\gamma)\in L_s\cap 4\N$.

The restriction  $\sigma|_{\bar\Gamma^s}$ will serve as coding functions in the construction of $\mathfrak{X}_s$, if $s\in\N$, and of the horizontal Argyros-Haydon structure, if $s=0$.

\item 
As in Section \ref{elementary BD example with U in Calkin} we denote
\[\mathcal{B} = \Bigg\{(\lambda_1,\dots,\lambda_k):k\in\mathbb{N},\;\lambda_n\in \frac{1}{4}D_n\;\text{ for }1\leq n\leq k,\text{ and }\Big\|\sum_{n=1}^k\lambda_nu_n^*\Big\|\leq 1\Bigg\}.\]

This will be used to achieve the partial lower $U$ condition. Recall that, by Remark \ref{pre two norming}, the subset
\[\left\{\sum_{n=1}^k\lambda_nu^*_n:\;k\in\mathbb{N}, (\lambda_n)_{n=1}^k\in\mathcal{B}\right\}\]
of the unit ball of $U^*$ is (1/8)-norming for $U$.

\end{enumerate}

Now we are ready to define inductively on $n\in\N$ the sets $\Delta_n\subset\bar\Delta_n$ alongside with $\Gamma_n=\bigcup_{i=0}^n\Delta_i$, $\Gamma_n^s=\Gamma_n\cap \bar\Gamma^s$, $s\in\N_0$. Let $\Delta_0= \Delta_0^0 = \bar\Delta_0 = \{(0)\}$. For any $n\in\N$ we define $\Delta_n$ to be the union of pairwise disjoints sets $\Delta_n^s$, $0\leq s\leq n$ defined below. 

\newcounter{XU}

For $s=n$ we let $\Delta_n^s = \{(n)\}\subset\bar\Gamma^s$ and for $s>n$ we let $\Delta_n^s=\emptyset$. For $1\leq s< n$ the set $\Delta_n^s\subset\bar\Gamma^s$ is the collection of nodes in $\bar\Delta_n$ of one of the following forms:
\begin{enumerate}[leftmargin=25pt,label=(\arabic*)]
    \item $\gamma=(n,m_{2j}^{-1},E,\pm n_{2j}^{-1/p},\eta)\in\bar\Delta_n$ with $2j\in L_s$ for some $s\in \N$, $\eta\in\Gamma_{n-1}^s$, 
    \item $\gamma=(n,m_{2j-1}^{-1},E,n_{2j-1}^{-1/p},\eta)\in\bar\Delta_n$ with $2j-1\in L_s$ for some $s\in \N$, $\eta\in\Gamma_{n-1}^s$ with $\mathrm{weight}(\eta)=m_{4i-2}^{-1}$ for some $4i-2\in L_s$ with $m_{4i-2}>n_{2j-1}^3$,
    \item $\gamma=(n,m_{2j}^{-1},\xi,E,\pm n_{2j}^{-1/p},\eta)\in\bar\Delta_n$ with $2j\in L_s$ for some $s\in \N$, $\xi,\eta\in\Gamma_{n-1}^s$  and $\mathrm{age}(\xi)<n_{2j}$,
    \item $\gamma=(n,m_{2j-1}^{-1},\xi, E,n_{2j-1}^{-1/p},\eta)\in\bar\Delta_n$ with $2j-1\in L_s$ for some $s\in \N$, $\xi,\eta\in\Gamma_{n-1}^s$ with  $\mathrm{age}(\xi)<n_{2j-1}$ and $\mathrm{weight}(\eta)=m_{\sigma(\xi)}^{-1}$.
\end{enumerate}

The set $\Delta_n^0\subset\bar\Gamma^0$ is the union of two sets $\Delta_n^{\mathrm{(hor)}}$ and $\Delta_n^{\mathrm{(gr)}}$. The definition of the former uses the notion of horizontal support (cf. Notation \ref{pre-bonding notation}) and the relation ``$\ll$'' (cf. Notation \ref{upsilon blocking}). Then, $\Delta_n^{\mathrm{(hor)}}\subset\bar\Gamma^0$ is the collection of nodes in $\bar\Delta_n$ of one of the following forms:
\begin{enumerate}[leftmargin=25pt,label=(\arabic*)]\setcounter{enumi}{4}
    \item\label{initial even h-node} $\gamma=(n,m_{2j}^{-1},E,\pm n_{2j}^{-1/p},\eta)\in\bar\Delta_n$ with $2j\in L_0\setminus\{0\}$ and $\eta\in\Gamma_{n-1}$,
    \item $\gamma=(n,m_{2j-1}^{-1},E,n_{2j-1}^{-1/p},\eta)\in\bar\Delta_n$ with $2j-1\in L_0$ and $\eta\in\Gamma_{n-1}$ with $\mathrm{weight}(\eta)=m_{4i-2}^{-1}$ for some $4i-2\in L_0$ with $m_{4i-2}>n_{2j-1}^3$,
    \item\label{successor even h-node} $\gamma=(n,m_{2j}^{-1},\xi, E,\pm n_{2j}^{-1/p},\eta)\in\bar\Delta_n$ with $2j\in L_0\setminus\{0\}$ and $\xi,\eta\in\Gamma_{n-1}$ with $\mathrm{supp}_\mathrm{h}(e^*_\xi)\ll\mathrm{supp}_\mathrm{h}(e^*_\eta\circ \bar P_E)$, $\rank(\xi) < \mathrm{supp}_\mathrm{h}(e_\eta^*\circ\bar P_E)$ and $\mathrm{age}(\xi)<n_{2j}$,
    \item $\gamma=(n,m_{2j-1}^{-1},\xi,E,n_{2j-1}^{-1/p},\eta)\in\bar\Delta_n$ with $2j-1\in L_0$ and $\xi,\eta\in\Gamma_{n-1}$ with $\mathrm{supp}_\mathrm{h}(e^*_\xi)\ll\mathrm{supp}_\mathrm{h}(e^*_\eta\circ \bar P_E)$, $\rank(\xi) < \mathrm{supp}_\mathrm{h}(e_\eta^*\circ\bar P_E)$,  $\mathrm{age}(\xi)<n_{2j-1}$, and $\mathrm{weight}(\eta)=m_{\sigma(\xi)}^{-1}$.
\end{enumerate}

Finally, the set $\Delta_n^{\mathrm{(gr)}}\subset\bar\Gamma^0$ is the collection of nodes in $\bar\Delta_n$ of one of the following forms:
\begin{enumerate}[leftmargin=25pt,label=(\arabic*)]\setcounter{enumi}{8}
    \item $\gamma=(n,m_0^{-1}, E, \lambda, \eta)\in \bar\Delta_n$ with $\eta\in \Gamma_{n-1}^1$, 
    $(\lambda)\in \mathcal{B}$ and $E=\{\mathrm{rank}(\eta)\}$, 
    \item $\gamma=(n,m_0^{-1},\xi, E, \lambda, \eta)\in \bar\Delta_n$ with 
    $\eta\in \Gamma_{n-1}^{s+1}$, 
    where $s=\mathrm{age}(\xi)$, 
    $\lambda(\xi)^\frown (\lambda)\in\mathcal{B}$,
    $E=\{\mathrm{rank}(\eta)\}$,
    and $\mathrm{age}(\xi) < \min\supp(e_\xi^*)$. 
\end{enumerate}

Let $\Gamma=\bigcup_{n=0}^\infty\Delta_n$ and  $\Gamma^s=\bar\Gamma^s\cap\Gamma=\bigcup_{n=s}^\infty\Delta_n^s$  for each $s\in\N_0$.

\begin{proposition}\label{Gamma's self-determined}
The set $\Gamma$ is a self-determined subset of $\bar\Gamma$. Therefore, by Proposition \ref{general s-d bd space} \ref{general s-d bd space1}, $\Gamma$ defines a Bourgain-Delbaen-$2$-$\mathscr{L}_\infty$-space $\mathfrak{X}_U = \mathfrak{X}_{(\Gamma_n,i_n)}$.

Furthermore, for $s\geq 1$, $\Gamma^s$ is a self-determined subset of $\Gamma$ defining a Bourgain-Delbaen-$2$-$\mathscr{L}_\infty$-space $\mathfrak{X}_{(\Gamma_n^s,i_n^s)}$.
\end{proposition}

\begin{proof}
It is enough to show that for any $n\in\mathbb{N}$ and $\gamma\in \Delta_{n}$ we have that $\bar c_\gamma^*:\ell_\infty(\bar\Gamma_{n-1})\to\mathbb{C}$ is in $\langle\{e^*_\eta\circ \bar P^{(n-1)}_E: \eta\in\Gamma_{n-1}, E\subset\N\}\rangle$, which is true by the definition of $\bar c^*_\gamma$ and $\Gamma$. The same reasoning applies to $\Gamma^s$ for each $s\in\N$, with $\Delta_{n}$, respectively $\Gamma_{n-1}$, replaced by $\Delta_{n}^s$, respectively $\Gamma^s_{n-1}$.
\end{proof}

\begin{notation}
Having defined the space $\mathfrak{X}_U = \mathfrak{X}_{(\Gamma_n,i_n)}$ we scrap the bar-notation for all objects associated to this space; we denote its basis by $(d_\gamma)_{\gamma\in\Gamma}$, its biorthogonal sequence by $(d_\gamma^*)_{\gamma\in\Gamma}$, the corresponding extension functionals by $(c^*_\gamma)_{\gamma\in\Gamma}$ and the FDD projections by $P_E$, $E\subset\N_0$. Also, for $s\in\mathbb{N}_0$, denote $\mathfrak{X}_s =[\{d_\gamma:\gamma\in\Gamma^s\}]$ and, for $s\geq 1$, $I_s = I_{\Gamma^s}$.
\end{notation}

\begin{remark}\label{horizontal in quotient}
Formally, horizontal support is defined with respect to the sets $\bar\Gamma^s$, however, its meaning carries over to the sets $\Gamma^s$. By this, we mean that if for $\gamma\in\Gamma$ and $E\subset\mathbb{N}_0$ we define $\mathrm{supp}_\mathrm{h}(e^*_\gamma\circ P_E) = \mathrm{supp}_\mathrm{h}(e^*_\gamma\circ \bar P_E)$, then it follows that
\[\mathrm{supp}_\mathrm{h}(e^*_\gamma\circ P_E) = \Big\{s\in\mathbb{N}:e^*_\gamma\circ P_E(d_\gamma)\neq 0\text{ for some }\gamma\in\Gamma^s\Big\}= \Big\{s\geq 1: e_\gamma^*\circ I_s\circ P_E\neq 0\Big\}.\]
Indeed, because $\Gamma$ is self-determined, $e_\gamma^*\circ \bar P_E$ can be expressed as a finite linear combination $\sum_{\eta\in\Gamma}\lambda_\eta\bar d_\eta^*$ and, thus, $\mathrm{supp}_\mathrm{h}(e_\gamma^*\circ\bar P_E) = \{s\in\mathbb{N}:\lambda_\eta\neq 0\text{ for some }\eta\in\Gamma^s\}$. By the preservation of formulae on self-determined sets, $e_\gamma^*\circ P_E = \sum_{\eta\in\Gamma}\lambda_\eta d_\eta^*$.

The notion of support with respect the the Bourgain-Delbaen FDD is also well-behaved, i.e., for $\gamma\in\Gamma$ and $E\subset\N_0$, $\supp(e_\gamma^*\circ P_E) = \supp(e_\gamma^*\circ\bar P_E)$.
\end{remark}

\subsection{Types of nodes in $\Gamma$} 
We discuss below different types of nodes and their role in the construction of $\mathfrak{X}_U$.  By the definition of $\Gamma$ and the preservation of formulae on self-determined sets, any node $\gamma\in\Gamma$, with $\mathrm{weight}(\gamma)\neq 0$, is of one of the following types:
\begin{enumerate}[leftmargin=12pt,label=\textbullet]
    
    \item It is simple, i.e. $\gamma\in\Gamma^s$, for some $s\in\N$. Then necessarily $s=\mathrm{index}(\gamma)$ and the evaluation analysis of $\gamma$ is of the form 
    \begin{equation}\label{eva analysis sim}
    e_\gamma^* = \sum_{r=1}^a d_{\xi_r}^* + \frac{1}{m_j}\sum_{r=1}^a \frac{1}{n_j^{1/p}} e_{\eta_r}^*\circ  P_{E_r}
    \end{equation} 
    with $j\in L_s$, $a\leq n_j$, and $\xi_r,\eta_r\in\Gamma^s$, for $r=1,\dots,a$.
   
    \item It is an h-node (horizontal node), i.e. $\gamma\in \Gamma^{\mathrm{(hor)}}:=\bigcup_{n=0}^\infty \Delta_n^{\mathrm{(hor)}}\subset\Gamma^0$, with evaluation analysis of the form 
    \begin{equation}\label{eva analysis hor}
    e_\gamma^* = \sum_{r=1}^a d_{\xi_r}^* + \frac{1}{m_j}\sum_{r=1}^a \frac{1}{n_j^{1/p}} e_{\eta_r}^*\circ  P_{E_r}
    \end{equation}
    with $j\in L_0\setminus\{0\}$, $a\leq n_j$,  $\xi_r\in\Gamma^{\mathrm{(hor)}}$ and $\eta_r\in \Gamma$, for $r=1,\dots,a$, and $\mathrm{supp}_\mathrm{h}(e^*_{\eta_1}\circ P_{E_1})\ll\cdots\ll\mathrm{supp}_\mathrm{h}(e^*_{\eta_a}\circ P_{E_a})$.
    
    \item It is a g-node (ground node), i.e. $\gamma\in \Gamma^{\mathrm{(gr)}}:=\bigcup_{n=0}^\infty \Delta_n^{\mathrm{(gr)}}\subset\Gamma^0$, with evaluation analysis of the form 
    \begin{equation}\label{eva analysis gr} e_\gamma^* = \sum_{r=1}^a d_{\xi_r}^* + \frac{1}{m_0}\sum_{r=1}^a \lambda_r d_{\eta_r}^*
    \end{equation}
    with  $\xi_r\in\Gamma^{\mathrm{(gr)}}$, $a\leq \min\supp(e_\gamma^*)$, $\eta_r\in \Gamma^r$, for $r=1,\dots,a$, and $(\lambda_1,\dots,\lambda_a)\in\mathcal{B}$.

\end{enumerate}
Also, for all types, $\rank(\eta_1)<\rank(\xi_1)<\rank(\eta_2)<\rank(\xi_2)<\cdots$, $E_1\subset [0,\rank(\xi_1))$, and $E_r\subset(\rank(\xi_{r-1}),\rank(\xi_r))$, for $2\leq r\leq \alpha$ (where for g-nodes, $E_r = \rank(\xi_r)$ and $d_{\xi_r}^* = e_{\eta_r}^*\circ P_{\{\rank(\eta_r)\}}$). Note that the form of the evaluation analysis of g-nodes follows from the fact that for  $\gamma=(n,m_0^{-1},E,\lambda,\eta)\in\Gamma^{\mathrm{(gr)}}$, $\mathrm{supp}_\mathrm{h}(e^*_\gamma)=\{1\}$, whereas  for $\gamma=(n,m_0^{-1},\xi,E,\lambda,\eta)\in\Gamma^{\mathrm{(gr)}}$,  $\max\mathrm{supp}_\mathrm{h}(e^*_\gamma)=\max\mathrm{supp}_\mathrm{h}(e^*_\xi)+1 = \mathrm{age}(\gamma)$.

The above discussion of the evaluation analysis of nodes in $\Gamma$ implies immediately that all assumptions of a bonding (Definition \ref{bonding definition}) are satisfied.
\begin{proposition}
\label{XU is bonding of Xs}
The space $\mathfrak{X}_U$ is a bonding of  $(\mathfrak{X}_s)_{s=1}^\infty$.
\end{proposition}

\begin{remark}
The three types of nodes discussed above refer to the different ingredients in the construction of $\mathfrak{X}_U$. For any fixed $s\in\N$ the self-determined set $\Gamma^s$, consisting only of simple nodes, defines the space $\mathfrak{X}_s$. The main difference between $\mathfrak{X}_s$ and the original Argyros-Haydon space, lies in coefficients of the form $n_j^{-1/p}m_j^{-1}$ instead of $m_j^{-1}$, which yield $q$-convex estimates in $\mathfrak{X}_s$. The definition of h-nodes repeats this scheme on the horizontal level, imbuing the horizontal structure of $\mathfrak{X}_U$ with tightness while achieving an upper $U$ bound. The collection of g-nodes provides the partial lower $U$ bound.
\end{remark}

\addtocontents{toc}{\SkipTocEntry}
\subsection*{Existence of nodes}
In the sequel, we will need to construct various types of nodes in abundance to fulfill their roles in proving the properties of $\mathfrak{X}_U$. In this context, we now outline how each type can be created by using other nodes as building blocks.

\begin{proposition}[Existence of h-nodes]
\label{existence of h-nodes}
Let $2j\in L_0\setminus\{0\}$, $1\leq a\leq n_{2j}$, and $n_{2j}\leq p_1<p_2<\cdots<p_a$. Let also $\eta_1,\ldots,\eta_a\in\Gamma$ and $E_1,\ldots,E_a$ be intervals of $\N_0$ such that the following are satisfied.
\begin{enumerate}[label=(\alph*),leftmargin=19pt]
    
    \item $\mathrm{rank}(\eta_1) \in [0, p_1)$ and $E_1\subset [0,p_1)$ and, for $2\leq r\leq a$, $\mathrm{rank}(\eta_r)\in (p_{r-1}, p_r)$ and $E_r\subset (p_{r-1}, p_r)$.

    \item for $2\leq r\leq a$, $\min\mathrm{supp}_\mathrm{h}(e_{\eta_r}^*\circ P_{E_r}) >  p_{r-1}$ and $\mathrm{supp}_\mathrm{h}(e_{\eta_{r-1}}^*\circ P_{E_{r-1}})\ll \mathrm{supp}_\mathrm{h}(e_{\eta_r}^*\circ P_{E_r})$.
    
\end{enumerate}
Then there exist h-nodes $\xi_1,\xi_2,\ldots,\xi_a = \gamma\in\Gamma^0$ with $\mathrm{rank}(\xi_r) = p_r$, for $1\leq r\leq a$, such that
\[e_\gamma^* = \sum_{r=1}^ad_{\xi_r}^* + \frac{1}{m_{2j}}\sum_{r=1}^a\frac{1}{n_{2j}^{1/p}}e_{\eta_r}^*\circ P_{E_r}.\]
\end{proposition}

\begin{proof}
Let $\xi_1 = (p_1,m_{2j}^{-1},E_1,n_{2j}^{-1/p},\eta_1)$ (as in item \ref{initial even h-node} on page \pageref{initial even h-node}) and, assuming that $\xi_1,\ldots,\xi_{r-1}$ have been defined, let $\xi_r = (p_r,m_{2j}^{-1},\xi_{r-1},E_r,n_{2j}^{-1/p},\eta_r)$ (as in item \ref{successor even h-node} on page \pageref{successor even h-node}). Then, $c_{\xi_1}^* = m_{2j}^{-1}n_{2j}^{-1/p} e_{\eta_1}^*\circ P_{E_1}$ and, for $2\leq r\leq a$, $c_{\xi_r}^* = e_{\xi_{r-1}}^* + m_{2j}^{-1}n_{2j}^{-1/p} e_{\eta_r}^*\circ P_{E_r}$. Use $d_{\xi_r}^* = e_{\xi_r}^* - c_{\xi_r}^*$ to deduce the desired conclusion.
\end{proof}

\begin{remark}
\label{plentitude of even h-nodes}
For every $2j\in L_0\setminus\{0\}$, $1\leq a\leq n_{2j}$ and $N\in\mathbb{N}$ there are many nodes $\gamma\in\Gamma^0$ of weight $m_{2j}^{-1}$ and $\mathrm{age}(\gamma) = a$, with $\mathrm{rank}(\gamma)$ arbitrarily large and $\min\mathrm{supp}_\mathrm{h}(e_\gamma^*)\geq N$. To see this, note, e.g., that there exist simple nodes $\eta$ with $\mathrm{rank}(\eta)$ and $\mathrm{index}(\eta)$ arbitrarily large, such as ones of zero weight; these can be used in Proposition \ref{existence of h-nodes}, because, for any such $\eta$ and interval $E$ of $\N_0$, $\mathrm{supp}_\mathrm{h}{(e_{\eta}^*\circ P_{E})}\subset\{\mathrm{index}(\eta)\}$.

Note that Proposition \ref{existence of h-nodes} can be applied to sequences of such nodes to create new ones of higher complexity.
\end{remark}

\begin{remark}
\label{plentitude of odd h-nodes}
A process analogous to the proof of Proposition \ref{existence of h-nodes} allows building, for every $2j-1\in L_0\setminus\{0\}$, h-nodes of odd weight $m_{2j-1}^{-1}$. Let $2j-1\leq p_1$ and $\eta_1\in\Gamma$ such that $\mathrm{weight}(\eta_1) = m^{-1}_{4i-2}$, for some $4i-2\in L_0\setminus\{0\}$ with $m_{4i-2} > n_{2j-1}^3$; by Remark \ref{plentitude of even h-nodes}, such $\eta_1$ can be found. Let $p_1>\mathrm{rank}(\eta_1)$, and an interval $E_1\subset [0,p_1)$ of $\N$, define $\xi_1 = (p_1,m_{2j-1}^{-1},E_1,n_{2j-1}^{-1/p},\eta_1)$. To proceed, we need $p_2$, $\eta_2$, and $E_2$ such that the same assumptions as in Proposition \ref{existence of h-nodes} are satisfied, but additionally also $\mathrm{weight}(\eta_2) = m^{-1}_{\sigma(\xi_1)}$; by Remark \ref{plentitude of even h-nodes} we can find such $\eta_2$. Define $\xi_2 = (p_2,m_{2j-1}^{-1},\xi_{1},E_2,n_{2j-1}^{-1/p},\eta_2)$ and proceed.

In practice (see Lemma \ref{existence h-dependent}), the nodes $\eta_1,\ldots,\eta_a$ need to be chosen with great care when we study bounded linear operators on $\mathfrak{X}_U$.
\end{remark}

The same method as in the proof of Proposition \ref{existence of h-nodes} yields the following.

\begin{proposition}[Existence of simple nodes]
\label{existence of simple even nodes}
Let $s\in\mathbb{N}$, $2j\in L_s$, $1\leq a\leq n_{2j}$, and $(s+1)\vee n_{2j}\leq p_1<p_2<\cdots<p_a$. Let also $\eta_1,\ldots,\eta_a\in\Gamma^s$ and $E_1,\ldots,E_a$ be intervals of $\N_0$ such that $\mathrm{rank}(\eta_1) \in [0, p_1)$ and $E_1\subset [0,p_1)$ and, for $2\leq r\leq a$, $\mathrm{rank}(\eta_r)\in (p_{r-1}, p_r)$ and $E_r\subset (p_{r-1}, p_r)$. Then there exist simple nodes $\xi_1,\xi_2,\ldots,\xi_a = \gamma\in\Gamma^s$ with $\mathrm{rank}(\xi_r) = p_r$, for $1\leq r\leq a$, such that
\[e_\gamma^* = \sum_{r=1}^ad_{\xi_r}^* + \frac{1}{m_{2j}}\sum_{r=1}^a\frac{1}{n_{2j}^{1/p}}e_{\eta_r}^*\circ P_{E_r}.\]
\end{proposition}

\begin{remark}
\label{plentitude of simple even nodes}
For every $s\in\N$, $2j\in L_s$, and $1\leq a\leq n_{2j}$ there are many simple nodes $\gamma$ in $\Gamma^s$ with $\mathrm{weight}(\gamma) = m_{2j}^{-1}$, $\mathrm{age}(\gamma) = a$, and $\mathrm{rank}(\gamma)$ arbitrarily large. This can, for example, be proved as follows. For $a = 1$, using that $\Gamma_s^s = \{(s,0)\}$ we can define, for any $n\geq 2j$ and interval $E$ of $\{0,\ldots,n-1\}$, the node $\gamma = (n,m_{2j}^{-1},E,n_{2j}^{-1/p},(s,0))$. For $1<a\leq n_{2j}$ we can apply Proposition \ref{existence of simple even nodes} to such nodes.

This process can be iterated to create simple nodes of increasing complexity.
\end{remark}

\begin{remark}
\label{plentitude of simple odd nodes}
The process of creating, for any $2j-1\in L_s$ and $1\leq a\leq n_{2j-1}$, simple nodes $\gamma$ in $\Gamma^s$ with $\mathrm{weight}(\gamma) = m_{2j-1}^{-1}$, $\mathrm{age}(\gamma) = a$, and $\mathrm{rank}(\gamma)$ arbitrarily large, is similar to the one outlined in Remark \ref{plentitude of odd h-nodes} and it uses Proposition \ref{existence of simple even nodes} and Remark \ref{plentitude of simple even nodes}.
\end{remark}

The existence of g-nodes is proved similarly to the existence of h-nodes.

\begin{proposition}[Existence of g-nodes]
\label{existence of g-nodes}
Let $a\in\N$, $\eta_1\in\Gamma^1,\ldots,\eta_a\in\Gamma^a$, $1 \leq p_1<p_2<\cdots<p_a$, and, for $1\leq r\leq a$, $\lambda_r\in{D}_{p_r}$ such that the following are satisfied.
\begin{enumerate}[label=(\alph*)]

\item $a\leq \mathrm{rank}(\xi_1)$, i.e., the Schreier condition is satisfied.

\item $\mathrm{rank}(\xi_1)\in [0,p_1)$ and, for $2\leq r\leq a$, $\mathrm{rank}(\xi_r)\in(p_{r-1},p_r)$.

\item $(\lambda_1,\ldots,\lambda_a)\in\mathcal{B}$.

\end{enumerate}
Then, there exist g-nodes $\xi_1,\ldots,\xi_a=\gamma\in\Gamma^0$ with $\mathrm{rank}(\xi_r) = p_r$, for $1\leq r\leq a$, such that
\[e_\gamma^* = \sum_{r=1}^ad_{\xi_r}^* + \frac{1}{m_0}\sum_{r=1}^a\lambda_rd_{\eta_r}^*.\]
\end{proposition}

\begin{remark}
\label{plentitude of g-nodes}
A wealth of g-nodes can be found by applying Proposition \ref{existence of g-nodes} to Remark \ref{plentitude of simple even nodes} or \ref{plentitude of simple odd nodes}.
\end{remark}

\section{The Calkin algebra of $\mathfrak{X}_U$}
\label{section calkin algebra XU}

In this section, we prove the main result, i.e., Theorem \ref{the theorem}. This is shown in two steps. In the first, we relatively easily obtain the representation of the unitization of $U$ in $\mathpzc{Cal}(\mathfrak{X}_U)$, by showing that $\mathfrak{X}_U$ satisfies the upper $U$ condition and partial lower $U$ condition, by design. The argument is in most regards identical to the one of Proposition \ref{BU projections are U} for the space $\mathfrak{B}_U$. The second step in the proof of the main theorem entails showing that this embedding of the unitization of $U$ is onto. This is performed by assuming two theorems about the paucity of operators defined on $\mathfrak{X}_U$ and its subspaces $\mathfrak{X}_s$. The proof of these theorems requires the use of HI methods, i.e., rapidly increasing sequences and the basic inequality and they are postponed until later.

\begin{theorem}
\label{XU projections are U}
The space $\mathfrak{X}_U$ is a bonding of $(\mathfrak{X}_s)_{s=1}^\infty$ that satisfies the upper $U$-condition and the partial lower $U$ condition with constant $1/8$. In particular, the following hold.
\begin{enumerate}[label=(\roman*),leftmargin=19pt]

\item For every $s\in\N$, $I_s$ is a norm-one projection onto $\mathfrak{X}_s$, which is isometrically isomorphic to $\mathfrak{X}_{(\Gamma_n^s,i_n^s)}$.

\item The sequence $(I_s)_{s=1}^\infty$, in both $\mathcal{L}(\mathfrak{X}_U)$ and $\mathpzc{Cal}(\mathfrak{X}_U)$, is a normalized basic sequence 128-equivalent to $(u_s)_{s=1}^\infty$ and the spaces $\mathcal{K}(\mathfrak{X}_U)$ and $[(I_s)_{s=1}^\infty]\oplus \mathbb{C}I$ form a direct sum in $\mathcal{L}(\mathfrak{X}_U)$.

\end{enumerate}

\end{theorem}

\begin{proof}
The first conclusion is a reiteration of Theorem \ref{bonding general statement} \ref{bonding general statement 1}, applicable because, by Proposition \ref{XU is bonding of Xs}, $\mathfrak{X}_U$ is a bonding of $(\mathfrak{X}_s)_{s=1}^\infty$. Once we show that it also satisfies the upper $U$ condition and the partial lower $U$ condition, then the second conclusion is a reiteration of Theorem \ref{bd condition on projections}.

The proof of the partial lower $U$ condition follows verbatim the proof for the space $\mathfrak{B}_U$ (Section \ref{elementary BD example with U in Calkin}, Proposition \ref{BU projections are U}) using ground nodes.

For the upper $U$ condition, pick $\gamma\in\Gamma^0$. If $\gamma\in\Gamma^{\mathrm{(gr)}}$, then follow word-for-word the argument for the upper $U$ condition of the space $\mathfrak{B}_U$ (see the paragraph before Proposition \ref{BU projections are U}). Assume, therefore, $\gamma\in\Gamma^{\mathrm{(hor)}}$ with
\[e_\gamma^* = \sum_{r=1}^a d_{\xi_r}^* + \frac{1}{m_j}\sum_{r=1}^a \frac{1}{n_j^{1/p}} e_{\eta_r}^*\circ  P_{E_r},\]
where $\vec\lambda(\gamma)=(m_j^{-1}n_j^{-1/p},\dots,m_j^{-1}n_j^{-1/p})$, and $\mathrm{supp}_\mathrm{h}(e_{\eta_1}^*\circ P_{E_1})\ll\cdots\ll\mathrm{supp}_\mathrm{h}(e_{\eta_a}^*\circ P_{E_a})$. Let $f_1,\dots,f_a$ be functionals in the unit ball of $U^*$ such that, for $1\leq r\leq a$, $\mathrm{supp}(f_r)\subset\mathrm{supp}_\mathrm{h}(e_{\eta_r}^*\circ  P_{E_r})$. In particular, $\supp(f_1)\ll\cdots\ll\mathrm{supp}(f_r)$. Therefore, by Assumption \ref{dual p assumption},
\[\Big\|\sum_{r=1}^a \frac{1}{m_j}\frac{1}{n_j^{1/p}}f_r\Big\|_{U^*} \leq \frac{1}{m_j}\theta_U^{-1}\frac{1}{n_j^{1/p}} a^{1/p}\leq 1,\]
because $a\leq n_j$.
\end{proof}

The first theorem that we assume concerns the scalar-plus-compact property of $\mathfrak{X}_s$ as a subspace of $\mathfrak{X}_U$.

\begin{theorem}
\label{coordinate scalar-plus-compact}
For every $s\in\mathbb{N}$ and bounded linear operator $T:\mathfrak{X}_s\to\mathfrak{X}_U$ there exists $\lambda\in\mathbb{C}$ such that $T-\lambda J$ is compact, where $J:\mathfrak{X}_s\to\mathfrak{X}_U$ denotes the inclusion operator.
\end{theorem}

The above is in fact equivalent to a combination of two statements. The first one is that $\mathfrak{X}_s$ has the scalar-plus-compact property, i.e., for every $T:\mathfrak{X}_s\to\mathfrak{X}_s$ there exists $\lambda\in\mathbb{C}$ such that $T-\lambda I$ 
is compact (Theorem \ref{xs scalar+compact}, proved in Section \ref{subsection xs-dependent} ). The second one is that for every $T:\mathfrak{X}_s\to\mathfrak{X}_U$, the operator $(I-I_s) \circ T$ is compact (Proposition \ref{same-component-reduction}, proved in Section \ref{incomparability of structures}). The equivalence of Theorem \ref{coordinate scalar-plus-compact} to the conjunction of these two statements immediately follows by writing a $T:\mathfrak{X}_s\to\mathfrak{X}_U$ as $I_s\circ T + (I-I_s)\circ T$ and identifying $I_s\circ T$ with an operator on $\mathfrak{X}_s$.

Although this is somewhat concealed, an important role in the proof of Theorem \ref{coordinate scalar-plus-compact} is played by the horizontal structure of $\mathfrak{X}_U$. To make this term more precise, we introduce, analogously to the case of functionals, for any $x\in\mathfrak{X}_U$, the horizontal support of $x$.

\begin{definition}~\phantom{A}
\label{horizontal vector support definition}
\begin{enumerate}[label=(\arabic*),leftmargin=19pt]

\item For $x\in\mathfrak{X}_U$ we define the horizontal support of $x$ as the subset of $\mathbb{N}$ given by
\[\mathrm{supp}_\mathrm{h}(x) = \Big\{s\in\mathbb{N}: d_\gamma^*(x)\neq 0\text{ for some }\gamma\in\Gamma^s\Big\}.\]

\item A sequence $(x_n)_{n=1}^\infty$ in $\mathfrak{X}_U$ is called horizontally block (h-block) if it is block with respect to the FDD $(Z_n)_{n=0}^\infty$ (cf. Proposition \ref{remark general BD spaces}) and the horizontal structure given by $(\mathfrak{X}_s)_{s=1}^\infty$, i.e., $\supp(x_n)<\supp(x_{n+1})$ and $\mathrm{supp}_\mathrm{h}(x_n)<\mathrm{supp}_\mathrm{h}(x_{n+1})$ for each $n\in\mathbb{N}$.

\end{enumerate}
\end{definition}

\begin{remark}
\label{remark h-block in X0}
As with functionals, the horizontal support of a vector does not account for its interaction with $\Gamma^0$. In particular, for any $x\in\mathfrak{X}_0$, $\mathrm{supp}_\mathrm{h}(x) = \emptyset$ and, thus, any block sequence $(x_n)_{n=1}^\infty$ in $\mathfrak{X}_0$ is automatically h-block.
\end{remark}

We define horizontally compact operators. Traditionally, this notion is defined with respect to an infinite dimensional Schauder decomposition of a Banach space and an operator with this property becomes ``small'' after erasing finitely many (infinite dimensional) components of the decomposition (see, e.g., {\cite[Definition 3]{argyros:raikoftsalis:2008})}, {\cite[Definition 7.1]{zisimopoulou:2014}}, and {\cite[Definition 5.2]{pelczar-barwacz:2023}}).

\begin{definition}
\label{h-compact definition}
A bounded linear operator $T:\mathfrak{X}_U\to Y$, where $Y$ is a Banach space, is called horizontally compact if for every bounded h-block sequence $(x_n)_{n=1}^\infty$ in $\mathfrak{X}_U$, $\lim_n\|Tx_n\| = 0$.
\end{definition}

Here, just as in {\cite{pelczar-barwacz:2023}}, the situation is more delicate than the traditional one because the space $\mathfrak{X}_U$ does not admit an infinite dimensional Schauder decomposition. Although the aforementioned intuition of ``smallness'' is still accurate, it is somewhat more delicate to handle. In fact, the horizontal compactness of a $T:\mathfrak{X}_U\to\mathfrak{X}_U$ is equivalent to the compactness of $T|_{\mathfrak{X}_0}$. We only use one direction of this equivalence which is rather elementary.

\begin{remark}
\label{h-compact on X0 compact}
If a bounded linear operator $T:\mathfrak{X}_U\to\mathfrak{X}_U$ is horizontally compact then, automatically, $T|_{\mathfrak{X}_0}:\mathfrak{X}_0\to \mathfrak{X}_U$ is compact. More specifically, $\lim_nT(I-P_n)|_{\mathfrak{X}_0}=0$ in operator norm. This is because every bounded block sequence in $\mathfrak{X}_0$ is automatically a bounded h-block sequence.
\end{remark}

The second theorem that we will assume in this section is that $\mathfrak{X}_U$ has the scalar-plus-horizontally compact property. This will be proved in Section \ref{subsection h-dependent sequences} (Theorem \ref{scalar+h-compact}).

\begin{theorem}
\label{scalar-plus-horizontally compact}
For every bounded linear operator $T:\mathfrak{X}_U\to\mathfrak{X}_U$ there exists $\lambda\in\mathbb{C}$ such that $T-\lambda I$ is horizontally compact.    
\end{theorem}

The aforementioned theorems, in tandem with Theorem \ref{XU projections are U} and the Bourgain-Delbaen-$\mathscr{L}_\infty$ structure of the space $\mathfrak{X}_U$, allow us to complete the proof of the main result (Theorem \ref{the theorem}).

\begin{theorem}
\label{main theorem modulo HI}
Every bounded linear operator $T$ on $\mathfrak{X}_U$ can be written as
\[T = \lambda_0 I + \sum_{s=1}^\infty\lambda_sI_s+K,\]
with $\lambda_s\in\mathbb{C}$, $s\in\mathbb{N}_0$, and $K$ compact. Therefore, $\mathcal{L}(\mathfrak{X}_U) = \mathcal{K}(\mathfrak{X}_U)\oplus[(I_s)_{s=1}^\infty]\oplus\mathbb{C}I$.
\end{theorem}

\begin{proof}
By Theorem \ref{XU projections are U}, we already know that $ \mathcal{K}(\mathfrak{X}_U)$ and $[(I_s)_{s=1}^\infty]\oplus\mathbb{C}I$ form a direct sum in $\mathcal{L}(\mathfrak{X}_U)$. By Theorem \ref{scalar-plus-horizontally compact} there exists $\lambda_0\in\mathbb{C}$ such that $S = T- \lambda_0 I$ is horizontally compact and, by Theorem \ref{coordinate scalar-plus-compact}, for every $s\in\mathbb{N}$, there exists $\lambda_s\in\mathbb{C}$ such that $S\circ I_s - \lambda_sI_s$ is compact. Thus, for every $N\in\mathbb{N}$, the operator
\[\Bigg(T - \Big(\lambda_0I + \sum_{s=1}^N\lambda_sI_s\Big)\Bigg) - \Bigg(S - \sum_{s=1}^N S\circ I_s\Bigg) = \sum_{s=1}^N\Big(S\circ I_s - \lambda_s I_s\Big)\]
is compact. From this, we obtain
\begin{align*}
    \dist\Big(T,\mathcal{K}(\mathfrak{X}_U)\oplus[(I_s)_{s=1}^\infty]\oplus\mathbb{C}I\Big)&\leq \dist\Big(T - \Big(\lambda_0I + \sum_{s=1}^N\lambda_sI_s\Big),\mathcal{K}(\mathfrak{X}_U)\Big)\\&=  \mathrm{dist}\Big(S - \sum_{s=1}^N S\circ I_s,\mathcal{K}(\mathfrak{X}_U)\Big).
\end{align*}
Therefore, it suffices to show that the last term goes to zero, as $N$ tends to infinity. Towards a contradiction assume that there exist $\varepsilon>0$ and an infinite $L\subset \mathbb{N}$ such that for each $N\in L$
\begin{equation}
\label{main theorem modulo HI eq1}
\mathrm{dist}\Big(S-\sum_{s=1}^NS\circ I_s,\mathcal{K}(\mathfrak{X}_U)\Big) > \varepsilon.
\end{equation}
We will construct, for arbitrarily large $m_0,N_0\in\mathbb{N}$, a vector $x\in\mathfrak{X}_U$ such that $\|x\|\leq 8$, $\min\mathrm{supp}(x) \geq m_0$, $\min\mathrm{supp}_\mathrm{h}(x) \geq N_0$, and $\|Sx\| \geq \varepsilon/2$. By the horizontal compactness of $S$, this would be absurd.

By perhaps taking a larger $N_0$ we may assume $N_0\in L$ and by then taking a larger $k_0$, by Remark \ref{h-compact on X0 compact}, we may assume that
\[\Big\|SP_{[k_0,\infty)}|_{\mathfrak{X}_0}\Big\| < \frac{\varepsilon}{8N_0+24}.\]
By \eqref{main theorem modulo HI eq1} there exist a norm-one finitely supported vector $u$ with
\[\Big\|\Big(S-\sum_{s=1}^{N_0}S\circ I_s\Big)P_{[k_0,\infty)}u\Big\| > \varepsilon\]
and, thus, $w = P_{[k_0,\infty)}u$ is a finitely supported vector of norm at most $4$ and $\min\supp(w)\geq k_0$ such that
\[\Big\|S \Big(I - \sum_{s=1}^{N_0}I_s\Big)w\Big\| > \varepsilon.\]
Because $w$ is finitely supported, for some $N_1\geq N_0$ we may write
\[\Big(I - \sum_{s=1}^{N_0}I_s\Big)w = w_0 + \sum_{s=N_0+1}^{N_1}I_sw,\]
with $w_0\in\mathfrak{X}_0$. Note that $\min\supp(w_0), \min\supp(I_{N_0+1}w),\ldots,\min\supp(I_{N_1}w)\geq k_0$. By Proposition \ref{cancellation of coordinates} there exists $y_0\in\mathfrak{X}_0$ such that $\min\supp(y_0) \geq k_0$ and
\begin{equation}
\label{main theorem modulo HI eq2}
\Big\|\sum_{s=N_0+1}^{N_1}I_sw - y_0\Big\| \leq 2\max_{N_0<s\leq N_1}\|I_sw\| \leq 8.
\end{equation}
Define $x = (I - \sum_{s=1}^{N_0}I_s)w - (w_0 + y_0)$. Then, $\min\supp(w)\geq k_0$, $\min\mathrm{supp}_\mathrm{h}(x) \geq N_0+1$ and, by \eqref{main theorem modulo HI eq2}, $\|x\| \leq 8$. Finally, to compute $\|Sx\|$, note that $\|w_0 + y_0\|\leq \|x\| + \|I - \sum_{s=1}^{N_0}I_s\|\|w\| \leq 4N_0 + 12$ and, therefore,
\[\|Sx\| \geq \Big\|S\Big(I - \sum_{s=1}^{N_0}I_s\Big)w\Big\| - \big\|SP_{[k_0,\infty)}|_{\mathfrak{X}_0}\big\|\|w_0+y_0\| \geq \varepsilon - \frac{\varepsilon}{8N_0 + 24}(4N_0 + 12)\geq \frac{\varepsilon}{2}.\]
\end{proof}

\part{Concepts from HI methods}
\label{HI part}

\section{Estimates on the basis in $q$-mixed-Tsirelson spaces}
\label{section mixed-tsirelson estimates}

We will require estimates on averages and repeated averages of the basis, and these are provided by the theory of mixed-Tsirelson spaces. These estimates, which are transferred into $\mathfrak{X}_U$ via the Basic Inequalities (Proposition \ref{h-basic inequality} and Proposition \ref{h-basic inequality refined} and their component versions), are at the core of estimating the norm of linear combinations of block sequences of special types, such as rapidly increasing sequences (see Section \ref{RIS section}) and dependent sequences (see Section \ref{dependent sequences}). These types of sequences are an established tool that allows controlling the behavior of bounded operators on the space $\mathfrak{X}_U$.

Recall that we have fixed a value $1<p\leq \infty$ and its conjugate exponent $1\leq q<\infty$ such that Assumption \ref{dual p assumption} holds, and a pair of sequences $(m_j)_{j=1}^\infty$, $(n_j)_{j=1}^\infty$ satisfying Assumption \ref{assumption integers}. We first make introductory comments for a general $q$-mixed-Tsirelson space $T_q[(\mathcal{A}_{\tilde n_j},\theta_j)_{j\in\N}]$, with $(\tilde n_j)_{j=1}^\infty$ arbitrary, and then proceed to estimates for specific such spaces with parameters associated to $(m_j)_{j=1}^\infty$ and $(n_j)_{j=1}^\infty$.

Recall that a finite rooted tree is a finite set $\mathcal{T}$ with a partial order ``$\leq$'' with a minimum element denoted $\emptyset$ such that, for each $t\in\mathcal{T}$, the set $\{s\in\mathcal{T}:s\leq t\}$ is well ordered. For a non-terminal (i.e., non-maximal) node $t\in\mathcal{T}$, we denote $\mathrm{succ}(t)$ the set of immediate successors of $t$. We recall the following useful tool. 

\begin{definition}\label{tree-analysis} Fix a sequence $(\tilde n_j)_{j=1}^\infty$ in $\N$ and a sequence $(\theta_j)_{j\in\N}$ in $(0,1]$. A family $(f_t)_{t\in\mathcal{T}}$ in $W_q[(\mathcal{A}_{\tilde n_j},\theta_j)_{j\in\N}]$, where $\mathcal{T}$ is a finite rooted tree, is called a tree-analysis of $f\in W_q[(\mathcal{A}_{\tilde n_j},\theta_j)_{j\in\N}]$, provided
\begin{enumerate}[leftmargin=19pt]

\item $f_\emptyset=f$,

\item $t\in \mathcal{T}$ is a terminal (i.e. maximal) node if and only if  $f_t\in \{\pm e_n^*: n\in\N\}$.

\item for every non-terminal $t\in \mathcal{T}$ we have $f_t=\theta_j\tilde n_j^{-1/p}\sum_{r\in succ(t)}f_s$, for some $j\in\N$ and a block sequence $(f_r)_{r\in succ(t)}$, for a suitable enumeration of $succ(t)$.

\end{enumerate}

\end{definition}

\begin{remark}
It is well known, and it follows from the minimality of $W_q[(\mathcal{A}_{\tilde n_j},\theta_j)_{j\in\N}]$ among all $W$ satisfying the conclusions of Definition \ref{q-convex mT}, that every functional in the norming set $W_q[(\mathcal{A}_{\tilde n_j},\theta_j)_{j\in\N}]$ admits a tree-analysis. Although this analysis is not necessarily unique, this does not present any obstacle.
\end{remark}

The following separation of weights $(\theta_j)_j$ and factors $(\tilde n_j^{-1/p})_j$ on the level of the tree-analysis will be useful for future estimates.

\begin{definition}
Let $f\in W_q[(\mathcal{A}_{\tilde n_j}, \theta_j)_{j\in\N}]$ with a tree-analysis $(f_t)_{t\in\mathcal{T}}$. To it, we associate a collection $(g_t)_{t\in\T}$ recursively as follows.
\begin{enumerate}[leftmargin=19pt]
    
    \item For a terminal node $t\in\mathcal{T}$, put  $g_t=f_t$ and

    \item For a non-terminal node $t\in\mathcal{T}$ with $f_t=\theta_j\tilde n_j^{-1/p}\sum_{r\in succ(t)}f_r$, let $g_t= \tilde n_j^{-1/p}\sum_{r\in succ(t)}g_r$.
    
\end{enumerate}
We also let $g_f = g_\emptyset$.
\end{definition}

\begin{remark}
A straightforward induction starting at the terminal nodes of $\mathcal{T}$ yields $g_f\in\ell_p$.
\end{remark}

\begin{lemma}
\label{trivial combinatorial estimates}
Let $f\in W_q[(\mathcal{A}_{\tilde n_j}, \theta_j)_{j\in\N}]$ with a tree-analysis $(f_t)_{t\in\mathcal{T}}$ and, for $t\in\mathcal{T}$, let $j_t\in\mathbb{N}$ such that $\mathrm{weight}(f_t) = \theta_{j_t}$. The following hold.
\begin{enumerate}[label=(\roman*)]

\item\label{trivial combinatorial estimates i} For a terminal node $t\in\T$, if $f_t = \pm e_k^*$,
\[f(e_k) = \Big(\prod_{s< t}\mathrm{weight}(f_s)\Big)g_f(e_k)\text{ and }|f(e_k)| = \prod_{s< t}\frac{\mathrm{weight}(f_s)}{\tilde n_{j_s}^{1/p}}.\]

\item\label{trivial combinatorial estimates ii} For an antichain $\mathcal A$ in $\mathcal{T}$ with $\#\mathcal{A}>1$,
\[\#\mathcal A\leq \max_{t\in\mathcal A}\prod_{s<t}\tilde n_{j_s}.\]

\item\label{trivial combinatorial estimates iii} For an antichain $\mathcal{B}$ of $\mathcal{T}$ and $x\in c_{00}(\mathbb{N})$,
\[|f(x)| \leq \Big(\sum_{t\in\mathcal{B}}|f_t(x)|^q\Big)^{1/q}.\]
    
\end{enumerate}
\end{lemma}

\begin{proof}
For the first assertion, using the definition of $(g_t)_{t\in T}$, if $t = t_0>t_1>\cdots>t_m$ is a maximal chain in $\mathcal{T}$, it can be easily shown by induction that, for $\ell=1,\ldots,m$, $f_{t_\ell}(e_k) = (\prod_{i=1}^\ell\mathrm{weight}(f_{t_i}))g_{t_\ell}(e_k)$ and $|f_{t_\ell}(e_k)| = \prod_{i=1}^\ell\mathrm{weight}(f_{t_i})\tilde n_{j_{t_i}}^{-1/p}$. For the second assertion, for $s\in\mathcal{T}$, let $\mathcal{A}_s = \{t\in\mathcal A:s \leq t\}$. Then, an induction starting at the terminal nodes easily yields that for $s\in\mathcal{T}$ with $\#\mathcal{A}_s>1$, $\#\mathcal{A}_s\leq \max_{t\in \mathcal{A}_s}\prod_{s\leq r<t}\tilde n_{t_r}$. For the last assertion, for $t\in\mathcal{T}$, let $\mathcal{B}_t = \{s\in\mathcal{B}:s\geq t\}$. Then, an induction starting at the terminal nodes and the H\"older inequality easily yield that for $s\in\mathcal{T}$, $|f_t(x)| \leq (\sum_{s\in \mathcal{B}_t}|f_s(x)|^q)^{1/q}$.
\end{proof}

\begin{fact}\label{estimate mixed tsirelson - small supremum norm}
For $f\in W_q[(\mathcal{A}_{\tilde n_j}, \theta_j)_{j\in\N}]$ with tree-analysis $(f_t)_{t\in\mathcal{T}}$ and $x\in c_{00}(\N)$,
\begin{equation*}
|f(x)| \leq \|x\|_q\max\Big\{\frac{|f(e_k)|}{|g_f(e_k)|}:k\in\mathrm{supp}(f)\cap \mathrm{supp}(x)\Big\}.
\end{equation*}
\end{fact}

\begin{proof}
Write $x = \sum_{k=1}^na_ke_k$ and for $A = \mathrm{supp}(f)\cap\mathrm{supp}(x)$, by the H\"older inequality, we have
\[ |f(x)|\leq \sum_{k\in A}|a_kg_f(e_k)|\frac{|f(e_k)|}{|g_f(e_k)|} \leq \|x\|_q\|g\|_p\max\Big\{\frac{|f(e_k)|}{|g_f(e_k)|}:k\in A\Big\},\]
which ends the proof, as $g_f$ belongs to the closed unit ball of $\ell_p$.
\end{proof}
    
For the rest of the section, we use the fixed parameters $(n_j)_{j=1}^\infty, (m_j)_{j=1}^\infty\subset\N$ satisfying Assumption \ref{assumption integers}. 

The estimates for averages of the basis in the proposition below are based on the scheme presented in \cite{argyros:todorcevic:2005} (see also \cite[Proposition 2.5]{argyros:haydon:2011}), with some adjustments in the spirit of \cite[Proposition 3.4]{deliyanni:manoussakis:2007}. They will be used to prove estimates on special sequences of vectors called rapidly increasing sequences (see Definitions \ref{def hris} and \ref{def component ris}) in $\mathfrak{X}_U$, following the usual scheme in Bourgain-Delbaen spaces of Argyros-Haydon type.

\begin{proposition}\label{estimates mixed tsirelson}
Fix $j_0\in \N$ and $\theta\in [1,4]$. Then for any $f\in W_q[(\mathcal{A}_{4n_j}, \frac{\theta}{m_j})_{j\in\N}]$ with $\mathrm{weight}(f)=m_i^{-1}\theta$, $i\in\N$, we have
\begin{equation}\label{1 estimate mixed tsirelson}
\left| f\left(\frac{1}{n_{j_0}^{1/q}}\sum_{k=1}^{n_{j_0}}e_k\right)\right|\leq \begin{cases} \frac{4}{m_i}  \ \ \ \text{if} \ \ i\geq j_0, \\
\frac{17}{m_{j_0}m_i} \ \ \ \text{if} \ \ i<j_0.
\end{cases}
\end{equation}
In particular, in $T_q[(\mathcal{A}_{4n_j},\frac{\theta}{m_j})_{j\in\N}]$,
\begin{equation}\label{2 estimate mixed tsirelson}
\Big\|\frac{1}{n_{j_0}^{1/q}}\sum_{k=1}^{n_{j_0}}e_k\Big\|\leq \frac{4}{m_{j_0}}.
\end{equation}
Moreover, for any $f\in W_q[(\mathcal{A}_{4n_j}, \frac{\theta}{m_j})_{j\neq j_0}]$, we have
\begin{equation}\label{3 estimate mixed tsirelson}
\left| f\left(\frac{1}{n_{j_0}^{1/q}}\sum_{k=1}^{n_{j_0}}e_k\right)\right|\leq \frac{3}{m_{j_0}^2}.
\end{equation}
\end{proposition}

\begin{proof}
We prove first the estimate \eqref{3 estimate mixed tsirelson}.  Take $f\in W_q[(\mathcal{A}_{4n_j}, \frac{\theta}{m_j})_{j\neq j_0}]$ with a tree-analysis $(f_t)_{t\in\mathcal{T}}$. Note that, in particular, for every $t\in\mathcal{T}$, $\mathrm{weight}(f_t)\neq j_0$. For any terminal node $t\in\mathcal{T}$, let $l(t)$ be the length of the branch linking $t$ and the root $\emptyset$ of $\mathcal{T}$. We also pick, for every $k\in\supp(f)$, a terminal node $t_k\in\mathcal{T}$ with $e^*_k=f_{t_k}$.

Let
\[A=\Big\{k\leq n_{j_0}: \ k\in\supp (f_t) \text{ for some }t\in\mathcal{T} \text{ with }\mathrm{weight}(f_t)\leq {m_{j_0+1}^{-1}}{\theta}\Big\}.\]
By Lemma \ref{trivial combinatorial estimates} \ref{trivial combinatorial estimates i}, it follows that, for any $k\in A$,
\[\frac{|f(e_k)|}{|g_f(e_k)|}\leq \frac{\theta}{m_{j_0+1}}.\]
Thus, by Fact \ref{estimate mixed tsirelson - small supremum norm} and  Assumption \ref{assumption integers} \ref{2 assumption integers} we have
\[\left| f\left(\frac{1}{n_{j_0}^{1/q}}\sum_{k\in A}e_k\right)\right|\leq \frac{\theta}{m_{j_0+1}}\leq \frac{1}{m_{j_0}^2}.\]

Because, for every $t\in\mathcal{T}$, $\mathrm{weight}(f_t)\neq j_0$ we conclude that for every $k\in \mathrm{supp}(f)\setminus A$ and $s<t_k$, $\mathrm{weight}(f_s) \geq m_{j_0-1}^{-1}\theta$. Let
\[B=\Big\{k\in\supp (f)\setminus A: \ l(t_k)\leq \log_2(m_{j_0})\Big\}.\]
By Lemma \ref{trivial combinatorial estimates} \ref{trivial combinatorial estimates ii}, $\#B = \#\{t_k:k\in B\}\leq (4n_{j_0-1})^{\log_2(m_{j_0})}$. Then, by Assumption \ref{assumption integers} \ref{4 assumption integers},
\[\left| f\left(\frac{1}{n_{j_0}^{1/q}}\sum_{k\in B}e_k\right)\right|\leq\frac{\# B}{n_{j_0}^{1/q}}\leq
\frac{(4n_{j_0-1})^{\log_2(m_{j_0})}}{n_{j_0}^{1/q}}\leq \frac{1}{m_{j_0}^2}.\]

Finally, let $C=\{1,\dots,n_{j_0}\}\setminus (A\cup B)$. 
Note that, for any $k\in C$, by Lemma \ref{trivial combinatorial estimates} \ref{trivial combinatorial estimates i} and Assumption \ref{assumption integers} \ref{1 assumption integers},
\[\frac{|f(e_k)|}{|g_f(e_k)|}\leq \Big(\frac{\theta}{m_1}\Big)^{\log_2(m_{j_0})}\leq \Big(\frac{1}{4}\Big)^{\log_2(m_{j_0})}=\frac{1}{m_{j_0}^2}.\]
Thus, by Fact \ref{estimate mixed tsirelson - small supremum norm}
\[\left| f\left(\frac{1}{n_{j_0}^{1/q}}\sum_{k\in C}e_k\right)\right|\leq \frac{1}{m_{j_0}^2},\]
which ends the proof of \eqref{3 estimate mixed tsirelson}.

Now we prove \eqref{1 estimate mixed tsirelson}. Take a norming functional $f\in W_q[(\mathcal{A}_{4n_j}, \frac{\theta}{m_j})_{j\in\mathbb{N}}]$ with a tree-analysis $(f_t)_{t\in\mathcal{T}}$ and $\mathrm{weight}(f)=m_i^{-1}{\theta}$. The case $i\geq i_0$ follows easily from Lemma \ref{trivial combinatorial estimates} \ref{trivial combinatorial estimates i} and Fact \ref{estimate mixed tsirelson - small supremum norm}. We therefore assume $i<j_0$ and set
\[A'=\Big\{k\leq n_j: k\in\supp (f_t) \text{ for some $t$ with } \mathrm{weight}(f_t)\leq \frac{\theta}{m_{j_0}}\Big\}.\]
Then, Lemma \ref{trivial combinatorial estimates} \ref{trivial combinatorial estimates i}, for any $k\in A'$,
\[\frac{|f(e_k)|}{|g_f(e_k)|}\leq \frac{\theta^2}{m_{j_0}m_i}\leq \frac{16}{m_{j_0}m_i}\]
and by Fact \ref{estimate mixed tsirelson - small supremum norm} we have
\[\left| f\left(\frac{1}{n_{j_0}^{1/q}}\sum_{k\in A'}e_k\right)\right|\leq \frac{16}{m_{j_0}m_i}.\]
Let $D=\{1,\dots,n_{j_0}\}\setminus A'$. By putting $\mathcal{T}' = \{t\in\mathcal{T}:f_t|_D\neq 0\}$, it is easy to see that $(f_t|_D)_{t\in\mathcal{T}'}$ is a tree-analysis of $f|_D$ and, in particular, $f|_D\in W_q[(\mathcal{A}_{4n_j},\frac{\theta}{m_j})_{j\neq j_0}]$. Thus by \eqref{3 estimate mixed tsirelson} and Assumption \ref{assumption integers} \ref{2 assumption integers}
\[\left| f\left(\frac{1}{n_{j_0}^{1/q}}\sum_{k\in D}e_k\right)\right|\leq \frac{3}{m_{j_0}^2}\leq \frac{1}{m_{j_0}m_i},\]
which ends the proof of \eqref{1 estimate mixed tsirelson}. Finally, \eqref{2 estimate mixed tsirelson} is an immediate consequence of the definition of the norm of $T_q[(\mathcal{A}_{4n_j},\frac{\theta}{m_j})_{j\in\N}]$ and \eqref{1 estimate mixed tsirelson}.
\end{proof}

As in most HI-type constructions, the estimates in Proposition \ref{estimates mixed tsirelson} are applicable, via the basic inequality (see Proposition \ref{h-basic inequality refined}), to special sequences, called rapidly increasing sequences (see Section \ref{RIS section}). Usually, they would also be directly applicable to a higher-complexity type of sequences, called dependent sequences (see Section \ref{dependent sequences}). In the present paper, the latter fails due to the $q$-convex nature of the Argyros-Haydon components present in the construction of $\mathfrak{X}_U$. We therefore need additional estimates on repeated averages of the basis of a $q$-mixed-Tsirelson space, provided in the proposition below, and later combined with a refined basic inequality (see Proposition \ref{h-basic inequality refined}).

\begin{proposition}\label{estimates 2 mixed tsirelson}
Let $j_0\in \N\setminus\{1\}$, $\theta\in [1,4]$, and  $j_1<\dots<j_{n_{j_0}}$ such that $4m_{j_0}^2n_{j_0}^2 < m_{j_1}$ and, for $1\leq i<n_{j_0}$, $4m_{j_0}^2n_{j_0}^2m_{j_i}\leq m_{j_{i+1}}$. Consider a block sequence $(y_i)_{i=1}^{n_{j_0}}$ of vectors of the form
\[y_i=\frac{m_{j_i}}{n_{j_i}^{1/q}}\sum_{k\in I_i}e_k,\]
where, for $1\leq i\leq n_{j_0}$, $\# I_i=n_{j_i}$. Let $f\in W_q[(\mathcal{A}_{4n_j},\frac{\theta}{m_j})_{j\in \N}]$ with a tree-analysis $(f_t)_{t\in\mathcal{T}}$ such that, for any $t\in\mathcal{T}$ with $\mathrm{weight}(f_t)=m_{j_0}^{-1}\theta$, $\supp(f_t)$ intersects at most two of sets $I_i$, $1\leq i\leq n_{j_0}$. Then,
\[\left|f\left(\frac{1}{n_{j_0}^{1/q}}\sum_{i=1}^{n_{j_0}}y_i\right)\right|\leq \frac{5}{m_{j_0}^2}.\]
\end{proposition}

\begin{proof}
We will show that there is a functional $g\in W_q[(\mathcal{A}_{4n_j}, \frac{\theta}{m_j})_{j\neq j_0}]$ with 
\begin{equation*}
\left|f\left(\frac{1}{n_{j_0}^{1/q}}\sum_{i=1}^{n_{j_0}}y_i\right)\right|\leq g\left(\frac{1}{n_{j_0}^{1/q}}\sum_{i=1}^{n_{j_0}}e_i\right)+\frac{2}{m_{j_0}^2}    
\end{equation*}
which by, Proposition \ref{estimates mixed tsirelson}, yields 
\begin{equation*}
\left|f\left(\frac{1}{n_{j_0}^{1/q}}\sum_{i=1}^{n_{j_0}}y_i\right)\right|\leq \frac{5}{m_{j_0}^2}.
\end{equation*}
By 1-unconditionality of the basis and appropriately modifying the tree-analysis of $f$, we can assume that $\supp(f)\subset\cup_{i=1}^{n_{j_0}}I_i$ and $f_t(e_k)\geq 0$ for all $k$'s and $t$'s. In the first part of this proof, we will, up to a point, simplify the tree-analysis of $f$ using Proposition \ref{estimates mixed tsirelson}. In the second part, we use the simplified tree-analysis as a model to define a tree-analysis of a functional $g$ acting on $n_{j_0}^{-1/q}\sum_{i=1}^{n_{j_0}}e_i$.

Define
\[\tilde{\mathcal{B}} = \Big\{t\in\mathcal{T}:\text{ for some }1\leq i\leq n_{j_0},\;\mathrm{weight}(f_t) = m_{j_i}^{-1}\theta\text{ and }\mathrm{supp}(f_t)\cap I_i\neq \emptyset\Big\}\]
and let $\mathcal B$ be the collection of minimal members of $\tilde{\mathcal{B}}$. Note that $\mathcal{B}$ is an antichain of $\mathcal{T}$ and, in particular, the sets $\mathrm{range}(f_t)$, $t\in\mathcal{B}$, are pairwise disjoint intervals of $\mathbb{N}$. Putting $K = \mathrm{supp}(f)\setminus \cup_{t\in\mathcal{B}}\mathrm{supp}(f_t)$, for each $1\leq i\leq n_{j_0}$, $f|_K$ has a tree-analysis not involving any of $m_{j_1}^{-1}\theta,\ldots,m_{j_{n_{j_0}}}^{-1}\theta$. Therefore, by Proposition \ref{estimates mixed tsirelson} \eqref{3 estimate mixed tsirelson},
\[\Big|f|_K\Big(\frac{1}{n_{j_0}^{1/q}}\sum_{i=1}^{n_{j_0}}y_i\Big)\Big| \leq \frac{1}{n_{j_0}^{1/q}}\sum_{i=1}^{n_{j_0}}\frac{1}{m_{j_i}^2}\leq \frac{1}{m_{j_0}^2},\]
by the growth assumption on $m_{j_1},\ldots,m_{j_{n_{j_0}}}$. For $1\leq i\leq n_{j_0}$, denote
\[\mathcal{B}_i = \Big\{t\in\mathcal B:\mathrm{weight}(f_t) = m_i^{-1}\theta\Big\}.\]
By Proposition \ref{estimates mixed tsirelson} \eqref{1 estimate mixed tsirelson}, for fixed $1\leq i\neq \ell\leq n_{j_0}$ and $t\in \mathcal{B}_\ell$,
\[|f_t(y_i)| \leq \begin{cases}\frac{4m_{j_i}}{m_{j_\ell}}&:\text{ if }i<\ell\\\frac{17}{m_{j_\ell}}&:\text{ if }i>\ell.\end{cases}\]
In either case, by the growth assumption on $m_{j_1},\ldots,m_{j_{n_{j_0}}}$, the above bound is at most $5m_{j_0}^{-2}n_{j_0}^{-2}$. Because, for each $1\leq \ell\leq n_{j_0}$ and $t\in \mathcal{B}_\ell$, $\mathrm{supp}(f_t)\cap I_\ell\neq\emptyset$, it easily follow that, for fixed $1\leq i\neq \ell\leq n_{j_0}$, there is at most one $t\in \mathcal{B}_\ell$ such that $\mathrm{supp}(f_t)\cap I_i\neq\emptyset$. Therefore,
\[\Big(\sum_{t\in\cup_{\ell\neq i}\mathcal{B}_\ell}|f_t(y_i)|^q\Big)^{1/q} \leq n_{j_0}^{1/q}\frac{5}{m_{j_0}^{2}n_{j_0}^{2}}.\]

Let $J = \cup_{i=1}^{n_{j_0}}\cup_{t\in \mathcal{B}_i}(\mathrm{supp}(f_t)\cap I_i)$ and $f' = f|_J$. Then, by Lemma \ref{trivial combinatorial estimates} \ref{trivial combinatorial estimates ii},
\[
\begin{split}
\Big|\Big(f-f'\Big)\Big(\frac{1}{n_{j_0}^{1/q}}\sum_{i=1}^{n_{j_0}}y_i\Big)\Big| &\leq \Big|f|_K\Big(\frac{1}{n_{j_0}^{1/q}}\sum_{i=1}^{n_{j_0}}y_i\Big)\Big| + \frac{1}{n_{j_0}^{1/q}}\sum_{i=1}^{n_{j_0}}\Big(\sum_{t\in\cup_{\ell\neq i}\mathcal{B}_\ell}|f_t(y_i)|^q    \Big)^{1/q}\\
&\leq \frac{1}{m_{j_0}^2} + \frac{5}{m^2_{j_0}n_{j_0}} \leq \frac{2}{m_{j_0}^2}.
\end{split}
\]

We will now construct $g\in W_q[(\mathcal{A}_{4n_j}, \frac{\theta}{m_j})_{j\neq j_0}]$ such that
\[f'\Big(\frac{1}{n_{j_0}^{1/q}}\sum_{i=1}^{n_{j_0}}y_i\Big) \leq g\Big(\frac{1}{n_{j_0}^{1/q}}\sum_{i=1}^{n_{j_0}}e_i\Big).\]
For $t\in\mathcal{T}$, denote $\mathrm{weight}(f_t)=m_{j_t}^{-1}\theta$. For $1\leq i\leq n_{j_0}$, pick $t^0_i\in\mathcal{B}_i$ such that
\[\prod_{s\leq t_i^0}\frac{\theta}{m_{j_s}n_{j_s}^{1/p}} = \max\Big\{\prod_{s\leq t}\frac{\theta}{m_{j_s}n_{j_s}^{1/p}}:t\in\mathcal{B}_i\Big\}.\]
Define
\[\mathcal{T}' = \Big\{t\in\mathcal{T}:\text{ for some }1\leq i\leq n_{j_0},\;t\leq t_{i}^0\Big\}.\]
Starting at the maximal elements of $\mathcal{T}'$, i.e., at the nodes $t_i^0$,$1\leq i\leq n_{j_0}$, we recursively construct functionals $(h_t)_{t\in\mathcal{T}'}$ and $(g_t)_{t\in\mathcal{T}'}$ as follows.
\begin{enumerate}[leftmargin=19pt,label=(\roman*)]
    
    \item If, for some $1\leq i\leq n_{j_0}$, $t = t_i^0$ then $h_t = m_{j_i}^{-1}n_{j_i}^{-1/p}\sum_{k\in I_i}e_k^*$ and $g_t = e_i^*$.

    \item For non-terminal $t\in\mathcal{T}'$, if $f_t = m_{j_s}^{-1}n_{j_s}^{-1/p}\sum_{s\in\mathrm{succ}(t)}f_t$, then,
    \begin{enumerate}[label=(\greek*)]
        
        \item if $j_s\neq j_0$,
        \[h_t = \frac{\theta}{m_{j_s}n_{j_s}^{1/p}}\sum_{s\in\mathrm{succ}(t)\cap \mathcal{T}'}h_s\text{ and }g_t = \frac{\theta}{m_{j_s}n_{j_s}^{1/p}}\sum_{s\in\mathrm{succ}(t)\cap \mathcal{T}'}g_s\]

        \item if $j_s = j_0$,
        \[h_t = \frac{\theta}{m_{1}n_{1}^{1/p}}\sum_{s\in\mathrm{succ}(t)\cap \mathcal{T}'}h_s\text{ and }g_t = \frac{\theta}{m_{1}n_{1}^{1/p}}\sum_{s\in\mathrm{succ}(t)\cap \mathcal{T}'}g_s.\]
        
    \end{enumerate}
    
\end{enumerate}

An induction starting at the maximal elements of $\mathcal{T}'$ easily yields that, for $t\in\mathcal{T}'$,
\[\mathrm{supp}(h_t) = \cup\{I_i:t\leq t_i^0\}.\]

The first main consequence of this is that for $t,s\in\mathcal{T}'$ incomparable, $\mathrm{supp}(h_t)\cap \mathrm{supp}(h_s)=\emptyset$. The second one is that, if $\mathrm{supp}(h_t)\cap I_i\neq \emptyset$ then $\mathrm{supp}(f_t)\cap I_i\neq\emptyset$. This implies in particular the following: if, for $t\in\mathcal{T}'$, $\mathrm{weight}(f_t) = m_{j_0}^{-1}\theta$, then $\#(\mathrm{succ}(f_t)\cap\mathcal{T}')\leq 2 \leq n_1$. Otherwise, there would be pairwise distinct $t_{i_1}^0$, $t^0_{i_2}$, $t^0_{i_3}$, all greater than $t$, meaning $\mathrm{supp}(f_t)$ intersects $I_{i_1}$, $I_{i_2}$, $I_{i_3}$. This contradicts the initial assumption.

Combining the observations of the above paragraph, it can be easily shown, by induction starting that the terminal nodes of $\mathcal{T}'$, that $h_t, g_t\in W_q[(\mathcal{A}_{4n_j}, \frac{\theta}{m_j})_{j\neq j_0}]$ and, for any $1\leq i\leq n_{j_0}$, $h_t(y_i) = g_t(e_i)$. In particular, for $h = h_\emptyset$ and $g = g_\emptyset$,
\[h\Big(\frac{1}{n_{j_0}^{1/q}}\sum_{i=1}^{n_{j_0}}y_i\Big) = g\Big(\frac{1}{n_{j_0}^{1/q}}\sum_{i=1}^{n_{j_0}}e_i\Big).\]
To complete the proof, it suffices to show that the estimate of $h$ dominates that of $f'$. We will show the following sufficient statement: for $k\in\mathbb{N}$, $f'(e_k)\leq h(e_k)$. If $k\in\mathbb{N}\setminus\cup_{i=1}^{n_{j_0}}\cup_{t\in\mathcal{B}_i}\mathrm{supp}(f_t)\cap I_i$, then $f'(e_k) = 0$. Otherwise, take the unique $1\leq i\leq n_{j_0}$ and $t\in\mathcal{B}_i$ such that $k\in\mathrm{supp}(f_t)\cap I_i$. By Proposition \ref{trivial combinatorial estimates} \ref{trivial combinatorial estimates i},
\[f'(e_k) \leq \prod_{s\leq t}\frac{\theta}{m_{j_s}n_{j_s}^{1/p}} \leq \prod_{s\leq t_i^0}\frac{\theta}{m_{j_s}n_{j_s}^{1/p}} \leq  \Big(\prod_{s\leq t_i^0\atop j_s\neq j_0}\frac{\theta}{m_{j_s}n_{j_s}^{1/p}}\Big) \Big(\frac{\theta}{m_2n_2^{1/p}}\Big)^{\#\{s\leq t_i^0:j_s = j_0\}} = h(e_k).\]
\end{proof}

\section{Rapidly increasing sequences  in $\mathfrak{X}_U$}
\label{RIS section}    

Rapidly increasing sequences (RIS) are a universal tool in the study of HI, and similar spaces (e.g., {\cite[Lemma 6]{schlumprecht:1991}}, {\cite{gowers:maurey:1993}}, {\cite{argyros:deliyanni:1997}}, {\cite{deliyanni:manoussakis:2007}}, and {\cite{argyros:haydon:2011}}). Traditionally, a $C$-RIS is a bounded sequence of vectors $(x_k)_{k=1}^\infty$ whose defining property is the existence of bounds $|f(x_k)| \leq Cm^{-1}_i$ where $f$ is a weighted functional of a weight $m_i^{-1}$ with $i$ less than a certain index $j_k$ associated to $x_k$. In practice, the most commonly found such sequences are $\ell_1$-averages of increasing length of block vectors. In some constructions, these are substituted with repeated $\ell_1$-averages of increasing complexity (e.g., in {\cite{argyros:deliyanni:1997}}) or with  $\ell_p$-averages of increasing length (e.g., in {\cite{deliyanni:manoussakis:2007}}). In most cases, RISs are used to detect strictly singular operators or compact operators (e.g., in {\cite{argyros:haydon:2011}}). In modern HI constructions a tool called the ``basic inequality'' provides high-precision bounds on linear combinations of such sequences. Typically, the proof of the basic inequality is an induction that uses the standard representation of a weighted functional (such as $f = m_i^{-1}\sum_{r=1}^af_r$ in \cite[Proposition 4.3]{argyros:tolias:2004}) or the evaluation analysis of a node (such as $e_\gamma^* = \sum_{r=1}^ad_{\xi_r}^* + m_i^{-1}\sum_{r=1}^ae_{\eta_r}^*\circ P_{E_r}$ in {\cite{argyros:haydon:2011}}) and it relies on the estimates of the form $|f(x_k)| \leq Cm_i$ in a substantial way. This process stumbles in the $q$-convex versions of the Argyros-Haydon spaces used in this paper. This is due to interference caused by the ``residual part'' $\sum_{r=1}^ad_{\xi_r}^*$ in the evaluation analysis of a node with $q$-convex estimates. To overcome this, we design RISs in such a way that this residual part contributes to a cancellation effect. This allows for a surprisingly clean estimate, called an ``evaluation refinement''. With this approach, the requirement for estimates of the form $Cm_i^{-1}$ becomes redundant, and it does not appear in the definition of RISs. The drawback of the cancellation effect design is that it is of an algebraic nature, and it cannot be directly used to study the subspace structure of $\mathfrak{X}_U$. This does not have a negative impact on the study of operators whose domain is $\mathfrak{X}_U$ or one of the spaces $\mathfrak{X}_s$, $s\in\N$.

In this paper, the definition of a RIS uses the notion of the local support of a vector introduced in {\cite[Definition 5.7]{argyros:haydon:2011}}.

\begin{notation}
Let $x\in\mathfrak{X}_U$ with $n = \max\supp(x)$.
\begin{enumerate}[label=(\arabic*)]

\item We define the local support of $x$ to be the subset of $\Gamma$ given by
\[\supp_\mathrm{loc}(x) = \big\{\gamma\in\Gamma_n: e_\gamma^*(x)\neq 0\big\}.\]

\item For $S\subset \supp_\mathrm{loc}(x)$ denote by $x{\upharpoonright}_S$ the vector
\[x{\upharpoonright}_S = i_n\Big(\sum_{\gamma\in S}e_\gamma^*(x)e_\gamma\Big),\]
where the $e_\gamma$ are seen as vectors in $\ell_\infty(\Gamma_n)$.
\end{enumerate}
\end{notation}

The local support is a rather capricious notion. For example, if $x$ and $y$ have disjoint local supports $S_1$ and $S_2$ it is not true that $x+y$ has local support $S_1\cup S_2$. This is because the local support depends on the maximum of the support of a vector. For similar reasons, if $x,y$ are two vectors with $supp_\mathrm{loc}(x) = \supp_\mathrm{loc}(y)$ and $x(\gamma) = y(\gamma)$, for $\gamma\in S$, then it is not necessarily true that $x = y$; it may happen that $\mathrm{max}\supp(x) \neq \max\mathrm{supp}(y)$. Finally, if $s\in\N$ and $x$ is a finitely supported vector in $\mathfrak{X}_s = [(d_\gamma)_{\gamma\in\Gamma^s}]$ then it is not necessarily true that $\supp_\mathrm{loc}(x)\subset\Gamma^s$. The following lemma outlines some regularity properties of this notion.

\begin{lemma}
\label{local support lemma}
Let $x\in\mathfrak{X}_U$ with $m=\min\supp(x)$ and $n = \max\supp(x)$. The following hold.
\begin{enumerate}[label=(\roman*),leftmargin=23pt]

\item\label{local starting point} $\supp_\mathrm{loc}(x) \subset \Gamma_n\setminus\Gamma_{m-1}$, making the convention $\Gamma_{-1}=\emptyset$.

\item\label{local partition} For a partition $S_1,\ldots,S_k$ of $\supp_\mathrm{loc}(x)$, $x = x{\upharpoonright}_{S_1}+\cdots+x{\upharpoonright}_{S_k}$.

\item\label{local quantities} For $S\subset\supp_\mathrm{loc}(x)$ the following hold for the vector $x{\upharpoonright}_S$.
\begin{enumerate}[label=(\greek*)]
    \item\label{local restriction norm} $\|x{{\upharpoonright}}_S\| \leq 2\|x\|$.
    \item\label{local restriction BD support} $\supp(x{\upharpoonright}_S)\subset[m,n]$.
    \item\label{local restriction local support} $\supp_\mathrm{loc}(x{\upharpoonright}_S)\subset S$ (it may happen that the inclusion is proper).
    \item\label{local restriction h-support}\label{local h-support} $\mathrm{supp}_\mathrm{h}(x{\upharpoonright}_S) \subset\mathrm{supp}_\mathrm{h}(x)$.
\end{enumerate}

\end{enumerate}

\end{lemma}

\begin{proof}
To show \ref{local starting point} let $\gamma\in\supp_\mathrm{loc}(x)$ with minimal rank $k$ and, thus, $x|_{\Gamma_{k-1}} = 0$. This yields $c_\gamma^*(x) = 0$ and $d_\gamma^*(x) = e_\gamma^*(x) - c_\gamma^*(x) \neq 0$. In other words, $m = \min\supp(x)\geq k$. For \ref{local partition} note
\[x = i_n( x|_{\Gamma_n}) = i_n\Big(\sum_{\gamma\in\Gamma_n}e_\gamma^*(x)e_\gamma\Big) = i_n\Big(\sum_{d=1}^k\sum_{\gamma\in S_d}e_\gamma^*(x)e_\gamma\Big) = \sum_{d=1}^kx{\upharpoonright}_{S_d}.\]

The first part of \ref{local quantities} follow easily; $\|x{\upharpoonright}_S\| = \|i_n(\sum_{\gamma\in S}e_\gamma^*(x)e_\gamma)\| \leq 2\|x|_{S}\|_\infty \leq 2\|x\|$. For the second part, put $k=\min\{\mathrm{rank}(\gamma):\gamma\in S\} \geq m$. The same argument used in the proof of \ref{local starting point} yields $\min\supp(x{\upharpoonright}_S) = k$. Although, somewhat surprisingly, it may not be true that $\max\supp(x{\upharpoonright}_S) =\max\{\mathrm{rank}(\gamma):\gamma\in S\}$, it is true that $\max\supp(x{\upharpoonright}_S) \leq n$. This follows from $x{\upharpoonright}_S\in i_n(\ell_\infty(\Gamma_n)) = \langle\{d_\gamma:\gamma\in\Gamma_n\}\rangle$. To show the third part of \ref{local quantities}, put $l = \max\supp(x{\upharpoonright}_S)$. Then, $\supp_\mathrm{loc}(x{\upharpoonright}_S)\subset\Gamma_l$ and, for $\gamma\in\Gamma_l$, $x{\upharpoonright}_S(\gamma) = (x{\upharpoonright}_S)|_{\Gamma_l}(\gamma) = (x{\upharpoonright}_S)|_{\Gamma_n}(\gamma)$. This yields $\mathrm{supp}_\mathrm{loc}(x{\upharpoonright}_S) = S\cap \Gamma_l\subset S$. For the fourth part of \ref{local quantities} let $s\in\mathrm{supp}_\mathrm{h}(x{\upharpoonright}_S)$, i.e., there exists $\gamma\in\Gamma^s_n$ such that $d_\gamma^*(x{\upharpoonright}_S) \neq 0$. Because $\Gamma^s$ is self-determined, by {\cite[Proposition 1.5 (b)]{argyros:motakis:2019}}, $d_\gamma^*\in\langle\{e_\eta^*:\eta\in\Gamma^s_n\}\rangle$ and, thus, there exists $\eta\in\Gamma^s_n$ with $e_\eta^*(x{\upharpoonright}_S) = e_\eta^*(x)\neq 0$. In conclusion, $s\in\mathrm{supp}_\mathrm{h}(x)$.
\end{proof}

\begin{notation}
\label{Xi and local supp}
For $\gamma\in\Gamma$ and a finitely supported $x\in\mathfrak{X}_U$ define
\begin{align*}
\Xi(\gamma) &= \big\{\xi\in\Gamma: |e_\gamma^*(d_\xi)| > \mathrm{weight}(\gamma)\big\}\text{ and}\\
\Xi(\gamma,x) &= \Xi(\gamma)\cap\supp_\mathrm{loc}(x) = \big\{\xi\in\Xi(\gamma): e^*_\xi(x)\neq 0\text{ and } \mathrm{rank}(\xi)\leq \max\supp(x)\big\}.
\end{align*}
\end{notation}

The set $\Xi(\gamma)$ can be explicitly described in terms of the evaluation analysis of $\gamma$. Furthermore, for certain block vectors $x_1,\ldots,x_N$ in $\mathfrak{X}_U$ and an interval $E$ of $\mathbb{N}$ we can write a useful formula for the evaluation of $e_\gamma^*\circ P_E$ on linear combinations of the $x_k$; if $\Xi(\gamma,x_k) = \emptyset$ then the vectors $x_1,\ldots,x_N$ behave like unit vector basis elements.

\begin{proposition}[Evaluation refinement]
\label{dxi elimination lemma}
Let $\gamma\in\Gamma$ have non-zero weight and evaluation analysis
\[    e_\gamma^* = \sum_{r=1}^a d_{\xi_r}^* + \frac{1}{m_j}\sum_{r=1}^a \lambda_r e_{\eta_r}^*\circ  P_{E_r}.\]
Then, $\Xi(\gamma) = \{\xi_1,\ldots,\xi_a\}$ and, for $1\leq r\leq a$, $\mathrm{age}(\xi_r) = r$. 

Furthermore, let $E$, $I$ be intervals of $\mathbb{N}$ and $(x_k)_{k\in I}$ be block vectors in $\mathfrak{X}_U$ such that, for $k\in I$,
\begin{enumerate}[label=(\alph*)]
    \item $\min(E)\leq \min\supp(x_k)$ and
    \item $\Xi(\gamma,x_k)=\emptyset$.
\end{enumerate}
Then, there exist subintervals $I_1<\cdots<I_a$ of $I$ such that for any scalars $(a_k)_{k\in I}$,
\begin{equation}
\label{dxi elimination lemma formula}
e_\gamma^*\circ P_E\Big(\sum_{k\in I}a_kx_k\Big) = \frac{1}{m_j}\sum_{r=1}^a\lambda_re_{\eta_r}^*\circ P_{E_r\cap E}\Big(\sum_{k\in I_r}a_kx_k\Big).
\end{equation}
\end{proposition}

\begin{proof}
The statement $\Xi(\gamma) = \{\xi_1,\ldots,\xi_a\}$ follows from the fact that $\|d_\xi\| = 1$, for all $\xi\in\Gamma$ (see Remark \ref{remark basis vector norm}). Next, fix intervals $E$, $I$ of $\N$ and block vectors $(x_k)_{k\in I}$ in $\mathfrak{X}_U$ such that, for $k\in I$, $\min(E)\leq \min\supp(x_k)$ and $\Xi(\gamma,x_k) = \emptyset$. For $1\leq r\leq a$ define
\[I_r =\Big\{k\in I: \min(E_r)\leq \max(E)\wedge\max\supp(x_k) < \mathrm{rank}(\xi_{r})\Big\}.\]
Although some of the $I_r$ may be empty, because $E_1<\mathrm{rank}(\xi_1)<\cdots<E_a<\mathrm{rank}(\xi_a)$, they are intervals and $I_1<\cdots<I_a$. Fix $k\in I$ and write
\[e_\gamma^*\circ P_E(x_k) = \sum_{r=1}^a\Big(\frac{1}{m_j}\lambda_r e_{\eta_r}^*\circ  P_{E_r\cap E}(x_k) + d_{\xi_r}^*(P_Ex_k)\Big).\]
We will prove that for $1\leq r\leq a$,
\[\frac{1}{m_j}\lambda_r e_{\eta_r}^*\circ  P_{E_r\cap E}(x_k) + d_{\xi_r}^*(P_Ex_k) = \left\{\begin{array}{ll}\frac{1}{m_j}\lambda_je^*_{\eta_{r}}\circ P_{E_{r}\cap E}(x_k)&\text{ if }k\in I_{r}\\0&\text{ if }k\notin I_{r}.\end{array}\right.\]
This easily yields the desired conclusion.

For $1\leq r\leq a$ such that $\mathrm{rank}({\xi_r }) < \min\supp(x_k)$, because $\max(E_r)<\mathrm{rank}(\xi_r)$,
\[\frac{1}{m_j}e_{\eta_r}^*\circ P_{E_r\cap E}(x_k)+ d_{\xi_r}^*( P_Ex_k) = 0.\]
For $1\leq r\leq a$ such that $\min\supp(x_k)\leq \mathrm{rank}({\xi_r}) \leq \max(E)\wedge\max\mathrm{supp}(x_k)$ we will show that $d_{\xi_r}^*(P_Ex_k) =- m_j^{-1}\lambda_re_{\eta_r}^*\circ P_{E_r\cap E}(x)$, which yields
\[\frac{1}{m_j}e_{\eta_r}^*\circ P_{E_r\cap E}(x_k)+ d_{\xi_r}^*( P_Ex_k) = 0.\]
Firstly, for such $1\leq r\leq a$, because $\min(E) \leq \min\supp(x_k)$ and $E_r<\mathrm{rank}(\xi_r)\leq \max(E)$, we have $P_{E_r}(x_k) = P_{E_r\cap E}(x_k)$. Next, if $r>1$, then
\begin{align*}
d_{\xi_r}^*(P_Ex_k) &= d_{\xi_r}^*(x_k) = e_{\xi_r}^*(x_k) - c_{\xi_r}^*(x_k) = \underbrace{\big(e_{\xi_r}^*(x_k) - e_{\xi_{r-1}}^*(x_k)\big)}_{=0,\text{ by }\Xi(\gamma,x_k)=\emptyset} - \frac{1}{m_j}\lambda_re_{\eta_r}^*\circ P_{E_r}(x_k)\\
&=- \frac{1}{m_j}\lambda_re_{\eta_r}^*\circ P_{E_r\cap E}(x_k),
\end{align*}
while if $r = 1$ the same equality holds by omitting the non-existing term $e_{\xi_{r-1}}^*$.

For $1\leq r\leq a$ such that $\mathrm{rank}(\xi_r) > \max(E)\wedge\max\supp(x_k)$ we have $d_{\xi_r}^*(P_Ex_k) = 0$, and thus,
\[\frac{1}{m_j}e_{\eta_r}^*\circ P_{E_r\cap E}(x_k)+ d_{\xi_r}^* (P_Ex_k) = \frac{1}{m_j}e_{\eta_r}^*\circ P_{E_r\cap E}(x_k).\]
Unless $\min(E_r) \leq \max(E)\wedge\max\supp(x_k) < \mathrm{rank}(\xi_r)$, we have ${m_j}^{-1}e_{\eta_r}^*\circ P_{E_r\cap E}(x_k) = 0$ and this condition is satisfied if and only if $k\in I_r$.
\end{proof}

\begin{remark}
The subintervals $I_1<\cdots<I_r$ given by Proposition \ref{dxi elimination lemma} need not cover $I$ and, in fact, gaps between them may be necessary for formula \eqref{dxi elimination lemma formula} to be correct.
\end{remark}

\subsection{Horizontally rapidly increasing sequences}
We define rapidly increasing sequences in the context of the horizontal structure of $\mathfrak{X}_U$. As already explained in the beginning of Section \ref{RIS section}, estimates depending on the weight of nodes are not part of this definition and instead we use the notion $\Xi(\gamma,x)$. This is done to better deal with the $q$-convex component of our version of the Argyros-Haydon construction.

\begin{definition}
\label{def hris}
Let $C > 0$ and $N=1$ or $N=2$. A horizontally block sequence $(x_k)_{k\in I}$, indexed over some interval $I$ of $\N$, is called a $(C,N)$-horizontally rapidly increasing sequence (or $(C,N)$-h-RIS) if
the following hold.
\begin{enumerate}[label=(\alph*)]

    \item\label{ris norm} For $k\in I$, $\|x_k\|\leq C$.

    \item\label{ris dxi collapse} For every $\gamma\in\Gamma$, $\#\{k\in I:\Xi(\gamma,x_k)\neq\emptyset\}\leq N$.

    \item\label{ris h-support growth} For $k\in I\setminus\{\min(I)\}$, $\max\supp(x_{k-1}) < \min\mathrm{supp}_\mathrm{h}(x_k)$.

\end{enumerate}
We will call $(x_k)_{k\in I}$ an h-RIS if it is a $(C,N)$-h-RIS for some $C>0$ and $N=1$ or $N=2$. 
\end{definition}

Requirement \ref{ris h-support growth} in the above definition is introduced to counteract the Schreier admissibility of ground nodes; imposing it makes the estimation of such a node on a linear combination of a h-RIS inconsequential.

\begin{remark}\label{interwined RIS}
Let $(x_k)_{k=1}^\infty$ be a $(C,1)$-h-RIS and $(y_k)_{k=1}^\infty$ be a $(D,1)$-h-RIS such that the sequence $x_1,y_1,x_2,y_2,\ldots$ is h-block and satisfies Assumption \ref{ris h-support growth} of Definition \ref{def hris}. Then, the sequence $x_1,y_1,x_2,y_2,\ldots$ is a $(C\vee D,2)$-h-RIS.
\end{remark}

The following statement provides upper bounds for estimates on h-RISs that follow directly and are used in the proof of the basic inequality for such sequences.
\begin{lemma}
\label{pre h-ris lemma}
Let $(x_k)_{k\in I}$ be a $(C,2)$-h-RIS. 
\begin{enumerate}[label=(\roman*)]

\item\label{pre h-ris lemma simple node} If $\gamma$ has weight zero or if $\gamma\in\Gamma^s$, for some $s\in\mathbb{N}$, then for every sequence $(a_k)_{k\in I}$,
\[\Big|e_\gamma^*\circ P_E\Big(\sum_{k\in I} a_kx_k\Big)\Big| \leq 4C\max_{k\in I}|a_k|.\]

\item\label{pre h-ris lemma ground node} If $\gamma\in\Gamma^0$, with $\mathrm{weight}(\gamma) = m_0^{-1}$, i.e., $\gamma$ is a g-node, then
\[\Big|e_\gamma^*\circ P_E\Big(\sum_{k\in I} a_kx_k\Big)\Big| \leq 12C\max_{k\in I}|a_k|.\]

\end{enumerate}
\end{lemma}

\begin{proof}
The proof of \ref{pre h-ris lemma simple node} is fairly straightforward. If $\mathrm{weight}(\gamma) = 0$ then $e_\gamma^* = d_\gamma^*$, which easily yields the result. If $\gamma\in\Gamma^s$, for some $s\in\mathbb{N}$, then, for every $\xi\in\Gamma\setminus\Gamma^s$, $e_\gamma^*(d_\xi) = 0$ (because $e_\gamma^*\in\langle\{d_\xi^*:\xi\in\Gamma^s\}\rangle)$. There is at most one $k_0\in I$ with $s\in\mathrm{supp}_\mathrm{h}(x_{k_0})$ and, therefore,
\[\Big|e_\gamma^*\circ P_E\Big(\sum_{k\in I} a_kx_k\Big)\Big| = \Big|e_\gamma^*\circ P_E(a_{k_0}x_{k_0})\Big|\leq 4C|a_{k_0}|.\]

Assume that $\mathrm{weight}(\gamma) = m_0^{-1}$. Without loss of generality we may assume that for all $k\in I$, $\max\supp(x_k)\geq \min(E)$ (because, otherwise, $P_E(x_k) = 0$). Similarly, we may also assume that for $k_0 = \min(I)$, $e_\gamma^*\circ P_E(x_{k_0})\neq 0$. In particular, for $k\in I\setminus\{k_0\}$, $\min(E)\leq \min\supp(x_k)$. Put $S = \{k_0\}\cup\{k\in I: \Xi(\gamma,x_k)\neq\emptyset\}$, which has at most $3$ members and, thus,
\[\Big|e_\gamma^*\circ P_E\Big(\sum_{k\in S}a_kx_k\Big)\Big|\leq 12C\max_{k\in I}|a_k|.\]
Write the evaluation analysis of $\gamma$
\[e_\gamma^* = \sum_{r=1}^ad_{\xi_r}^*+\frac{1}{m_0}\sum_{r=1}^a \lambda_rd^*_{\eta_r}\]
By the definition of ground nodes, as $e_\gamma^*\circ P_E(x_{k_0})\neq 0$, $a = \age(\gamma) \leq \min\supp(e_\gamma^*)\leq    \max\supp(x_{k_0})$ and for $1\leq r\leq a$,  $\mathrm{index}(\eta_r) = r$. Apply Proposition \ref{dxi elimination lemma} to find subsets $F_1<\cdots<F_a$ of $I\setminus S$ such that
\begin{equation*}
e_\gamma^*\circ P_E\Big(\sum_{k\in I\setminus S}a_kx_k\Big) = \frac{1}{m_0}\sum_{r=1}^a\lambda_rd_{\eta_r}^*\circ P_E\Big(\sum_{k\in F_r}a_kx_k\Big)
\end{equation*}
The above quantity is in fact zero and this yields the desired estimate. Indeed, for $1\leq r\leq a$ and $k\in F_r$, because $k > k_0$, by Definition \ref{def hris} \ref{ris h-support growth}, $\min\mathrm{supp}_\mathrm{h}(x_k)  > \max\supp(x_{k_0}) \geq a$ and, therefore, $\mathrm{index}(\xi_r) = r\notin\mathrm{supp}_\mathrm{h}(P_Ex_k)\subset \mathrm{supp}_\mathrm{h}(x_k)$.
\end{proof}

\begin{proposition}[Basic inequality for h-RIS]
\label{h-basic inequality}
Let $(x_k)_{k\in I}$ be a $(C,2)$-h-RIS, indexed over an interval $I$ of $\N$, and let $(a_k)_{k\in I}$ be a sequence of complex numbers. Let $\gamma\in\Gamma$ and $E$ be an interval of $\N_0$. Then there exist $k_0\in I$ and $g\in W_q[(\mathcal{A}_{4n_j},4^{1/p}m_j^{-1})_{j\in L_0\setminus\{0\}}]$ such that
\begin{enumerate}[label=(\roman*)]

    \item either $g = 0$ or $\mathrm{weight}(g) = 4^{1/p}\mathrm{weight}(\gamma)$ and $\supp(g) \subset I\setminus\{k_0\}$ and

    \item\[\Big|e_\gamma^*\circ P_E\Big(\sum_{k\in I}a_kx_k\Big)\Big| \leq 12C\Big|(e_{k_0}^*+g)\Big(\sum_{k\in I}|a_k|e_k\Big)\Big|.\]

\end{enumerate}
\end{proposition}

\begin{proof}
The proof is performed by induction on the rank of $\gamma\in\Gamma$. The inductive hypothesis states that for a given rank $n$, the conclusion holds for all $\gamma\in\Gamma$ of rank $n$, all intervals $E$ of $\mathbb{N}_0$, and all subintervals $J$ of $I$. For the base case, assume that $\gamma\in\Gamma$ has rank zero, and thus also weight zero, $E$ is an interval of $\N_0$, and $J$ is a subinterval of $I$.  Then $e_\gamma^* = d_\gamma^*$ and, thus, either $e_\gamma^*\circ P_E(\sum_{k\in J}a_kx_k) = 0$ or there is $k_0\in J$ for which
\[\Big|e_\gamma^*\circ P_E\Big(\sum_{k\in J}a_kx_k\Big)\Big| = \big|e_\gamma^*\big(a_{k_0}x_{k_0}\big)\big| \leq  C|a_{k_0}|.\]
In the first case, take $k_0 = \min(J)$ and in both cases, put $g= 0$. For the inductive step, let $\gamma\in\Gamma$, $E$ be a given interval of $\mathbb{N}$, and $J$ be a given subinterval of $I$. Assume that the conclusion holds for all $\eta\in\Gamma$ with $\mathrm{rank}(\eta)<\mathrm{rank}(\gamma)$, intervals $G$ of $\N_0$, and subintervals $J'$ of $J$. If $\gamma$ is in some $\Gamma^s$, $s\in\N$, or of weight $m_0^{-1}$, or of weight zero then, by Lemma \ref{pre h-ris lemma}, we can easily find $k_0$ in $J$ such that by setting $g=0$ the conclusion holds. We treat the remaining case, i.e., when $\gamma\in\Gamma^0$ is of weight $m_j^{-1}$, for some $j\in L_0\setminus\{0\}$. Assume, without loss of generality, that for all $k\in J$, $\min(E)\leq \max\supp(x_k)$ (because, otherwise, $P_E(x_k) = 0$). Note that, in particular, for $k\in J\setminus \{\min(J)\}$ we have $\min(E)\leq \min\supp(x_k)$. Put $S = \{\min(J)\}\cup\{k\in J:\Xi(\gamma,x_k)\neq\emptyset\}$ and note $\#S\leq 3$. Pick $k_0\in S$ such that $|a_{k_0}| = \max_{k\in S}|a_k|$ and note
\begin{equation}
\label{h-basic inequality eq1}
\Big|e_\gamma^*\circ P_E\Big(\sum_{k\in S}a_kx_k\Big)\Big| \leq 4C(N+1)\Big|e_{k_0}^*\Big(\sum_{k\in F}|a_k|e_k\Big)\Big|.
\end{equation}

Because $\min(S) = \min(J)$ and $\#S\leq 3$ there are intervals $J^{(1)}<J^{(2)}<J^{(3)}$ of $\mathbb{N}_0$ (some of which may be empty) such that $J\setminus S = J^{(1)}\cup J^{(2)}\cup J^{(3)}$. Write the evaluation analysis of $\gamma$ as in \eqref{eva analysis hor}:
\[e_\gamma^* = \sum_{r=1}^a d_{\xi_r}^* + \frac{1}{m_j}\sum_{r=1}^a \frac{1}{n_j^{1/p}} e_{\eta_r}^*\circ  P_{E_r},\]
with $a\leq n_j$, $\eta_r\in \Gamma$ with $\mathrm{rank}(\eta_r)<\rank(\gamma)$, for $r=1,\dots,a$, and $E_1 <\cdots<E_a$. Apply Proposition \ref{dxi elimination lemma} to find, for $i=1,2,3$, subintervals $J^{(i)}_1<\cdots< J^{(i)}_a$ of $J^{(i)}$ with
\begin{equation*}
e_\gamma^*\circ P_E\Big(\sum_{k\in J^{(i)}}a_kx_k\Big) = \frac{1}{m_j}\sum_{r=1}^a\frac{1}{n_j^{1/p}}e_{\eta_r}^*\circ P_{E_r\cap E}\Big(\sum_{k\in J^{(i)}_r}a_kx_k\Big).
\end{equation*}
Put, for $i=1,2,3$ and $1\leq r\leq a$,
\begin{align*}
\tilde J^{(i)}_r = \Big\{k\in J^{(i)}_r:E_r\cap\mathrm{range}(x_k)\neq\emptyset\Big\}\text{ and }R^{(i)} = \Big\{1\leq r\leq a:\tilde J^{(i)}_r\neq\emptyset\Big\},
\end{align*}
such that
\begin{equation*}
e_\gamma^*\circ P_E\Big(\sum_{k\in J^{(i)}}a_kx_k\Big) = \frac{1}{m_j}\sum_{r\in R^{(i)}}^a\frac{1}{n_j^{1/p}}e_{\eta_r}^*\circ P_{E_r\cap E}\Big(\sum_{k\in \tilde J^{(i)}_r}a_kx_k\Big).
\end{equation*}

By the inductive hypothesis, for $i=1,2,3$ and $r\in R^{(i)}$, we can find $k^{(i)}_r\in \tilde J^{(i)}_r$ and $g^{(i)}_r\in W_q[(\mathcal{A}_{4n_j},4^{1/p}m_j^{-1})_{j\in L_0\setminus\{0\}}]$ with $\supp(g^{(i)}_r)\subset \tilde J^{(i)}_r\setminus\{k^{(i)}_r\}$, such that
\begin{equation}
\label{h-basic inequality eq2}
\Big|e_\gamma^*\circ P_E\Big(\sum_{k\in \tilde J^{(i)}_r}a_kx_k\Big)\Big| \leq \frac{12C}{m_jn_j^{1/p}}(e^*_{k^{(i)}_r} + g^{(i)}_r)\Big(\sum_{k\in \tilde J^{(i)}_r}|a_k|e_k\Big).
\end{equation}
Define the following functional:
\begin{equation*}
g = \frac{4^{1/p}}{m_j(4n_j)^{1/p}}\sum_{i=1}^3\sum_{r\in R^{(i)}}\Big(g^{(i)}_r|_{[0,k^{(i)}_r)}+e_{k^{(i)}_r}^*+g^{(i)}_r|_{(k_r,\infty)}\Big)
\end{equation*}
which has, at most, $3\#R^{(1)}+3\#R^{(2)}+3\#R^{(3)}$ summands.
\begin{claim}
$\#R^{(1)}+\#R^{(2)}+\#R^{(3)}\leq a +2$.
\end{claim}
Given the claim, $3\#R^{(1)}+3\#R^{(2)}+3\#R^{(3)} \leq 3n_j+6 \leq 4n_j$ and thus, $g$ is a member of $W_q[(\mathcal{A}_{4n_j},4^{1/p}m_j^{-1})_{j\in L_0\setminus\{0\}}]$. Note that the support of $g$ is in $J^{(1)}\cup J^{(2)}\cup J^{(3)}\subset J\setminus \{k_0\}$ and, by \eqref{h-basic inequality eq1} and \eqref{h-basic inequality eq2},
\[\Big|e_\gamma^*\circ P_E\Big(\sum_{k\in J}a_kx_k\Big)\Big| \leq 12C\Big|(e_{k_0}^*+g)\Big(\sum_{k\in J}|a_k|e_k\Big)\Big|.\]
To prove the claim we verify the sufficient condition $\max(R^{(i)})\leq \min(R^{(i+1)})$, for $i=1,2$. Let $r \in R^{(i)}$, $s\in R^{(i+1)}$ for which we will show $r\leq s$. Pick $k\in J^{(i)}_r$ with $E_r\cap\mathrm{range}(x_k)\neq \emptyset$ and $l\in J^{(i+1)}_s$ with $E_s\cap\mathrm{range}(x_l)\neq\emptyset$. Because $J^{(i)}<J^{(i+1)}$, $k\in J^{(i)}_r\subset J^{(i)}$, and $l\in J^{(i+1)}_s\subset J^{(i+1)}$ we deduce $k<l$ and thus $\min(E_r)\leq \max(\mathrm{range}(x_k)) \leq \min(\mathrm{range}(x_l)) \leq\max(E_s)$. It follows that $r\leq s$.
\end{proof}

\begin{corollary}
\label{h-ris norm estimates}
Let $j_0\in\N$ and $(x_k)_{k=1}^{n_{j_0}}$ be a $(C,2)$-h-RIS in $\mathfrak{X}_U$. Then,
\[\Big\|\frac{1}{n_{j_0}^{1/q}}\sum_{k=1}^{n_{j_0}}x_k\Big\| \leq \left\{\begin{array}{ll}\frac{96C}{m_{j_0}}&\text{ if }j_0\in L_0\setminus\{0\}\\ \frac{72C}{m_{j_0}^2}&\text{ if }j_0\notin L_0\setminus\{0\}.\end{array}\right.\]
\end{corollary}

\begin{proof}
By the basic inequality for h-RISs and Proposition \ref{estimates mixed tsirelson} \eqref{2 estimate mixed tsirelson}, for any $\gamma\in\Gamma$ there exist $k_0\in\{1,\ldots,n_{j_0}\}$ and $g\in W_q[(\mathcal{A}_{4n_j},4^{1/p}m_j^{-1})_{j\in L_0\setminus\{0\}}]$ such that
\begin{align*}
\Big|e_\gamma^*\Big(\frac{1}{n_{j_0}^{1/q}}\sum_{k=1}^{n_{j_0}}x_k\Big)\Big| &\leq 12C(e_{k_0}^*+g)\Big(\frac{1}{n_{j_0}^{1/p}}\sum_{k=1}^{n_{j_0}}e_k\Big)\\
&\leq 24C\Big\|\frac{1}{n_{j_0}^{1/q}}\sum_{k=1}^{n_{j_0}}e_k\Big\|.
\end{align*}
The conclusion follows from Proposition \ref{estimates mixed tsirelson}.
\end{proof}

There are certain assumptions on the local support of a bounded horizontally block sequence $(x_k)_{k=1}^\infty$ that imply that it has a subsequence that is a h-RIS.

\begin{notation}
Let $(x_k)_{k=1}^\infty$ be a block sequence in $\mathfrak{X}_U$.

\begin{enumerate}[label=(\arabic*),leftmargin=19pt]
\item We say that $(x_k)_{k=1}^\infty$ has uniform local age $a_0$, where $a_0\in\mathbb{N}$, if for every $k\in\mathbb{N}$ and $\gamma\in\supp_\mathrm{loc}(x_{k})$, $\mathrm{age}(\gamma) = a_0$.

\item We say that $(x_k)_{k=1}^\infty$ has rapidly increasing local ages if for every $k\in\mathbb{N}$ and $\gamma\in\supp_\mathrm{loc}(x_{k+1})$, if we put $j_k = \max\supp(x_k)$, then $\mathrm{age}(\gamma)> n_{j_k}$.

\end{enumerate}
\end{notation}

\begin{proposition}
\label{standard hris}
Let $(x_k)_{k=1}^\infty$ be a horizontally block sequence in $\mathfrak{X}_U$ with $C = \sup_k\|x_k\|<\infty$ that either has uniform local age $a_0$, for some $a_0\in\N$, or rapidly increasing local ages. Then $(x_k)_{k=1}^\infty$ has a subsequence that is a $(C,1)$-h-RIS.
\end{proposition}

\begin{proof}
Note that in Definition \ref{def hris} condition \ref{ris h-support growth} can always be achieved easily by passing to subsequences so we only need to focus on achieving condition \ref{ris dxi collapse}.

We first treat the case in which $(x_k)_{k=1}^\infty$ has uniform local ages $a_0$. 
Let $\gamma\in\Gamma$ and 
write $\Xi(\gamma) = \{\xi_1,\ldots,\xi_a\}$ with $\mathrm{age}(\xi_r) = r$, for $1\leq r\leq a$. If $a<a_0$ then, for every $k\in\N$, $\Xi(\gamma,x_k) = \Xi(\gamma)\cap\supp_\mathrm{loc}(x_k) = \emptyset$. If $a\geq a_0$ then, for every $k\in\mathbb{N}$, $\Xi(\gamma,x_k)\subset \{\xi_{a_0}\}\cap\supp_\mathrm{loc}(x_k)$, which is non-empty for at most one $k\in\mathbb{N}$.

Next, let $(x_k)_{k=1}^\infty$ have 
rapidly increasing local ages.
Take $\gamma\in\Gamma$ and assume that there is $k$ such that $\Xi(\gamma,x_k)\neq\emptyset$. We will show that for $l>k$, $\Xi(\gamma,x_l) = \emptyset$. We distinguish between different cases for the weight of $\gamma$. If it is of zero weight, then $\Xi(\gamma) = \{\gamma\}$ and thus $\Xi(\gamma)\subset\supp_\mathrm{loc}(x_k)$, which is disjoint from $\supp_\mathrm{loc}(x_l)$. If $\gamma$ has weight $m_j^{-1}$, for some $j\geq 2$, write $\Xi(\gamma) = \{\xi_1,\ldots,\xi_a\}$, with $a\leq n_j$ and  $\mathrm{weight}(\xi_r) = m_j^{-1}$, for $1\leq r\leq a$. Take $\xi_r\in\Xi(\gamma)\cap\supp_\mathrm{loc}\supp(x_k)$. By the definition of $\Gamma$, the weight of each node is controlled by its rank and, thus, $j\leq \max\mathrm{supp}(x_k)$. Because $(x_k)_{k=1}^\infty$ has rapidly increasing local ages, we have that $\supp_\mathrm{loc}(x_l)\cap\Xi(\gamma) = \emptyset$. Finally, if $\mathrm{weight}(\gamma) = m_0^{-1}$, write $\Xi(\gamma) = \{\xi_1,\ldots,\xi_a\}$ with $a\leq\min\supp(e^*_\gamma)$ and $\mathrm{age}(\xi_r) = r$, $1\leq r\leq a$. Because there is $1\leq r\leq a$ such that $\xi_r\in\supp_\mathrm{loc}(x_k)$, $\min\supp(e_\gamma^*)\leq\max\supp(x_k)$. Therefore, $\Xi(\gamma,x_l) = \Xi(\gamma)\cap \supp_\mathrm{loc}(x_l) = \emptyset$.
\end{proof}

Every bounded horizontally block sequence can be apropriately decomposed into h-RIS components and, thus, h-RISs can be used to test horizontal compactness of a bounded linear operator. 

\begin{proposition}
\label{horizontal compactness via h-ris}
Let $Y$ be a Banach space and $T:\mathfrak{X}_U\to Y$ be a bounded linear operator. If for every $C>0$ and $(C,1)$-h-RIS $(y_k)_{k=1}^\infty$ in $\mathfrak{X}_U$, $\lim_k\|Ty_k\| = 0$ then $T$ is horizontally compact.
\end{proposition}

\begin{proof}
We start with an arbitrary bounded horizontally block sequence $(x_k)_k$ and will prove that it has a subsequence with $\lim_k\|Tx_k\| = 0$. For each $k\in\mathbb{N}$ denote $S_k = \supp_\mathrm{loc}(x_k)$. Importantly, by Lemma \ref{local support lemma} \ref{local quantities} \ref{local restriction BD support} and \ref{local restriction h-support}, for any $S'_k\subset S_k$, the sequence $(x{\upharpoonright}_{S'_k})_{k=1}^\infty$ is horizontally block. For $k,a\in\mathbb{N}$ denote $S_k^a = \{\gamma\in S_k:\mathrm{age}(\gamma) = a\}$. By Lemma \ref{local support lemma} \ref{local quantities} \ref{local restriction local support}, for every $a\in\mathbb{N}$ the sequence $(x_k{\upharpoonright}_{S_k^a})_k$ has uniform local ages $a$ and, by Proposition \ref{standard hris}, $\lim_k\|T(x_k{\upharpoonright}_{S_k^a})\|\to 0$.
We may pass to a subsequence of $(x_k)_{k=1}^\infty$, again denoted by $(x_k)_{k=1}^\infty$, such that, for all $k\geq 2$, putting $j_k = \max\supp(x_{k-1})$ we have
\[\Big|\big\|Tx_k\big\|-\big\|Tx_k{\upharpoonright}_{S_k\setminus \cup_{a=1}^{j_k}S_k^a}\big\|\Big| \leq \big\|Tx_k{\upharpoonright}_{\cup_{a=1}^{j_k-1}S_k^a}\big\| < k^{-1}.\]
By Lemma \ref{local support lemma} \ref{local quantities} \ref{local restriction local support}, the sequence $(x_k{\upharpoonright}_{S_k\setminus \cup_{a=1}^{j_k}S_k^a})_{k=1}^\infty$ has rapidly increasing local ages and by Proposition \ref{standard hris} we deduce that $\lim_k\|Tx_k{\upharpoonright}_{S_k\setminus \cup_{a=1}^{j_k}S_k^a}\| = 0$.
\end{proof}

\begin{corollary}
\label{h-block w-null}
Every bounded horizontally block sequence in $\mathfrak{X}_U$ is weakly null.
\end{corollary}

\begin{proof}
If the conclusion is false, there exists a bounded linear functional $f:\mathfrak{X}_U\to\mathbb{C}$ and a bounded horizontally block sequence $(x_k)_{k=1}^\infty$ such that $\limsup_k|f(x_k)| > 0$. By Proposition \ref{horizontal compactness via h-ris}, there exists a $(C,1)$-h-RIS $(y_k)_{k=1}^\infty$ and $\varepsilon > 0$  such that, for all $k\in\N$, $\mathfrak{Re}(f(y_k)) \geq \varepsilon$. Then, for arbitrary $j_0\in L_0\setminus\{0\}$,
\[\mathfrak{Re}\Big(f\Big(\frac{1}{n_{j_0}^{1/q}}\sum_{k=1}^{n_{j_0}}y_k\Big)\Big) \geq\varepsilon n_{j_0}^{1/p}\]
and by Proposition \ref{h-basic inequality}
\[\Big|f\Big(\frac{1}{n_{j_0}^{1/q}}\sum_{k=1}^{n_{j_0}}y_k\Big)\Big| \leq \|f\|\Big\|\frac{1}{n_{j_0}^{1/q}}\sum_{k=1}^{n_{j_0}}c_ky_k\Big\|\leq \|f\|96C\frac{1}{m_{j_0}}.\]
We deduce $n_{j_0}^{1/p}/m_{j_0}\leq 96C\|f\|/\varepsilon$. This contradicts Assumption \ref{assumption integers} \ref{4 assumption integers}, which implies that for any $\alpha, \beta>0$, $\lim_jn_{j}^\alpha/m_j^\beta = \infty$.
\end{proof}

\subsection{$\mathfrak{X}_s$ rapidly increasing sequences}\label{xs-ris}
We need to redefine RISs in the context of each of the subspaces $\mathfrak{X}_s$, $s\in\N$, of $\mathfrak{X}_U$. Because each $\Gamma^s$ is a self-determined subset of $\Gamma$, we can treat $\mathfrak{X}_s$ as the Bourgain-Delbaen space $\mathfrak{X}_{(\Gamma_n^s,i_n^s)}$ to which it is isometrically isomorphic.  The definitions and proofs are almost identical and practically carry over unaltered. For reasons of formal soundness, we point out the technical details of the identification $\mathfrak{X}_s\equiv\mathfrak{X}_{(\Gamma_n^s,i_n^s)}$ and rephrase the required definitions, but don't repeat any steps of the proofs, such as the evaluation refinement.

It follows from Theorem \ref{XU projections are U} that, if $r_{\Gamma^s}:\ell_\infty(\Gamma)\to\ell_\infty(\Gamma^s)$ denotes the usual restriction map and $R_s = r_{\Gamma^s}|_{\mathfrak{X}_U}$, then $\mathrm{Im}(R_s) =\mathfrak{X}_{(\Gamma_n^s,i_n^s)}$ and  $R_s|_{\mathfrak{X}_s}$ is an isometry onto $\mathfrak{X}_{(\Gamma_n^s,i_n^s)}$. Denote by $d_\gamma^s$, $c_\gamma^{s*}$, $d_\gamma^{s*}$, and $P_E^s$ the usual Bourgain-Delbaen objects associated to $\mathfrak{X}_{(\Gamma_n^s,i_n^s)}$ and maintain the notation $e_\gamma^*$, for $\gamma\in\Gamma^s$.

\begin{fact}
\label{component fact}
Let $s\in\N$.
    \begin{enumerate}[label=(\arabic*)]
        
        \item\label{component fact 1} For $\gamma\in\Gamma^s$, $d_\gamma^s = R_sd_\gamma$.

        \item\label{component fact 2} Let $\gamma\in\Gamma^s$ be a simple node with the evaluation analysis
        \[e_\gamma^* = \sum_{r=1}d_{\xi_r}^* + \frac{1}{m_j}\sum_{r=1}^a\frac{1}{n_j^{1/p}}e_{\eta_r}^*\circ P_E.\]
    Then, for $1\leq r\leq a$, $\xi_r$, $\eta_r\in\Gamma^s$ and $e_\gamma^*$, in $\ell_\infty(\Gamma_s)^*$, has evaluation analysis
    \[e_\gamma^* = \sum_{r=1}d_{\xi_r}^{s*} + \frac{1}{m_j}\sum_{r=1}^a\frac{1}{n_j^{1/p}}e_{\eta_r}^*\circ P^s_E.\]

    \item\label{component fact 3} $I_s = (R_s|_{\mathfrak{X}_s})^{-1}R_s$.
    
    \end{enumerate}
\end{fact}

\begin{proof}
Fact \ref{component fact 1} follows from the definition of $i_n^s$, fact \ref{component fact 2} follows from the preservation of formulae for self-determined sets (Fact \ref{preservation of formulae}) and fact \ref{component fact 3} follows from fact \ref{component fact 1} and that $\langle\{e_\gamma^*:\gamma\in\Gamma^s\}\rangle = \langle\{d_\gamma^*:\gamma\in\Gamma^s\}\rangle$ and, thus, for $\eta\notin\Gamma^s$, $R_s(d_\eta) = 0$.
\end{proof}

We redefine the notion of local support of a finitely supported vector $x$ in $\mathfrak{X}_{(\Gamma_n^s,i_n^s)}$, with $\max\supp(x) = n$, as $\supp^s_\mathrm{loc}(x) = \{\gamma\in\Gamma^s:\mathrm{rank}(\gamma)\leq n\text{ and }e_\gamma^*(x)\neq 0\}$. For $\gamma\in\Gamma^s$ we define $\Xi^s(\gamma) = \{\eta\in\Gamma^s:|e_\gamma^*(d_\xi)| > \mathrm{weight}(\gamma)\}$ and for a finitely supported vector $x\in\mathfrak{X}_{(\Gamma_n^s,i_n^s)}$ denote $\Xi^s(\gamma,x) = \Xi^s(\gamma)\cap\supp_\mathrm{loc}^s(x)$. Although these are not required, it can be seen that for $\gamma\in\Gamma^s$, $\Xi^s(\gamma) = \Xi(\gamma)$ and for $x\in\mathfrak{X}_U$, $\supp_\mathrm{loc}^s(R_sx) = \supp_\mathrm{loc}(x)\cap\Gamma^s$.

\begin{definition}
\label{def component ris}
Let $s\in\N$, $C > 0$ and $N = 1$ or $N = 2$. A block sequence $(x_k)_{k\in I}$ in $\mathfrak{X}_{(\Gamma_n^s,i_n^s)}$ indexed over some interval $I$ of $\N$, is called a $(C,N)$-$\mathfrak{X}_s$ rapidly increasing sequence (or a $(C,N)$-$\mathfrak{X}_s$-RIS) if
the following hold.
\begin{enumerate}[label=(\alph*)]
\item\label{Xs ris norm} For $k\in I$, $\|x_k\|\leq C$.
\item\label{Xs ris dxi collapse} For every $\gamma\in\Gamma^s$, $\#\{k\in I:\Xi(\gamma,x_k)\neq\emptyset\}\leq N$.

\end{enumerate}
We call a sequence $(x_k)_{k\in I}$ in $\mathfrak{X}_s$ a $(C,N)$-$\mathfrak{X}_s$-RIS if its isometric image $(R_sx_k)_{k\in I}$ in $\mathfrak{X}_{(\Gamma_n^s,i_n^s)}$ is a $(C,N)$-$\mathfrak{X}_s$-RIS. We will call a sequence $(x_k)_{k\in I}$ in $\mathfrak{X}_s$ an $\mathfrak{X}_s$-RIS if it is a $(C,N)$-$\mathfrak{X}_s$-RIS for some $C>0$ and $N=1$ or $N=2$. 
\end{definition}

The only difference between Definition \ref{def hris} and the above lies in the omission of condition \ref{ris h-support growth} due to the lack of ground nodes in $\Gamma^s$.

\begin{proposition}[Basic inequality for $\mathfrak{X}_s$-RIS]
\label{Xs-basic inequality}
Let $s\in\N$, $(x_k)_{k\in I}$ be a $(C,2)$-$\mathfrak{X}_s$-RIS in $\mathfrak{X}_s$, indexed over an interval $I$ of $\N$, and let $(a_k)_{k\in I}$ be a sequence of complex numbers. Let $\gamma\in\Gamma^s$ and $E$ be an interval of $\N_0$. Then there exist $k_0\in I$ and $g\in W_q[(\mathcal{A}_{4n_j},4^{1/p}m_j^{-1})_{j\in L_s}]$ such that
\begin{enumerate}[label=(\roman*)]

    \item either $g = 0$ or $\mathrm{weight}(g) = \mathrm{weight}(\gamma)$ and $\supp(g) \subset I\setminus\{k_0\}$ and

    \item\[\Big|e_\gamma^*\circ P_E\Big(\sum_{k\in I}a_kx_k\Big)\Big| \leq 12C\Big|(e_{k_0}^*+g)\Big(\sum_{k\in I}|a_k|e_k\Big)\Big|.\]

\end{enumerate}
\end{proposition}

\begin{corollary}
\label{Xs-ris norm estimates}
Let $j_0\in\N$ and $(x_k)_{k=1}^{n_{j_0}}$ be a $(C,2)$-$\mathfrak{X}_s$-RIS. Then,
\[\Big\|\frac{1}{n_{j_0}^{1/q}}\sum_{k=1}^{n_{j_0}}x_k\Big\| \leq \left\{\begin{array}{ll}\frac{96C}{m_{j_0}}&\text{ if }j_0\in L_s\\ \frac{72C}{m_{j_0}^2}&\text{ if }j_0\notin L_s.\end{array}\right.\]
\end{corollary}

\begin{proposition}
\label{compactness via Xs-ris}
Let $s\in\mathbb{N}$, $Y$ be a Banach space and $T:\mathfrak{X}_s\to Y$ be a bounded linear operator. If for every $C>0$ and $(C,1)$-$\mathfrak{X}_s$-RIS $(y_k)_{k=1}^\infty$ in $\mathfrak{X}_s$, $\lim_k\|Ty_k\| = 0$ then for every bounded horizontally block sequence $(x_k)_{k=1}^\infty$ in $\mathfrak{X}_s$, $\lim_k\|Tx_k\| = 0$.
\end{proposition}

\begin{corollary}
\label{Xs-block w-null}
Let $s\in\N$. Every bounded block sequence in $\mathfrak{X}_s$ is weakly null.
\end{corollary}

\begin{lemma}
\label{close to h-block}
Let $(x_k)_{k=1}^\infty$ be a block sequence in $\mathfrak{X}_U$ such that, for all $s\in\N$, $\lim_k\|I_sx_k\| = 0$. Then, $(x_k)_{k=1}^\infty$ has a subsequence that is arbitrarily close to a horizontally block sequence. That is, for every sequence $(\varepsilon_k)_{k=1}^\infty$ of positive real numbers there exist a subsequence $(x'_k)_{k=1}^\infty$ of $(x_k)_{k=1}^\infty$ and a horizontally block sequence $(y_k)_{k=1}^\infty$ such that, for all $k\in\N$, $\|x_k'-y_k\| \leq \varepsilon_k$.
\end{lemma}

\begin{proof}
By a diagonalization argument, we may pass to a subsequence of $(x_k)_{k=1}^\infty$, again denoted by $(x_k)_{k=1}^\infty$, such that the following hols: for all $k\geq 2$, putting $s_k = \max\mathrm{supp}_\mathrm{h}(x_{k-1})$, for all $s\leq s_k$, $\|I_sx_k\| \leq \varepsilon_ks_k^{-1}$. Define $y_1 = x_1$ and, for $k\geq 2$, $y_k = (I - \sum_{s=1}^{s_k}I_s)x_k$. Recall that, for all $s\in\N$, $I_s = I_{\Gamma^s}$ which yields that, for all $k\in\N$, $\supp(y_k)\subset \supp(x_k)$ and, if $k\geq 2$, $\min\mathrm{supp}_\mathrm{h}(y_k) > s_k \geq \max\mathrm{supp}_\mathrm{h}(x_{k-1}) \geq \max\mathrm{supp}_\mathrm{h}(y_{k-1})$. Therefore, $(y_k)_{k=1}^\infty$ is horizontally block and, for all $k\in\N$, $\|x_k - y_k\| \leq \sum_{s=1}^{s_k}\|I_sx_k\| \leq \varepsilon_k$.
\end{proof}

\begin{proposition}\label{block w-null}
Every bounded block sequence in $\mathfrak{X}_U$ is weakly null and, thus, $\mathfrak{X}_U^*\simeq\ell_1(\N)$.
\end{proposition}

\begin{proof}
If every bounded block sequence is weakly null, then the basis $(d_\gamma)_{\gamma\in\Gamma}$, with an appropriate enumeration as in Remark \ref{basis enumeration remark}, is shrinking and, therefore, $\mathfrak{X}_U^* = [d^*_\gamma:\gamma\in\Gamma] = [ e^*_\gamma:\gamma\in\Gamma]$. By Proposition \ref{dual equivalence}, $(e_\gamma^*)_{\gamma\in\Gamma}$ is $2$-equivalent to the unit vector basis of $\ell_1(\Gamma)$.

Let $(x_k)_{k=1}^\infty$ be a bounded block sequence and assume, towards contradiction, that for some $\varepsilon > 0$ and $f\in\mathfrak{X}_U^*$, $\mathfrak{Re}(f(x_k)) \geq \varepsilon$. In particular, for every $y\in\mathrm{conv}\{x_k:k\in\N\}$, $\mathfrak{Re}(f(y)) \geq \varepsilon$. Because, for all $s\in\N$, $I_s = I_{\Gamma^s}$, the sequence $(I_sx_k)_{k=1}^\infty$ is a block sequence of $(d_\gamma)_{\gamma\in\Gamma^s}$ and, therefore, by Corollary \ref{Xs-block w-null}, it is weakly null. By Mazur's theorem and a diagonal argument, we may pass to a block sequence $(y_k)_{k=1}^\infty$ in $\mathrm{conv}\{x_k:k\in\N\}$ such that, for all $s\in\N$, $\lim_k\|I_sy_k\| = 0$. Lemma \ref{close to h-block} and Corollary \ref{h-block w-null} yield that $(y_k)_{k=1}^\infty$ is weakly null, which contradicts the fact $\mathfrak{Re}(f(y_k)) \geq \varepsilon$, for all $k\in\N$.
\end{proof}

\subsection{Incomparability of structures}
\label{incomparability of structures}
In this section we prove, for every $s\in\N$, smallness properties of operators $T:\mathfrak{X}_U\to\mathfrak{X}_s$ and of operators $T:\mathfrak{X}_s\to\mathfrak{X}_U$.

The following is a standard functional analysis perturbation argument that makes the proofs in this section tidier.

\begin{fact}
\label{standard perturbation}
Let $X$, $Y$ be Banach spaces, $T:X\to Y$ be a bounded linear operator, $(x_k)_{k=1}^\infty$ be a seminormalized Schauder basic sequence in $X$, and $(y_k)_{k=1}^\infty$ be a sequence in $Y$ such that $\lim_k\|Tx_k-y_k\| = 0$. Then, for every $\varepsilon > 0$, there exists a nuclear, and thus compact, operator $S:X\to Y$ such that $\|S\|\leq \varepsilon$ and $(T+S)x_k = y_k$ for infinitely many $k\in 2\N$ and infinitely many $k\in 2\N-1$. 
\end{fact}

\begin{proof}
Let $(f_k)_{k=1}^\infty$ denote the biorthogonal sequence of $(x_k)_{k=1}^\infty$ in $[(x_k)_{k=1}^\infty]^*$ and, for $k\in\N$, let $x_k^*\in X^*$ be a norm preserving extension of $f_k$. In particular, $(x_k^*)_{k=1}^\infty$ is bounded. Let $z_k=Tx_k-y_k$, $k\in\N$, and pick $J\subset\N$ with $J\cap 2\N$ and $J\cap(2\N-1)$ infinite, such that that $\sum_{k\in J}\|z_k\|<\varepsilon\inf_k\|x^*_k\|^{-1}$. Define
\[S = \sum_{k\in J}x_k^*\otimes z_k\]
which satisfies the desired properties.
\end{proof}

\begin{proposition}
\label{horizontal compactness of projections}
For every $s\in\mathbb{N}$, every bounded linear operator $T:\mathfrak{X}_U\to\mathfrak{X}_s$ is horizontally compact.    
\end{proposition}

\begin{proof}
Recall that $\mathfrak{X}_s\equiv\mathfrak{X}_{(\Gamma_n^s,i_n^s)}$. It is sufficient to prove that an arbitrary bounded linear operator $T:\mathfrak{X}_U\to\mathfrak{X}_{(\Gamma_n^s,i_n^s)}$ sends bounded horizontally block sequences to norm-null ones. If the conclusion fails, by Proposition \ref{horizontal compactness via h-ris} there exist a $(C,1)$-h-RIS $(x_k)_{k=1}^\infty$ and $\varepsilon>0$ such that $\liminf_k\|Tx_k\| > \varepsilon$. By Corollary \ref{h-block w-null}, $(Tx_k)_{k=1}^\infty$ is weakly null and, by passing to a subsequence of $(x_k)_{k=1}^\infty$ and using Fact \ref{standard perturbation}, we may assume that it is a block sequence in $\mathfrak{X}_{(\Gamma_n^s,i_n^s)}$. For each $k\in\N$ pick $\eta_k\in\Gamma^s$ such that $\mathfrak{Re}(e_{\eta_k}^*(Tx_k)) > \varepsilon$. By passing to a further subsequence, we may assume that there are $(p_k)_{k=1}^\infty$ in $\N$ such that, for all $k\in\N$, $\min\supp(Tx_k) \leq \mathrm{rank}(\eta_k)\vee\max\supp(Tx_k)<p_k<\min\supp(Tx_{k+1})$ and, in particular, for arbitrary $2j\in L_s$, by Proposition \ref{existence of simple even nodes}, there exists $\gamma\in\Gamma^s$ such that
\[\mathfrak{Re}\Big(e_\gamma^*\Big(T\Big(\frac{1}{n^{1/q}_{2j}}\sum_{k=2j+1}^{n_{2j}+2j}c_kx_k\Big)\Big)\Big) \geq  \frac{\varepsilon}{m_{2j}}.\]
However, by Corollary \ref{h-ris norm estimates}, we have
\[\Big\|\frac{1}{n^{1/q}_{2j}}\sum_{k=2j+1}^{n_{2j}+2j}x_k\Big\| \leq \frac{48C}{m_{2j}^2}\]
and, therefore, $\|T\| \geq \varepsilon(48C)^{-1}m_{2j}$. As $2j$ was arbitrary, this is absurd.
\end{proof}

\begin{proposition}
\label{incomparable coordinates}
Let $s\neq t\in\N$. Every bounded linear operator $T:\mathfrak{X}_s\to\mathfrak{X}_t$ is compact.    
\end{proposition}

\begin{proof}[Comment on proof]
In the proof of Proposition \ref{horizontal compactness of projections} we showed that every bounded linear operator $T:\mathfrak{X}_U\to\mathfrak{X}_s$ maps bounded horizontally block sequences to norm-null ones. The same argument can be repeated to prove that a bounded linear map $T:\mathfrak{X}_s\to\mathfrak{X}_t$ sends bounded block sequences to zero and is thus compact.
\end{proof}

Although finding lower bounds for linear combinations of appropriate block sequences in each $\mathfrak{X}_s$ is relatively straightforward, the situation is slightly more delicate for h-block sequences in $\mathfrak{X}_U$. The difficulty lies in the fact that the definition of h-nodes is restricted by a growth condition on the h-supports of functionals of the form $e_\eta^*\circ P_E$. The following lemma addresses this hurdle.

\begin{lemma}
\label{h-block norming}
Let $(x_k)_{k=1}^\infty$ be a bounded horizontally block sequence in $\mathfrak{X}_U$ with $\liminf\|x_k\| > 0$. Then, there exist an infinite subset $\mathbb{M}$ of $\N$, a sequence $(\eta_k)_{k\in\mathbb{M}}$ in $\Gamma$, and a sequence of successive intervals $(F_k)_{k\in\mathbb{M}}$ such that  $\lim_{k\in \mathbb{M}}\min\mathrm{supp}_\mathrm{h}(e_{\gamma_k}^*\circ P_{F_k}) = \infty$ and $\inf_{k\in\mathbb{M}}|e_{\gamma_k}^*\circ P_{F_k}(x_k)| > 0$.
\end{lemma}

\begin{proof}
We will require the following statement.
    \begin{claim}
        Let $x\in\mathfrak{X}_U$, $m\leq \min\mathrm{supp}_\mathrm{h}(x)\wedge\min\supp(x)$, $\gamma\in\Gamma$, and $E$ be an interval of $\N_0$ such that $e_\gamma^*\circ P_E(x)\neq 0$. Then, there exists $\xi\in\Gamma$ and an interval $F$ of $E$ such that:
        \begin{enumerate}[label=(\roman*)]
            
            \item $\displaystyle{\Big|e_\xi^*\circ P_F(x)\Big| \geq \frac{1}{2}\Big|e_\gamma^*\circ P_E(x)\Big| - m\sup_{\eta\in\Gamma^0}|d_\eta^*(x)|}$ and

            \item $\min\mathrm{supp}_\mathrm{h}(e_\xi^*\circ P_F)\geq m$.
        
        \end{enumerate}
    \end{claim}
First, given the claim, we deduce the conclusion. If, for some $\varepsilon > 0$, the set
\[\mathbb{M} = \Big\{k\in\N: \sup_{\eta\in\Gamma^0}|d_\eta^*(x_k)| > \varepsilon\Big\}\]
is infinite then we may choose, for each $k\in\mathbb{M}$, $\eta_k\in\Gamma^0$ with $\rank(\eta_k)\in\supp(x_k)$ such that $|e_{\eta_k}^*\circ P_{\{\mathrm{rank}(\eta_k)\}}(x_k)| = |d_{\eta_k}^*(x_k)| > \varepsilon$. Putting, for each $k\in\mathbb{M}$, $F_k = \{\rank(\eta_k)\}$ we have $\mathrm{supp}_\mathrm{h}(e_{\eta_k}^*\circ P_{F_k}) = \emptyset$; this case is settled. Otherwise, $\lim_k\supp_{\eta\in\Gamma^0}|d_\eta^*(x_k)| = 0$. For every $k\in\mathbb{N}$, put $E_k = \mathrm{range}(x_k)$ 
and pick $\gamma_k\in\Gamma$ such that $|e_{\gamma_k}^*\circ P_{E_k}(x_k)| > (1/2)\|x_k\|$. We may pass to a subsequence of $(x_k)_{k=1}^\infty$, again denoted by $(x_k)_{k=1}^\infty$, such that, for some $\varepsilon>0$, $\inf_k\|x_k\| >\varepsilon$ and, for all $k\in\N$, $\min\mathrm{supp}_\mathrm{h}(x_k)\wedge\min\supp(x_k)\geq k$ and $\sup_{\eta\in\Gamma^0}|d_{\eta}^*(x_k)| \leq \varepsilon/(4k)$. By the claim, for each $k\in\mathbb{N}$, there exists an interval $F_k$ of $E_k$ and $\xi_k\in\Gamma$ such that $|e_{\xi_k}^*\circ P_{F_k}(x_k)| \geq \varepsilon/4$ and $\min\mathrm{supp}_\mathrm{h}(e_{\xi^*}\circ P_{F_k})\geq k$;  the conclusion follows.

We proceed to prove the Claim. Assume first that $\gamma$ is a simple node of index $s$. Because $\Gamma^s$ is self-determined, $e_\gamma^*\in\langle\{d^*_\eta:\eta\in\Gamma^s\}\rangle$ and, thus, $e_\gamma^*\circ P_E\in\langle\{d^*_\eta:\eta\in\Gamma^s,\mathrm{rank}(\eta)\in E\}\rangle$. Because $e_\gamma^*\circ P_E(x) \not= 0$ we obtain $s\in\mathrm{supp}_\mathrm{h}(x_k)$ and, thus, $s\geq m$. We may simply choose $\xi = \gamma$ and $F = E$. Next, assume that $\gamma$ is a g-node and write its evaluation analysis
\[e_\gamma^* = \sum_{r=1}^ad_{\xi_r}^* + \frac{1}{m_0}\sum_{r=1}^a\lambda_rd_{\eta_r}^*.\]
For $r=1,\ldots,m-1$, because $\mathrm{index}(\xi_r) = r$, we have $d_{\xi_r}^*(P_Ex) = 0$. Put $\xi=\gamma$ and $F = E\cap [\mathrm{rank}(\eta_m),\infty)$. Then, $\min\mathrm{supp}_\mathrm{h}(e_\xi^*\circ P_F)\geq m$ and
\[|e_{\xi}^*\circ P_F(x)| = \Big|e_\gamma^*\circ P_E(x) - \sum_{1\leq r\leq (m-1)\atop \rank(\xi_r)\in E}d_{\xi_r}^*(x)\Big| \geq \Big|e_\gamma^*\circ P_E(x)\Big| - (m-1)\sup_{\eta\in\Gamma^0}|d_\eta^*(x)|.\]

To handle h-nodes we need to perform an induction on $\mathrm{rank}(\gamma) = \ell=0,1,\ldots$. Having fixed $x$ and $m$ the inductive hypothesis states that the Claim is true for every $\gamma$ (of any type) up to rank $\ell$ and all intervals $E$ of $\N_0$. The base case follows easily from the fact that the unique $\gamma$ with rank zero has empty horizontal support. For the $\ell$'th inductive step let $\gamma\in\Gamma$ be a h-node of rank $\ell$ and $E$ be an interval $\N$. Write the evaluation analysis of $\gamma$
\[e_\gamma^* = \sum_{r=1}^ad_{\xi_r}^* + \frac{1}{m_j}\sum_{r=1}^a\frac{1}{n_j^{1/p}}e_{\eta_r}^*\circ P_{E_r},\]
let $r_0 = \min\{1\leq r\leq m: \mathrm{rank}(\xi_r) \geq m\}$ and note
\[e_\gamma^*\circ P_E(x) = \frac{1}{m_jn_j^{1/p}}e_{\eta_{r_0}}^*\circ P_{E_{r_0}}(P_Ex) + \sum_{r=r_0}^ad_{\xi_r}^*(P_Ex) + \frac{1}{m_j}\sum_{r=r_0+1}^a\frac{1}{n_j^{1/p}}e^*_{\eta_r}\circ P_{E_r}(P_Ex)\]
We distinguish two cases.

\noindent{\bf Case 1:} $\displaystyle{\Big|\frac{1}{m_jn_j^{1/p}}e_{\eta_{r_0}}^*\circ P_{E_{r_0}}(P_Ex)\Big| < \frac{1}{2}\Big|e_\gamma^*\circ P_E(x)\Big|.}$

Put $\xi = \gamma$ and $F = E\cap (\mathrm{rank}(\xi_{r_0}),\infty)$ to obtain
\[\min\mathrm{supp}_\mathrm{h}(e_\xi^*\circ P_F) = \min\Big(\bigcup_{r=r_0+1}^a\mathrm{supp}_\mathrm{h}(e_{\eta_r}^*\circ P_{E_r})\Big) > \mathrm{rank}(\xi_{r_0}) \geq m.\]
and
\begin{align*}\Big|e_\xi^*\circ P_F(x)\Big|  &= \Big|\sum_{r=r_0}^ad_{\xi_r}^*(P_Ex) + \frac{1}{m_j}\sum_{r=r_0+1}^a\frac{1}{n_j^{1/p}}e^*_{\eta_r}\circ P_{E_r}(P_Ex)\Big|\\
&= \Big|e_\gamma^*\circ P_E(x) - \frac{1}{m_jn_j^{1/p}}e_{\eta_{r_0}}^*\circ P_{E_{r_0}}(P_Ex)\Big| \geq \frac{1}{2}\Big|e_\gamma^*\circ P_E(x)\Big|.
\end{align*}

\noindent{\bf Case 2:} $\displaystyle{\Big|\frac{1}{m_jn_j^{1/p}}e_{\eta_{r_0}}^*\circ P_{E_{r_0}}(P_Ex)\Big| \geq \frac{1}{2}\Big|e_\gamma^*\circ P_E(x)\Big|.}$

We deduce that $|e_{\eta_{r_0}}^*\circ P_{E_{r_0}\cap E}(x)| \geq (m_jn_j^{1/p}/2)|e_\gamma^*\circ P_E(x)| \geq |e_\gamma^*\circ P_E(x)|$. By the inductive hypothesis, there exists $\xi\in\Gamma$ and an interval $F$ of $E_{r_0}\cap E$ such that
\begin{align*}
|e_\xi^*\circ P_F(x)| \geq \frac{1}{2}|e_{\eta_{r_0}}^*\circ P_{E_{r_0}\cap E}(x)| - m\sup_{\eta\in\Gamma^0}|d_\eta^*(x)| \geq \frac{1}{2}|e_\gamma^*\circ P_E(x)| - m\sup_{\eta\in\Gamma^0}|d_\eta^*(x)|.
\end{align*}
\end{proof}

\begin{proposition}
\label{same-component-reduction}
For every $s\in\N$ and bounded $T:\mathfrak{X}_s\to\mathfrak{X}_U$, $(I-I_s)\circ T$ is compact.
\end{proposition}

\begin{proof}
Assume that $S = (I-I_s)\circ T$ is not compact. By Proposition \ref{compactness via Xs-ris} there exist a $(C,1)$-$\mathfrak{X}_s$-RIS $(x_k)_{k=1}^\infty$ in $\mathfrak{X}_s$ such that $\inf_k\|Sx_k\| > 0$. Obviously, $I_s\circ S = 0$ and, by Proposition \ref{incomparable coordinates}, for $t\neq s$, $I_t\circ S$ is compact. By Corollary \ref{Xs-block w-null}, $(x_k)_{k=1}^\infty$ is weakly null and, thus, for all $t\in\N$, $\lim_kI_t\circ S(x_k) = 0$. By Lemma \ref{close to h-block} and Fact \ref{standard perturbation}, passing to a subsequence of $(x_k)_{k=1}^\infty$, again denoted by $(x_k)_{k=1}^\infty$, we may assume that $(Sx_k)_{k=1}^\infty$ is a horizontally block sequence with $\inf_k\|Sx_k\| > 0$. By Lemma \ref{h-block norming} we can find $\varepsilon >0$, $(\eta_k)_{k=1}^\infty$ in $\Gamma$, intervals $(F_k)_{k=1}^\infty$ of $\N_0$ and $(p_k)_{k=1}^\infty$ in $\N$ such that, for all $k\in\N$,
\begin{enumerate}[label=(\arabic*)]

    \item $|e_{\eta_k}^*\circ P_{F_k}(Sx_k)| \geq \varepsilon$ and $F_k$ is an interval of $\mathrm{range}(Sx_k)$. 

    \item $\mathrm{rank}(\eta_1) \in [0, p_1)$ and $\mathrm{range}(x_1)\cup E_1\subset [0,p_1)$ and, for $k\geq 2$, $\mathrm{rank}(\eta_k)\in (p_{k-1}, p_k)$ and $\mathrm{range}(x_k)\cup E_k\subset (p_{k-1}, p_k)$.

    \item for $k\geq 2$, $\min\mathrm{supp}_\mathrm{h}(e_{\eta_k}^*\circ P_{F_k}) > p_{k-1}$ and $\mathrm{supp}_\mathrm{h}(e_{\eta_{k-1}}^*\circ P_{E_{k-1}})\ll \mathrm{supp}_\mathrm{h}(e^*_{\eta_k}\circ P_{E_k})$.

\end{enumerate}
In particular, for arbitrary $2j\in L_0\setminus\{0\}$, by Proposition \ref{existence of h-nodes}, there exists $\gamma\in\Gamma^0$ such that, for $c_k = \overline{e_{\eta_k}^*\circ P_{F_k}(Sx_k)}/|e_{\eta_k}^*\circ P_{F_k}(Sx_k)|$,
\[\Big|e_\gamma^*\Big(S\Big(\frac{1}{n^{1/q}_{2j}}\sum_{k=2j+1}^{n_{2j}+2j}c_kx_k\Big)\Big)\Big| \geq  \frac{\varepsilon}{m_{2j}}.\]
However, by Corollary \ref{Xs-ris norm estimates}, we have
\[\Big\|\frac{1}{n^{1/q}_{2j}}\sum_{k=2j+1}^{n_{2j}+2j}c_kx_k\Big\| \leq \frac{48C}{m_{2j}^2}\]
and, therefore, $\|S\| \geq \varepsilon(48C)^{-1}m_{2j}$. As $2j$ was arbitrary, this is absurd.
\end{proof}

\section{Dependent sequences  in $\mathfrak{X}_U$}\label{dependent sequences}

Dependent sequences of weighted averages of increasing lengths, dictated by parameters of some odd-weight node in $\Gamma$, form the  second level of special sequences and are the next canonical tool in the study of Argyros-Haydon-type spaces. We apply a canonical procedure based on the properties of dependent sequences  to prove the scalar-plus-compact type properties of $\mathfrak{X}_U$ both with respect to the horizontal structure of $\mathfrak{X}_U$ and in each component $\mathfrak{X}_s$, $s\in\N$, separately. More precisely, we prove that every bounded operator $T:\mathfrak{X}_U\to\mathfrak{X}_U$ is a scalar multiple of the identity plus a horizontally compact operator and that, for $s\in\N$,  every bounded operator $T:\mathfrak{X}_s\to\mathfrak{X}_s$ is a scalar multiple of the identity plus a compact operator.

Unlike the typical situation, due to the $q$-convex component in the definition of the space $\mathfrak{X}_U$, dependent sequences are not RISs; thus, we need refined tools to handle them, presented in the sequel. 

We define and thoroughly study horizontally dependent sequences. We formulate the corresponding $\mathfrak{X}_s$ counterparts but do not give proofs of their properties and application as they are analogous to the horizontal case.

\subsection{Horizontally dependent sequences}
\label{subsection h-dependent sequences}

We define h-exact pairs and h-dependent sequences. We describe their construction process, prove a refined version of the Basic Inequality (Proposition \ref{h-basic inequality}) and conclude the subsection with the proof that every bounded operator  $T:\mathfrak{X}_U\to\mathfrak{X}_U$ is a scalar multiple of the identity plus a horizontally compact operator. 

\begin{definition}[h-exact pair]\label{def h-exact pair} Let $C>0$, $j\in L_0\setminus\{0\}$, $\varepsilon\in\{0,1\}$, and $(x_k)_{k=1}^\infty$ be a $(C,2)$-h-RIS. A pair $(x,\eta)\in\mathfrak{X}_U\times \Gamma$ is called a $(j,\varepsilon)$-h-exact pair of $(x_k)_{k=1}^\infty$, provided
\begin{enumerate}[label=(\alph*)]

\item\label{h-exact-1} $x=m_jn_j^{-1/q}\sum_{k\in A}x_k$, for some $A\subset\N$ with $\# A=n_j$,

\item\label{h-exact-2} $\mathrm{weight}(\eta)=m_j^{-1}$,

\item\label{h-exact-3} $e^*_\eta(x)=1$, if $\varepsilon=1$, and $|e^*_\eta(x)|\leq m_j^{-1}$, if $\varepsilon=0$.

\end{enumerate}
\end{definition}

\begin{lemma}\label{h-exact estimates} Fix $C>0$,  $j\in L_0\setminus\{0\}$ and $\varepsilon\in\{0,1\}$. Then any $(j,\varepsilon)$-h-exact pair $(x,\eta)$ of a $(C,2)$-h-RIS $(x_k)_{k=1}^\infty$ satisfies the following. 
\begin{enumerate}[label=(\alph*)]

\item\label{h-exact estimate 1} $|d^*_\xi(x)|\leq 3Cm_j^{-1}$, for any $\xi\in\Gamma$.

\item\label{h-exact estimate 2} $\|x\|\leq 96C$.

\item\label{h-exact estimate 3} For any $\gamma\in\Gamma$ with $\mathrm{weight}(\gamma)=m_i^{-1}$,  $i\neq j$, and any interval $E$ we have
\[|e^*_{\gamma}\circ P_E(x) |\leq \left\{
\begin{array}{ll}
24Cm_j^{-1} & \mbox{if } i>j, \\
216Cm_i^{-1} & \mbox{if } i<j.
\end{array}
\right.   \]

\end{enumerate}
\end{lemma}

\begin{proof}
For property \ref{h-exact estimate 1}, by Definition \ref{def hris} , $|d_\xi^*(x)| \leq m_jn_j^{-1/q}\|d_\xi^*\|\sup_k\|x_k\|$. Recall that $\|d_\xi^*\|\leq 3$,  
$\sup_k\|x_k\|\leq C$, and by Assumption \ref{assumption integers} \ref{4 assumption integers}, $m_jn_j^{-1/q}\leq m_j^{-1}$. Property \ref{h-exact estimate 2} follows from Corollary \ref{h-ris norm estimates}.   For \ref{h-exact estimate 3}  take $\gamma\in\Gamma$ with $\mathrm{weight}(\gamma)=m_i^{-1}$, $i\neq j$, and interval $E$. By Proposition \ref{h-basic inequality}, for some $k_0$ and 
$g\in W_q[(\mathcal{A}_{4n_j},4^{1/p}m_j^{-1})_{j\in L_0\setminus\{0\}}]$ with  either $g = 0$ or $\mathrm{weight}(g) = m_i^{-1}$, we have
\[|e_\gamma^*\circ P_E(x)| \leq 12C\Big|(e_{k_0}^*+g)\Big(\frac{m_j}{n_j^{1/q}}\sum_{k=1}^{n_j}e_k\Big)\Big|.\]
We combine Proposition \ref{estimates mixed tsirelson} \eqref{1 estimate mixed tsirelson} with Assumption \ref{assumption integers} to finish the proof. 
\end{proof}

Proposition \ref{existence of h-nodes} immediately yields the following.

\begin{lemma}[Existence of h-exact pairs]\label{existence exact pairs} Let $2j\in L_0\setminus\{0\}$ and $\varepsilon\in\{0,1\}$. Fix, for some $C>0$,  a $(C,2)$-h-RIS $(x_k)_{k=1}^{n_{2j}}$. 
Let $(\mu_k)_{k=1}^{n_{2j}}$ be nodes in $\Gamma$ and  $(F_k)_{k=1}^{n_{2j}}$ be intervals of $\N_0$ that  satisfy, for some  indices $n_{2j}\leq p_1<p_2<\cdots<p_{n_{2j}}$, the following.
\begin{enumerate}[label=(\alph*),leftmargin=19pt]
   
    \item\label{existence exact 1} $\mathrm{rank}(\mu_1) \in [0, p_1)$,  $\supp(x_1)\cup F_1\subset[0,p_1)$ and, for $2\leq k\leq n_{2j}$, $\mathrm{rank}(\mu_k)\in (p_{k-1}, p_k)$, $\supp(x_k)\cup F_k\subset (p_{k-1},p_k)$.
   
    \item\label{existence exact 2}  
    For $2\leq k\leq n_{2j}$, $\min\mathrm{supp}_\mathrm{h}(e_{\mu_k}^*\circ P_{F_k}) >  p_{k-1}$ and $\mathrm{supp}_\mathrm{h}(e_{\mu_{k-1}}^*\circ P_{F_{k-1}})\ll \mathrm{supp}_\mathrm{h}(e_{\mu_k}\circ P_{F_k})$.
    
    \item\label{existence exact 3} $e^*_{\mu_k}\circ P_{F_k}(x_k)=1$, if $\varepsilon=1$, and $|e^*_{\mu_k}\circ P_{F_k}(x_k)|\leq m_{2j}^{-1}$, if $\varepsilon=0$. 
\end{enumerate}
Let $x=m_{2j}n_{2j}^{-1/q}\sum_{r=1}^{n_{2j}}x_k$. Then there is $\eta\in\Gamma$ such that $e_\eta^*$ has evaluation analysis
\[e_\eta^* = \sum_{k=1}^{n_{2j}}d_{\xi_k}^* + \frac{1}{m_{2j}}\sum_{k=1}^{n_{2j}}\frac{1}{n_{2j}^{1/p}}e_{\mu_k}^*\circ P_{F_k},\]
for some nodes $(\xi_k)_{k=1}^{n_{2j}}$ in $\Gamma$, with $\mathrm{rank}(\xi_k) = p_k$,  $1\leq k\leq n_{2j}$, and such that $(x,\eta)$ is a $(2j,\varepsilon)$-h-exact pair of $(x_k)_{k=1}^\infty$.  
\end{lemma}

\begin{definition}[h-dependent sequence]\label{def h-dependent} Let $C>0$, $2j_0-1\in L_0$, $\varepsilon\in\{0,1\}$, and $(x_k)_{k=1}^\infty$ be a $(C,2)$-h-RIS. An h-block sequence $(y_i)_{i=1}^{n_{2j_0-1}}$ is called a $(2j_0-1,\varepsilon)$-h-dependent sequence of $(x_k)_k$, provided there 
is an h-node $\zeta\in\Gamma$ with evaluation analysis
\[e_\zeta^* = \sum_{i=1}^{n_{2j_0-1}}d_{\zeta_i}^* + \frac{1}{m_{2j_0-1}}\sum_{i=1}^{n_{2j_0-1}}\frac{1}{n_{2j_0-1}^{1/p}}e_{\eta_i}^*\circ P_{E_i}\]
satisfying the following.
\begin{enumerate}[label=(\alph*),leftmargin=19pt]

\item \label{h-dependent sequence 1} $\mathrm{weight}(\eta_1)=m_{4j_1-2}^{-1}$, where $4j_1-2\in L_0$ satisfies $m_{4j_1-2}\geq n_{2j_0-1}^3$,
and $\mathrm{weight}(\eta_i)=m_{4j_i}^{-1}$, with  $4j_i\in L_0\setminus \{0\}$, for $2\leq i\leq n_{2j_0-1}$.

\item \label{h-dependent sequence 2} $\supp(y_i)\subset E_i$ for $1\leq i\leq n_{2j_0-1}$.

\item \label{h-dependent sequence 3} $(y_1,\eta_1)$ is a $(4j_1-2,\varepsilon)$-h-exact pair of $(x_k)_{k=1}^\infty$ and $(y_i,\eta_i)$ is a $(4j_i,\varepsilon)$-h-exact pair of $(x_k)_{k=1}^\infty$ for $2\leq i\leq n_{2j_0-1}$. 

\end{enumerate}
\end{definition}

It is important that the exact pair components $(y_1\eta_1),\ldots,(y_{n_{2j_0-1}},\eta_{n_{2j_0-1}})$ of an h-dependent sequence are built on a common underlying RIS $(x_k)_{k=1}^\infty$. This is further elaborated on in the paragraph preceding Proposition \ref{h-basic inequality refined}.

\begin{lemma}[Existence of h-dependent sequences]\label{existence h-dependent}
Let $C>0$,  $2j_0-1\in L_0$, $\varepsilon\in\{0,1\}$, and $(x_k)_{k=1}^\infty$ be a $(C,2)$-h-RIS. Let also $(z_k)_{k=1}^\infty$ be an h-block sequence, $(\mu_k)_{k=1}^\infty$ be nodes in $\Gamma$,  and $(F_k)_{k=1}^\infty$ be successive intervals of $\N_0$ with $\lim_k\min\mathrm{supp}_\mathrm{h}(e^*_{\mu_k}\circ P_{F_k})=\infty$. Assume, furthermore, that $\lim_ke^*_{\mu_k}\circ P_{F_k}(x_k)=0$, if $\varepsilon=0$, and $e^*_{\mu_k}\circ P_{F_k}(x_k)=1$, $k\in\N$, if $\varepsilon=1$. 

Then, there is a $(2j_0-1,\varepsilon)$-h-dependent sequence $(y_i)_{i=1}^{n_{2j_0-1}}$ of $(x_k)_{k=1}^\infty$ satisfying the following.
\begin{enumerate}[label=(\roman*)]

\item For $1\leq i\leq n_{2j_0-1}$, $y_i=m_{l_i}n_{l_i}^{-1/q}\sum_{k\in A_i}x_k$, for some $A_i\subset\N$ with $\# A_i=n_{l_i}$.

\item The associated h-node $\zeta\in\Gamma$ satisfies, for $1\leq i\leq n_{2j_0-1}$,
\[m_{2j_0-1}n_{2j_0-1}^{1/p}e^*_\zeta\left(\frac{m_{l_i}}{n_{l_i}^{1/q}}\sum_{k\in A_i}z_k\right)=\frac{1}{n_{l_i}}\sum_{k\in A_i}e^*_{\mu_k}\circ P_{F_k}(z_k).\]

\end{enumerate}

Moreover, for $1\leq i\leq n_{2j_0-1}$, we can choose $A_{i}\subset 2\N$, if $i$ is even, and $A_{i}\subset 2\N-1$, if $i$ is odd.
\end{lemma}

\begin{proof} 
Passing to a subsequence, we can assume that, for some $(p_k)_{k=1}^\infty$, we have  
\begin{enumerate}[label=(\alph*)]
    \item $\rank(\mu_1)\in [0,p_1)$, $\supp(x_1)\cup\supp(z_1)\cup F_1\subset [0,p_1)$, and, for $k\geq 2$, $\rank(\mu_k)\in (p_{k-1},p_k)$, $\supp(x_k)\cup \supp(z_k)\cup F_k\subset (p_{k-1},p_k)$.
    \item For $k\geq 2$, $\min\mathrm{supp}_\mathrm{h}(e^*_{\mu_k}\circ P_{F_k})> p_{k-1}$ and $\mathrm{supp}_\mathrm{h}(e^*_{\mu_{k-1}}\circ P_{F_{k-1}})\ll\mathrm{supp}_\mathrm{h}(e^*_{\mu_k}\circ P_{F_k})$. 
\end{enumerate}
Concerning the ``moreover'' part, notice that we can pass to a subsequence as above containing infinitely many even and infinitely many odd indices.

We treat first the case of $\varepsilon=1$. We build a $(2j_0-1,1)$-h-dependent sequence $(y_i)_{i=1}^{n_{2j_0-1}}$ of the $(C,2)$-h-RIS $(x_k)_{k=1}^\infty$ inductively in the following way. Pick $4i_1-2\in L_0$ with $m_{4i_1-2}>n_{2j_0-1}^3$ and let $k_1=0$. Applying Lemma \ref{existence exact pairs}, pick a suitable $(4j_1-2,1)$-h-exact pair $(y_1,\eta_1)$ of $(x_k)_{k\in A_1}$, for any $A_1\subset\N$, $\# A_1=n_{4j_1-2}$, with $e^*_{\eta_1}=\sum_{k\in A_1}d^*_{\xi_k}+m_{4j_1-2}^{-1}n_{4j_1-2}^{-1/p}\sum_{k\in A_1}e^*_{\mu_k}\circ  P_{F_k}$, $\rank(\xi_k)\in (p_l)_{l=1}^\infty$ for all $k\in A_1$. For the ``moreover'' part, we choose $A_1$ from $2\N-1$.  
Pick $q_1>\supp(y_1)\cup\supp(w_1)$, where $w_1=m_{4j_1-2}n_{4j_1-2}^{-1/q}\sum_{k\in A_1}z_k$. Let $E_1=(0,q_1)$ and $\zeta_1=(q_1,m_{2j_0-1}^{-1},E_1,n_{2j_0-1}^{-1/p}, \eta_1)$.  

Let $4j_2=\sigma(\xi_1)$. Pick $k_2$ with $\supp(x_{k_2})>q_1$ and $\supp(z_{k_2})>q_1$.  Applying Lemma \ref{existence exact pairs}, pick a suitable $(4j_2,1)$-h-exact pair $(y_2,\eta_2)$ of $(x_k)_{k\in A_2}$, for any $A>k_2$ with $\# A_2=n_{4j_2}$, such that $e^*_{\eta_2}=\sum_{k\in A_2}d^*_{\xi_k}+m_{4j_2}^{-1}n_{4j_2}^{-1/p}\sum_{k\in A_2}e^*_{\mu_k}\circ  P_{F_k}$, $\rank(\xi_k)\in (p_l)_{l=1}^\infty$ for all $k\in A_2$. Again, to achieve the ``moreover'' part of the statement, choose $A_2$ from $2\N$. Pick $q_2>\supp(y_2)\cup\supp(w_2)$, where $w_2=m_{4j_2}n_{4j_2}^{-1/q}\sum_{k\in A_2}z_k$. Let $E_2=(q_1,q_2)$ and $\xi_2=(q_2,m_{2j_0-1}^{-1},\xi_1,E_2,n_{2j_0-1}^{-1/p},\eta_2)$. 

Continuing the process described above, as in Proposition \ref{existence of h-nodes}, according to Remark \ref{plentitude of even h-nodes}, we obtain a $(2j_0-1,1)$-h-dependent sequence $(y_i)_{i=1}^{n_{2j_0-1}}$ of $(x_k)_{k=1}^\infty$ with associated $\zeta=\zeta_{n_{2j_0-1}}\in\Gamma$. To prove the desired estimation of $e^*_\zeta$ on weighted averages $(w_i)_{i=1}^{n_{2j_0-1}}$ of $(z_k)_{k=1}^\infty$, fix $1\leq i\leq n_{2j_0-1}$. Note that $\supp(e^*_\zeta)\cap \supp(w_i)\subset \supp(e^*_{\eta_i}\circ P_{E_i})\cap\supp(w_i)\subset \cup_{k\in A_i}\supp(e^*_{\mu_k}\circ P_{F_k})$. Moreover, $m_{l_i}^{-1}=\mathrm{weight}(\eta_i)$. Thus, 
\begin{align*}
m_{2j_0-1}n_{2j_0-1}^{1/q}e^*_{\zeta}(w_i)=\frac{1}{m_{l_i}n_{l_i}^{1/p}}\sum_{k\in A_i}e^*_{\mu_k}\circ P_{F_k}\Big(\frac{m_{l_i}}{n_{l_i}^{1/q}}\sum_{k\in A_i}z_k\Big)=\frac{1}{n_{l_i}}
\sum_{k\in A_l}e^*_{\mu_k}\circ P_{F_k}(z_k).
\end{align*}
which ends the proof for $\varepsilon=1$. 

Now assume that $\varepsilon=0$ and $\lim_ke^*_{\mu_k}\circ P_{F_k}(x_k)=0$. We proceed as in the previous case, with the following  additional condition on the choice of $k_1,k_2,\dots$. For any $1\leq i\leq n_{2j_0-1}$, pick $k_i$ such that, additionally, for any $k\geq k_i$, $|e^*_{\mu_k}\circ P_{F_k}(x_k)|\leq m_{l_i}^{-1}$, where $l_1=4j_1-2$ and $l_i=4j_i$ for any $2\leq i\leq n_{2j_0-1}$. Then the resulting pairs $(y_i,\eta_i)$, $1\leq i\leq n_{2}$ satisfy
\[|e^*_{\eta_i}\circ P_{E_i}(y_i)|=\Big|\frac{1}{m_{l_i}n_{l_i}^{1/p}}\sum_{k\in A_i}e^*_{\mu_k}\circ P_{F_k}\Big(\frac{m_{l_i}}{n_{l_i}^{1/q}}\sum_{k\in A_i}x_k\Big)\Big|\leq \frac{1}{n_{l_i}}\sum_{k\in A_i}|e^*_{\mu_k}\circ P_{F_k}(x_k)|\leq \frac{1}{m_{l_i}}.\]
As the other properties of the chosen nodes are preserved, we finish the proof of the case of $\varepsilon=0$.
\end{proof}

As already discussed, h-dependent sequences are not h-RISs; thus, the Basic Inequality for h-RISs (Proposition \ref{h-basic inequality}) cannot directly provide norm estimates on linear combinations of h-dependent sequences. Instead, we require that h-dependent sequences are sequences of weighted averages of a common RIS. To exploit this, we state and prove a refined version of the Basic Inequality, which provides upper bounds on linear combinations of h-dependent sequences in terms of estimates on nested averages in the auxiliary $q$-mixed-Tsirelson space.

\begin{proposition}[Second basic inequality for h-RIS]
\label{h-basic inequality refined} Fix $C>0$ and $2j_0-1\in L_0$. 
Let $(x_k)_{k\in I}$ be a $(C,2)$-h-RIS, indexed over an interval $I$ of $\N$, and let $(a_k)_{k\in I}$ be a sequence of complex numbers, satisfying the following, for some partition $I=\cup_{i=1}^b I_i$ into intervals $I_1<\dots<I_b$.
\begin{enumerate}[label=(\alph*)]

\item\label{basic inequality refined 1} For $1\leq i\leq b$, $\|\sum_{k\in I_i}a_ke_k\|\geq 1$, in  $T_q[(\mathcal{A}_{4n_j},4^{1/p}m_j^{-1})_{j\in L_0\setminus\{0\}}]$.

\item\label{basic inequality refined 2} For any $\gamma\in\Gamma$ with $\mathrm{weight}(\gamma)=m_{2j_0-1}^{-1}$ and any intervals $E\subset \N$, $F\subset I$, 
\[\Big|e^*_\gamma\circ P_E\Big(\sum_{k\in F}a_kx_k\Big)\Big|\leq C.\]

\end{enumerate}
Then, for any node   $\gamma\in\Gamma$ and interval $E$ of $\N_0$, there exist $g\in W_q[(\mathcal{A}_{4n_j},4^{1/p}m_j^{-1})_{j\in L_0\setminus\{0\}}]$ and $k_0\in I$ such that
\begin{enumerate}[label=(\roman*)]

\item\label{h-dep basic inequality support} $\supp(g) \subset I\setminus\{k_0\}$,

\item\label{h-dep basic inequality tree-analysis} $g$ has a tree-analysis $(g_t)_{t\in\mathcal{T}}$ such that for any $t\in\mathcal{T}$ with $\mathrm{weight}(g_t)=4^{1/p}m_{2j_0-1}^{-1}$,  $\supp(g_t)$ intersects at most two of the intervals $I_1,\ldots,I_b$ and

\item\label{h-dep basic inequality estimate}
\[\Big|e_\gamma^*\circ P_E\Big(\sum_{k\in I}a_kx_k\Big)\Big| \leq 12C\Big|(e_{k_0}^*+g)\Big(\sum_{k\in I}|a_k|e_k\Big)\Big|.\]

\end{enumerate}

\end{proposition}

\begin{proof}

The proof if performed by induction on the rank of $\gamma\in\Gamma$. 
The inductive hypothesis states that for a given rank $n$ the conclusion holds for all $\gamma\in\Gamma$ of rank $n$, all intervals $E$ of $\mathbb{N}_0$, and all subintervals $F\subset I$. We repeat the inductive proof of Proposition \ref{h-basic inequality} line for line with one exception, treating the case $\mathrm{weight}(\gamma)=m_{2j_0-1}^{-1}$ separately. 
Note that in the proof of Proposition \ref{h-basic inequality}, for a fixed $\gamma\in \Gamma$ with $\mathrm{weight}(\gamma)\neq m_{2j_0-1}^{-1}$, a non-zero $g$ is constructed as a suitable weighted average of functionals provided by the inductive hypothesis and basic functionals $(e_k^*)_k$. Thus, assuming the functionals building $g$ satisfy \ref{h-dep basic inequality tree-analysis}, we obtain \ref{h-dep basic inequality tree-analysis} also for $g$. 

Take now $\gamma\in \Gamma$ with $\mathrm{weight}(g)=m_{2j_0-1}^{-1}$, an interval $E\subset \N_0$ and  a subinterval $F\subset I$, and assume that the statement holds for all $\gamma'\in \Gamma$ with $\mathrm{rank}(\gamma')<\mathrm{rank}(\gamma)$, all intervals $E'\subset \N_0$ and $F'\subset F$. 

\noindent\textbf{Case 1:} $F$ intersects at most two of the intervals $I_1,\ldots,I_b$. Then we 
proceed as in the proof of Proposition \ref{h-basic inequality} for $\gamma$ with weight $m_j^{-1}$, $j\in L_0\setminus\{0\}$, obtaining suitable  $k_0\in F$ and $g$ with the support intersecting at most two of intervals $(I_i)_{i=1}^b$, thus satisfying conditions \ref{h-dep basic inequality support}, \ref{h-dep basic inequality tree-analysis}, \ref{h-dep basic inequality estimate}. 

\noindent\textbf{Case 2:} $F$ intersects at least three of the intervals $I_1,\ldots,I_b$. Pick $1\leq i_0\leq b$ with $I_{i_0}\subset  F$ and $k_0\in F\setminus I_{i_0}$. Take $g\in W_q[(\mathcal{A}_{4n_j},4^{1/p}m_j^{-1})_{j\in L_0\setminus\{0\}}]$ supported on $I_{i_0}$ with $g(\sum_{k\in I_{i_0}}|a_k|x_k)\geq 1$.  Then, 
by  assumption,
\[\Big|e^*_\gamma\circ P_E\Big(\sum_{k\in F}a_ke_k\Big)\Big|\leq C\leq C(e_{k_0}^*+g)\Big(\sum_{k\in F}|a_k|e_k\Big)\]
and conditions \ref{h-dep basic inequality support}, \ref{h-dep basic inequality tree-analysis} and \ref{h-dep basic inequality estimate} follow immediately.
\end{proof}

\begin{corollary}\label{sum of h-dependent sequence} Fix $C\geq 1$, $2j_0-1\in L_0$ and $\varepsilon\in\{0,1\}$. 
Let $(y_i)_{i=1}^{n_{2j_0-1}}$ be a $(2j_0-1, \varepsilon)$-h-dependent sequence of some $(C,2)$-h-RIS. Then 
\[\Big\|
\frac{1}{n_{2j_0-1}^{1/q}}\sum_{i=1}^{n_{2j_0-1}}(-1)^{\varepsilon i}y_i\Big\|\leq \frac{2^{20}C}{m_{2j_0-1}^2}.
\]
\end{corollary}

\begin{proof} 
This proof consists of a two-step process. In the first one, we show that the assumptions of Proposition \ref{h-basic inequality refined} are satisfied, for an appropriate h-RIS and scalar coefficients. In the second step we invoke an estimate on nested averages in the auxiliary $q$-mixed-Tsirelson space, namely Proposition \ref{estimates 2 mixed tsirelson}.

The reasoning of the first step follows the lines of the  proof of {\cite[Lemma 6.5]{argyros:haydon:2011}}. Let each $y_i$, $1\leq i\leq n_{2j_0-1}$, be of the form $y_i=m_{l_i}n_{l_i}^{-1/q}\sum_{k\in A_i}x_k$, for some $(C,2)$-h-RIS $(x_k)_{k=1}^\infty$, where $l_1=4j_1-2$ and $l_i=4j_i$, for $2\leq i\leq n_{2j_0-1}$. Let the associated h-node $\zeta$ have evaluation analysis
\[e_\zeta^* = \sum_{i=1}^{n_{2j_0-1}}d_{\zeta_i}^* + \frac{1}{m_{2j_0-1}}\sum_{i=1}^{n_{2j_0-1}}\frac{1}{n_{2j_0-1}^{1/p}}e_{\eta_i}^*\circ P_{E_i}.\]

We define $(a_k)_{k\in A}$, $A=\cup_{i=1}^{n_{2j_0-1}}A_i$, by letting $a_k=(-1)^{\varepsilon i}m_{l_i}n_{l_i}^{-1/q}$ whenever $k\in A_i$, $1\leq i\leq n_{2j_0-1}$. Assumption \ref{basic inequality refined 1} of Proposition \ref{h-basic inequality refined} is satisfied, as witnessed by taking, for $1\leq i\leq n_{2j_0-1}$, $g = 4^{1/q}m^{-1}_{l_i}n^{-1}_{l_i}\sum_{k\in A_i}(-1)^{\varepsilon i}e_k^*$. We will prove assumption \ref{basic inequality refined 2} with the constant $2^{11}C$.

\noindent\textbf{Case 1:} $\varepsilon=1$. Fix $\zeta'\in\Gamma$ with $\mathrm{weight}(\zeta')=m_{2j_0-1}^{-1}$, intervals $E\subset\N$ and $F\subset A$. Let $e^*_{\zeta'}$ have evaluation analysis
\[e_{\zeta'}^* = \sum_{i=1}^{a}d_{\zeta'_i}^* + \frac{1}{m_{2j_0-1}}\sum_{i=1}^{a}\frac{1}{n_{2j_0-1}^{1/p}}e_{\eta'_i}^*\circ P_{E'_i},\]
for some $a\leq n_{2j_0-1}$ and with $\mathrm{weight}(\eta'_i)=m_{l'_i}^{-1}$, $1\leq i\leq a$, where $l'_1=4j'_1-2$ and $l'_i=4j'_i$, for $2\leq i\leq a$. By the definition of special h-nodes, for some $1\leq l\leq a$, we have  
\begin{enumerate}[label=(\roman*)]
    \item\label{sum of dependent tech 1} for $r<l$, $\eta'_r=\eta_r$, $E'_r=E_r$, $\zeta'_r=\zeta_r$,
    \item\label{sum of dependent tech 2} for $i>l$ and $r\leq a$, $j'_r\neq j_i$. 
\end{enumerate}

For any $i<l$, as $\supp(y_i)\subset E_i=E'_i$ by \ref{sum of dependent tech 1}, we have
\begin{equation}\label{sum of dependent est 1}
e^*_{\zeta'}(y_i)=e^*_\zeta(y_i)=\frac{1}{m_{2j_0-1}n_{2j_0-1}^{1/p}}e^*_{\eta_i}(y_i)=\frac{1}{m_{2j_0-1}n_{2j_0-1}^{1/p}}.
\end{equation}
Also, for any $i>l$ and $1\leq r\leq a$, by \ref{sum of dependent tech 2}, Lemma \ref{h-exact estimates}\ref{h-exact estimate 3} and item \ref{h-dependent sequence 1} of Definition \ref{def h-dependent}, we have
\begin{equation}\label{sum of dependent est 2} 
|e_{\eta'_r}\circ P_{E'_r}(y_i)|\leq 216C\frac{1}{\min\{m_{l_1},m_{l_1'}\}}\leq 216C\frac{1}{n_{2j_0-1}^2}.
\end{equation}

Let $E'=E\cap \cup_{k\in F}\supp(x_k)$. Pick $1\leq i_1\leq i_2\leq n_{2j_0-1}$ with $\min E'\in \supp(y_{i_1})$ and $\max E'\in\supp(y_{i_2})$. We will assume $i_1<l<i_2$ as the other cases have simpler estimates. Estimate, using Lemma \ref{h-exact estimates}, item \ref{h-dependent sequence 1} of Definition \ref{def h-dependent}, \eqref{sum of dependent est 1}, and \eqref{sum of dependent est 2},
\begin{align*}
    \Big|e^*_{\gamma'}\circ P_E\Big(\sum_{k\in F}a_kx_k\Big)\Big|&= \Big|e^*_{\gamma'}\circ P_{E'}\Big(\sum_{i=i_1}^{i_2}(-1)^iy_i\Big)\Big|\\
    &\leq |e^*_{\gamma'}\circ P_{E'}(y_{i_1})|+\Big|e^*_{\gamma'}\Big(\sum_{i_1<i<l}(-1)^iy_i\Big)\Big|+  |e^*_{\gamma'}\circ P_{E'}(y_l)|\\
    &+\sum_{r=1}^a\sum_{i>l} |d^*_{\zeta'_r}(y_i)|+\sum_{r=1}^a\sum_{i>l}|e^*_{\eta'_r}\circ P_{E_r'}(y_i)|+|e^*_{\gamma'}\circ P_{E'}(y_{i_2})|\\
    &\leq 384C+1+384C+3C\frac{n_{2j_0-1}^2}{m_{l_1}}+216C\frac{n_{2j_0-1}^2}{m_{l_1}}+384C\\
    &\leq 2^{11}C.
\end{align*}

\noindent\textbf{Case 2:} $\varepsilon=0$. We repeat the estimates of Case 1, with one modification, i.e. instead of \eqref{sum of dependent est 1} we have 
\begin{equation}
|e^*_{\zeta'}(y_i)|=|e^*_\zeta(y_i)|=\frac{1}{m_{2j_0-1}n_{2j_0-1}^{1/p}}|e^*_{\eta_i}(y_i)|\leq \frac{1}{m_{2j_0-1}n_{2j_0-1}^{1/p}}\cdot \frac{1}{m_{l_i}},
\end{equation}
 which also ensures
 \[ \Big|e^*_{\gamma'}\Big(\sum_{i_1<i<l}y_i\Big)\Big|\leq 1.\]
The other estimates of Case 1 follows readily.

Note that the parameters $2j_0-1,l_1,\dots,l_{n_{2j_0-1}}\in\N$ satisfy the assumption of Proposition \ref{estimates 2 mixed tsirelson} by the definition of odd-weight h-nodes, the definition of $\sigma$, and Assumption \ref{2 assumption integers}. Therefore, in both cases, Proposition \ref{h-basic inequality refined} and Proposition \ref{estimates 2 mixed tsirelson}  yield 
\begin{align*}
    \Big\|\sum_{i=1}^{n_{2j_0-1}}(-1)^{\varepsilon i}y_i\Big\|\leq 12\cdot 2^{11}C\Big(\frac{m_{l_1}}{n_{l_1}^{1/q}}+28\frac{n_{2j_0-1}^{1/q}}{m_{2j_0-1}^2}\Big)\leq 2^{20}C\frac{n_{2j_0-1}^{1/q}}{m_{2j_0-1}^2},
\end{align*}
which ends the proof. 
\end{proof}

Having the above standard corollary, we proceed to the study of bounded operators on $\mathfrak{X}_U$ with respect to its horizontal structure. The experienced reader will notice that we have not encoded the rational convex combinations of functionals on $\mathfrak{X}_U$ in $\Gamma$, which typically allow the direct proof of the scalar-plus-horizontally-compact property via the Hahn-Banach theorem. Instead, we prove the desired behaviour of a bounded operator first on the basis of $\mathfrak{X}_U$ and then extend it to any h-RIS, in both cases by using h-dependent sequences and their properties. 

\begin{proposition}\label{scalar + h-compact on basis} Let $T\in\mathcal{L}(\mathfrak{X}_U)$. Then 
\[\lim_{\rank(\gamma)\to\infty, \gamma\in\Gamma^0} \mathrm{dist}(Td_\gamma,\mathbb{C}d_\gamma)=0.\]  
\end{proposition}
\begin{proof}
The proof follows the reasoning of {\cite[Proposition 7.4]{motakis:2024}}, we present it here for the sake of completeness. We proceed in two steps. 

\noindent\textbf{Step 1.} We prove that 
\[\lim_{\rank(\gamma)\to\infty, \gamma\in\Gamma^0}\|P_{\{\rank(\gamma)\}}Td_\gamma-Td_\gamma\|=0.\]
Suppose 
$\inf_k \|P_{\{\rank(\gamma_k)\}}Td_{\gamma_k} - Td_{\gamma_k}\|>0$ for some $(\gamma_k)_{k=1}^\infty$ in $\Gamma^0$ with $\rank(\gamma_k)$ strictly increasing.
By Proposition, \ref{horizontal compactness of projections} $\lim_k\|I_s(Td_{\gamma_k})\| = 0$, for any $s\in\N$. By Lemma \ref{close to h-block} and  Fact \ref{standard perturbation}, passing to a subsequence if necessary, we can assume that 
$(Td_{\gamma_k})_{k=1}^\infty$ is an h-block sequence. 

By our assumption, passing to a subsequence, we have either $\inf_k\|P_{(\rank(\gamma_k),\infty)}Td_{\gamma_k}\|>0$ or $\inf_k\|P_{[0,\rank(\gamma_k))}Td_{\gamma_k}\|>0$. We shall treat the first case; the second one is proved analogously. By Lemma \ref{h-block norming}, passing to a subsequence, we can pick  nodes $(\mu_k)_{k=1}^\infty\subset\Gamma$ and successive intervals $(F_k)_{k=1}^\infty$, with $F_k>\rank(\gamma_k)$, for all $k\in\N$,  such that $\lim_k\min\mathrm{supp}_\mathrm{h}(e^*_{\mu_k}\circ P_{F_k})=\infty$ and  $\delta:=\inf_{k\in\N}\mathfrak{Re}(e^*_{\mu_k}\circ P_{F_k}(Td_{\gamma_k}))>0$. Note that $e^*_{\mu_k}\circ P_{F_k}(d_{\gamma_k})=0$ for any $k\in\N$. 
 
Now we apply Lemma \ref{existence h-dependent} for $\varepsilon=0$, the $(1,1)$-h-RIS  $(d_{\gamma_k})_{k=1}^\infty$,  the h-block sequence $(Td_{\gamma_k})_{k=1}^\infty$, and the nodes $(\mu_k)_{k=1}^\infty$ and intervals  $(F_k)_{k=1}^\infty$ chosen above. We, thus, obtain for any $2j_0-1\in L_0$, a suitable $(2j_0-1,0)$-h-dependent sequence $(y_i)_{i=1}^{n_{2j_0-1}}$ with associated $\zeta\in\Gamma^0$. For $y=\frac{m_{2j_0-1}}{n_{2j_0-1}^{1/q}}\sum_{i=1}^{n_{2j_0-1}}y_i$ we have 
\[\mathfrak{Re}\Big(e^*_\zeta(Ty)\Big)=\frac{1}{n_{2j_0-1}}\sum_{i=1}^{n_{2j_0-1}}\mathfrak{Re}\Big(e^*_\zeta(y_i)\Big) = \frac{1}{n_{2j_0-1}}\sum_{i=1}^{n_{2j_0-1}}\frac{1}{n_{l_i}}\sum_{k\in A_i}\mathfrak{Re}\Big(e^*_{\mu_k}\circ P_{F_k}(Td_{\gamma_k})\Big)\geq\delta\]
(for the last inequality recall that $\# A_i=n_{l_i}$ for any $1\leq i\leq n_{2j_0-1}$.) 
On the other hand, by Corollary \ref{sum of h-dependent sequence}, $\|y\|\leq 2^{20} m_{2j_0-1}^{-1}$, which for sufficiently big $j_0$ contradicts the boundedness of $T$. 

\noindent\textbf{Step 2.}
We assume now that $\supp(Td_\gamma)=\{\rank(\gamma)\}$ for all $\gamma\in\Gamma^0$ and  we show that 
\[
\lim_{\rank(\gamma)\to\infty,\gamma\in\Gamma^0}\|Td_\gamma-d^*_\gamma(Td_\gamma)d_\gamma\|=0.
\]
Suppose now that $\inf_k \|Td_{\gamma_k}-d^*_{\gamma_k}(Td_{\gamma_k})d_{\gamma_k}\|>0$, for some $(\gamma_k)_{k=1}^\infty$ in $\Gamma^0$, with $\lim_k\rank(\gamma_k) = \infty$. As before, passing to a subsequence if necessary, we may assume that $(Td_{\gamma_k})_{k=1}^\infty$ is an h-block sequence.  Recall that, for $n\in\N_0$, $(d_{\gamma'})_{\gamma'\in\Delta_n}$ is 2-equivalent to the unit vector basis of $\ell_\infty(\Delta_n)$. Thus, the assumption on $(\gamma_k)_{k=1}^\infty$ yields a sequence $(\gamma_k')_{k=1}^\infty$, $\gamma'_k\neq\gamma_k$, with $\rank(\gamma'_k)=\rank(\gamma_k)$, $k\in\N$, such that $\inf_k|d^*_{\gamma_k'}(Td_{\gamma_k})|\geq \delta>0$. As $(Td_{\gamma_k})_{k=1}^\infty$ is h-block, $(d_{\gamma'_k})_{k=1}^\infty$ is also h-block. 

As in Step 1, we apply   Lemma \ref{existence h-dependent} for $\varepsilon=0$, the $(1,1)$-h-RIS  $(d_{\gamma_k})_{k=1}^\infty$,  the h-block sequence $(Td_{\gamma_k})_{k=1}^\infty$, the nodes $(\gamma_k')_{k=1}^\infty$ and intervals  $(\{\rank(\gamma_k')\})_{k=1^\infty}$. We, thus, obtain for any $2j_0-1\in L_0$ a suitable $(2j_0-1,0)$-h-dependent sequence $(y_i)_{i=1}^{n_{2j_0-1}}$ with associated $\zeta\in\Gamma^0$. As in Step 1, for $y=m_{2j_0-1}^{-1}n_{2j_0-1}^{1/q}\sum_{i=1}^{n_{2j_0-1}}y_i$ we have $e^*_\zeta(Ty)\geq \delta$, whereas by Corollary \ref{sum of h-dependent sequence}, $\|y\|\leq 2^{20} m_{2j_0-1}^{-1}$, which gives contradiction for sufficiently large $j_0$. \end{proof}

\begin{proposition}\label{null on basis to null on h-RIS}
Let $T\in\mathcal{L}(\mathfrak{X}_U)$ satisfy $\lim_k\|Td_{\gamma_k}\|=0$, for some sequence $(\gamma_k)_{k=1}^\infty$ of pairwise distinct members of $\Gamma^0$. Then,  $\lim_k\|Tx_k\|=0$, for any $(C,1)$-h-RIS $(x_k)_{k=1}^\infty$.     
\end{proposition}

\begin{proof}
Fix a $(C,1)$-h-RIS $(x_k)_{k=1}^\infty$. By Proposition \ref{horizontal compactness of projections}, $\lim_k\|I_s(Tx_k)\| = 0$, for any $s\in\N$. By  Lemma \ref{close to h-block} and  Fact \ref{standard perturbation}, passing to a subsequence if necessary, we can assume that, for all $k\in\N$, $Td_{\gamma_k} = 0$ and that  $(Tx_k)_{k=1}^\infty$ is an h-block sequence. Assume, towards contradiction, that $\inf_{k\in\N}\|Tx_k\|>0$. By Lemma \ref{h-block norming}, passing to a subsequence, we can pick  nodes $(\mu_k)_{k=1}^\infty$ in $\Gamma$ and successive intervals $(F_k)_{k=1}^\infty$ of $\N_0$  such that $\lim_k\min\mathrm{supp}_\mathrm{h}(e^*_{\mu_k}\circ P_{F_k}) = \infty$ and  $\delta:=\inf_{k\in\N}\mathfrak{Re}(e^*_{\mu_k}\circ P_{F_k}(Tx_k))>0$. We consider two cases. 

\noindent\textbf{Case 1:} $\liminf_k|e^*_{\mu_k}(P_{F_k}x_k)| =  0$. We will not use the sequence $(d_{\gamma_k})_{k=1}^\infty$ in this case. Passing to a subsequence, we can assume that $\lim_ke^*_{\mu_k}(P_{F_k}x_k)=0$. 

We use Lemma \ref{existence h-dependent} for $\varepsilon=0$, the  $(C,2)$-h-RIS  $(x_k)_{k=1}^\infty$,  the h-block sequence $(Tx_k)_{k=1}^\infty$, the nodes $(\mu_k)_{k=1}^\infty$ and the intervals  $(F_k)_{k=1}^\infty$. Thus, for any $2j_0-1\in L_0$ we obtain a $(2j_0-1,0)$-h-dependent sequence $(y_i)_{i=1}^{n_{2j_0-1}}$, where, for $1\leq i\leq n_{2j_0-1}$, $y_{i}=m_{l_{i}}n_{l_{i}}^{-1/q}\sum_{k\in A_{i}}x_k$ and an associated $\zeta\in\Gamma^0$. For $y=m_{2j_0-1}n_{2j_0-1}^{-1/q}\sum_{i=1}^{n_{2j_0-1}}y_i$ we estimate, by Lemma \ref{existence h-dependent},
\begin{align*} 
\mathfrak{Re}\Big(e^*_\zeta(Ty)\Big)&= \frac{m_{2j_0-1}}{n^{1/q}_{2j_0-1}}\sum_{i=1}^{n_{2j_0-1}}\mathfrak{Re}\Big(e^*_\zeta(Ty_i)\Big) = \frac{1}{n_{2j_0-1}}\sum_{i=1}^{n_{2j_0-1}}\frac{1}{n_{l_{i}}}\sum_{k\in A_{i}}\mathfrak{Re}\Big(e^*_{\mu_k}(P_{F_k}Tx_k)\Big)>\delta.\\
\end{align*}
On the other hand, by Corollary \ref{sum of h-dependent sequence}, $\|y\|\leq 2^{20}C m_{2j_0-1}^{-1}$, which for sufficiently large $j_0$, contradicts the boundedness of $T$.

\noindent\textbf{Case 2:} $\liminf_k|e^*_{\mu_k}(P_{F_k}x_k)| > 0$. By passing to a subsequence, there is $\delta_1>0$ such that, for all $k\in\N$, $|e^*_{\mu_k}(P_{F_k}x_k)| > \delta_1$. The sequence $(x_k')_{k=1}^\infty$, where $x_k'=(e^*_{\mu_k}(P_{F_k}x_k))^{-1}x_k$, $k\in\N$, is a $(C',1)$-h-RIS, where $C'=C\delta_1^{-1}$, such that  $|e^*_{\mu_k}(P_{F_k}Tx_k')|\geq \delta/(4C)$, $k\in\N$. Thus, by passing to a subsequence, we either have that, for all $k\in\N$, $\mathfrak{Re}(e^*_{\mu_k}(P_{F_k}Tx_k')) \geq \sqrt 2\delta/(8C)$, or, for all $k\in\N$, $\mathfrak{Im}(e^*_{\mu_k}(P_{F_k}Tx_k')) \geq \sqrt 2\delta/(8C)$. We  will assume that the first alternative holds.

We will use the full statement of Lemma \ref{existence h-dependent} for $\varepsilon=1$,  the $(C,2)$-h-RIS  $(d_{\gamma_1},x_1',d_{\gamma_2},x_2',\dots)$ (see Remark \ref{interwined RIS}),  the horizontally block sequence $(0, Tx_1',0, Tx_2',\dots)$,  the sequence of nodes $(\gamma_1,\mu_1, \gamma_2,\mu_2,\dots)$, and the sequence of intervals  $(\{\rank(\gamma_1)\}, F_1,\{\rank(\gamma_2)\}, F_2, \dots,)$. Thus, for any $2j_0-1\in L_0$, we obtain a $(2j_0-1,1)$-h-dependent sequence $(y_i)_{i=1}^{n_{2j_0-1}}$, where, for $1\leq i\leq n_{2j_0 - 1}$ odd, $y_{i}=m_{l_{i}}n_{l_{i}}^{-1/q}\sum_{k\in A_{i}}d_{\gamma_k}$, and, for $1\leq i\leq n_{2j_0-1}$ even,  $y_{i}=m_{l_{i}}n_{l_{i}}^{-1/q}\sum_{k\in A_{i}}x_k'$, and an associated $\zeta\in\Gamma^0$. For $y=m_{2j_0-1}n_{2j_0-1}^{1/q}\sum_{i=1}^{n_{2j_0-1}}(-1)^iy_i$ we estimate, by Lemma \ref{existence h-dependent},
\begin{align*} 
\mathfrak{Re}\Big(e^*_\zeta(Ty)\Big)&= \frac{m_{2j_0-1}}{n^{1/q}_{2j_0-1}}\sum_{i=1}^{\lfloor n_{2j_0-1}/2\rfloor}\mathfrak{Re}\Big(e^*_\zeta(Ty_{2i})\Big) = \frac{1}{n_{2j_0-1}}\sum_{i=1}^{\lfloor n_{2j_0-1}/2\rfloor}\frac{1}{n_{l_{2i}}}\sum_{k\in A_{2i}}\mathfrak{Re}\Big(e^*_{\mu_k}(P_{F_k}Tx'_k)\Big)\\
&\geq \frac{\lfloor n_{2j_0-1}/2\rfloor}{n_{2j_0-1}}\frac{\sqrt 2\delta}{8C}.\\
\end{align*}
On the other hand, by Corollary \ref{sum of h-dependent sequence}, $\|y\|\leq 2^{20}C m_{2j_0-1}^{-1}$, which for sufficiently large $j_0$, contradicts the boundedness of $T$.
\end{proof}

Now we are ready to prove Theorem \ref{scalar-plus-horizontally compact} of Section \ref{section calkin algebra XU}, that is, the following.

\begin{theorem}\label{scalar+h-compact}
For every bounded linear operator $T:\mathfrak{X}_U\to\mathfrak{X}_U$ there exists $\lambda\in\mathbb{C}$ such that $T-\lambda I$ is horizontally compact.    
\end{theorem}
\begin{proof}
By Proposition \ref{scalar + h-compact on basis}, there exist $\lambda\in\mathbb{C}$ and a sequence of nodes $(\gamma_k)_{k=1}^\infty$ in $\Gamma^0$ with $\lim_k\|(T-\lambda I)(d_{\gamma_k})\|=0$. By Proposition \ref{null on basis to null on h-RIS} $\lim_k\|(T-\lambda I)(x_k)\|=0$ for any $(C,1)$-h-RIS $(x_k)_{k=1}^\infty$, thus, by Proposition \ref{horizontal compactness via h-ris}, $(T-\lambda I)$ is horizontally compact.  
\end{proof}

\subsection{$\mathfrak{X}_s$-dependent sequences}\label{subsection xs-dependent} We present now, as in Section \ref{xs-ris}, variants of results on dependent sequences for the component spaces $\mathfrak{X}_s$, $s\in\N$. We only state the results omitting proofs, as they are a direct translation of the reasoning from the previous subsection into the $\mathfrak{X}_s$ setting, with h-RISs replaced by $\mathfrak{X}_s$-RISs, $\Gamma^0$ by $\Gamma^s$, and $L_0$ by $L_s$, respectively. In what follows, we will use the notation of Section \ref{xs-ris}.

\begin{definition}[$\mathfrak{X}_s$-exact pair] Let $s\in\N$, $C>0$, $j\in L_s$ and $\varepsilon\in\{0,1\}$. Fix a $(C,2)$-$\mathfrak{X}_s$-RIS $(x_k)_{k=1}^\infty$. A pair $(x,\eta)\in\mathfrak{X}_U\times \Gamma$ is called a $(j,\varepsilon)$-$\mathfrak{X}_s$-exact pair of $(x_k)_{k=1}^\infty$, provided
\begin{enumerate}[label=(\alph*)]
\item $x=m_jn_j^{-1/q}\sum_{k\in A}x_k$, for some $A\subset\N$, with $\# A=n_j$,
\item $\mathrm{weight}(\eta)=m_j^{-1}$,
\item $e^*_\eta(x)=1$, if $\varepsilon=1$, and $|e^*_\eta(x)|\leq m_j^{-1}$, if $\varepsilon=0$.
\end{enumerate}
\end{definition}

\begin{definition}[$\mathfrak{X}_s$-dependent sequence]\label{def xs-dependent} Let $s\in\N$, $C>0$, $2j_0-1\in L_s$ and $\varepsilon\in\{0,1\}$. Fix a $(C,2)$-$\mathfrak{X}_s$-RIS $(x_k)_{k=1}^\infty$. A block sequence $(y_i)_{i=1}^{n_{2j_0-1}}$ in $\mathfrak{X}_s$ is called a $(2j_0-1,\varepsilon)$-$\mathfrak{X}_s$-dependent sequence of $(x_k)_k$, provided there  
is an node $\zeta\in\Gamma^s$ with evaluation analysis
\[e_\zeta^* = \sum_{i=1}^{n_{2j_0-1}}d_{\zeta_i}^* + \frac{1}{m_{2j_0-1}}\sum_{i=1}^{n_{2j_0-1}}\frac{1}{n_{2j_0-1}^{1/p}}e_{\eta_i}^*\circ P_{E_i}\]
satisfying the following.
\begin{enumerate}[label=(\alph*),leftmargin=19pt]
\item $\mathrm{weight}(\eta_1)=m_{4j_1-2}^{-1}$, where $4j_1-2\in L_s$ satisfies $m_{4j_1-2}\geq n_{2j_0-1}^3$,
and $\mathrm{weight}(\eta_i)=m_{4j_i}^{-1}$, with  $4j_i\in L_s$, for $2\leq i\leq n_{2j_0-1}$.
\item $\supp(y_i)\subset E_i$, for $1\leq i\leq n_{2j_0-1}$.
\item $(y_1,\eta_1)$ is a $(4j_1-2,\varepsilon)$-$\mathfrak{X}_s$-exact pair of $(x_k)_{k=1}^\infty$ and $(y_i,\eta_i)$ is a $(4j_i,\varepsilon)$-$\mathfrak{X}_s$-exact pair of $(x_k)_{k=1}^\infty$, for $2\leq i\leq n_{2j_0-1}$. 
\end{enumerate}
\end{definition}

The following is proved by the $\mathfrak{X}_s$-variant of the refinement of the Basic Inequality \ref{h-basic inequality refined}. 
\begin{proposition} Fix $s\in\N$, $C\geq 1$, $2j_0-1\in L_0$ and $\varepsilon\in\{0,1\}$. 
Let $(y_i)_{i=1}^{n_{2j_0-1}}$ be a $(2j_0-1, \varepsilon)$-$\mathfrak{X}_s$-dependent sequence of a $(C,2)$-$\mathfrak{X}_s$-RIS. Then 
\[\Big\|
\frac{1}{n_{2j_0-1}^{1/q}}\sum_{i=1}^{n_{2j_0-1}}(-1)^{\varepsilon i}y_i\Big\|\leq \frac{2^{20}C}{m_{2j_0-1}^2}
\]
\end{proposition}

Having the above and repeating line for line  the reasoning of Section \ref{subsection h-dependent sequences}, we obtain the following. 

\begin{proposition} Let $T\in\mathcal{L}(\mathfrak{X}_s)$, $s\in\N$. Then 
\[\lim_{\rank(\gamma)\to\infty, \gamma\in\Gamma^s} \mathrm{dist}(Td_\gamma,\mathbb{C}d_\gamma)=0.\]  
\end{proposition}

\begin{proposition}\label{null on basis to null on xs-RIS}
Let $T\in\mathcal{L}(\mathfrak{X}_s)$, $s\in\N$, satisfy $\lim_k\|Td_{\gamma_k}\|=0$ for some sequence $(\gamma_k)_{k=1}^\infty$ of pairwise distinct members of $\Gamma^s$. Then, $\lim_k\|Tx_k\|=0$, for any $(C,1)$-$\mathfrak{X}_s$-RIS $(x_k)_{k=1}^\infty$.     
\end{proposition}

\begin{theorem}\label{xs scalar+compact}
For every bounded linear operator $T:\mathfrak{X}_s\to\mathfrak{X}_s$, $s\in\N$, there exists $\lambda\in\mathbb{C}$ such that $T-\lambda  I$ is compact.    
\end{theorem}

The above theorem, combined with Proposition \ref{same-component-reduction}, proves  Theorem \ref{coordinate scalar-plus-compact} of Section \ref{section calkin algebra XU} and thus ends the justification of the main Theorem \ref{main theorem modulo HI}.

\part*{Appendix}

\section{Spaces with unconditional bases without a $c_0$ asymptotic version}
\label{asymptotics section}
In this section we provide the promised proof of Theorem \ref{no asymptotic space c0 equivalent}. It is a combination of standard techniques with tools about asymptotic spaces from {\cite{maurey:milman:tomczak-jaegermann:1995}} and {\cite{knaust:odell:schlumprecht:1999}}.
\label{asymptotic proofs section}
\begin{definition}
\label{def asymptotic quantification}
For a Banach space $X$ with an FDD $\boldsymbol{E}$ and $n\in\mathbb{N}$ define
\begin{align*}
\theta_n^{\boldsymbol{E}} &= \inf\Big\{\Big\|\sum_{k=1}^ny_k\Big\|:\; (y_k)_{k=1}^n\text{ is an ${\boldsymbol{E}}$-block sequence of vectors with norm at least one}\Big\}\text{ and}\\
q^{\boldsymbol{E}}_\infty &= \limsup_n\Big(\frac{\log(n)}{\log(\theta_n^{\boldsymbol{E}})}\Big).
\end{align*}
\end{definition}

\begin{lemma}
\label{probably standard}
For a Banach space $X$ with an FDD ${\boldsymbol{E}}$, either $q_\infty^{\boldsymbol{E}}<\infty$ or $\sup_n\theta_n^{\boldsymbol{E}}<\infty$ and these alternatives are mutually exclusive.
\end{lemma}

\begin{proof}
The exclusivity of the alternatives is obvious. By elementary properties of the logarithm, the quantity $\limsup_n\log(n)/\log(\theta_n^{\boldsymbol{E}})$ is invariant under passing to an equivalent norm of $X$. We may therefore assume that the FDD of $X$ is monotone. This immediately yields that for all $n\in\N$, $1\leq\theta_n^{\boldsymbol{E}}\leq \theta_{n+1}^{\boldsymbol{E}}$. A standard computation also yields that for all $n,m\in\N$, $\theta^{\boldsymbol{E}}_{nm} \geq \theta_n^{\boldsymbol{E}}\theta_m^{\boldsymbol{E}}$ and, in particular, for all $n,k\in\N$, $\theta_{n^k}^{\boldsymbol{E}}\geq (\theta_n^{\boldsymbol{E}})^k$. 

If $q_\infty^{\boldsymbol{E}} = \infty$, under the assumption of monotonicity, we will show that for all $n\in\N$, $\theta_n^{\boldsymbol{E}} = 1$. Obviously, $\theta_1^{\boldsymbol{E}} = 1$ so fix $n_0\geq 2$. We will show that, for arbitrary $0<\delta<1$, $\theta_{n_0}^{\boldsymbol{E}}\leq n_0^{2\delta}$, which yields the desired conclusion. We may then take $n\geq n_0$ such that $\theta_{n}^{\boldsymbol{E}}\leq n^{\delta}$. Choose $k\in\N$ with $n/n_0\leq n_0^k\leq n$ and, thus,
$(\theta_{n_0}^{\boldsymbol{E}})^k \leq \theta_{n_0^k}^{\boldsymbol{E}} \leq \theta_n^{\boldsymbol{E}} \leq n^{\delta} \leq (n_0^{(k+1)})^\delta.$
\end{proof}

\begin{remark}
Because $\theta_n^{\boldsymbol{E}}\leq n$, for all $n\in\N$, $q_\infty^{\boldsymbol{E}}\geq 1$. If $q_\infty^{\boldsymbol{E}}<\infty$ then for $q>q_\infty^{\boldsymbol{E}}$ there exists $n_0$ such that for all $n\geq n_0$, $\theta_n^{\boldsymbol{E}}\geq n^{1/q}$. Let $\theta = (\lambda\cdot n_0^{1/q})^{-1}$, where $\lambda$ is the FDD projection constant. Then for all $n\in\N$ and ${\boldsymbol{E}}$-block sequences $(y_k)_{k=1}^n$ in $X$ of vectors of norm at least one
\[\Big\|\sum_{k=1}^ny_k\Big\|\geq \theta n^{1/q}.\]
\end{remark}

\begin{remark}
For a Banach space $X$ with an unconditional FDD, $c_0$ is finitely ${\boldsymbol{E}}$-block representable in $X$ if and only if $\sup_n\theta_n^{\boldsymbol{E}}<\infty$.
\end{remark}

The following is well known, it follows, e.g., by applying arguments similar to those in {\cite{johnson:1974}  or by tracking the definitions of {\cite[Section 3]{knaust:odell:schlumprecht:1999}}. However, we could not find a direct reference.
\begin{proposition}
\label{no c0 block equivalent}
Let $X$ be a Banach space with an unconditional FDD ${\boldsymbol{E}}$. Assume that $c_0$ is not finitely ${\boldsymbol{E}}$-block represented in $X$ and, thus, $1\leq q_\infty^{\boldsymbol{E}}<\infty$. Then, for every $q>q_\infty^{\boldsymbol{E}}+1$, there exists $\theta>0$ such that for all ${\boldsymbol{E}}$-block vectors $(y_k)_{k=1}^n$,
\[\Big\|\sum_{k=1}^ny_k\Big\| \geq \theta\Big(\sum_{k=1}^n\|y_k\|^q\Big)^{1/q}.\]
\end{proposition}

\begin{proof}
    Without loss of generality, we may assume that the FDD ${\boldsymbol{E}} = (E_s)_{s=1}^\infty$ of $X$ is $1$-unconditional. First, take $q_0>q_\infty^{\boldsymbol{E}}$ and $\theta_0>0$ such that for all ${\boldsymbol{E}}$-block vectors $(y_k)_{k=1}^n$ in $X$ or norm at least one we have
    \[\Big\|\sum_{k=1}^ny_k\Big\| \geq\theta_0 n^{1/q_0}.\]

    Recall that for $f\in X^*$ we define $\mathrm{supp}(f) = \{k\in\N: E_k\not\subset\mathrm{ker}(f)\}$. In the first step of our proof we will show that for $\Theta = (q_0\theta_0)^{-q_0/(q_0+1)}+(q_0/\theta_0^{q_0})^{1/(q_0+1)}$  and $r = (q_0+1)/q_0$ every ${\boldsymbol{E}}$-block sequence $(f_k)_{k=1}^n$ in $X^*$ of functionals of norm at most one satisfies
\begin{equation}
\label{no c0 block equivalent eq1}
    \Big\|\sum_{k=1}^nf_k\Big\| \leq \Theta n^{1/r}.
\end{equation}
Let $x = \sum_{s=1}^\infty x_s$, where, for $s\in\N$, $x_s\in E_s$, with $\|x\|\leq 1$ and for $1\leq k\leq n$ define
\[y_k = \sum_{s\in\supp(f_k)}x_s.\]
By $1$-unconditionality, for every $F\subset\{1,\ldots,n\}$, $\|\sum_{k\in F}y_k\|\leq 1$. For $\delta > 0$ to be determined later put $F_\delta = \{1\leq k\leq n:\|y_k\| \geq \delta\}$ and note that
$\|\sum_{k\in F_\delta}y_k\| \geq \delta\theta_0(\#F_\delta)^{1/q_0}$ and, thus, $\#F_\delta\leq (\delta\theta_0)^{-q_0} $. We immediately obtain
\[\Big|\sum_{k=1}^nf_k(x)\Big|\leq \sum_{k=1}^n|f_k(y_k)| \leq (\delta\theta_0)^{-q_0} + n\delta.\]
Taking $\delta = q_0^{1/(q_0+1)}\cdot \theta_0^{-q_0/(q_0+1)}\cdot n^{-1/(q_0+1)}$ yields \ref{no c0 block equivalent eq1}.

Next take $q > q_0 + 1$ and let $p = q/(q-1)$ be the conjugate exponent of $q$. Note $p<r$ and, therefore, for $n_0 =\lceil \Theta^{rp/(r-p)}\rceil$, $\Theta n_0^{1/r}\leq n_0^{1/p}$. We deduce that for every ${\boldsymbol{E}}$-block sequence $(f_k)_{k=1}^{n_0}$ of functionals in $B_{X^*}$  we have $n_0^{-1/p}\sum_{k=1}^{n_0}f_k\in B_{X^*}$, and this is iterable. Let $K$ denote the subset of $c_{00}$ that is built as follows.
\begin{equation*}
K_0 = \{\pm e^*_i: i\in\N\},\; K_1 = K_0\cup\Big\{n_{0}^{-1/p}\sum_{k = 1}^n\pm e_{i_k}^*:1\leq n \leq n_0\text{ and }i_1<\cdots<i_n\Big\},\\
\end{equation*}
and if $K_m$ has been defined let
\begin{equation*}
K_{m+1} = K_m\cup \Big\{n_{0}^{-1/p}\sum_{k=1}^ng_k:1\leq n\leq n_0\text{ and } (g_k)_{k=1}^n \text{ are block vectors in } K_m\Big\}.
\end{equation*}
It follows directly by induction that for every $g = (a_1,a_2,\ldots,a_n,0,0,\ldots)\in K=\cup_mK_m$ and ${\boldsymbol{E}}$-block vectors $(f_k)_{k=1}^n$ in $B_{X^*}$ we have $\sum_{k=1}^na_kf_k\in B_{X^*}$. By {\cite[Theorem 1.2]{bellenot:1986}} (see also {\cite[Theorem 3.a]{argyros:deliyanni:1992}}), $(2 n_0^{1/q})^{-1}(B_{\ell_p}\cap c_{00})\subset \mathrm{conv}(K)$. This yields that for every $n\in\N$ and ${\boldsymbol{E}}$-block sequence $(y_k)_{k=1}^n$ in $X$,
\[\Big\|\sum_{k=1}^ny_k\Big\| \geq \frac{1}{2 n_0^{1/q}}\Big(\sum_{k=1}^n\|y_k\|^q\Big)^{1/q}.\]
\end{proof}

The proof of Theorem \ref{no asymptotic space c0 equivalent} relied on Proposition \ref{no c0 block equivalent} above in combination with certain central theorems about asymptotic versions of a Banach space. We list them for completeness.

\begin{definition}[{\cite[Definition 1.1.3]{maurey:milman:tomczak-jaegermann:1995}}]
    Let $X$ be a Banach space with an FDD $(x_s)_{s=1}^\infty$ and let $E$ be an $n$-dimensional normed space with a distinguished normalized basis $(e_i)_{i=1}^n$. We say that $E$ is in the asymptotic structure of $X$ if for every $\varepsilon>0$
    \begin{align*}
    \text{for every }k_1\in\mathbb{N}\text{ there exists }&y_1\in \langle\{x_s:s>k_1\}\rangle,\\
    \text{for every }k_2\in\mathbb{N}\text{ there exists }&y_2\in \langle\{x_s:s>k_2\}\rangle,\\
    &\vdots\\
    \text{for every }k_n\in\mathbb{N}\text{ there exists }&y_n\in \langle\{x_s:s>k_n\}\rangle,
\end{align*}
such that $(y_i)_{i=1}^n$ is $(1+\varepsilon)$-equivalent to $(e_i)_{i=1}^n$.

A Banach space $Z$ with a Schauder basis $(z_i)_{i=1}^\infty$ is called an asymptotic version of $X$ if, for every $n\in\N$, $Z_n = \langle\{z_i:1\leq i\leq n\}\rangle$ is an asymptotic space of $X$.
\end{definition}

\begin{theorem}[{\cite[Theorem 2.3]{maurey:milman:tomczak-jaegermann:1995}}]
\label{universal asymptotic version}
Let $X$ be a Banach space with an FDD ${\boldsymbol{E}}$. Then there exists a universal asymptotic version $Z$ of $X$. That is, any $n$-dimensional vector space $E$ with a distinguished normalized basis $(e_i)_{i=1}^n$ is in the asymptotic structure of $X$ if an only if, for every $\varepsilon>0$, $(e_i)_{i=1}^n$ is $(1+\varepsilon)$-equivalent to a normalized ${\boldsymbol{E}}$-block sequence in $Z$.
\end{theorem}

\begin{notation}
Let $X$ be a Banach space with an FDD ${\boldsymbol{E}}$. A sequence $(x_k)_{k=1}^n$ in $X$ is called a skipped ${\boldsymbol{E}}$-block sequence if, for $1\leq k<n$,  $\max(\supp_{{\boldsymbol{E}}}(x_k)) +1 < \min(\supp_{\boldsymbol{E}}(x_{k+1}))$.
\end{notation}

\begin{theorem}[{\cite[Proposition 1.9]{knaust:odell:schlumprecht:1999}}]
\label{KOS Theorem}
Let $X$ be a Banach space with a monotone FDD ${\boldsymbol{E}}$ and let $Z$ be a universal asymptotic version of $X$. Then, for every sequence $(\varepsilon_k)_{k=1}^\infty$ in $(0,1)$, there exists a blocking $\boldsymbol{\Upsilon}$ of ${\boldsymbol{E}}$ such that the following holds. Any normalized $\boldsymbol{\Upsilon}$-skipped block sequence $(y_k)_{k=1}^n$ in $X$ with $n \leq \min(\supp_{\boldsymbol{\Upsilon}}(y_1))$, is $(1+\varepsilon_n)$-equivalent to a normalized block sequence in $Z$.
\end{theorem}

\begin{proof}[Proof of Theorem \ref{no asymptotic space c0 equivalent}]
We will use a slightly weaker assumption, namely that $U$ has an unconditional FDD ${\boldsymbol{E}}$ and $c_0$ is not an asymptotic version of $U$. One way to proceed is the following. Let $Z$ be a universal asymptotic version of $U$ and let $\boldsymbol{\Upsilon} $ be given by Theorem \ref{KOS Theorem}. Consider $q_\infty^Z$, as in Definition \ref{def asymptotic quantification}, using block sequences of the basis of $Z$. The conclusion of Theorem \ref{KOS Theorem}, together with the unconditionality assumption,  fairly easily yields that $q_\infty^{\boldsymbol{\Upsilon}} = q_\infty^Z$. By the definition of a universal asymptotic version, $c_0$ is not finitely block represented in $Z$ and, thus, $q_\infty^{\boldsymbol{\Upsilon}} = q_\infty^Z <\infty$.

Alternatively, it is possible to replace the use of $q_\infty^{\boldsymbol{\Upsilon}}$ with {\cite[Proposition 3.5]{knaust:odell:schlumprecht:1999}}.
\end{proof}

\bibliographystyle{plain}
\bibliography{bibliography}

\end{document}